\documentclass[times,final]{elsarticle}
\usepackage[margin=2.5cm]{geometry}

\usepackage{framed,multirow}

\usepackage{amsmath}
\usepackage{mathrsfs}
\usepackage{breqn}
\usepackage{bm}
\usepackage{amssymb}
\usepackage{latexsym}
\usepackage{subcaption}
\usepackage{algorithm}
\usepackage{algcompatible}
\usepackage[export]{adjustbox}
\usepackage{placeins}
\usepackage{url}
\usepackage{xcolor}
\usepackage{graphicx}
\usepackage{sidecap}
\definecolor{newcolor}{rgb}{.8,.349,.1}
\usepackage{hyperref}
\usepackage{cleveref}

\allowdisplaybreaks
\begin{document}



\title{A variational volume-of-fluid approach for front propagation}


\author[1]{Ali Fakhreddine}
\author[1]{Karim Alam\' e}
\author[1,2]{Krishnan Mahesh\corref{cor1}}

\cortext[cor1]{Corresponding author: 
  Tel.: +1-612-624-4175;
  Email: krmahesh@umich.edu}

\address[1]{Department of Aerospace Engineering and Mechanics, University of Minnesota, MN 55455, USA}
\address[2]{Department of Naval Architecture and Marine Engineering, University of Michigan, MI 48109, USA}


\begin{abstract}  A variational volume-of-fluid (VVOF) methodology for evolving interfaces under curvature-dependent speed is devised. The interface is reconstructed geometrically using the analytic relations of Scardovelli and Zaleski \cite{Scardovelli2000} and the advection of the volume fraction is performed using the algorithm of Weymouth and Yue (WY) \cite{Weymouth2010} with a technique to incorporate a volume conservation constraint. The proposed approach has the advantage of simple implementation and straightforward extension to more complex systems. Canonical curves and surfaces traditionally investigated by the level set (LS) method are tested with the VVOF approach and results are compared with existing work in LS. \end{abstract}



\maketitle

\section{Introduction}
Traditional interface-capturing (IC) methods, namely the Volume-of-Fluid (VOF), Level-Set (LS), and Phase-Field (PF) methods, have been used to study a wide range of multiphase flow problems over the past decade. Some of the applications include bubbly flows, spray atomization, breaking waves, and phase change (i.e. boiling, cavitation, solidification).  While all three interface-capturing methods mentioned achieve the common task of advecting a fluid interface given a prescribed velocity field, they differ at the level of accuracy, efficiency, and implementation. This consequently makes one method more favorable than the other based on the level of fidelity that the physical problem imposes on the numerical method. IC methods can be categorized by the type of approach used to track the interface. VOF and LS methods fall under sharp-interface approaches while PF methods fall under diffuse-interface approaches. The primary emphasis of this work will be on the former category.
\par VOF and LS methods have long existed in the two-phase modeling community, starting with the early work of Hirt and Nichols \cite{Hirt1981} who introduced the fractional VOF technique to study problems involving free-surface boundaries, and Osher and Sethian \cite{Osher1988} -- the developers of the traditional LS formulation. Further developments in that area gave rise to improved VOF \cite{Zaleski2007,Weymouth2010} and LS formulations that emphasize discrete mass conservation and accurate interface representation. A hybrid method originally developed by Sussman and Pucket \cite{Sussman2000} also currently exists and is known as the coupled-level set and volume-of-fluid method (CLSVOF). The hybrid method utilizes level sets to accurately compute curvature and interface normals and uses the VOF framework to improve mass conservation. Although CLSVOF achieves higher accuracy than the use of either VOF or LS separately, its employment is still less common due to higher computational cost and the inherent difficulty in achieving a balanced load between VOF and LS operations which often leads to poor parallel scalability. 
\par
\begin{figure}[H]
\centering
\includegraphics[width=.7\textwidth]{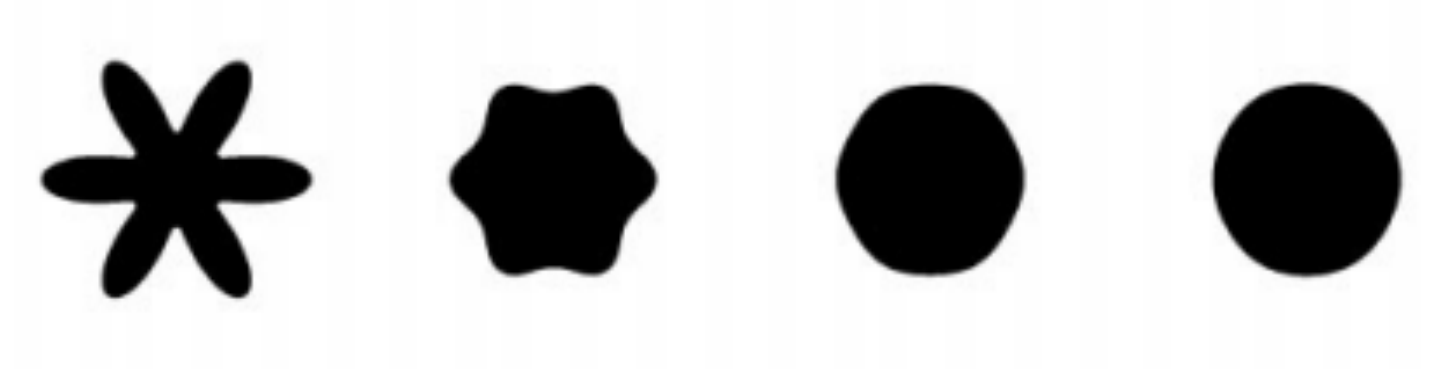}
\caption{Evolution of a star shape under curvature-driven motion \cite{Chopp1999}.}
\label{example_curvM}
\end{figure}

 At the most fundamental level, both VOF and LS present attractive frameworks for modeling two-phase phenomena. However, it is well established that the LS method yields more precise interface normals and curvature estimation which is owed to the mathematical properties of the level-set itself. This naturally made the study of propagating fronts with curvature-dependent speed an area that is unique to the LS framework. Curvature-driven motion can be described as the evolution of an interface due to a self-generated velocity field $\overrightarrow{V}$ that is normal to the interface and proportional to its curvature. An example of this type of motion is shown in Figure \ref{example_curvM} where a 6-petal star shape evolves into a circle due to its curvature; the troughs propagate outwards and the peaks propagate inwards with an interfacial speed proportional to $\kappa$. Chopp and Sethian \cite{Chopp1993} studied hypersurfaces moving under flow that depends on mean curvature using the LS method of Osher and Sethian \cite{Osher1988}. In a later work, Chopp and Sethian \cite{Chopp1999} presented a discussion on various numerical schemes to model the motion of surfaces under the intrinsic Laplacian of curvature. In their work, they test the curvature-driven motion of different shapes (convex and non-convex) using their proposed algorithm. 
While using the LS method is suitable for the study of curvature-driven motion, the presence of high curvature regions imposes a more stringent requirement to solve the reinitialization equation (see \cite{Chopp1993_2,Sussman1994,Peng1999} for details on reinitialization). As the level-set evolves in time, its values begin to drift away from the original signed distance function (SDF); the reinitialization step serves to retain that property of the level-set before the proceeding advection. The frequency of the reinitialization step increases when regions of high curvature are present since the deviation of the level-set from an SDF is more pronounced in those areas of the interface. This has an adverse effect on computational cost depending on the complexity of interface topology. Initial solutions in that direction evolved in the context of variational LS formulations, starting with the work of Li et al. \cite{Li2010} in the image processing community. They proposed an energy functional that regularizes the distance function as the level set deviates from an SDF. This is more commonly known as the distance regularized level set evolution (DRLSE) where the regularization term penalizes the deviation of the level-set from an SDF, therefore, bypassing the need for an explicit reinitialization step. More recently, Alam\'e et al. \cite{Alame2020} devised a robust numerical methodology to predict equilibrium interfaces by utilizing a DRLSE approach. Their variational LS formulation comes in the context of Gibbs free energy minimization where they demonstrate comparisons between analytical solutions of liquid-air interface equilibria and the proposed variational LS method. This work on variational LS provides a stepping stone for a VOF-based formulation that is accurate, efficient, intuitive, and easy to implement. A variational VOF approach does not currently exist in the literature to the best of the authors' knowledge. Such an approach opens the door to a wider range of applications that involve modeling interfacial phenomena using VOF (e.g. bubble dynamics, phase change, etc.) with the additional advantage of discrete mass conservation and interface sharpness if the advection scheme permits these features. The methodology is inherently applicable to geometric and algebraic VOF since the advection equation is invariant with respect to the type of VOF, and since the principles through which the approach is derived do not change based on the method of interface representation. In the current work, the variational VOF (VVOF) approach will use a conservative geometric volume-of-fluid method following the work of Weymouth and Yue (WY) \cite{Weymouth2010}. We investigate constrained and unconstrained curvature-driven motion, with more emphasis on the mathematical definition of the constraint in the following sections. The evolution of curves and surfaces traditionally presented in the LS community is studied using VVOF and qualitative comparisons are drawn.
\par The paper is organized as follows: \S 2 will provide an overview of the volume-of-fluid method with canonical validation problems, \S 3 will present the proposed variational VOF approach and its details, \S 4 will demonstrate results of different geometries evolving under curvature-dependent speed and the relevant discussions, and the paper is summarized in \S 5 with general remarks on the future extension of the method.
\section{Volume-of-Fluid method}
Given a fluid boundary consisting of two phases that share an interface at a particular location $\mathbf{x}$, the phase indicator function $H(\mathbf{x},t)$ is defined as
\begin{align}
H(\mathbf{x},t) =
\begin{cases}
1 \quad \text{if $\mathbf{x}$ is in fluid 1}\\ 0 \quad \text{if $\mathbf{x}$ is in fluid 2}
\end{cases}
\end{align}
As the interface moves and deforms, the shape of each fluid changes, however, each fluid parcel is assumed to be immiscible in the absence of phase change so that the material derivative of $H$ remains zero i.e
 \begin{align}\label{mat_H}
 \frac{DH}{Dt}=\frac{\partial H}{\partial t} + \mathbf{u}\cdot\nabla H=0 \quad .
 \end{align}
This is equivalent to 
 \begin{align}\label{H_advection}
 \frac{\partial H}{\partial t} + \nabla\cdot(\mathbf{u}H) = H(\nabla\cdot\mathbf{u})
 \end{align}
 using $\nabla\cdot(\mathbf{u}H)=\mathbf{u}\cdot\nabla H+H(\nabla\cdot\mathbf{u})$. The volume fraction $C$, often referred to as the color function, is used to define the spatial average of $H$ in each computational cell such that 
 \begin{align}\label{c_avg}
 \displaystyle{C_{i,j,k} = \frac{1}{\delta x \delta y \delta z}\int_{\Omega} H(\mathbf{x},t)}  \ dV
 \end{align}
 where $\Omega$ represents the computational domain, and the subscripts $i$, $j$, and $k$ denote the index locations of the control volume on the Cartesian grid in the $x$, $y$, and $z$ directions, respectively. This consequently results in a continuous distribution of volume fraction values that range between $0<C<1$. We first integrate Eq. \ref{H_advection} over the computational cell and use the definition Eq. \ref{c_avg} of the color function $C$ to obtain
 \begin{align}\label{adv_int}
 \delta x \delta y \delta z \frac{\partial C_{i,j,k}(t)}{\partial t} + \oint_{\partial \Omega} (\mathbf{u}\cdot\mathbf{n}) \ H(\mathbf{x},t) \ dA = \int_{\Omega}H(\mathbf{x},t)(\nabla\cdot\mathbf{u}) \ dV
 \end{align}
 where $\partial \Omega$ is the cell boundary line and $\mathbf{n}$ is the outgoing unit normal. Note that for the case of incompressible flow i.e. $\nabla\cdot\mathbf{u}$ equal to zero, the presence or absence of the righthand side in Eq. \ref{adv_int} is related to the implementation details of the chosen advection algorithm. We also note that our proposed method uses geometric VOF, therefore the description of algebraic VOF will remain brief in what follows as it is included for completeness.
 \subsection{Algebraic VOF}
 In algebraic VOF (aVOF) methods, the color function $C$ is obtained with a numerical approximation (volume-averaged, polynomial or hyperbolic-tangent representation) of $H$. The key difference between algebraic VOF and geometric VOF is the lack of the need to geometrically reconstruct the interface in order to compute the fluxes. Alternatively, the fluxes are calculated algebraically which has the advantage of straightforward implementation on arbitrary meshes and computational efficiency. As aVOF methods evolved, different approaches to compute fluxes were developed, e.g. compressive schemes and THINC (tangent of hyperbola for interface capturing) schemes.
 \par The main idea behind compressive schemes is the assignment of a donor (upwind) and an acceptor (downwind) cell for every cell face with respect to the underlying flow field. Then, based on the angle between the local interface normal and the cell face normal, donor-acceptor schemes switch between compressive downwind and diffusive upwind schemes. Among some of the compressive schemes are: the high-resolution interface capturing scheme (HRIC) \cite{Muzaferija1998}, compressive interface capturing scheme for arbitrary meshes (CICSAM) \cite{Ubbink1999}, switching technique for advection and capturing surfaces (STACS) \cite{Darwish2006}, and high-resolution artificial compressive formulation (HiRAC) \cite{Heyns2013}.
 \par The THINC scheme was first introduced by Xiao et al. \cite{Xiao2005} in 1D. It involves the piecewise reconstruction of the volume fraction inside each cell by making use of the hyperbolic tangent function which is continuous with a controllable thickness. The primary advantage of the THINC scheme is that it eliminates the numerical diffusion of the interface without the need for compression techniques or explicit interface reconstruction. 
 \par While both flux computation approaches make aVOF a more attractive method over the geometric counterpart, aVOF methods have yet to be extensively tested for large-scale simulations. Compressive aVOF methods are in general less accurate than geometric VOF, and the cost of three-dimensional THINC schemes is still a topic of further research.
 
 \subsection{Geometric VOF}
In geometric VOF (gVOF) methods, the interface is approximated geometrically in a computational cell. This type of VOF is summarized in two fundamental steps: the reconstruction step and the advection step.
 \begin{enumerate}[I]
 \item \textbf{Reconstruction step}: The most commonly used reconstruction approach is the piecewise linear interface calculation (PLIC) scheme \cite{Debar}. For PLIC methods, the reconstruction is a two-step procedure: the normal $\mathbf{m}$ (note that $\mathbf{n}=\mathbf{m}/|\mathbf{m}|$) is first determined using the values of $C$ in a given cell and its neighborhood. The equation of the interface segment (a plane in three dimensions) is then written as
 \begin{align}
 	\mathbf{m}\cdot\mathbf{x} = m_xx+m_yy+m_zz=\alpha
 \end{align}
 where $\alpha$ is adjusted until the area under the interface is equal to $\delta_x \delta_y \delta_z C_{i,j,k}$ and $\mathbf{x}$ is defined as the position vector extending from the origin of a computational cell to the interface. \\
 \item \textbf{Advection step}: Given the velocity field of the reconstructed interface, $C$ of the reference phase needs to be shared across neighboring cells in a manner that guarantees volume conservation. Two advection algorithms currently exist to achieve that task: (1) split methods where a series of one-dimensional advections are performed in each spatial direction via operator-splitting of Eq. \ref{adv_int}, and (2) unsplit methods where a single complex flux calculation is done in cells containing the interface.
 \end{enumerate}
 
\subsubsection{Reconstruction method}
The interface reconstruction employed uses the PLIC methodology, and Young's finite difference method to compute the local normal. Young's finite difference method was first developed by Youngs \cite{Youngs1984} and independently by Li \cite{Li1995}. The normal $\mathbf{m}$ is estimated as a gradient using finite difference approximation such that $\mathbf{m}=-\nabla_{\delta x_i}C$. In 2D, Young's method can be summarized by the following steps:
\begin{enumerate}
\item The normal $\mathbf{m}$ is evaluated at the four corners of the central cell $(i,j)$; e.g. the x-- and y--components of $\mathbf{m}$ on the top-right corner are given by
\begin{equation}
	m_{x:i+1/2,j+1/2}=-\frac{1}{2\delta x}(C_{i+1,j+1}+C_{i+1,j}-C_{i,j+1}-C_{i,j})
\end{equation}
\begin{equation}
	m_{y:i+1/2,j+1/2}=-\frac{1}{2\delta y}(C_{i+1,j+1}-C_{i+1,j}+C_{i,j+1}-C_{i,j})
\end{equation}
and similarly for the other three corners. Note that the cell is assumed to be uniform.

\item The cell-centered value of $\mathbf{m}$ is then obtained by averaging the four cell-corners such that 
\begin{equation}
\mathbf{m}_{i,j}=\frac{1}{4}(\mathbf{m}_{i+1/2,j+1/2}+\mathbf{m}_{i+1/2,j-1/2}+\mathbf{m}_{i-1/2,j+1/2}+\mathbf{m}_{i-1/2,j-1/2})
\end{equation}
\end{enumerate}
After $\mathbf{m}_{i,j}$ is found, we use the analytical relations of Scardovellie and Zaleski \cite{Scardovelli2000} to find the slope and intercept of the interface segment.
\subsubsection{Advection method}
For advection, the WY method provides an improvement to the traditional split-direction VOF advection algorithm by making it conservative. Consider the advection equation Eq. \ref{adv_int} written in discrete form after integrating in space and time such that
\begin{equation}
C^{n+1}_{i,j,k}-C^n_{i,j,k} = -\sum_f F_f + \int^{t^{n+1}}_{t^n} \delta t\int_{\Omega}H(\nabla\cdot\mathbf{u})dV \quad ,
\end{equation}
where the first term in the RHS is the sum of the flux term $F_f$ over all faces of $\mathbf{u}H$ defined by
\begin{equation}
F_f=\int^{t^{n+1}}_{t^n}\delta t\int_f u_f(\mathbf{x},t)H(\mathbf{x},t)dx \quad ,
\end{equation}
where $u_f=\mathbf{u}\cdot\mathbf{n_f}$ and $n_f$ is the unit normal vector pointing out of the face $f$. Therefore, the general 3D directionally split advection scheme can be written as
\begin{equation}
C^{(n,d+1)}_{i,j,k}-C^{(n,d)}_{i,j,k} = \Delta F_d + \int^{t^{n+1}}_{t^n}\delta t\int_{\Omega}H\frac{\partial u_d}{\partial x_d}dV = \Delta F_d + \underbrace{c_d\frac{\partial u_d}{\partial x_d}}_{\text{Dilatation term}}
\end{equation}
where $d$ is the Cartesian index that indicates the direction of advection, $\Delta F_d$ is the net flux in that direction whose sign determines if the faces are donating or receiving regions, and $c_d$ is the compression coefficient. Note that when alternating the directions, the sweeps are rotated cyclically to avoid preferential direction. The improvement to the existing split advection method comes in the context of meeting the following set of concurrent requirements -- if for a given algorithm 1) the flux terms are conservative, 2) the divergence term sums to zero, and 3) no clipping or filling of a cell is needed due to the violation of $0<C<1$, then the algorithm must conserve $C$ exactly. Previous VOF advection methods estimate the integral in the dilation term using the volume fraction $C$; this does not necessarily lead to obtaining zero dilatation for incompressible flow since the value of $c_d$ may not be the same along the three Cartesian directions. Alternatively, Weymouth and Yue \cite{Weymouth2010} suggested defining $c_d$ as
\begin{equation}
c_d=H(C^n_{i,j,k}-1/2)
\end{equation}
such that $c_d=0$ if $C^n_{i,j,k}<1/2$ and $c_d=1$ if $C^n_{i,j,k}\ge1/2$ -- this definition guarantees directional invariance. Hence, as long as $c_d$ is treated fully explicitly i.e. using $C$ from the previous time step, the simple advection method:
 \begin{equation}
 \Delta C = \frac{\delta t}{\delta x \delta y \delta z}\Bigg{(}\Delta F_d + c_d\frac{\partial u_d}{\partial x_d}\Bigg{)} \quad \text{for $d=1,..,N$}
 \end{equation}
 satisfies flux conservation and zero velocity divergence. As for the requirement of no overfilling or over-emptying, this can be enforced with the following Courant number restriction
 \begin{equation}
 \delta t\sum_{d=1}^N\bigg{|}\frac{u_d}{\delta {x_d}}\bigg{|}< \frac{1}{2}
 \end{equation}
More details can be found in Tryggvason et al. \cite{Tryggvason2011}, Alam\'e \cite{AlameThesis}, and  Mirjalili et al. \cite{ctrreport}.
\subsection{VOF validation}
This subsection presents validation problems for VOF (2D and 3D) under rigid body rotation and periodic vortex flow. These test cases were performed to determine the extent to which interface integrity is maintained before attempting to propagate the interface with curvature-dependent speed. This represents a more systematic way to separate possible limitations of the variational approach (Sec. \ref{vvof}) from the limitations of VOF itself.
\subsubsection{Zalesak disc and Star shape (2D)}  
\begin{figure}[H]
\centering
\includegraphics[width=.325\textwidth]{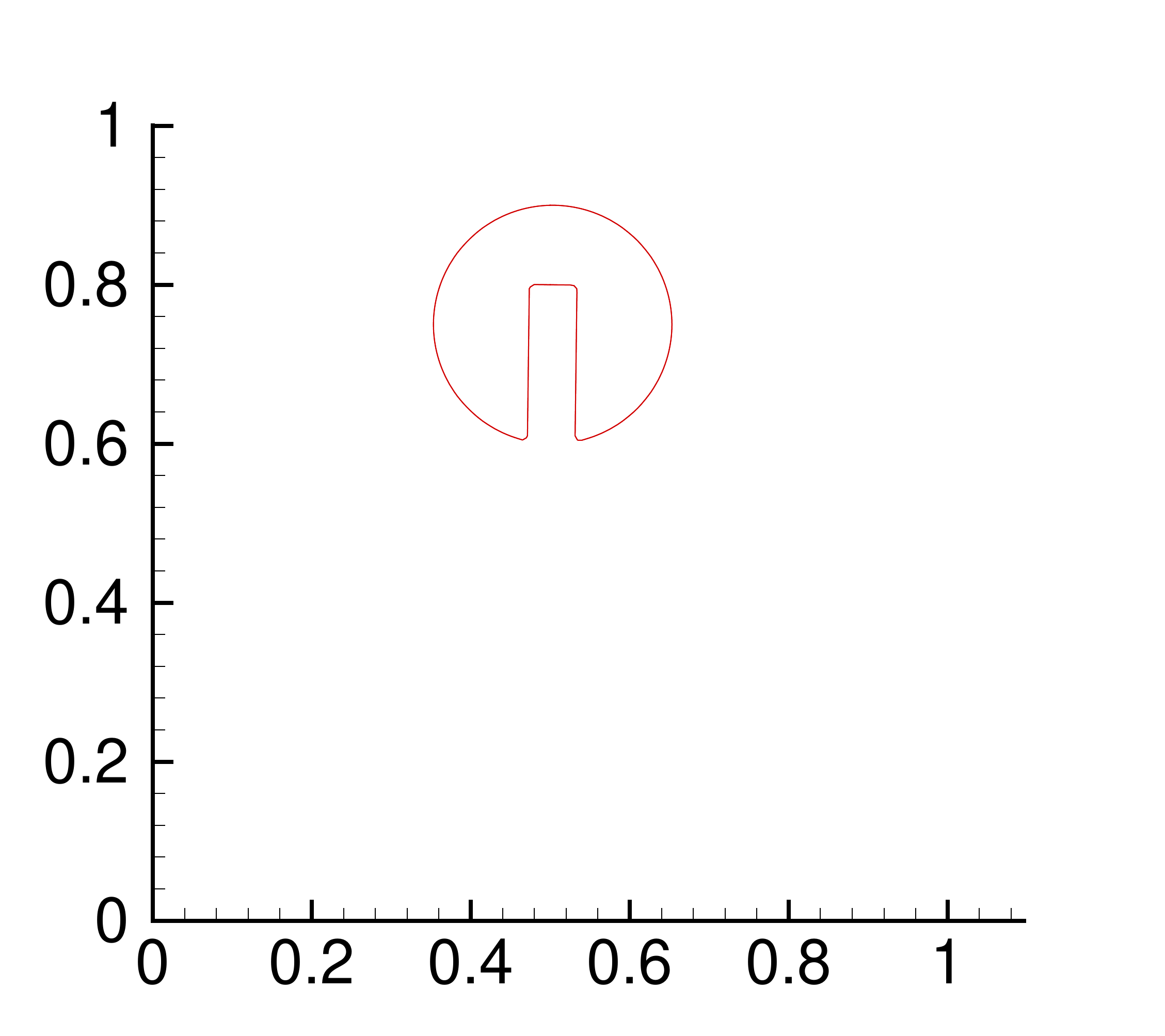}
\includegraphics[width=.325\textwidth]{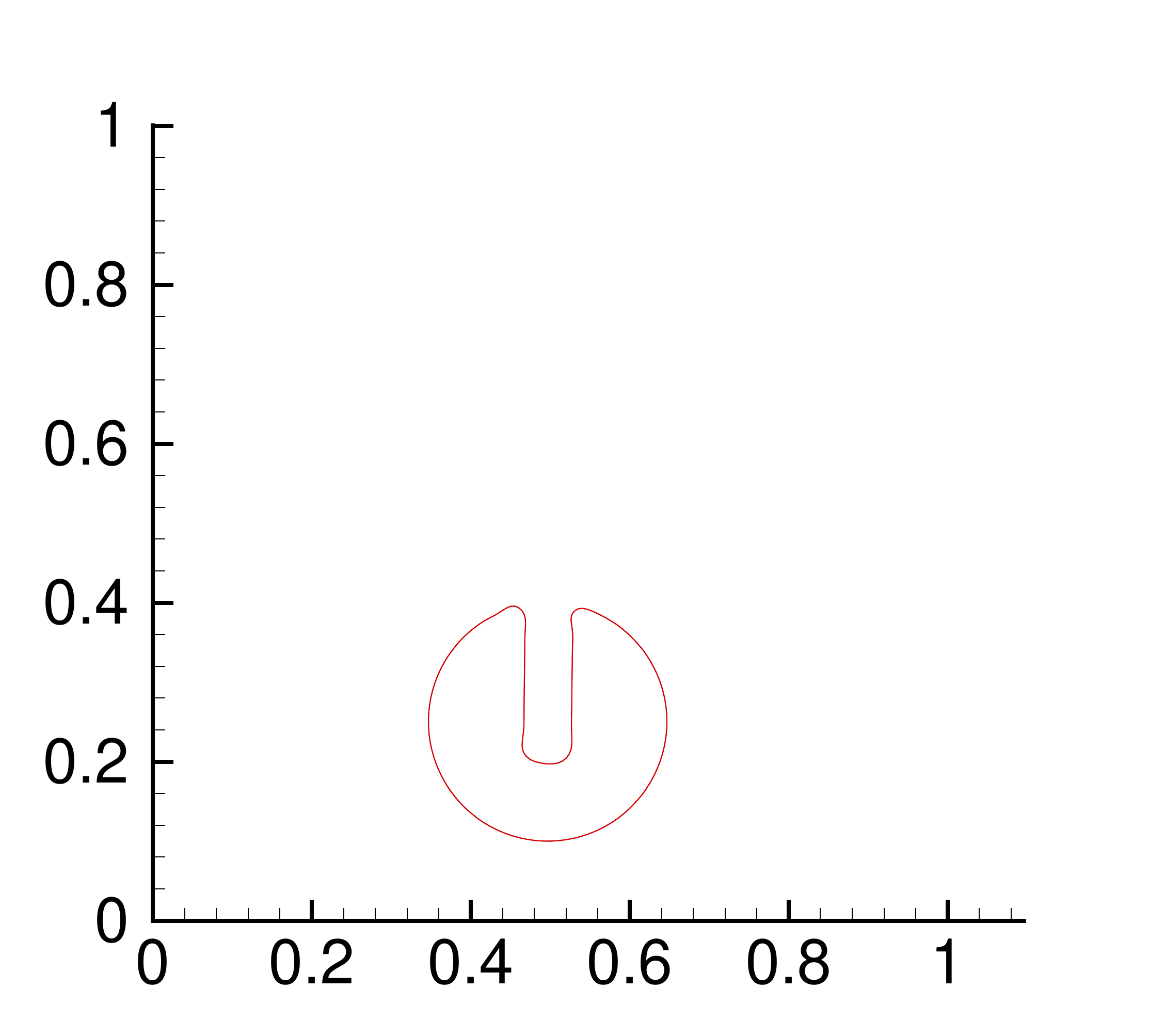}
\includegraphics[width=.325\textwidth]{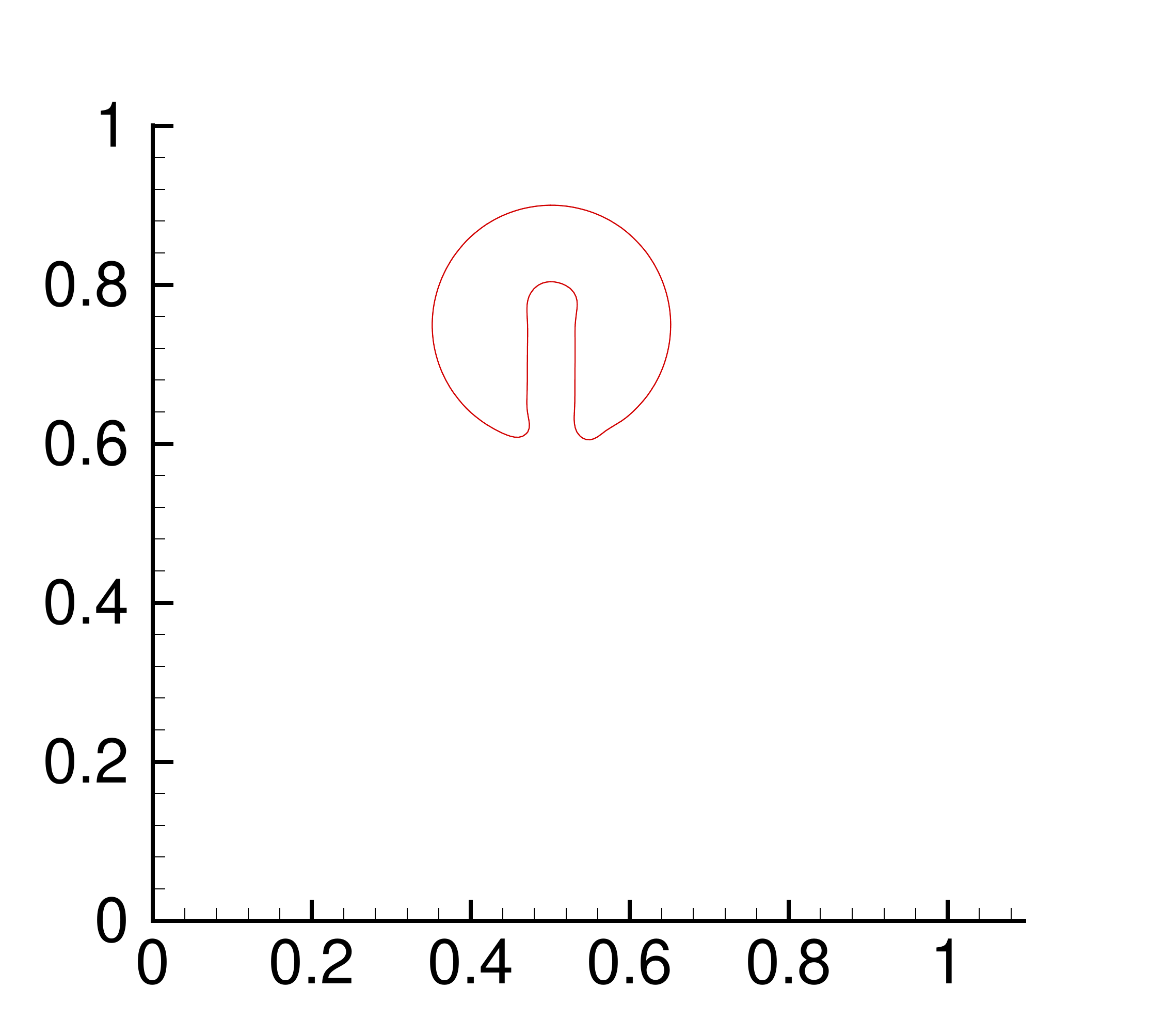} \\
\caption{Evolution of the color function of the Zalesak disc under rigid body rotation. The images are snapshots of the interface in time from left to right where the disc rotates clockwise starting from its initial position at $(0.5,0.75)$ and returning to the same point. Snapshots shown are for  $t=0,0.5,$ and $1$.}
\label{zalesak_disc}
\end{figure}

\begin{figure}[H]
\centering
\includegraphics[width=.325\textwidth]{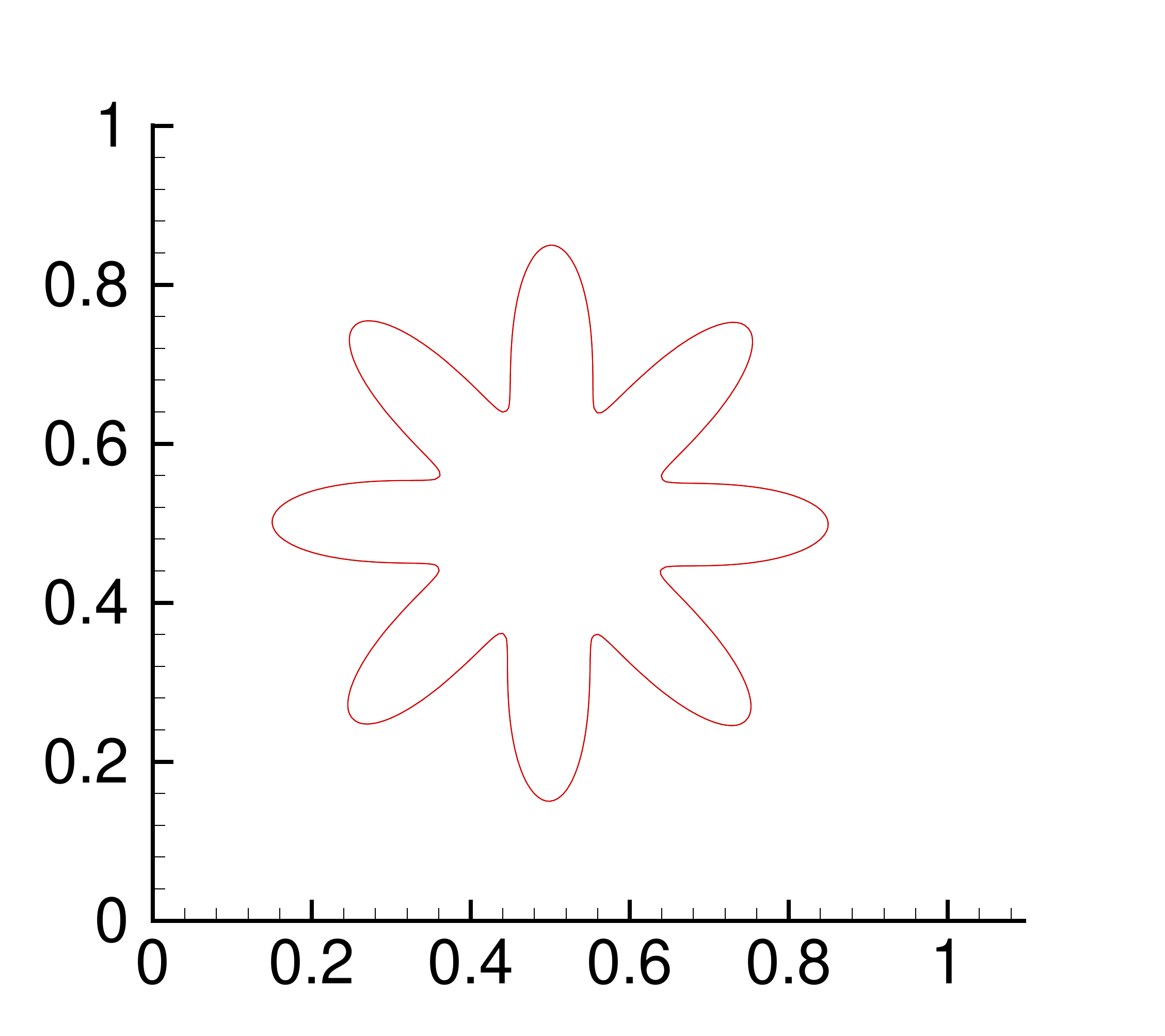}
\includegraphics[width=.325\textwidth]{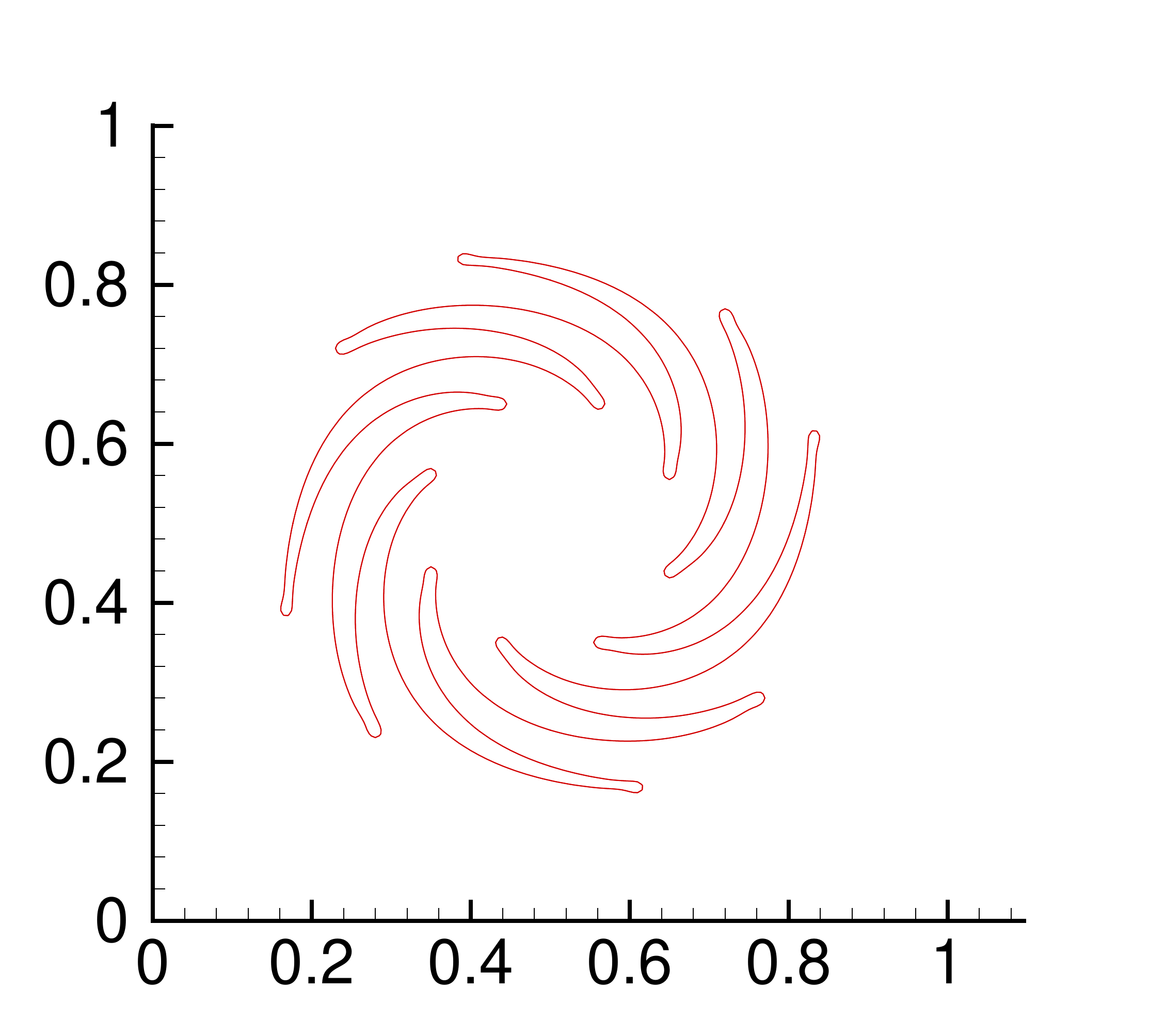}
\includegraphics[width=.325\textwidth]{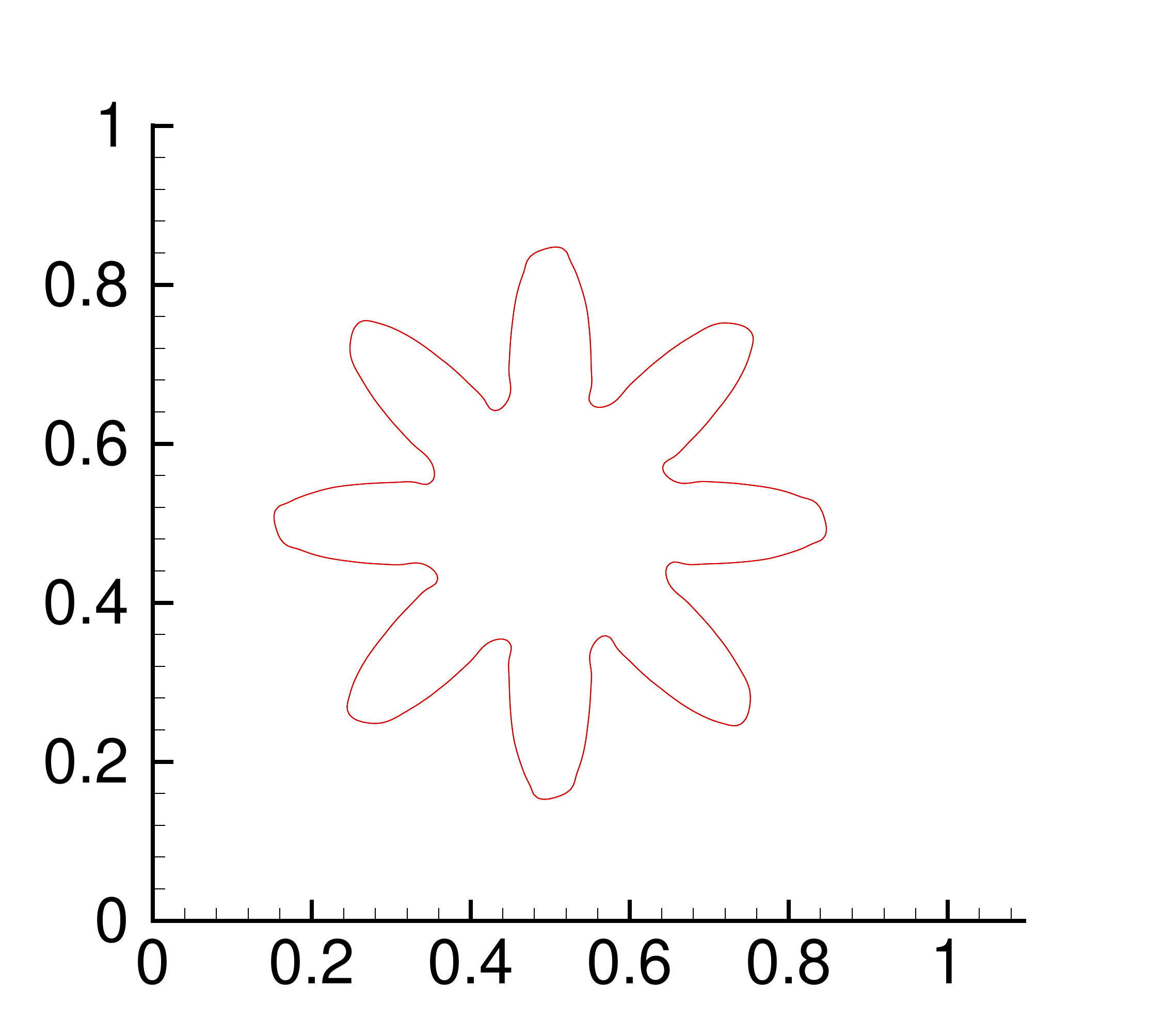} \\
\caption{Evolution of the color function of a star-shape under vortex flow. The images are snapshots of the interface in time from left to right where the star is stretched clockwise until $T/2$ at which the flow is reversed and the initial shape is recovered. Snapshots shown are for  t = 0.0, 1.2, 1.0.}
\label{star_vortex}
\end{figure}

Consider a unit square domain with an initial slotted circular fluid parcel centered at $(0.5,0.75)$. The fluid parcel has a radius of $r=0.15$ and its slot has dimensions $l=0.2$ and $w=0.06$. The slotted disc is allowed to cycle once inside the computational domain under the following velocity field
\begin{equation}
u(y) = 2\pi(y-y_c); \quad
v(x) = -2\pi(x-x_c)
\end{equation}
which describes rigid body rotation, where $x_c$ and $y_c$ are the coordinates of the center of rotation. Fig. \ref{zalesak_disc} shows the interface of the disc ($C\approx 0.5$) at three different instants of the cycle. Note that the results in Fig. \ref{zalesak_disc} are for a $256\times 256$ grid. The $CFL$ number was fixed at $0.192$, and the time-step chosen was equal to $\Delta t=0.001,0.0005,0.00025,$ and $0.000125$ for grid sizes ($32\times 32$), ($64\times 64$), ($128\times 128$), and ($256\times 256$), respectively.
For the second 2D validation, consider a unit square domain with an 8-pointed star shape centered at $(0.5,0.5)$. The star shape is allowed to stretch clockwise until $t=T/2$ at which the flow is reversed and the interface unstretches back to its initial shape at $t=T$. The motion is induced by the following velocity field:
\begin{equation}
u(x,y) = -2\text{sin}^2(\pi x)\ \text{sin}(\pi y) \ \text{cos}(\pi y)\ \text{cos}\Bigg(\frac{\pi}{T}t\Bigg); \quad
v(x,y) = 2\text{sin}^2(\pi y) \ \text{sin}(\pi x) \ \text{cos}(\pi x)\ \text{cos}\Bigg(\frac{\pi}{T}t\Bigg)
\end{equation}
Note that the parametrization of the star shape will be provided in Sec. \ref{results}. Fig. \ref{star_vortex} shows interface evolution ($C\approx 0.5$) of the prescribed geometry at $t=0$, $1.2$, and $2.0$ where $t=1.2$ represents the time at which maximum stretching of the star shape petals occurs. Note that the results in Fig. \ref{star_vortex} are for a $256\times 256$ grid. The $CFL$ number was fixed at $0.064$, and the time-step chosen was equal to $\Delta t=0.001,0.0005,0.00025,$ and $0.000125$ for grid sizes ($32\times 32$), ($64\times 64$), ($128\times 128$), and ($256\times 256$), respectively.

\begin{figure}[H]
\centering
        \begin{subfigure}[c]{0.32\textwidth} 
        \includegraphics[width=\textwidth]{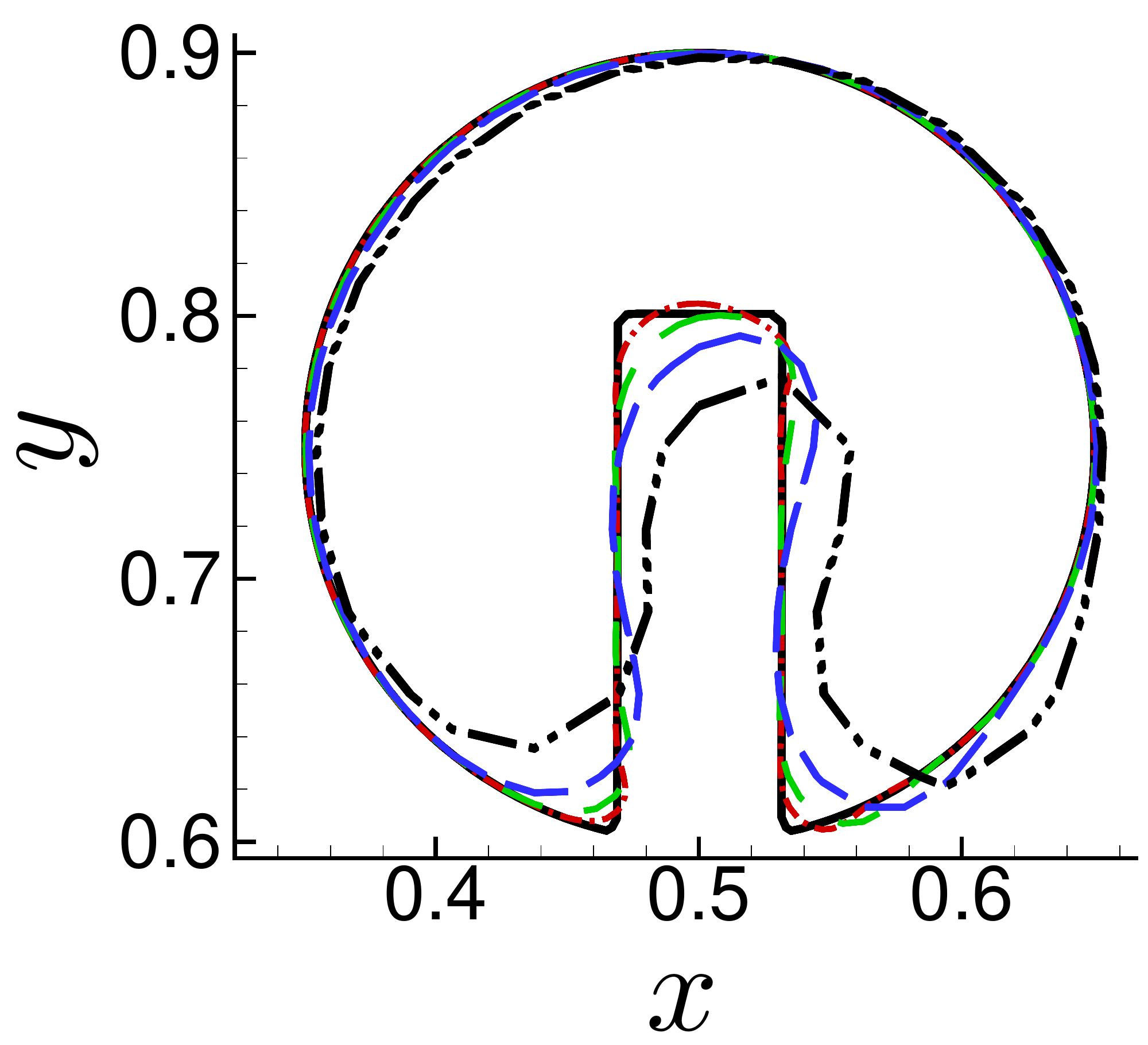}
        \end{subfigure}
         \begin{subfigure}[c]{0.325\textwidth} 
        \includegraphics[width=\textwidth]{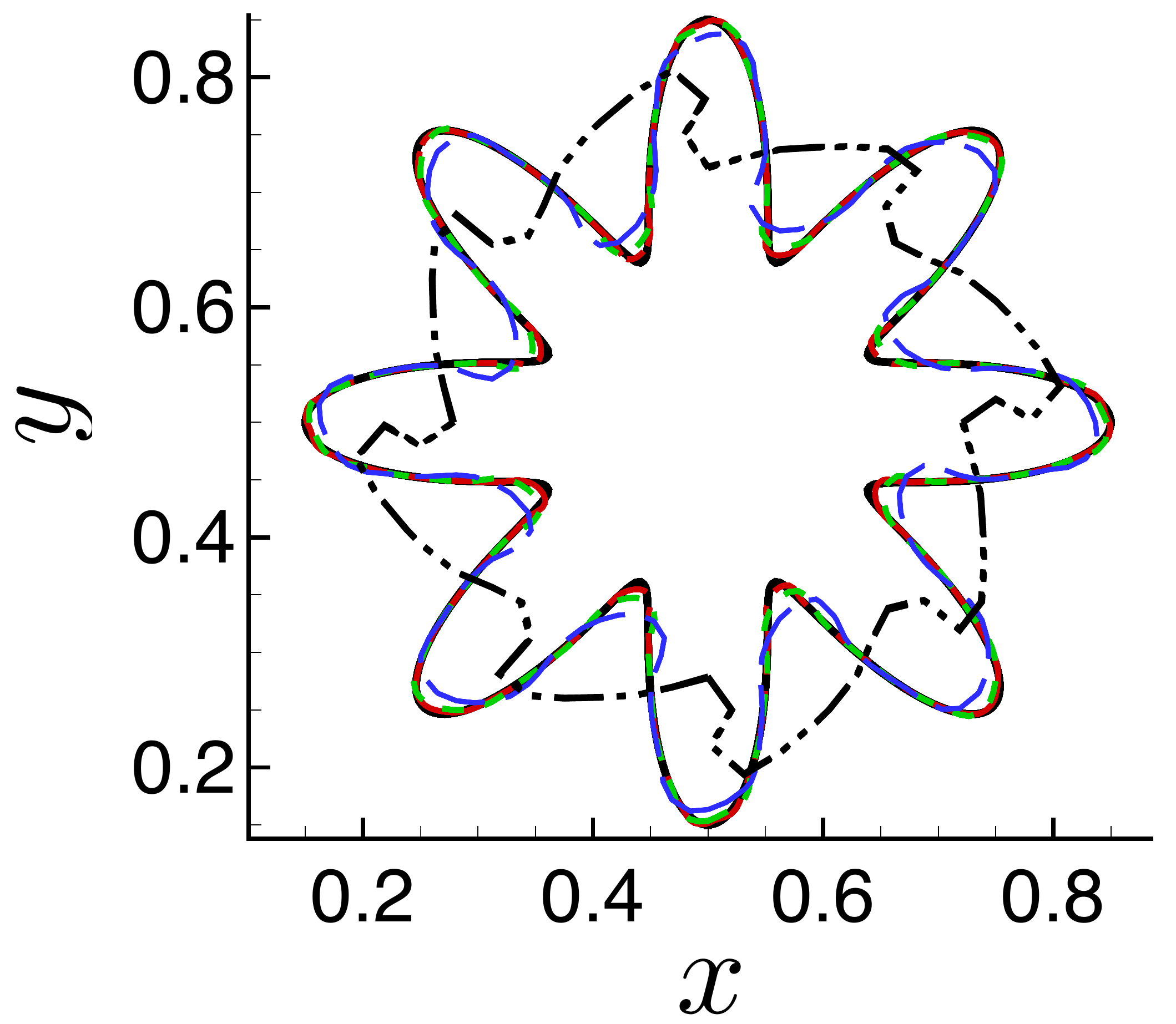}
        \end{subfigure} \hspace{0.5cm}
        \begin{subfigure}[c]{0.25\textwidth} 
        \includegraphics[width=\textwidth]{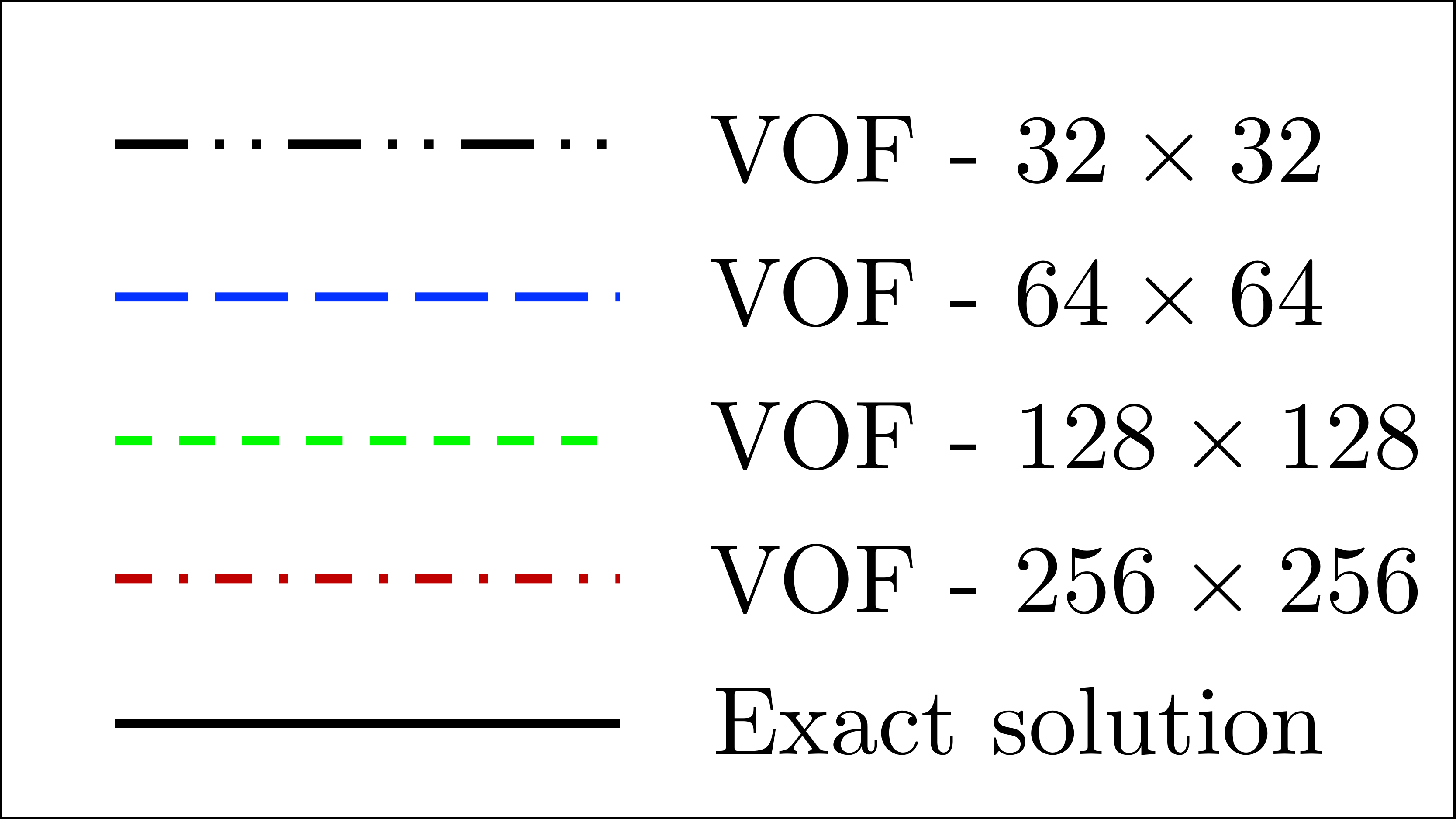}
        \end{subfigure}
    \caption{Comparison between initial and final interface position for different grid sizes for a) Zalesak disc and b) Star shape.}
   \label{comparison}
\end{figure}

\begin{figure}[H]
    \centering
    \begin{subfigure}[c]{0.47\textwidth}
    \includegraphics[width=\textwidth]{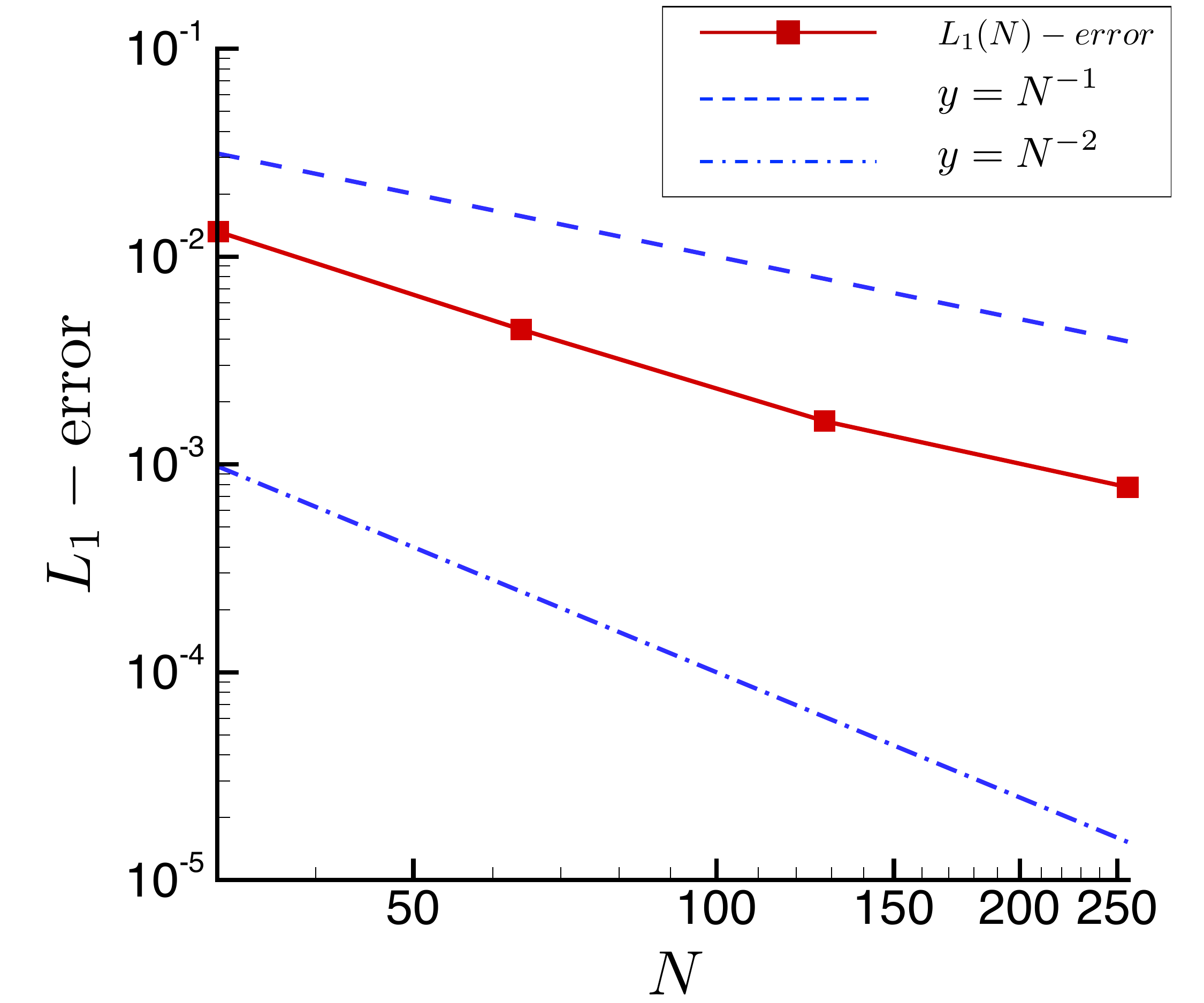}
    \caption{}
    \label{zale_conv}
    \end{subfigure}
    \quad
    \begin{subfigure}[c]{0.47\textwidth}
    \includegraphics[width=\textwidth]{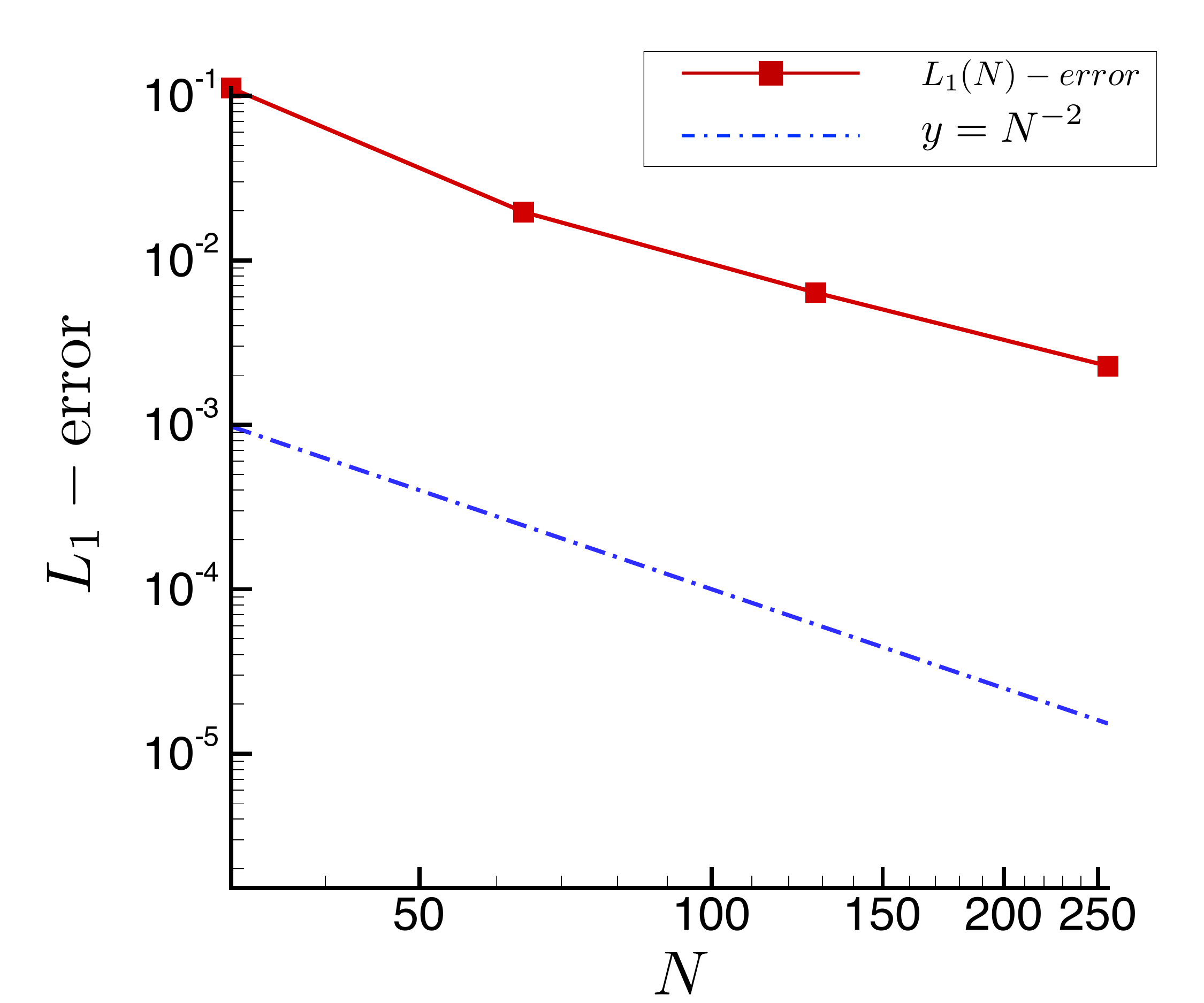}
    \caption{}
    \label{starV_conv}
    \end{subfigure}
    \caption{Error convergence upon grid refinement for a) Zalesak disc and b) Star shape.}
    \label{L1_2D}
\end{figure}

A comparison between the initial and final interface shape is shown in Fig. \ref{comparison} for both 2D geometries and the deviation from the initial shape is measured via the $L_1$ norm error. The $L_1$ norm is defined as follows:
\begin{equation}
L_1 = \frac{1}{NxNy}\sum_{i=1}^{N_x} \sum_{j=1}^{N_y} \bigg|C^f(i,j)-C^0(i,j)\bigg|
\end{equation}
where the superscripts $f$ and $0$ denote final and initial states, respectively. The order of convergence is given by the following relation,
\begin{equation}
n_{L_1}=\displaystyle{\frac{log(L_{1,2N})-log(L_{1,N})}{log(2)}}
\end{equation}
where $L_{1,N}$ is the $L_1$ error using an $N_x\times N_y$ grid and $L_{1,2N}$ is the $L_1$ error for a $2N_x\times 2N_y$ grid. The order of convergence was found to be $n_{L_1}\approx 1.6$ for the Zalesak disc and $n_{L_1}=2.5$ for the star shape as indicated in Fig \ref{L1_2D}. The effect of the first grid refinement is more apparent in Fig. \ref{star_cont} where the final shape deviates significantly from the initial shape using the $32\times 32$ grid.

\vspace{1cm}

\subsection{Sphere (3D)}
\begin{figure}[H]
\centering
\includegraphics[width=.325\textwidth]{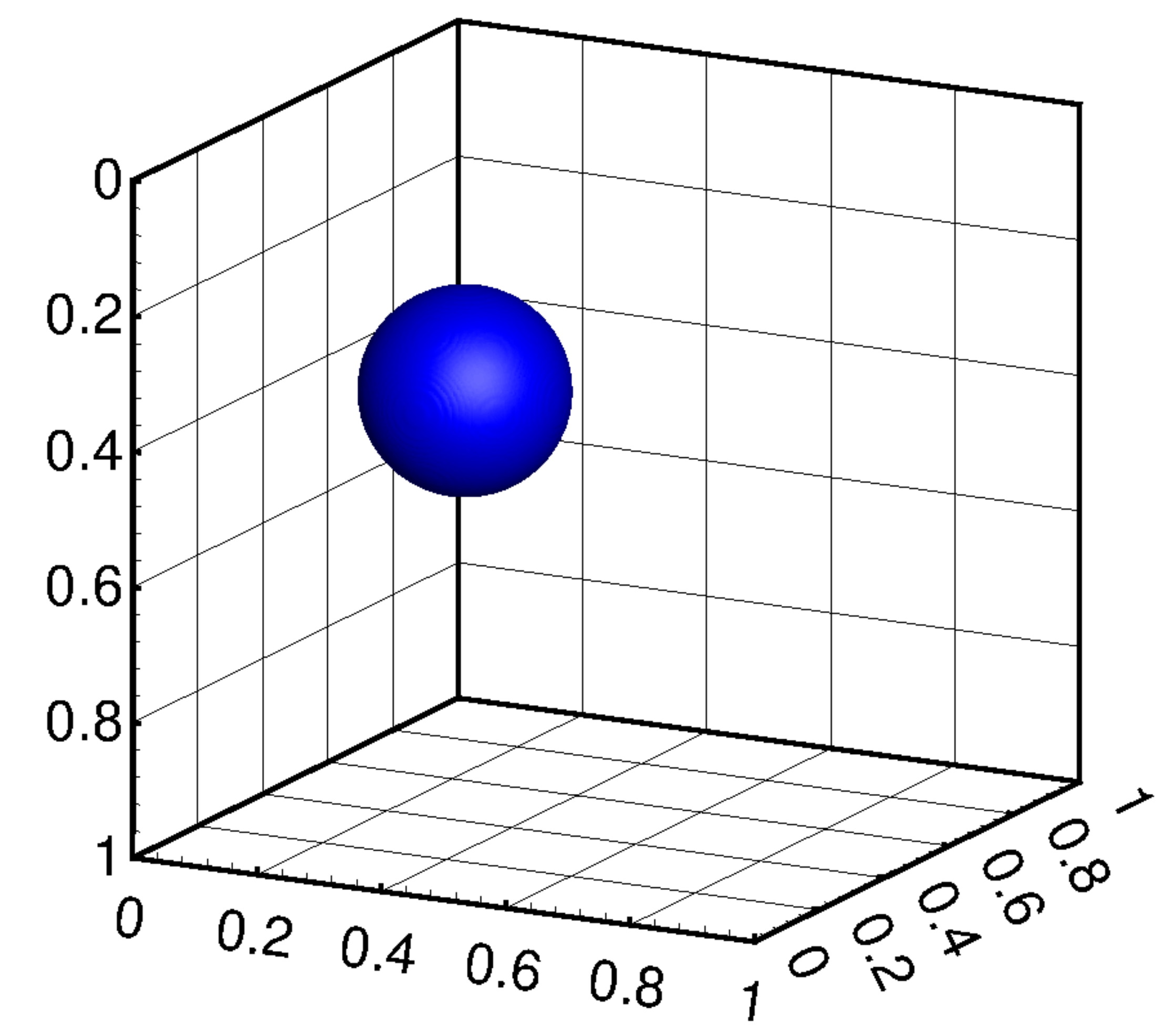}
\includegraphics[width=.325\textwidth]{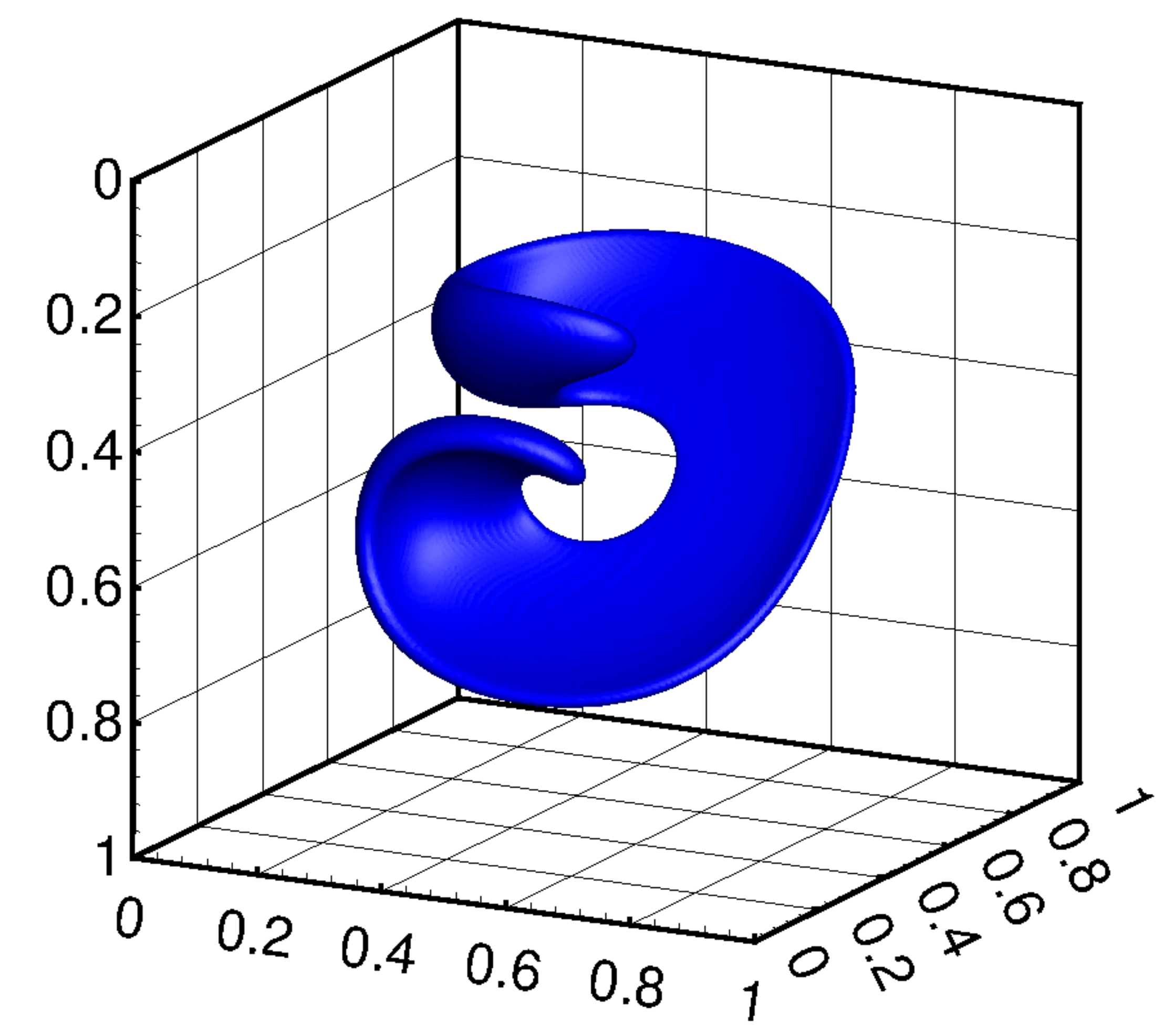}
\includegraphics[width=.325\textwidth]{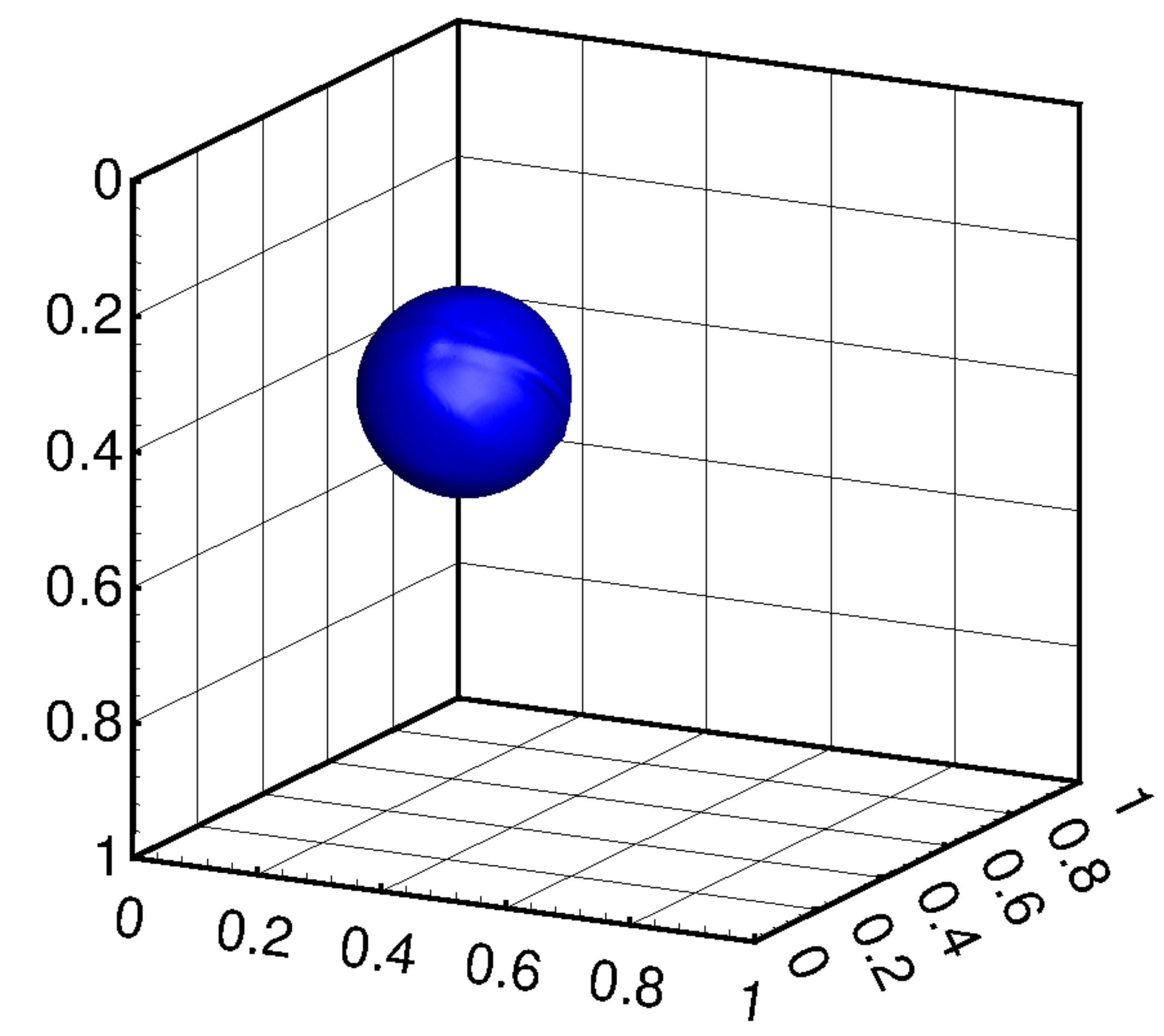} \\
\caption{Evolution of a 3D sphere in a vortex field for a grid size $256^3$ and a period $T=3$. The images are snapshots of the interface in time going from left to right where the middle image is a snapshot at the instant of maximum strain.}
\label{vortex_box}
\end{figure}

The evolution of a 3D sphere is examined. A similar setup is investigated using a $256\times 256\times 256$ grid where the radius of the sphere is $r=0.15$,  its center is at $(0.35,0.35,0.35)$ in a unit cube, and the time period of advection is $T=3$. The velocity field is directly given by:
\begin{align}
u(x,y,z) &= 2\text{sin}^2(\pi x) \ \text{sin}(2\pi y) \ \text{sin}(2\pi z) \ \text{cos}\Bigg(\frac{\pi}{T}t\Bigg) \\
v(x,y,z) &= -\text{sin}^2(\pi y) \ \text{sin}(2\pi x) \ \text{sin}(2\pi z)\ \text{cos}\Bigg(\frac{\pi}{T}t\Bigg) \\
w(x,y,z)&= -\text{sin}^2(\pi z) \ \text{sin}(2\pi x) \ \text{sin}(2\pi y)\ \text{cos}\Bigg(\frac{\pi}{T}t\Bigg) 
\end{align}
The interface stretches significantly (see Fig. \ref{vortex_box}) before returning to its initial position, however, mass is conserved to machine precision.

\section{Variational VOF}\label{vvof}
From the variational point of view, constrained curvature-driven motion of an interface can be achieved via the minimization of interfacial surface area whilst keeping the volume constant. Mathematically, we are trying to minimize the following auxiliary energy functional

\begin{equation}\label{laplacian}
	\mathscr{L}(C) = \mathcal{E}(C) + \lambda G(C)
\end{equation}

\noindent where $\mathcal{E}(C)$ is the interfacial surface energy, $G(C)$ is the volume conservation constraint, and $\lambda$ is the Lagrange multiplier dictating the strength of the volume conservation constraint. The interfacial surface energy and the volume constraint can be expressed as
\begin{center}
\begin{tabular}{c|c}
\textbf{Interfacial surface energy} & \textbf{Volume conservation constraint} \\
\hline \\
$\begin{aligned}[t] 
\mathcal{E}(C) = \int_{\Omega} \delta(C) \ |\nabla C| \ d\mathbf{x} \\
\end{aligned}$ & 
$\begin{aligned}[t]
G(C) = \int_{\Omega} C \ d\mathbf{x} \\
\end{aligned}$
\end{tabular}
\end{center}
whereby 
\begin{equation}
\mathscr{L}(C) = \int_{\Omega} \delta(C) \ |\nabla C| \ d\mathbf{x} +  \int_{\Omega} \lambda C \ d\mathbf{x} \quad .
\end{equation}
The evolution equation of $C$ can be written as a gradient flow that minimizes $\mathscr{L}(C)$ such that 
\begin{equation}\label{var_adv}
\frac{\partial C}{\partial t} = -\frac{\partial \mathscr{L}}{\partial C} \quad .
\end{equation}
We will later show that the advection equation of $C$ can be retained using Eq. \ref{var_adv}. Upon taking the Fr\'echet derivative in the $L^2$ norm of $\mathscr{L}(C)$,
\begin{equation}
	\lim_{\epsilon \to 0}\frac{1}{\epsilon}\bigg{[}\mathscr{L}(C+\epsilon\chi)-\mathscr{L}(C)\bigg{]} = \Bigg{<}\frac{\partial\mathscr{L}}{\partial C},\chi\Bigg{>}=\Bigg{<}\frac{\partial \mathcal{E}}{\partial C},\chi\Bigg{>}+\lambda\Bigg{<}\frac{\partial G}{\partial C},\chi\Bigg{>}
\end{equation}
where $\chi$ is a test function in $L^2$. Therefore, $\frac{\partial \mathscr{L}}{\partial C}= \delta(C)\Bigg{[}\lambda-\nabla\cdot\bigg{(}\frac{\nabla C}{|\nabla C|}\bigg{)}\Bigg{]}$ which leads to the following form of the advection equation
\begin{equation}\label{advection_1}
\frac{\partial C}{\partial t} = \delta(C)\Bigg{[}\nabla\cdot\bigg{(}\frac{\nabla C}{|\nabla C|}\bigg{)}-\lambda\Bigg{]}\quad .
\end{equation}

\noindent To determine the value of the Lagrange multiplier, we enforce the fact that the constraint does not vary in time i.e. volume conservation is maintained at every instant. Mathematically,
\begin{equation*}
\frac{dG(C)}{dt} = 0 \ \forall  \ t \in [0,t_f]
\end{equation*}
\noindent where $t_f$ is the final time. For $dG(C)/dt=0$, $\lambda$ has the following expression
\begin{equation}
	\lambda = \frac{\displaystyle{\int_{\Omega}\nabla\cdot\Bigg{(}\frac{\nabla C}{|\nabla C|}\Bigg{)}\ d\mathbf{x}}}{\displaystyle{\int_{\Omega}\delta(C) \ d\mathbf{x}}} =\frac{\displaystyle{\int_{\Omega}\kappa\ d\mathbf{x}}}{\displaystyle{\int_{\Omega}\delta(C) \ d\mathbf{x}}} = \bar{\kappa} \quad .
\end{equation}
In other words, the Lagrange multiplier dictating the applied constraint is the mean interface curvature. This in itself is an interesting result; it implies that at every point on the interface, the local curvature will have to balance with mean interface curvature to get a minimal energy solution. Geometrically, the minimal energy solution is a circle (a sphere in 3D) for a closed interface.  \par The Dirac delta $\delta(C)$ is typically approximated by $|\nabla C|$ in VOF methods which in this case is a requirement in order to ensure that the motion due to mean curvature is gradient descent of Eq. \ref{laplacian}. Hence, Eq. \ref{advection_1} can be re-written in the following form,
\begin{align}\label{adv_norm}
	\frac{\partial C}{\partial t} &= (\kappa-\bar{\kappa})\ |\nabla C| \quad . 
\end{align}
The more general form of Eq. \ref{adv_norm} is given by 
\begin{equation}\label{ls_eqn}
	\frac{\partial C}{\partial t} + (V_n\mathbf{n}+V_t\mathbf{T})\cdot\nabla C=0
\end{equation}
where $V_n$ is the component of velocity in the normal direction, otherwise known as the normal velocity, and $V_t$ is the tangential component. Since the tangential component is identically zero for motion in the normal direction,  Eq. \ref{ls_eqn} reduces to
\begin{equation} \label{advection2}
	\frac{\partial C}{\partial t} + V_n|\nabla C|=0 \quad .
\end{equation}
Determining an expression for $V_n$ is purely dependent on the physical phenomenon under study. For example, $V_n$ will be proportional to interface curvature $\kappa$ for curvature-driven motion, it will depend on $\Delta P$ or $\Delta T$ for multiphase flows, and it can also be a constant $a$ where $a>0$ or $a<0$.
Since maintaining interface sharpness is paramount, using the form of the advection equation in Eq. \ref{advection2} is disadvantageous as we will later show. Further, the form in Eq. \ref{advection2} does not represent an intuitive starting point for directionally-split advection with VOF. Using the fact that 

\begin{equation}
\displaystyle{\mathbf{n}\cdot\nabla C=\frac{\nabla C}{|\nabla C|}\cdot\nabla C=|\nabla C|} \quad ,
\end{equation}
Eq. \ref{adv_norm} can be written as 
\begin{equation}
\frac{\partial C}{\partial t} - (\kappa-\bar{\kappa}) \ \mathbf{n}\cdot\nabla C = 0
\end{equation}

\noindent leading to the original advection equation of $C$ for incompressible flow
\begin{equation}
\frac{\partial C}{\partial t} + \mathbf{u}\cdot\nabla C = 0
\end{equation}
\noindent where $\mathbf{u}=-(\kappa-\bar{\kappa}) \ \mathbf{n}$ and the unit normal to the interface is $\mathbf{n}=\nabla C/|\nabla C|$. Since the advection algorithm we are following is directionally split, the unit normal vector can be used to project the normal interface velocity in the three Cartesian directions such that
\begin{equation}\label{vel_proj}
u_{\Gamma}=-(\kappa-\bar{\kappa})\frac{\partial C}{\partial x}\frac{1}{|\nabla C|}; \quad v_{\Gamma}=-(\kappa-\bar{\kappa})\frac{\partial C}{\partial y}\frac{1}{|\nabla C|}; \quad w_{\Gamma}=-(\kappa-\bar{\kappa})\frac{\partial C}{\partial z}\frac{1}{|\nabla C|}
\end{equation}
where the subscript $\Gamma$ designates the interface. The computation of the fluxes uses the velocities at the cell face $f$ which we calculate as a simple average based on the following stencil representation
\begin{figure}[H]
\centering
\includegraphics[width=.4\textwidth]{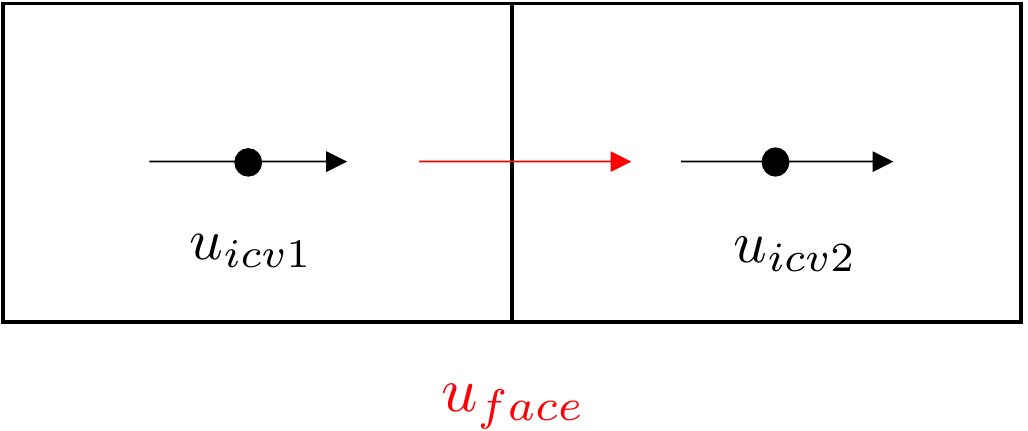}
\end{figure}
such that $\displaystyle{u_{face}=(u_{icv1}+u_{icv2}})/2$.

Note the following:
\begin{enumerate}
\item The sign of $u_{\Gamma}$, $v_{\Gamma}$ and $w_{\Gamma}$ in Eq. \ref{vel_proj} is based on the convention used to define the normal direction i.e. inwards or outgoing. Therefore, omitting the negative sign that multiplies $(\kappa-\bar{\kappa})$ in each term is also valid as long as it is consistent with the direction of the normal.
\item For unconstrained curvature-driven interface motion ($\lambda=0$), $\bar{\kappa}$ vanishes for all components in Eq. \ref{vel_proj}.
\end{enumerate}

\subsection{Curvature estimation}
\begin{figure}[H]
\centering
\includegraphics[width=0.45\textwidth]{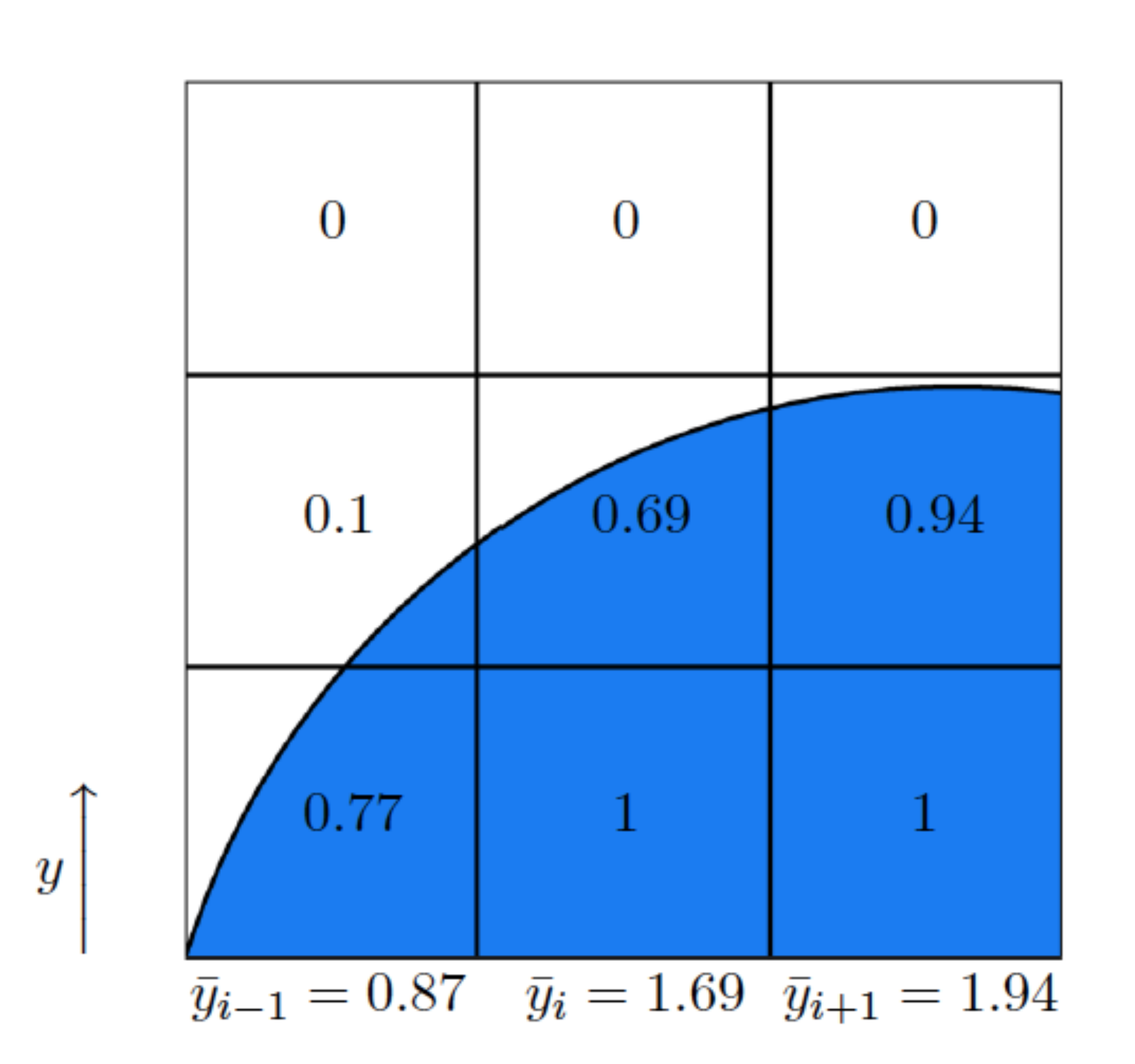}
\caption{Example of the HF method. The normals satisfy $|n_y|>|n_x|$, therefore the summation is in the dominant vertical direction. Note that in this example, the assumption is that $\delta x=\delta y=1$.}
\label{example2d}
\end{figure}

The curvature is defined by $\kappa = \nabla\cdot\mathbf{n}$ and estimated using the height function (HF). Similar to the reconstruction method of WY, where the normal is estimated using the local height function, the HF method is extended to compute curvature \cite{Cummins2005}. The following is a 2D example to demonstrate a sample calculation. \\ As shown in Fig. \ref{example2d}, we have $|n_y|>|n_x|$, therefore the summation of the height function is computed in the vertical direction. Similar to the definition used in the HF method for estimating normals, define:
\begin{equation}
\bar{y}_{i,j}=\sum^1_{l=-1}C_{i,j}\delta y_{i,j+l}
\end{equation}
such that the curvature is given by 
\begin{equation}
\kappa=\frac{n_y}{|n_y|}\frac{\bar{y}_{xx}}{(1+\bar{y}_x^2)^{3/2}}
\end{equation}
where $n_y/|n_y|$ gives the sign of curvature. The derivatives $\bar{y}_{xx}$ and $\bar{y}_x$ are computed using a central difference scheme which leads to second-order method convergence.

\noindent The height function in three dimensions is a straightforward extension, such that
\begin{equation}
\bar{y}_{i,j,k}=\sum^1_{l=-1}C_{i,j+l,k}\delta y_{i,j+l,k}\quad .
\end{equation}
However, the curvature expression is slightly different whereby
\begin{equation}
\kappa=\frac{n_y}{|n_y|}\left(\frac{\bar{y}_{xx}+\bar{y}_{zz}+\bar{y}_{xx}\bar{y}_z^2+\bar{y}_{zz}\bar{y}_x^2-2\bar{y}_{xz}\bar{y}_x\bar{y}_z}{(1+\bar{y}^2_x+\bar{y}_z)^{3/2}}\right)
\end{equation}

\section{Results and Discussion}\label{results}
\begin{figure}[H]
    \centering
    \begin{subfigure}{0.33\textwidth}
        \includegraphics[width=\textwidth]{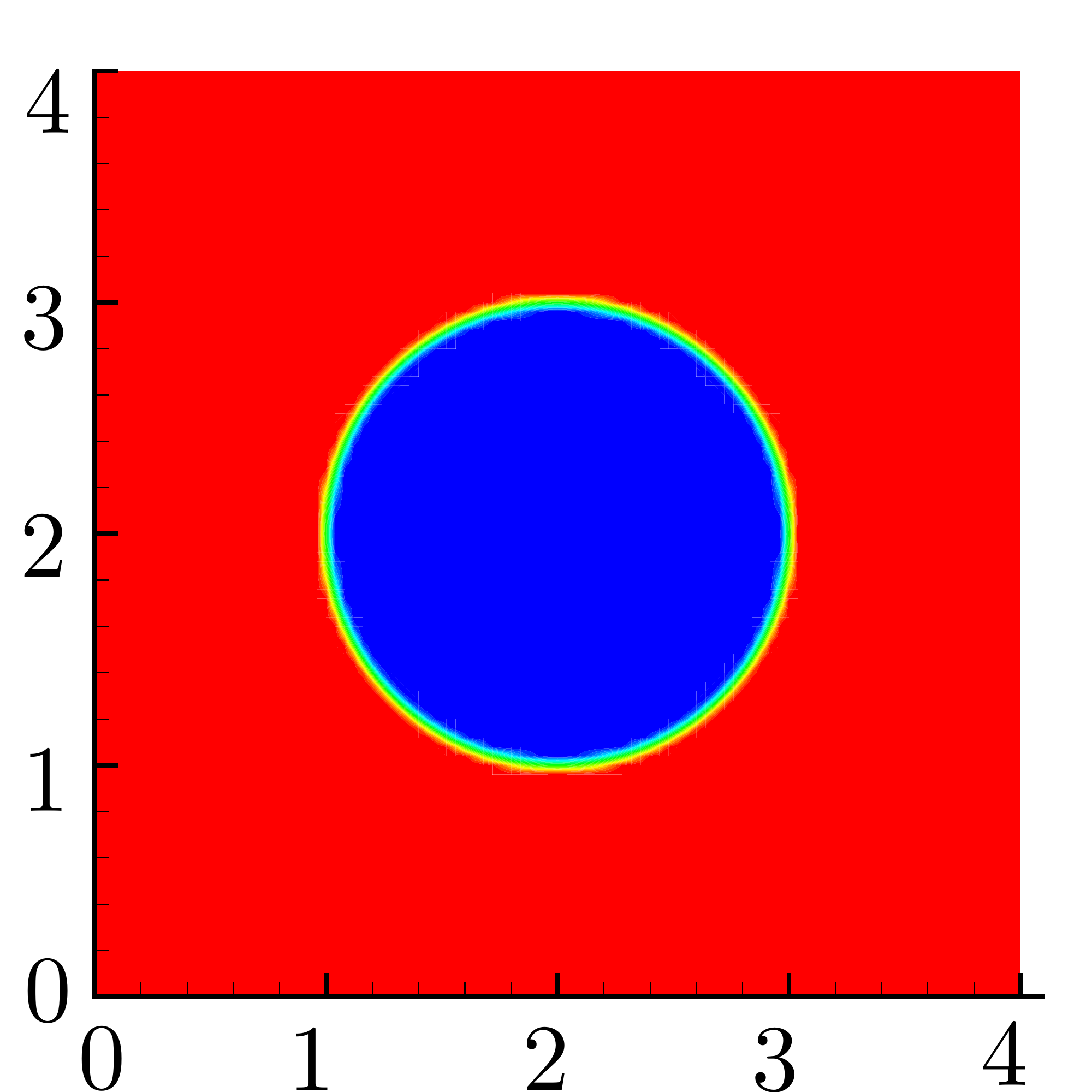}
    \end{subfigure}
    \begin{subfigure}{0.33\textwidth}
        \includegraphics[width=\textwidth]{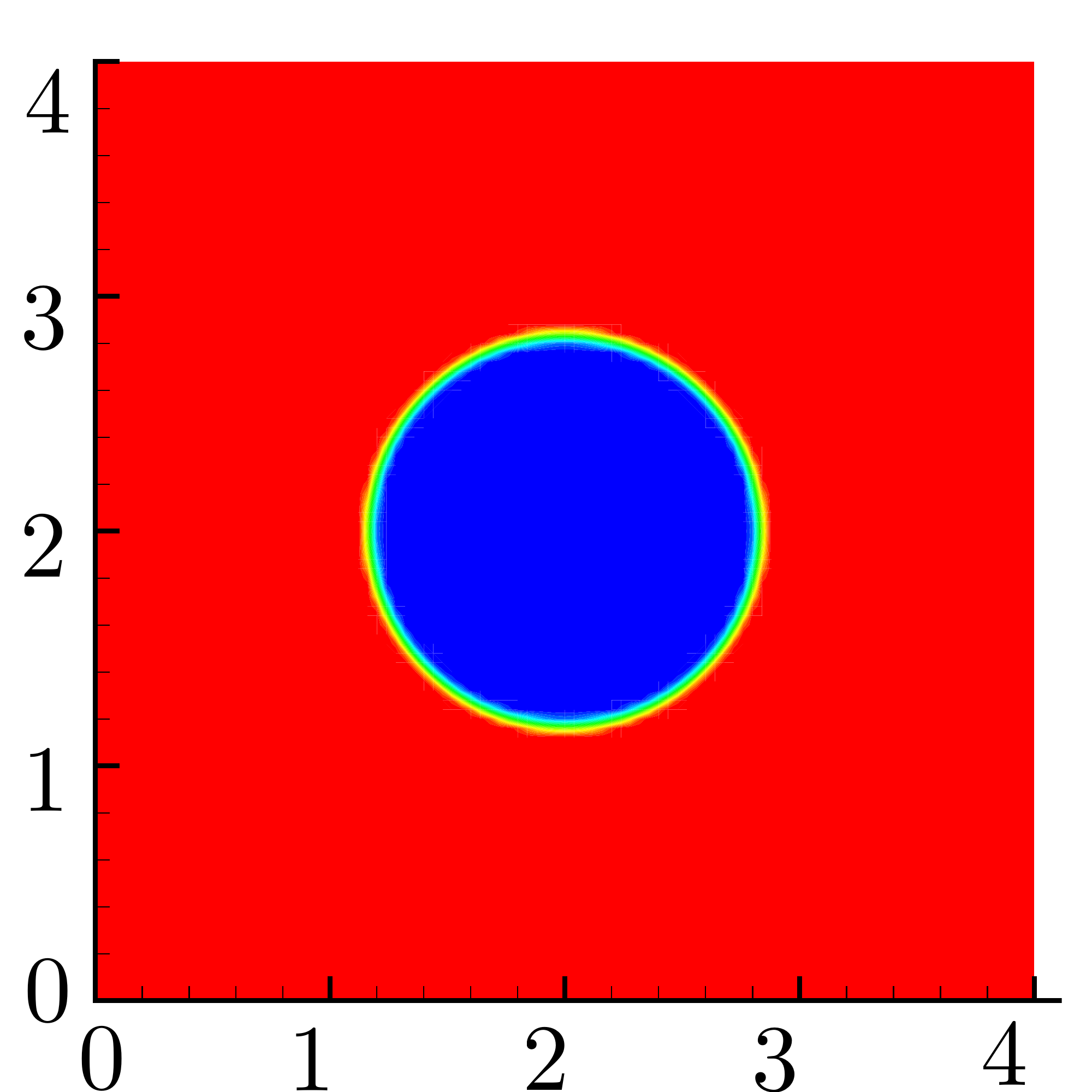}
    \end{subfigure}
    \begin{subfigure}{0.33\textwidth}
        \includegraphics[width=\textwidth]{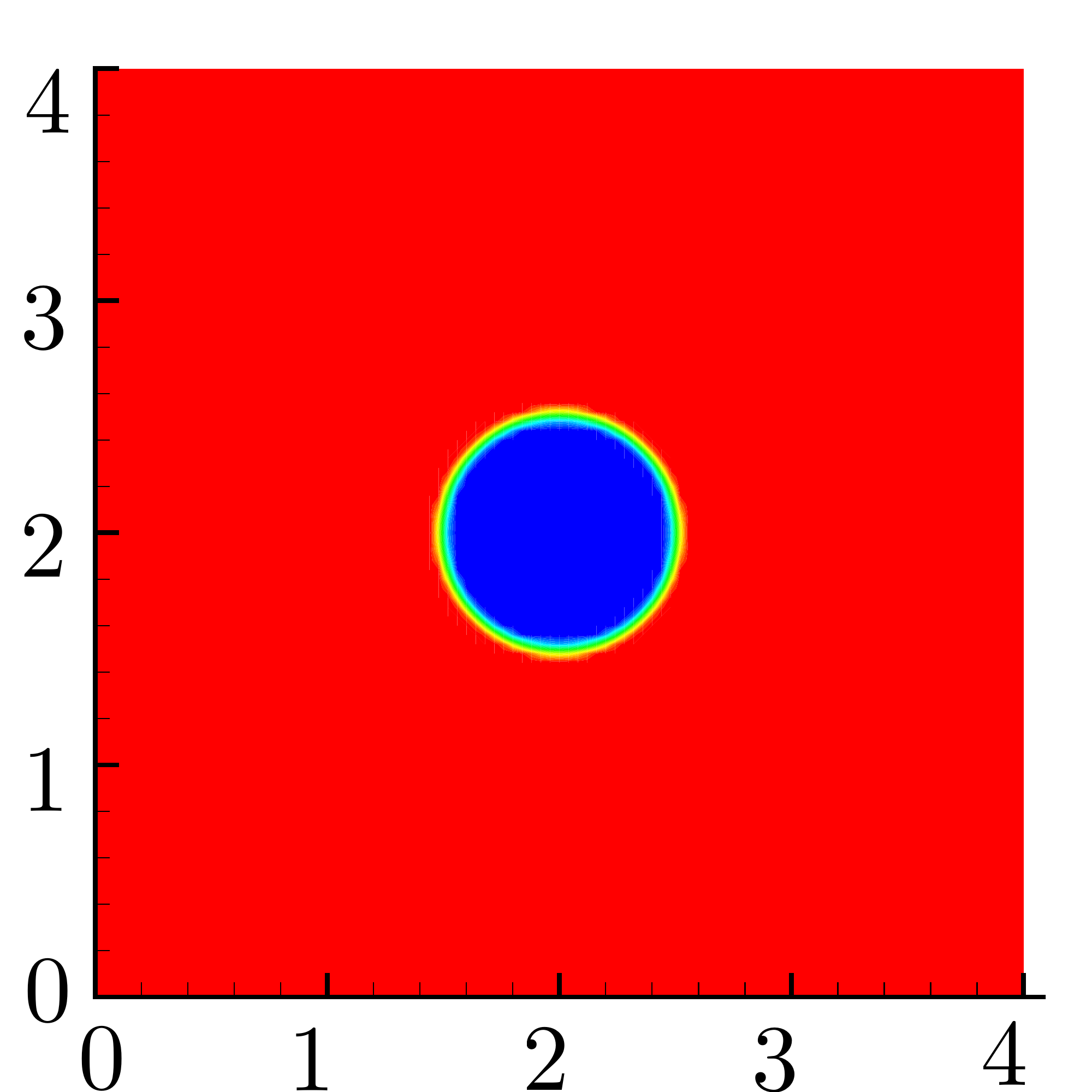}
    \end{subfigure} 
    \caption{Snapshots of the color function as the bubble collapse using directionally-split VOF.}  
    \label{cfunc_snaps} 
\end{figure}

In order to test the validity of the suggested projection in Eq. \ref{vel_proj}, Rayleigh-Plesset (RP) bubble collapse was simulated. The velocity field driving the problem is interfacial in nature and serves as an adequate check given the presence of an analytical solution. In the absence of inertial effects, surface tension, and non-condensable gas (NCG), the RP equation reduces to
\begin{equation}\label{RP}
R\ddot{R} + \frac{3}{2}(\dot{R})^2 = \frac{P_B(t)-P_\infty(t)}{\rho_L}
\end{equation}
where $\Delta P = P_B(t)-P_\infty(t)$ and $\rho_L$ is the liquid density. A constant $\Delta P$ was chosen and $\rho_L$ was set to 1. Note that the value of $\Delta P$ is not of particular importance to the current problem since the comparison with the analytical solution is done on a non-dimensional basis. Also, note that the solutions $R$ and $\dot{R}$ are retrieved after solving Eq. \ref{RP} using the Runge-Kutta RK4 method. Since $u(R,t)=\dot{R}$ acts normal to the interface, the velocity field can be written as 
\begin{equation}\label{rp_split}
u_{\Gamma} = \dot{R}\frac{\partial C}{\partial x}\frac{1}{|\nabla C|}; \quad v_{\Gamma} = \dot{R}\frac{\partial C}{\partial y}\frac{1}{|\nabla C|}; \quad w_{\Gamma} = \dot{R}\frac{\partial C}{\partial z}\frac{1}{|\nabla C|}
\end{equation}
similar to Eq. \ref{vel_proj}. A bubble of radius $r=1$ is initialized at the center of a cubic domain with $L=4$, and its interface evolves via $\dot{R}$ for different grid resolutions.
\begin{figure}[H]
\centering
\begin{subfigure}{0.495\textwidth}
\includegraphics[width=\textwidth]{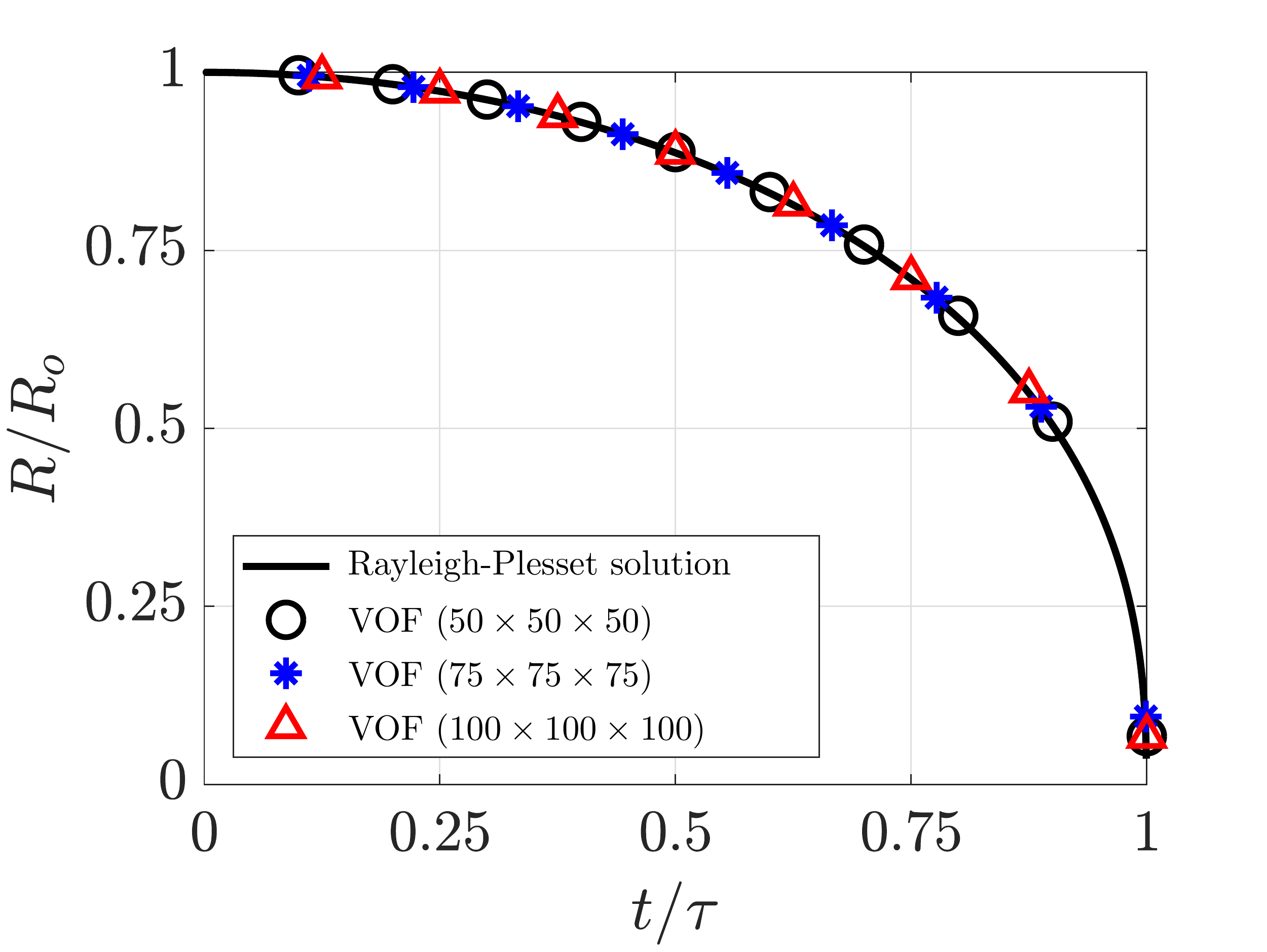}
\end{subfigure}
\begin{subfigure}{0.495\textwidth}
    \includegraphics[width=\textwidth]{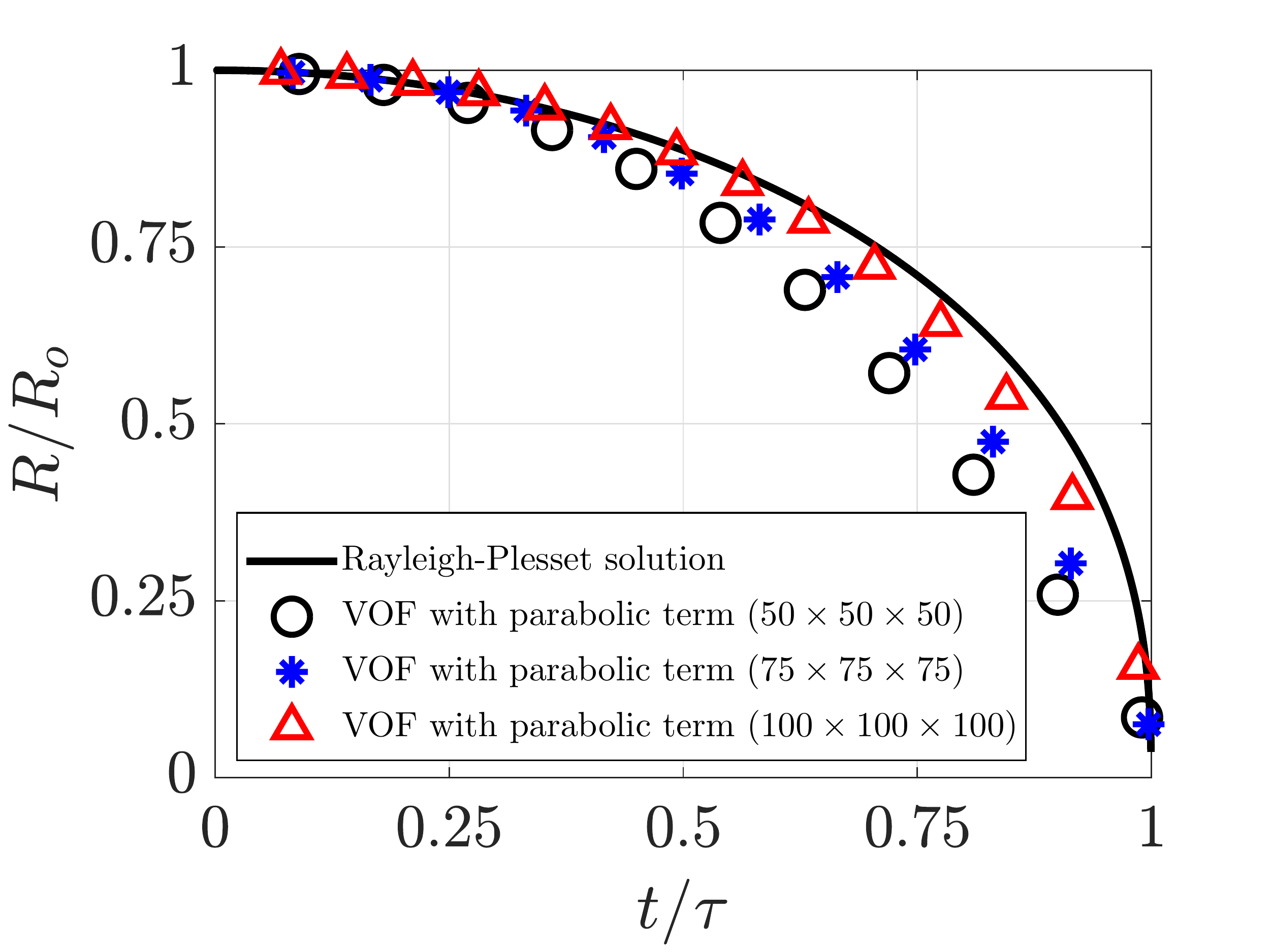}
    \end{subfigure}
\caption{Comparison of bubble collapse profile for three different grid sizes using dir-split VOF (left) and nondir-split VOF (right).}
\label{cRP}
\end{figure}
VOF advection demonstrates good agreement with the RP solution even for the coarsest grid as shown in Fig. \ref{cRP} (left). The result signifies that modeling the prescribed velocities adequately has a stronger effect on the accuracy of the result than the excessive refinement of the interfacial area. This note is of particular relevance to problems involving phase change modeling (e.g. cavitation, boiling, etc.). Fig. \ref{cfunc_snaps} shows central slices of the color function at different time instants as the bubble shrinks where it can be seen that interface sharpness is maintained.
\noindent However, upon driving the interface with $\dot{R}$ via a parabolic term on the right-hand side (same form as Eq. \ref{adv_norm}) such that
\begin{equation}\label{rp_source}
	\frac{\partial C}{\partial t} = -\dot{R}|\nabla C|\quad ,
\end{equation}
it is evident that agreement with the RP solution is significantly compromised as shown in Fig. \ref{cRP} (right). The mismatch is reduced upon grid refinement, however, agreement with RP is still unsatisfactory for the finest grid. Fig. \ref{cfunc_snaps} shows central slices of the color function at different time instants where interface diffusion increases consistently as the bubble shrinks in size. The smearing of the interface has a direct effect on the accuracy of the calculated velocity components which leads to the mismatch with RP. Therefore, advection using Eq. \ref{rp_source} is not recommended. Note that the scheme employed throughout this work is that demonstrated in Eq. \ref{rp_split} and not that in Eq. \ref{rp_source}.
\begin{figure}[H]
    \centering
    \begin{subfigure}{0.33\textwidth}
        \includegraphics[width=\textwidth]{Figures/colorfunc_nx100_dirsplit_00000001.pdf}
    \end{subfigure}
    \begin{subfigure}{0.33\textwidth}
        \includegraphics[width=\textwidth]{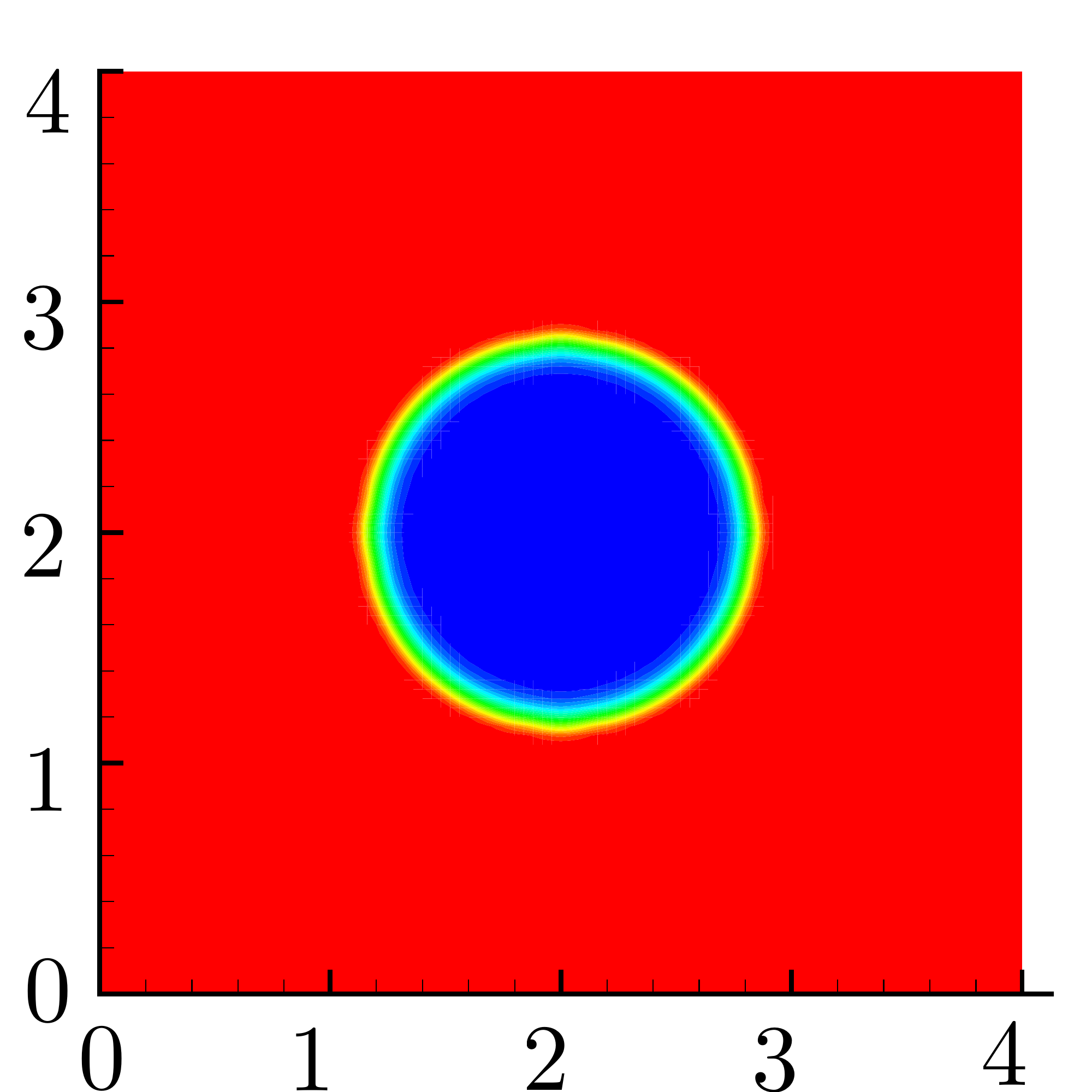}
    \end{subfigure}
    \begin{subfigure}{0.33\textwidth}
        \includegraphics[width=\textwidth]{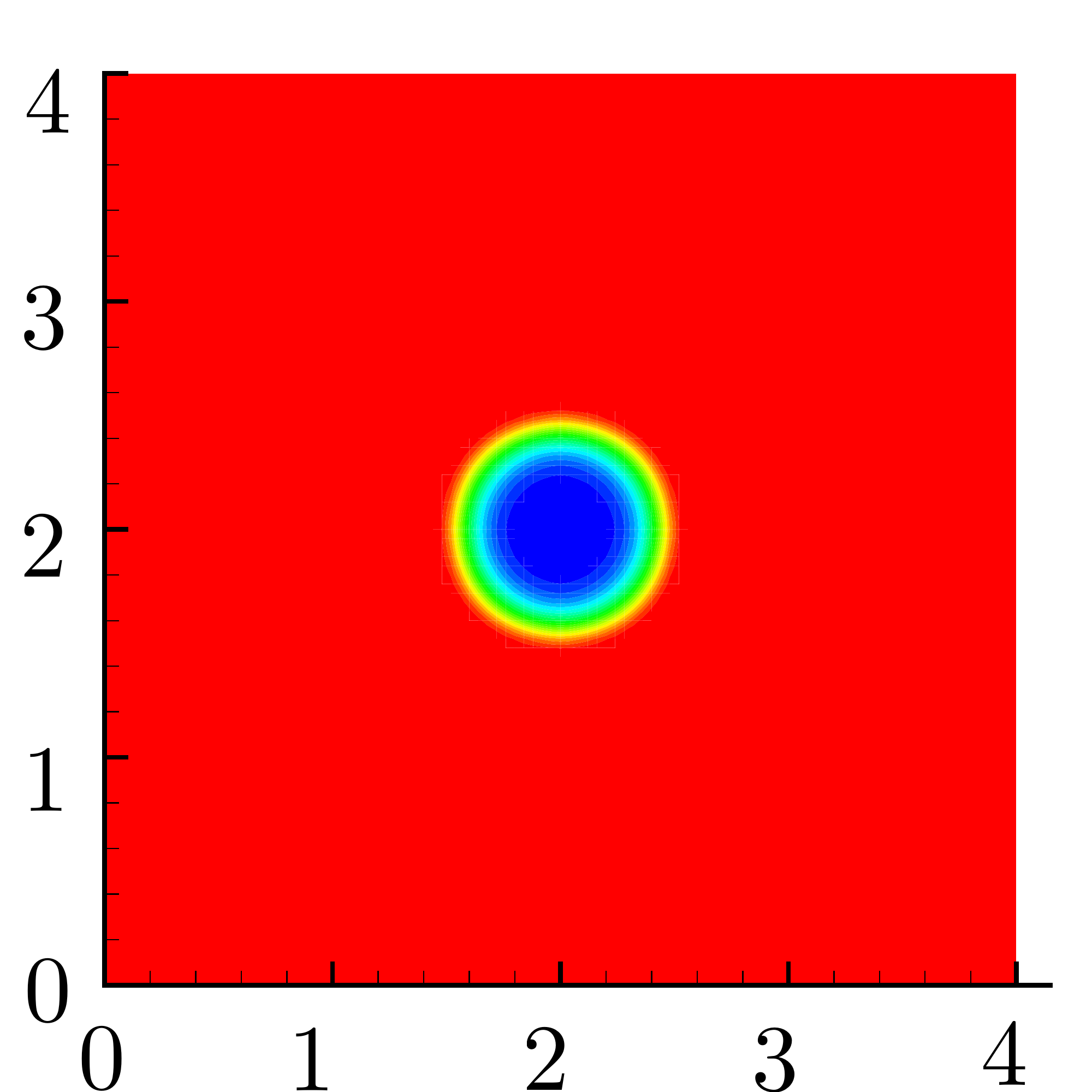}
    \end{subfigure}    
    \caption{Snapshots of the color function as the bubble collapse using nondir-split VOF.}  
\end{figure}

\subsection{Curvature-driven motion}
In what follows, we investigate the curvature-driven motion of various geometries using the proposed approach. The interface is defined using the zeroth level-set $\phi(\mathbf{x,0})$ which is then used to construct the color function for VOF. We also make use of boolean operators to initialize more complex geometries which represent a union of multiple level-set functions.
\subsubsection{Pointed star} 
\begin{figure}[H]
\centering
\includegraphics[width=.33\textwidth]{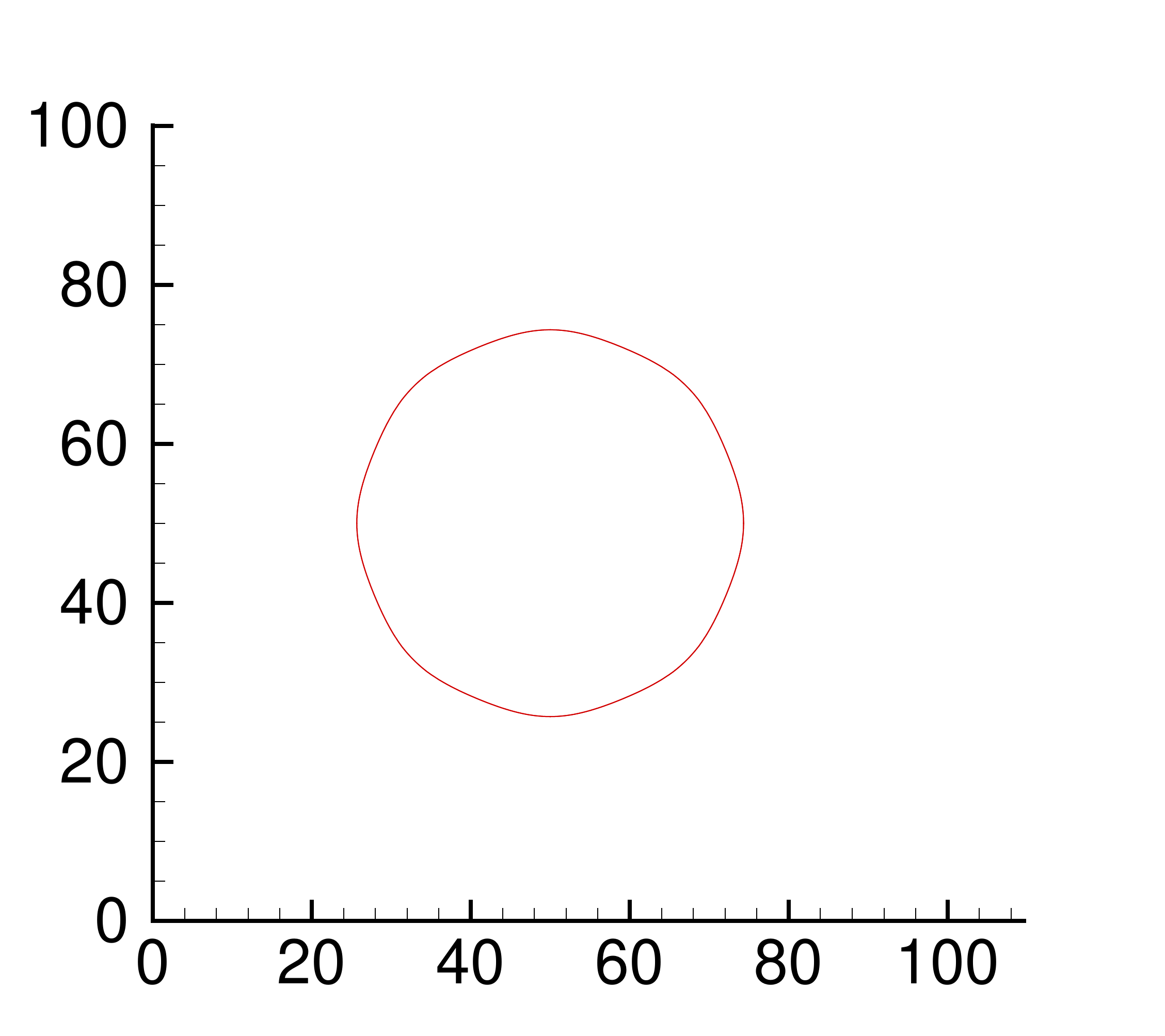}
\includegraphics[width=.33\textwidth]{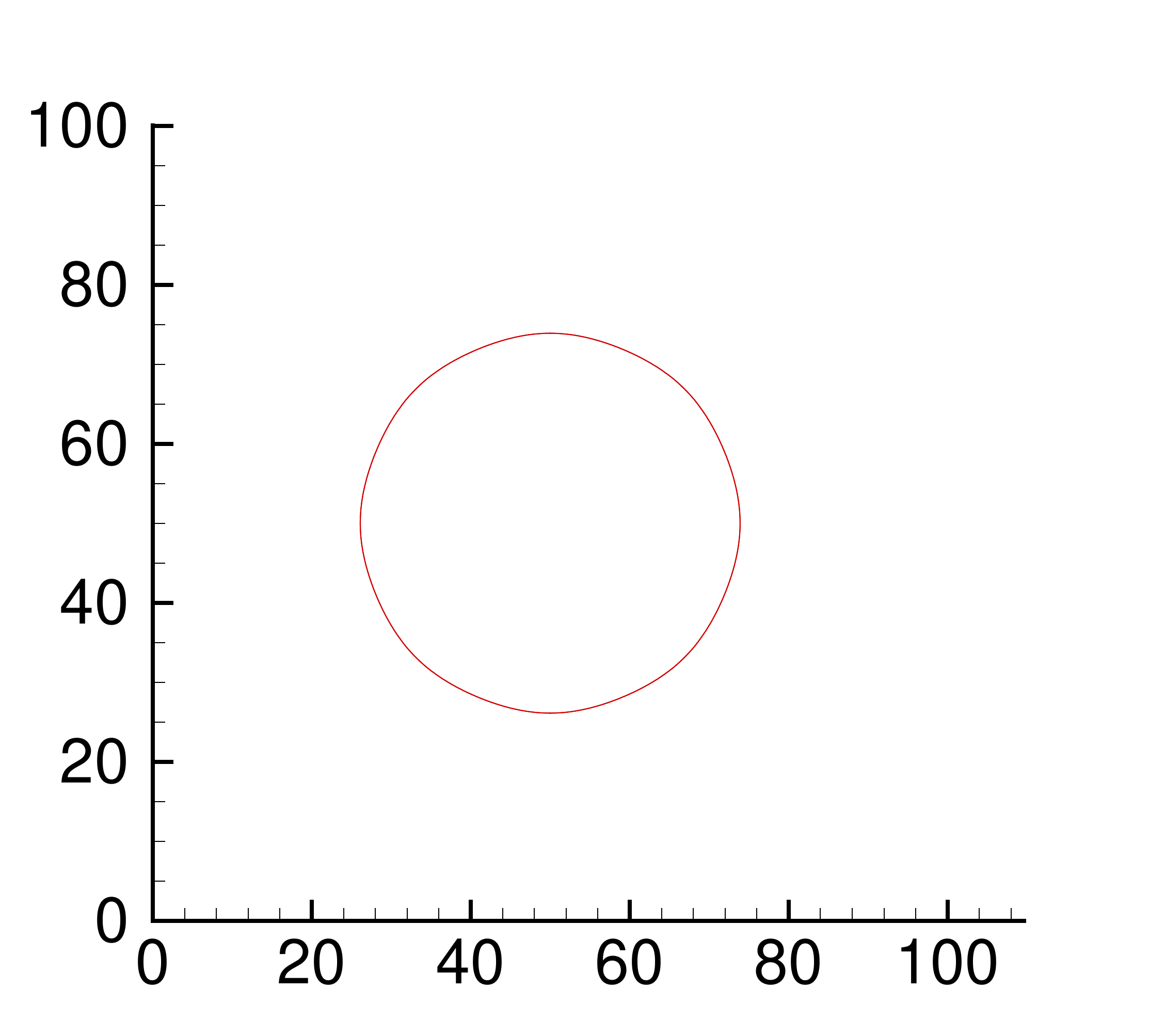} 
\includegraphics[width=.33\textwidth]{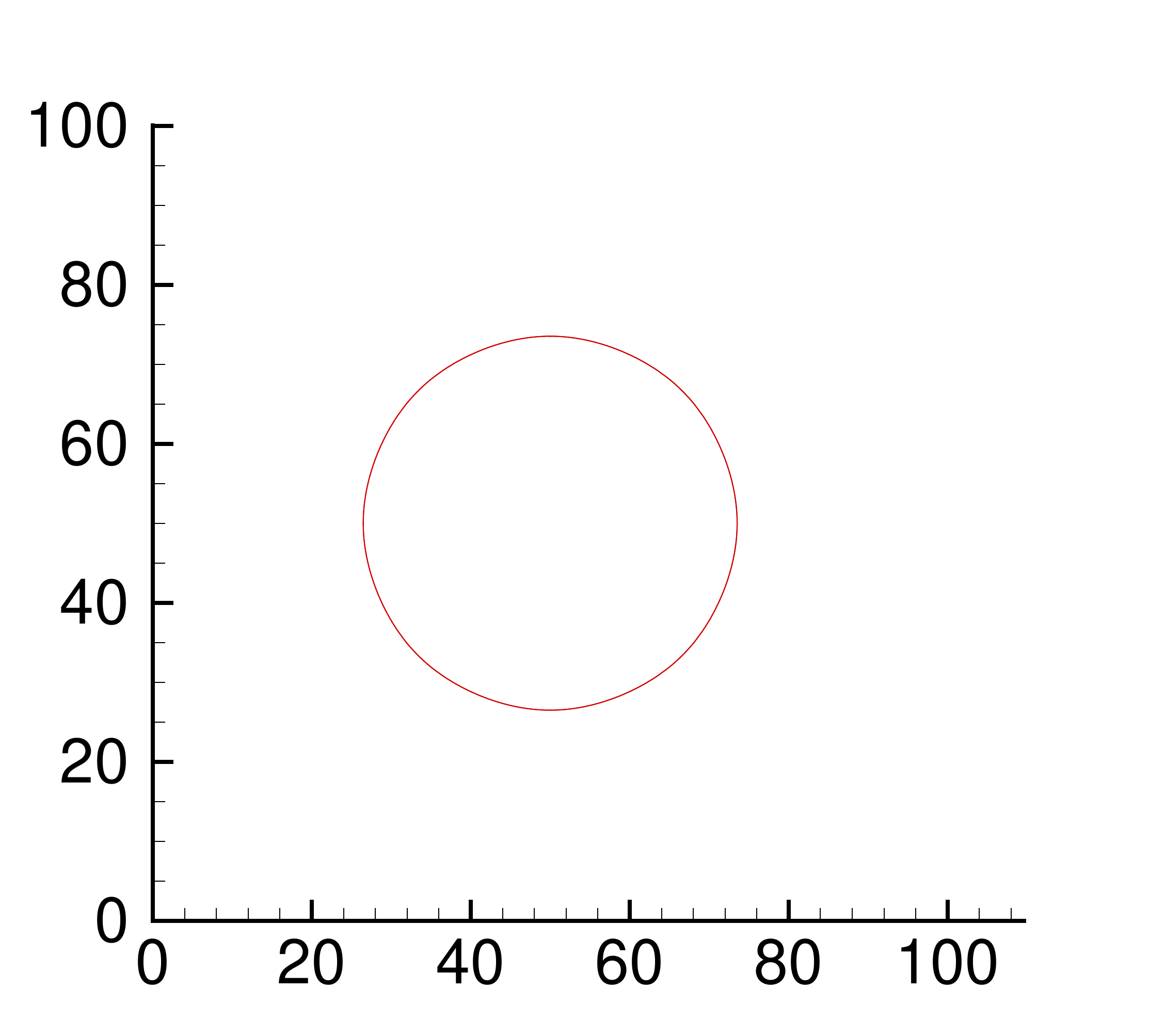} 
\caption{Evolution of the color function of a star-shape under curvature-driven flow. The images are snapshots of the interface in time from left to right, top to bottom where the star collapses unto a circle for  t = 0.0, 7.5, 15.0, 22.5, 30.0, 37.5, 45.0, 52.5, 60.0.}
\label{pointed_star}
\end{figure}

\begin{figure}[H]
\centering
\begin{subfigure}[c]{0.45\textwidth}
\includegraphics[width=\textwidth]{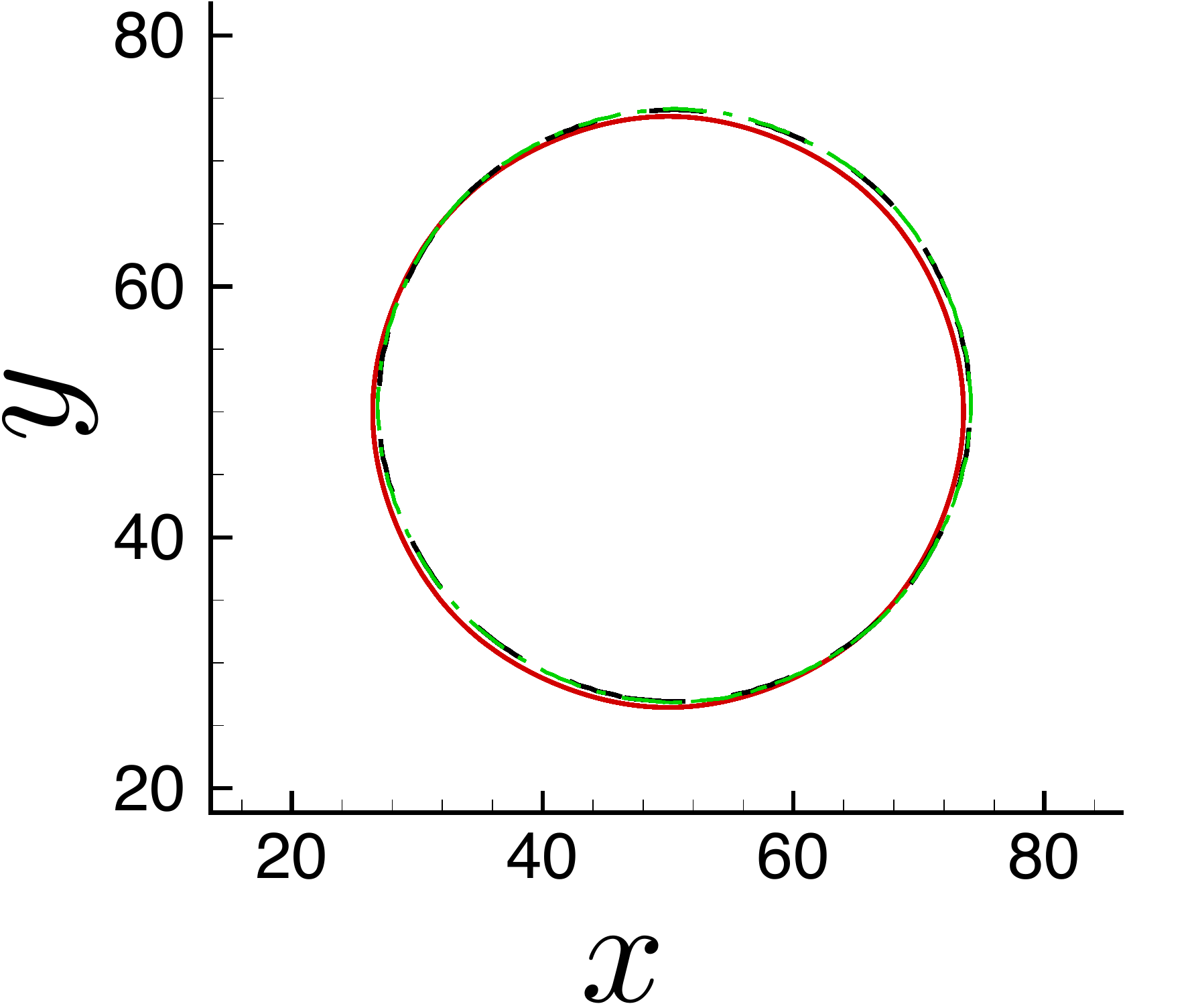}
\end{subfigure}
\qquad
\begin{subfigure}[c]{0.4\textwidth}
\includegraphics[width=\textwidth]{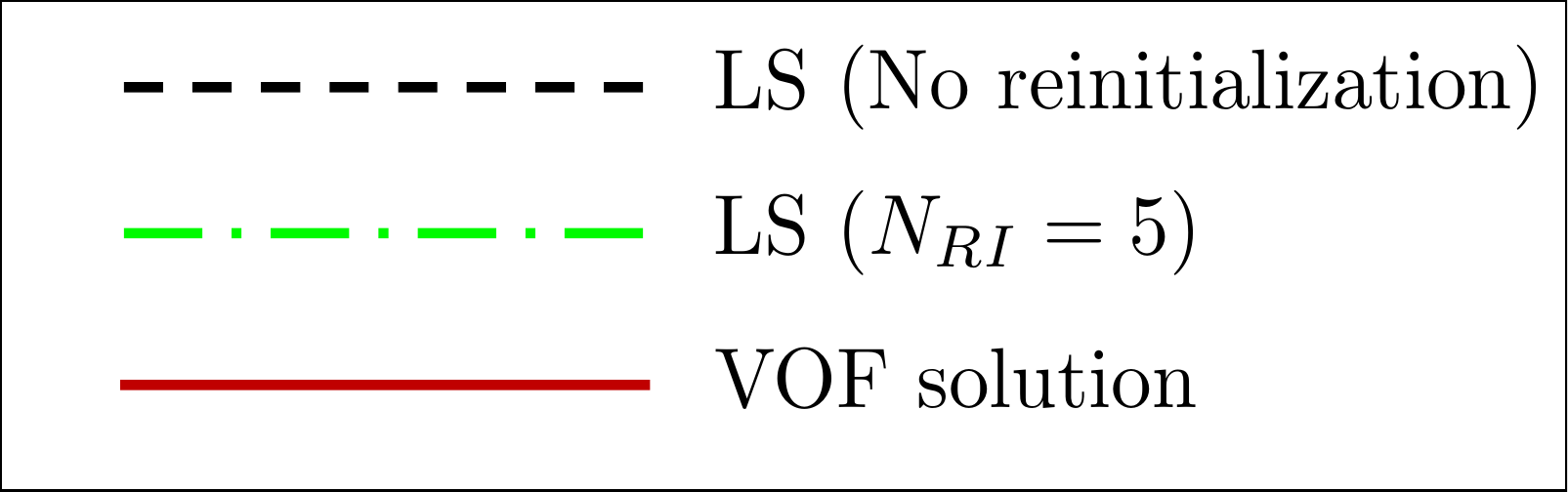}
\end{subfigure}
\caption{Comparison of final contour shape between VVOF and LS with and without reinitialization for the case of the pointed star.}
    \label{petal_comp}
\end{figure}

\noindent Consider a pointed star given by the following parameterized curve
\begin{equation}
\gamma(s) = \big{[}\mathcal{A} + \mathcal{B} \ \text{cos}(\mathcal{K} \cdot2\pi s)\big{]}\big{[}\text{cos}(2\pi s),\text{sin}(2\pi s)\big{]} \ \ \ \text{for} \ \ \ s\in[0,1],
\end{equation}
such that the zeroth level set, defined from polar coordinates,  is given by
\begin{equation}
    \phi(\mathbf{x},0) = ||\mathbf{x}-\mathbf{x_c}||-\bigg{\{}\mathcal{A}+\mathcal{B} \ \text{cos}\bigg{[}\mathcal{K}\cdot \text{arctan}\bigg{(}\frac{y-y_c}{x-x_c}\bigg{)}\bigg{]}\bigg{\}}
\end{equation}

\noindent where $\mathcal{A}$ is the shift of the star petals from the center, $\mathcal{B}$ is a scaling factor for the star size, $\mathcal{K}$ is the number of petals, and  $\mathbf{x_c}$ is the offset parameter. The petals of the star are shifted from its center by 25 units of length ($\mathcal{A}=25$) which effectively increases the relative size of the petal shape in comparison to the main body, the scaling factor $\mathcal{B}$ is chosen to be 10, and the number of petals is 8 ($\mathcal{K}=8$). The computational domain is a square with a side length of 100 units and a grid resolution of $200 \times 200$. The time-step is $\Delta t=0.05$ and the star-shaped level curve is centered at $(50,50)$. Fig. \ref{pointed_star} shows the progression of the color function front under curvature-driven flow for $t=0.0, 7.5, 15, 22.5, 30.0, 37.5, 45.0, 52.5, 60.0.$. The interface collapses to a circle under the curvature-driven motion due to the collapse of the peaks ($+\kappa$) inwards and the propagation of the troughs outwards ($-\kappa$). Similar results have been reported in \cite{Smereka2003,Alame2020} for level set methods. We note that for geometries involving sharp corners (e.g. pointed star), the smoothness of the interface at those locations is paramount to the physical evolution of the full geometry. This is particularly important in the context of VOF since errors in normal estimation will grow with every directional sweep leading eventually to the rupture of the interface if the corners are not well-resolved.
\begin{figure}[H]
\centering
\includegraphics[width=.33\textwidth]{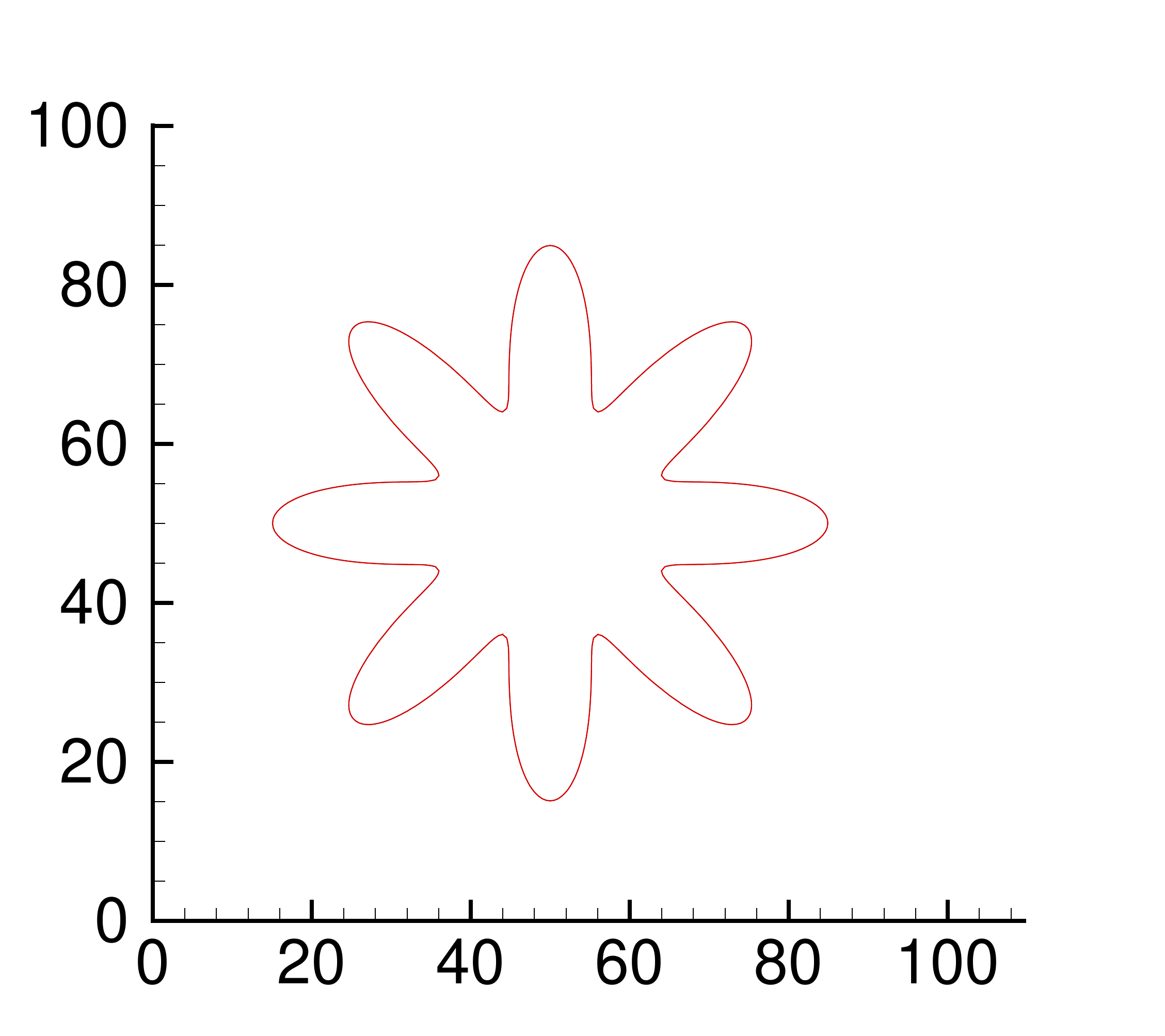}
\includegraphics[width=.33\textwidth]{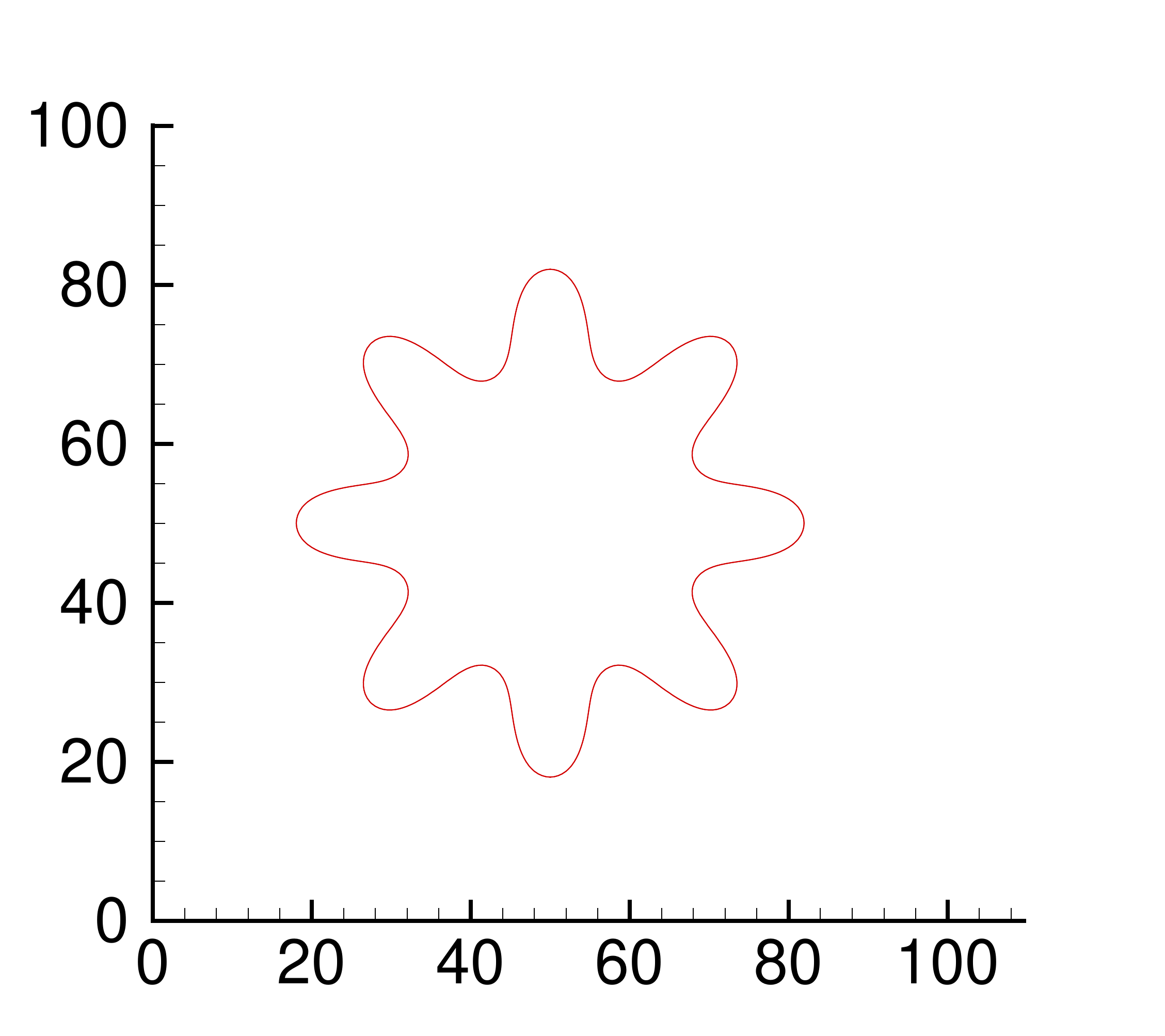} 
\includegraphics[width=.33\textwidth]{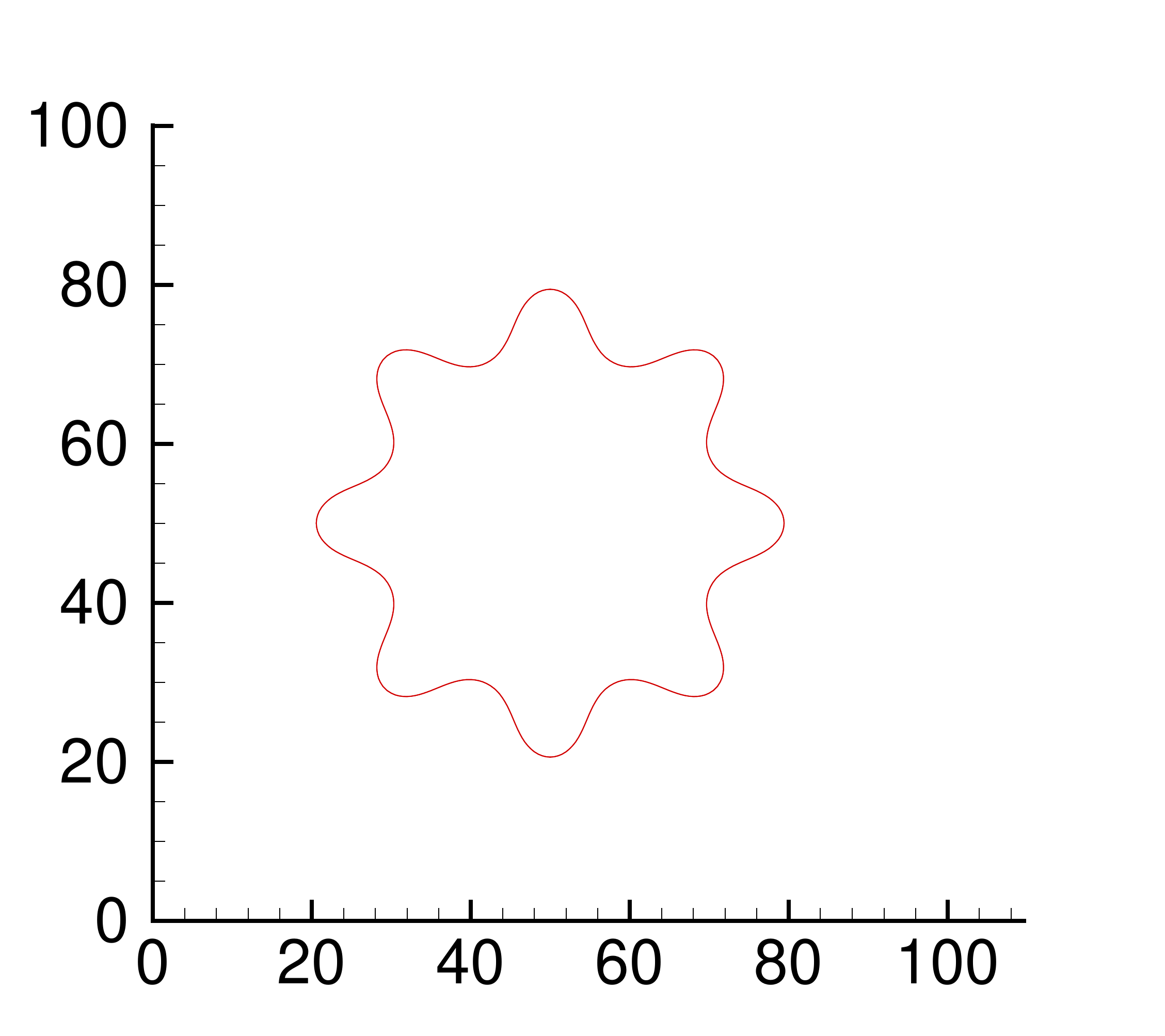}
\end{figure}
\vspace*{1cm}
\begin{figure}[H]
\centering
\includegraphics[width=.33\textwidth]{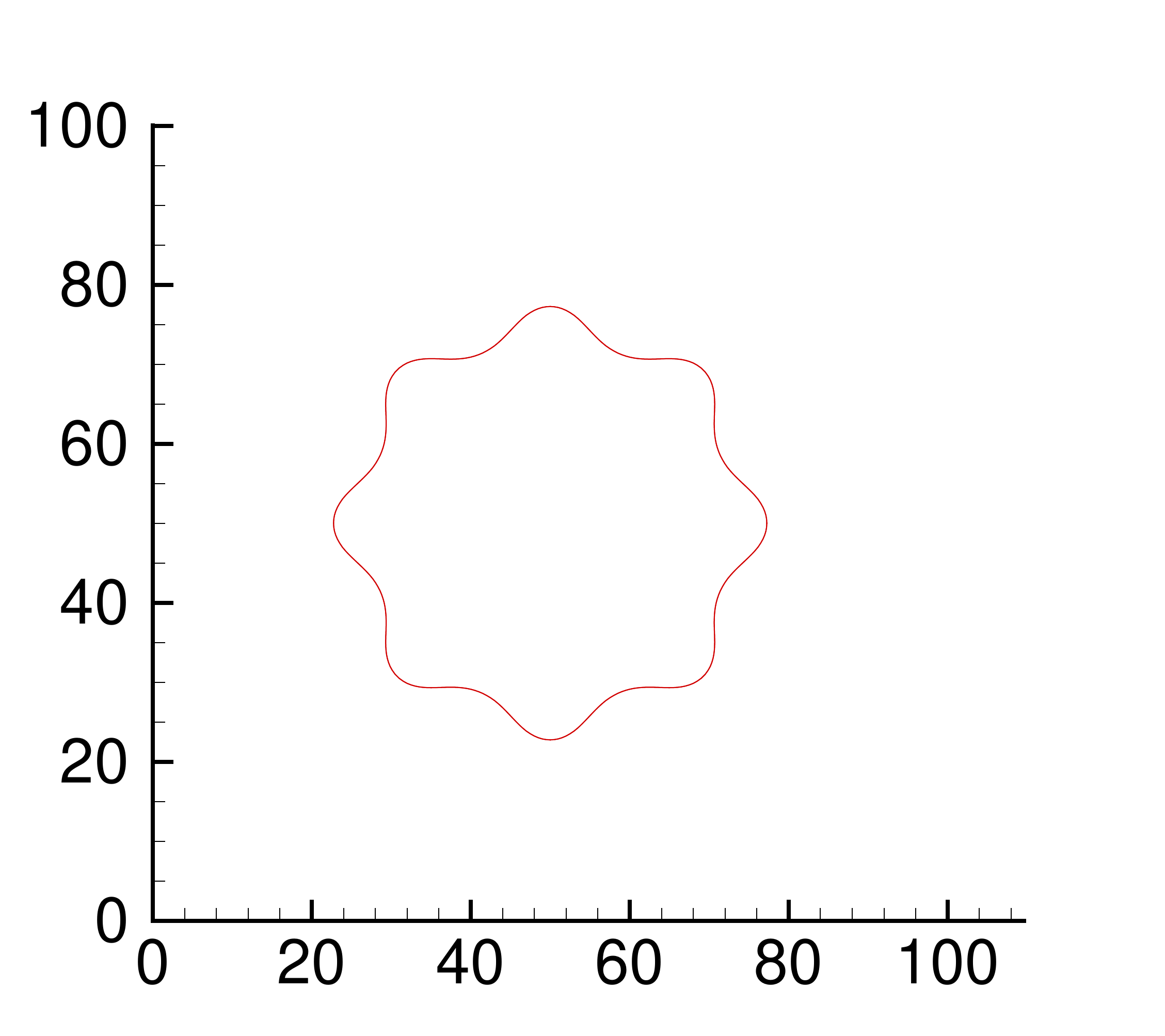} 
\includegraphics[width=.33\textwidth]{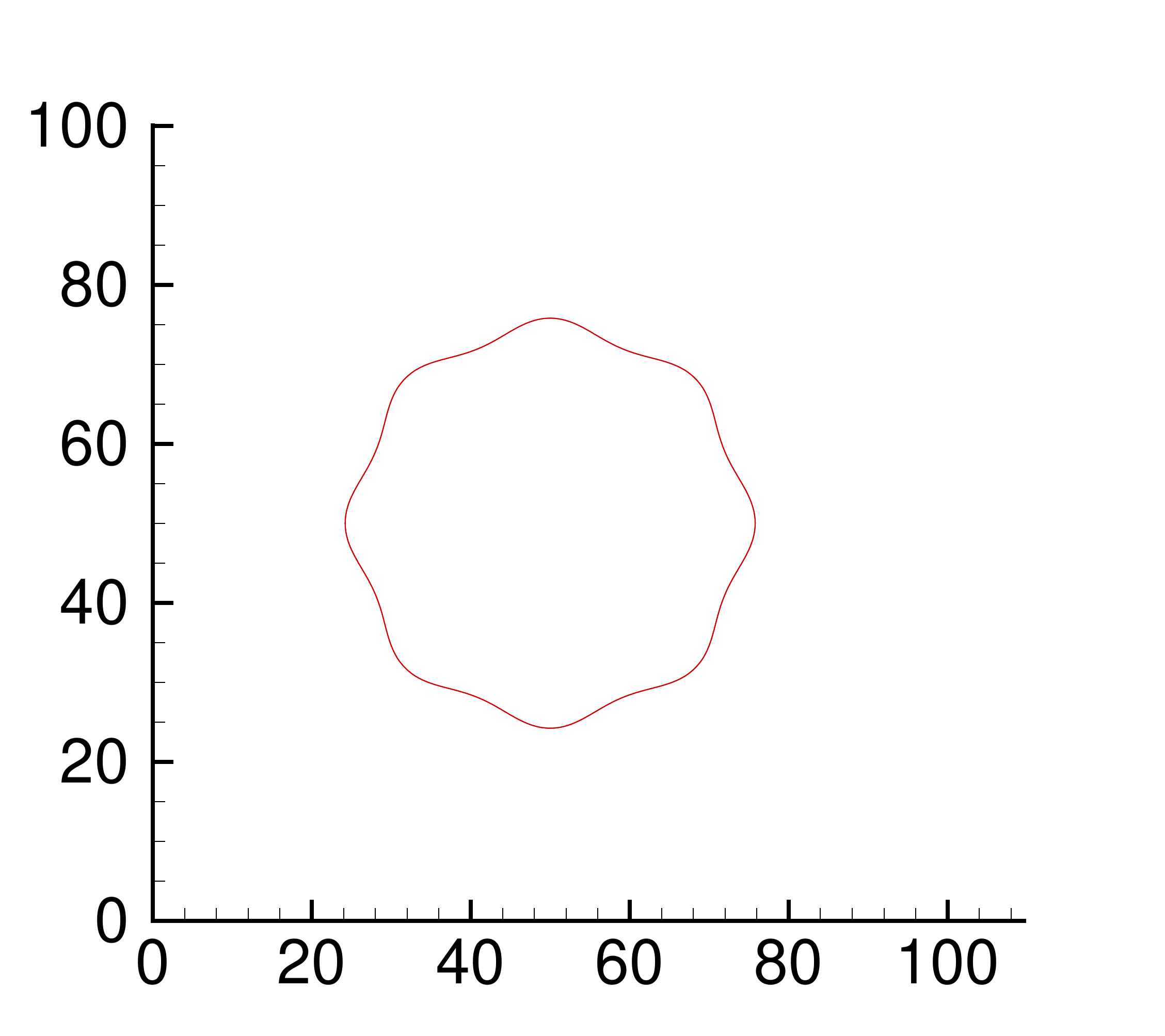}
\includegraphics[width=.33\textwidth]{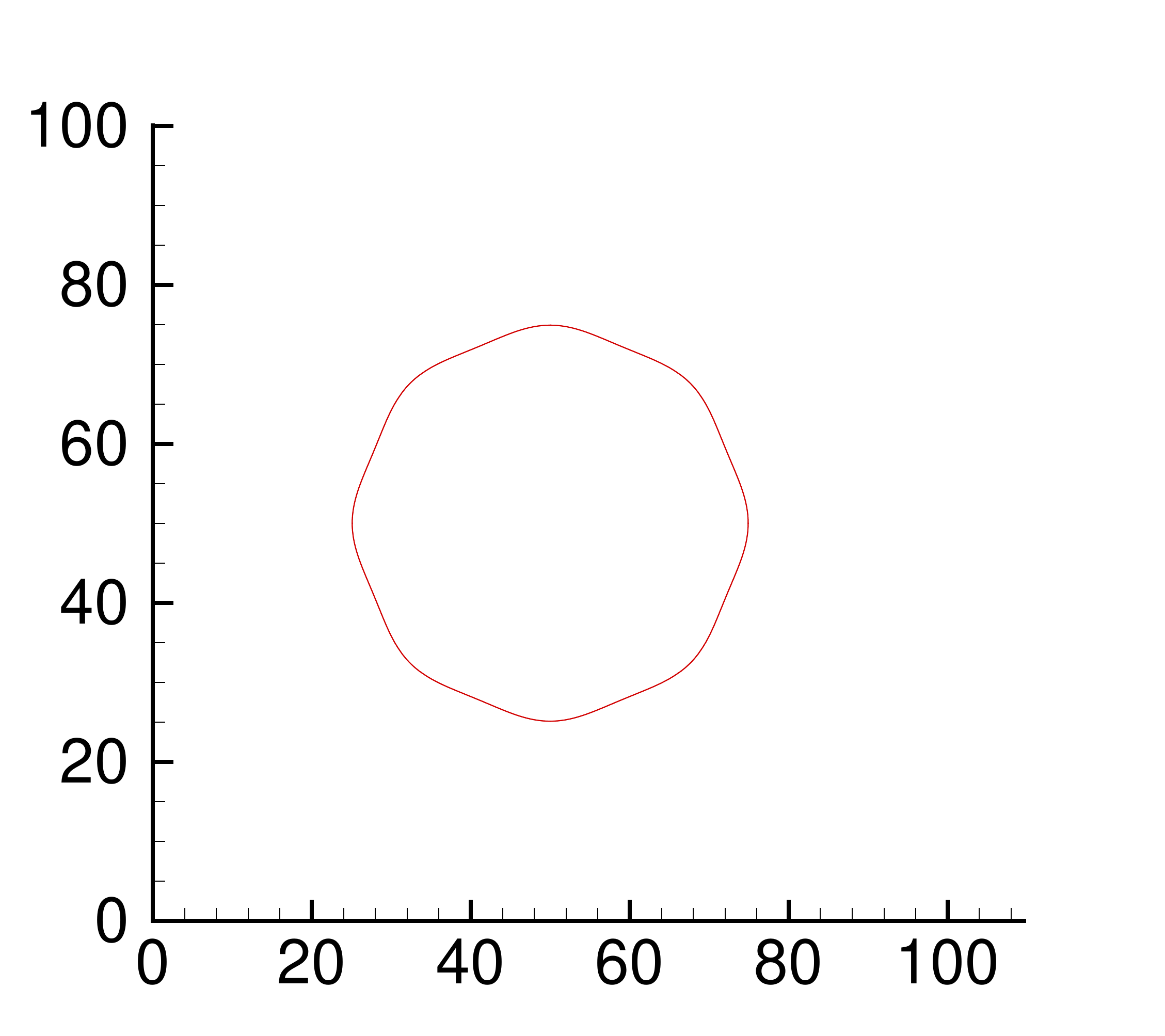}
\end{figure}
The results in Fig. \ref{pointed_star} were compared with standard LS with reinitialization for the same problem parameters. Fig. \ref{petal_comp} shows a comparison between the interface contour of VVOF and LS with and without reinitialization. Reinitialization has no observable effect on the final shape of the level-set and the final result is consistent for both interface capturing methods. \par Since the pointed star should evolve into a perfect circle under curvature-driven motion, the roundness (sphericity in 3D) of the final contour was calculated for both methods. Roundness is the measure of how closely the shape of an object approaches that of a mathematically perfect circle, and among the most common roundness error definitions is circularity which is often used in digital image processing. Circularity is the ratio of the perimeter squared of shape to its area such that
\begin{equation}
\mathcal{S} = \frac{Perimeter^2}{4\pi \times Area} \quad .
\end{equation}
The ratio between the $\mathcal{S}_{LS}$ and $\mathcal{S}_{VVOF}$ was found to be less than 1 which indicates higher roundness with VOF and a contour shape that more closely approximates a perfect circle. Note that the variable $N_{RI}$ is Fig. \ref{petal_comp} refers to the number of reinitialization steps.

\subsubsection{Wound spiral}
\begin{figure}[H]
\centering
\includegraphics[width=.33\textwidth]{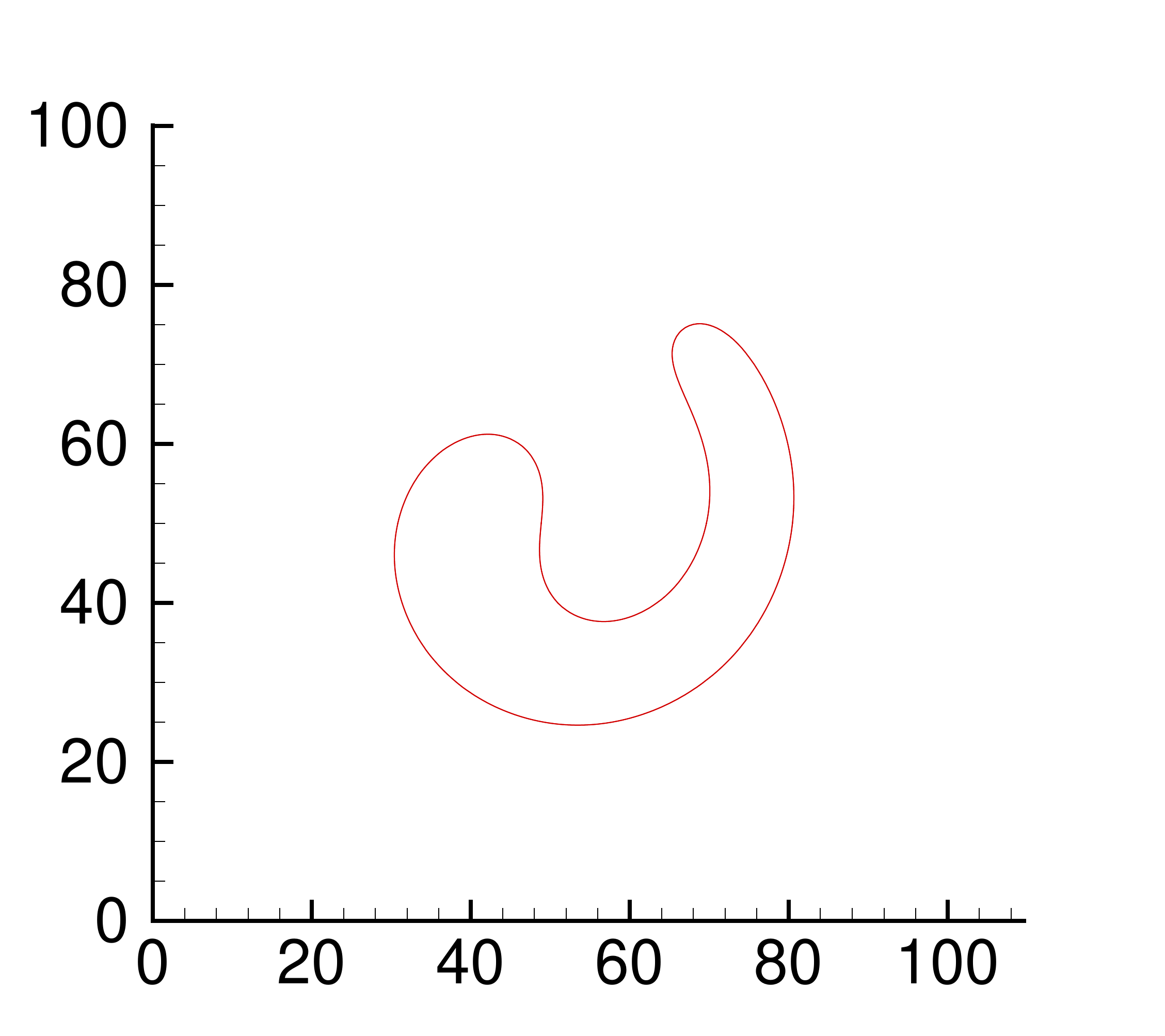}
\includegraphics[width=.33\textwidth]{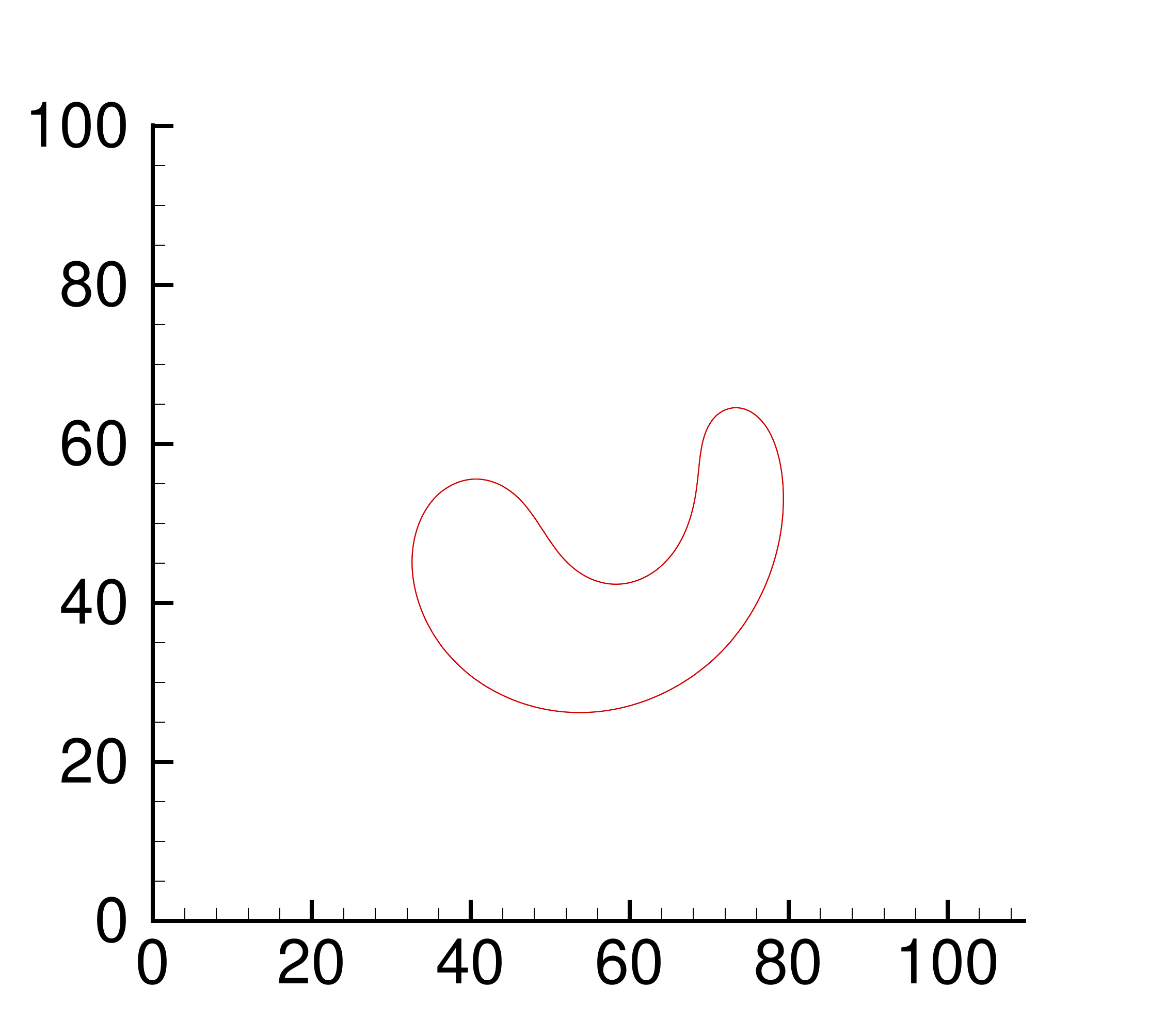} 
\includegraphics[width=.33\textwidth]{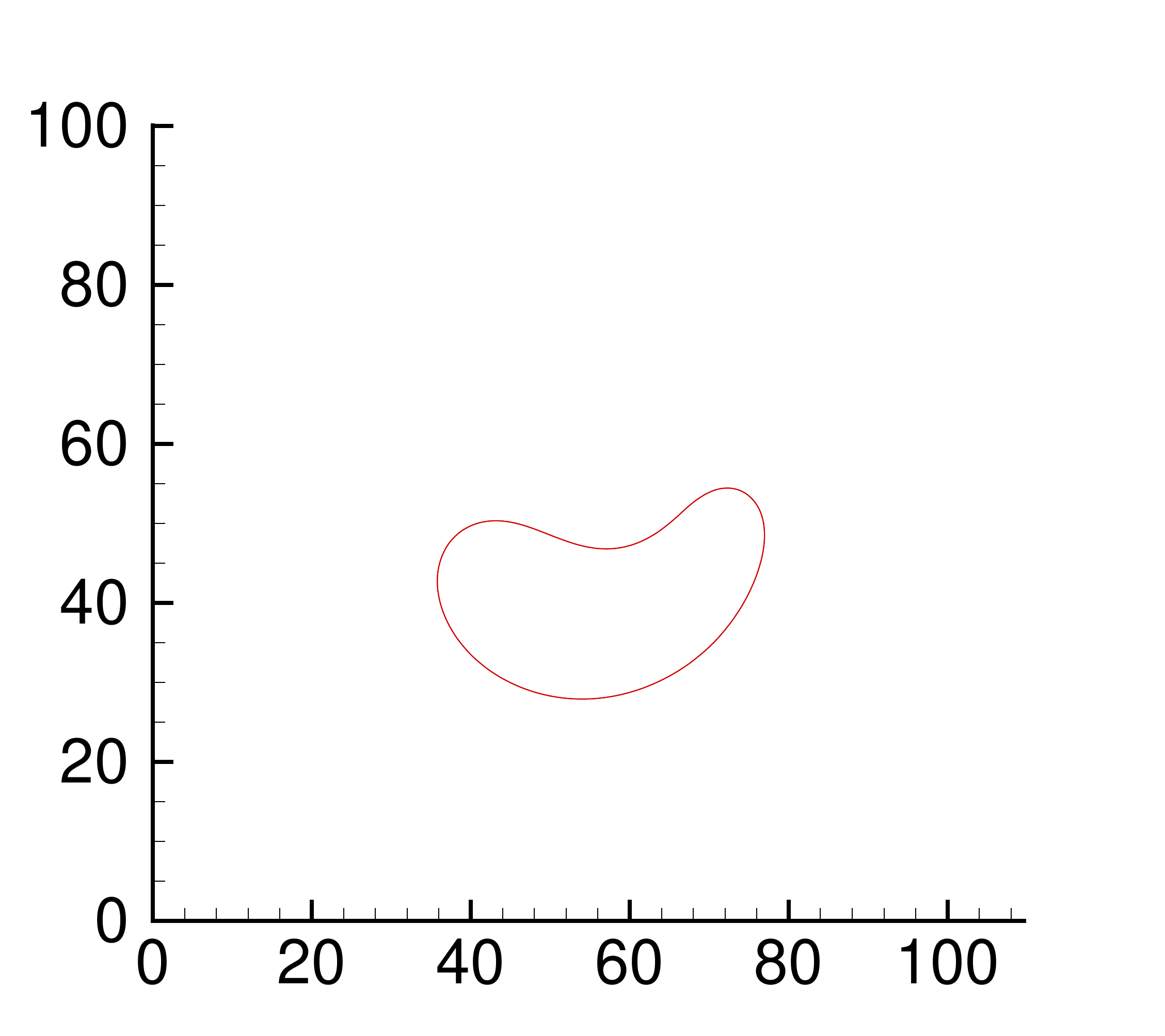}
\caption{Evolution of the color function of a spiral under curvature-driven flow. The images are snapshots of the interface in time from left to right, top to bottom where the spiral unwinds for t = 0.0, 37.5, 75.0, 112.5, 150.0, 187.5, 225.0, 262.5, and 300.0.}
\label{wound_spiral}
\end{figure}

\noindent Consider the curvature-driven motion of a wound spiral given by the following parameterization
\begin{equation}
    \theta = 2\pi D \sqrt{s} \ \ \ \text{for} \ \ \ s\in[0,1],
\end{equation}
such that $s=(k+a)/(n_p+a)$, where $k$ is an integer that loops over the total number of points $n_p$, and the constant $a$ determines the shape of the spiral head in the center. The location of the points describing the spiral is defined by
\begin{equation}
    x_p = x_c + s_f \bigg{[}\frac{D\sqrt{s} \ \text{cos}(\theta)}{1+D}\bigg{]} \ \ \ \ \ \text{and} \ \ \ \ \ y_p = y_c + s_f \bigg{[}\frac{D\sqrt{s} \ \text{sin}(\theta)}{1+D}\bigg{]}
\end{equation}
where $\mathbf{x_c}$ is the location of the domain center, $s_f$ is the scaling factor for the spiral size, and $D$ is the number of spirals. We define a distance function $d(\mathbf{x})=||\mathbf{x}-\mathbf{x_p}||$ such that the zeroth level set has the following form
\begin{equation}
    \phi(\mathbf{x},0) = d(\mathbf{x}) - w
\end{equation}
where $w$ is the width of the spirals. The computational domain is a square of side length equal to 100 units and a grid size of $100 \times 100$. The time-step is $\Delta t = 0.05$, the level curve center $\mathbf{x_c}=(50,50)$, $n_p=400$, $D=2.5$, $a=3$, and $s_f=50$.  The simulation is stopped after 6000 iterations ($t_{final}=300.0$) and snapshots of the interface evolution are shown in Fig \ref{wound_spiral} for $t=0.0, 37.5, 75.0, 112.5, 150.0, 187.5,225.0,262.5,$ and $ 300.0$. Problem parameters describing spiral geometry were chosen to match \cite{Alame2020} for comparison. The spiral unwinds under curvature-driven motion, and its size decreases in the process until it eventually acquires a bean shape that later collapses rapidly. The interface maintains its sharpness even at the last stages of collapse which is consistent with the bubble collapse problem discussed earlier.

\begin{figure}[H]
\centering
\includegraphics[width=.33\textwidth]{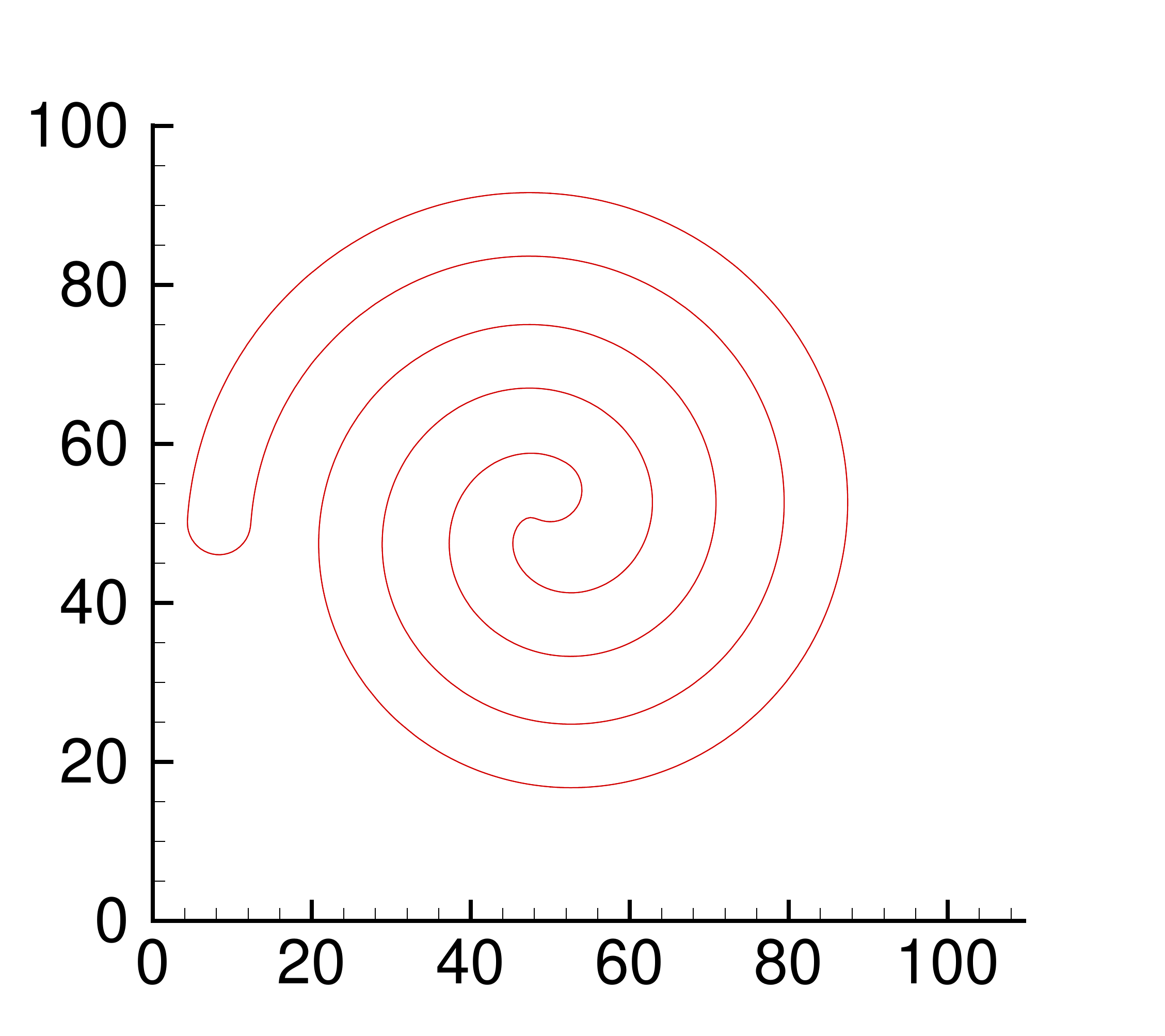}
\includegraphics[width=.33\textwidth]{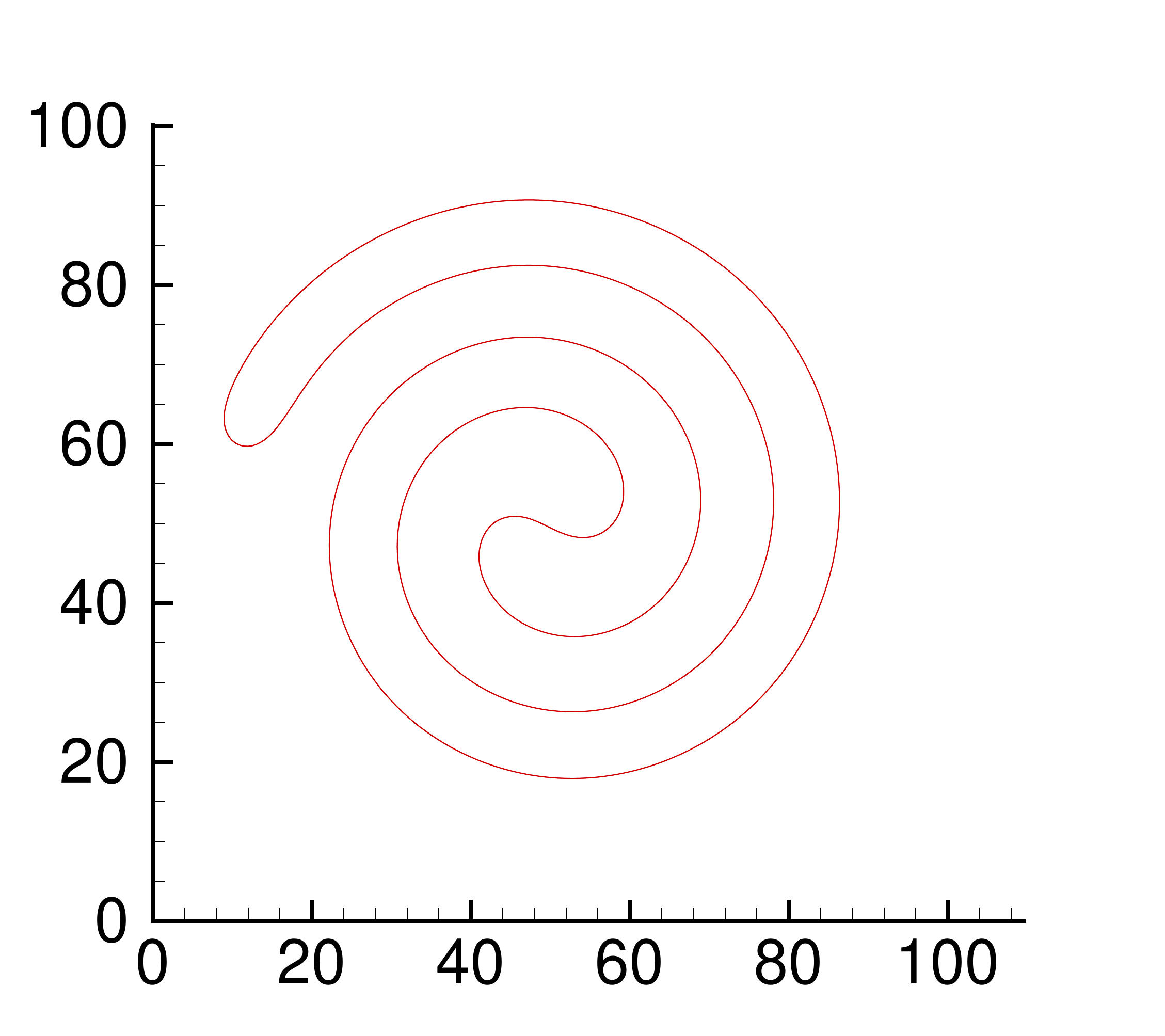} 
\includegraphics[width=.33\textwidth]{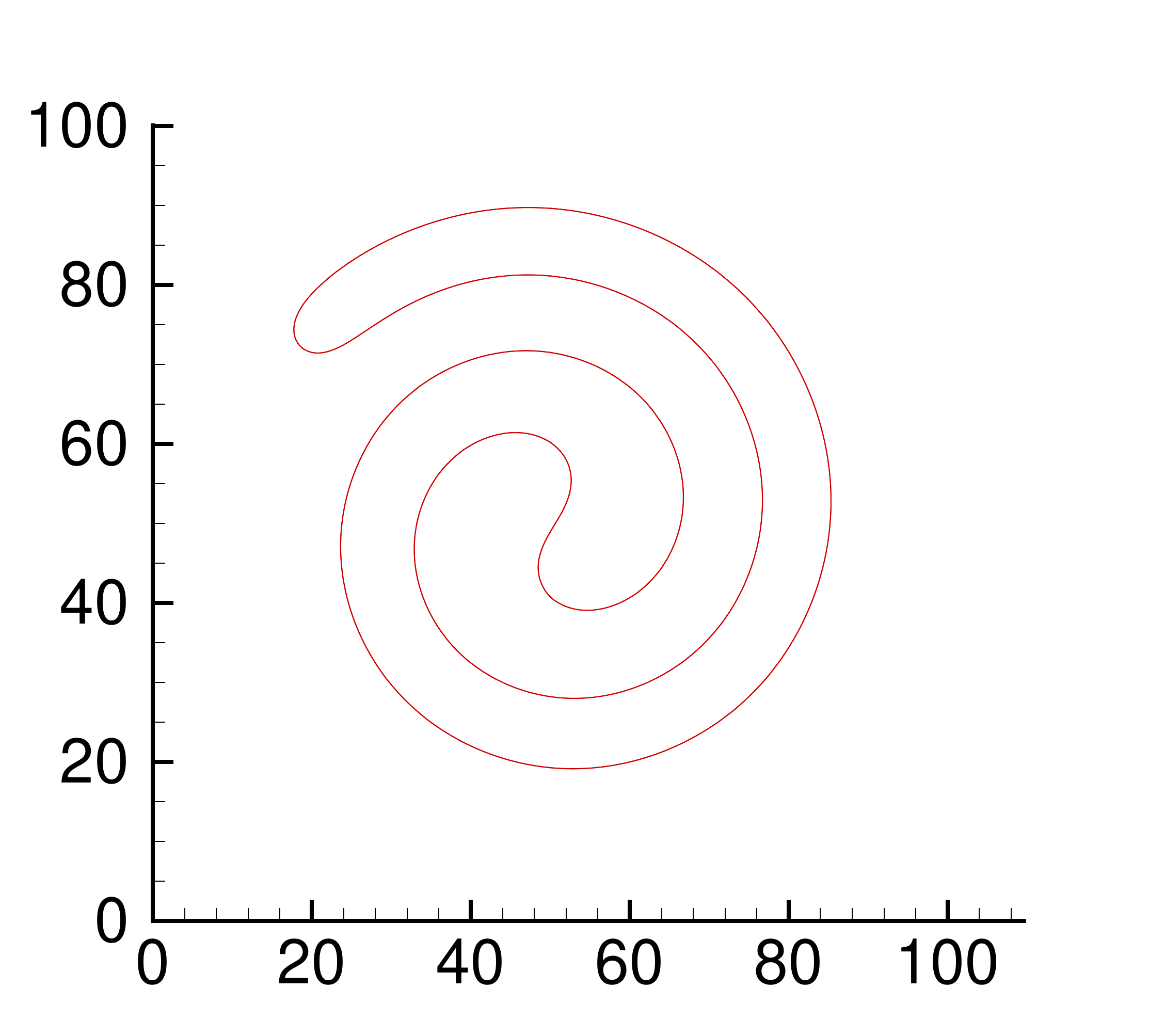}
\end{figure}
\begin{figure}[H]
\centering
\includegraphics[width=.33\textwidth]{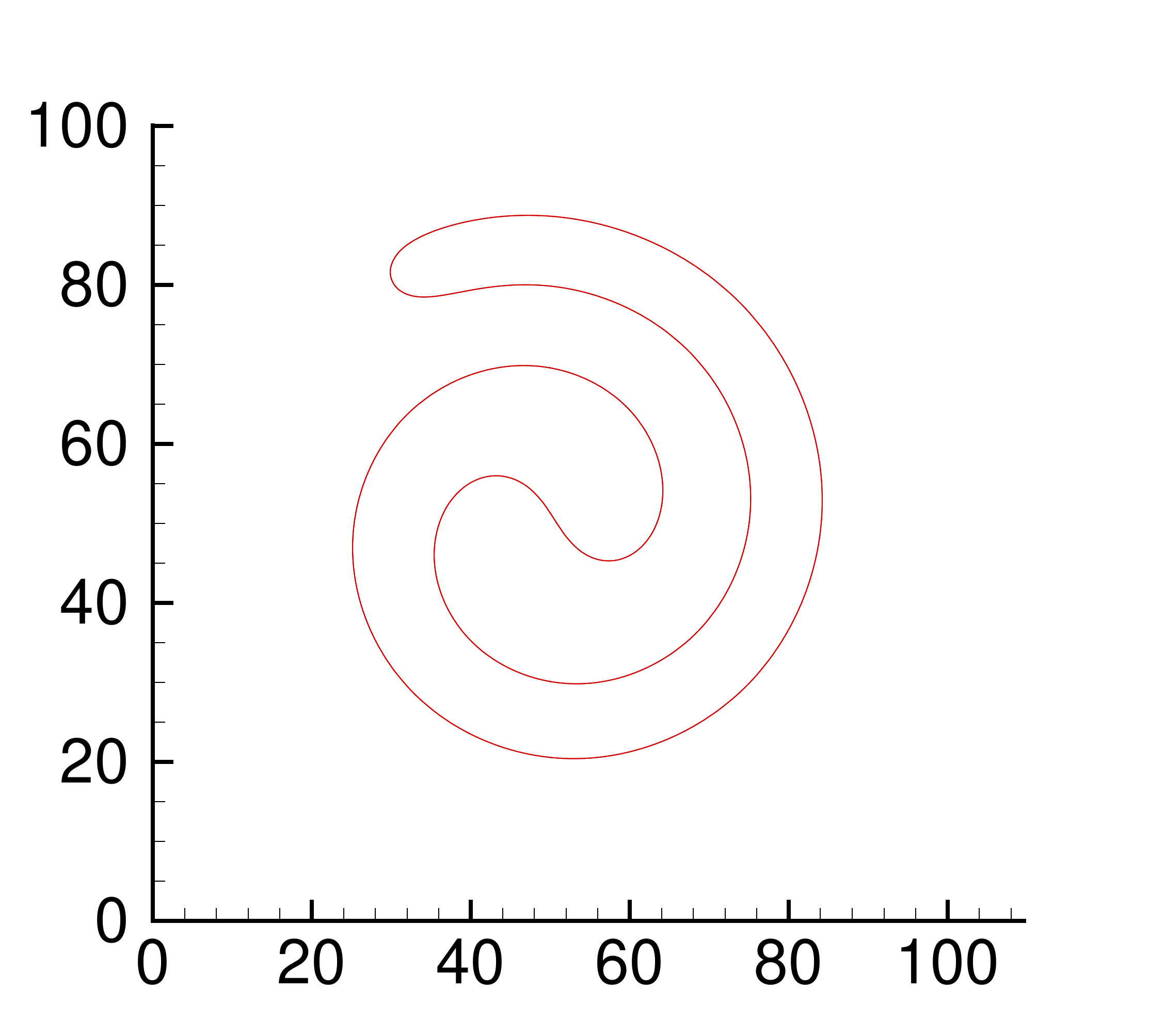}
\includegraphics[width=.33\textwidth]{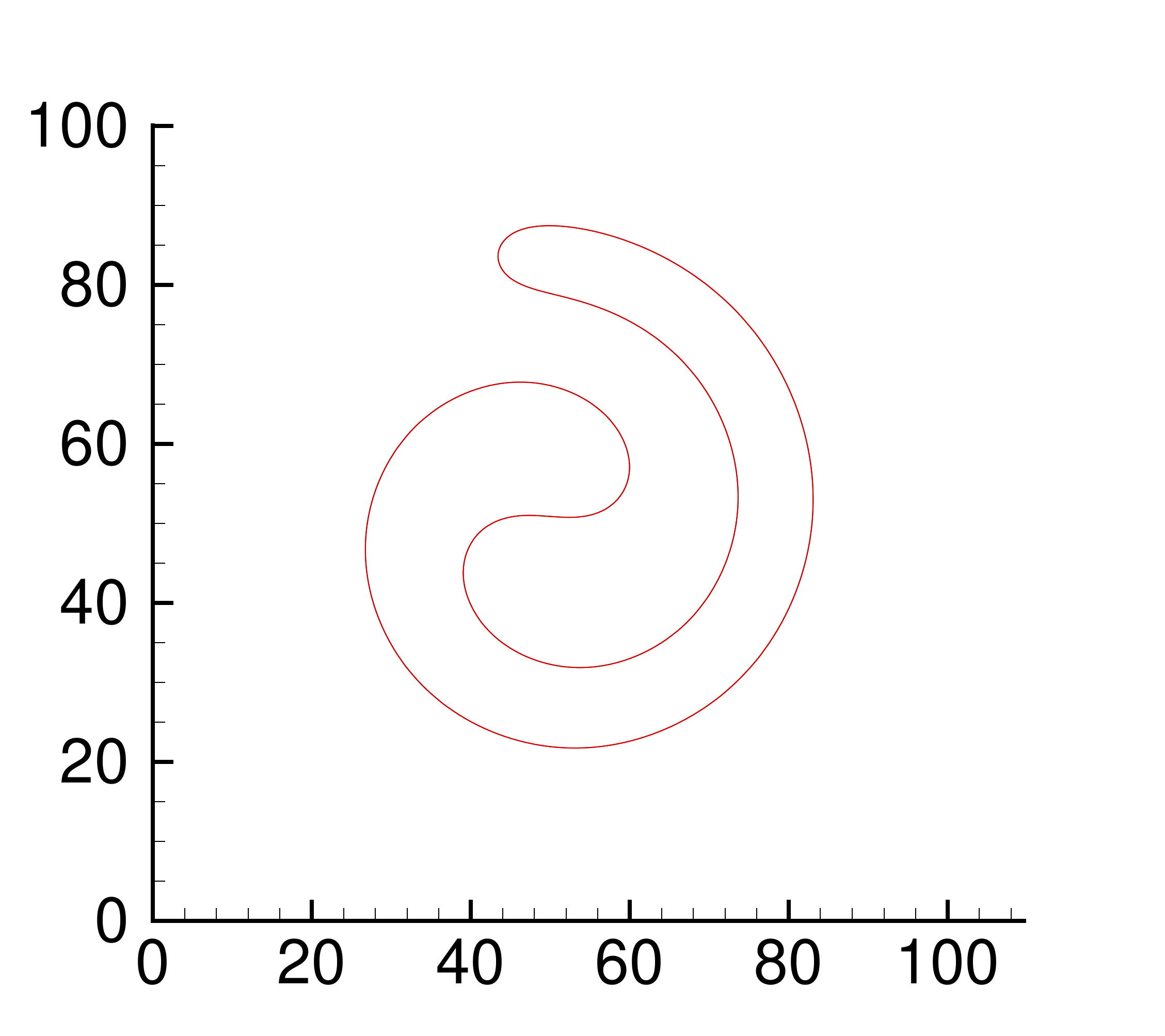} 
\includegraphics[width=.33\textwidth]{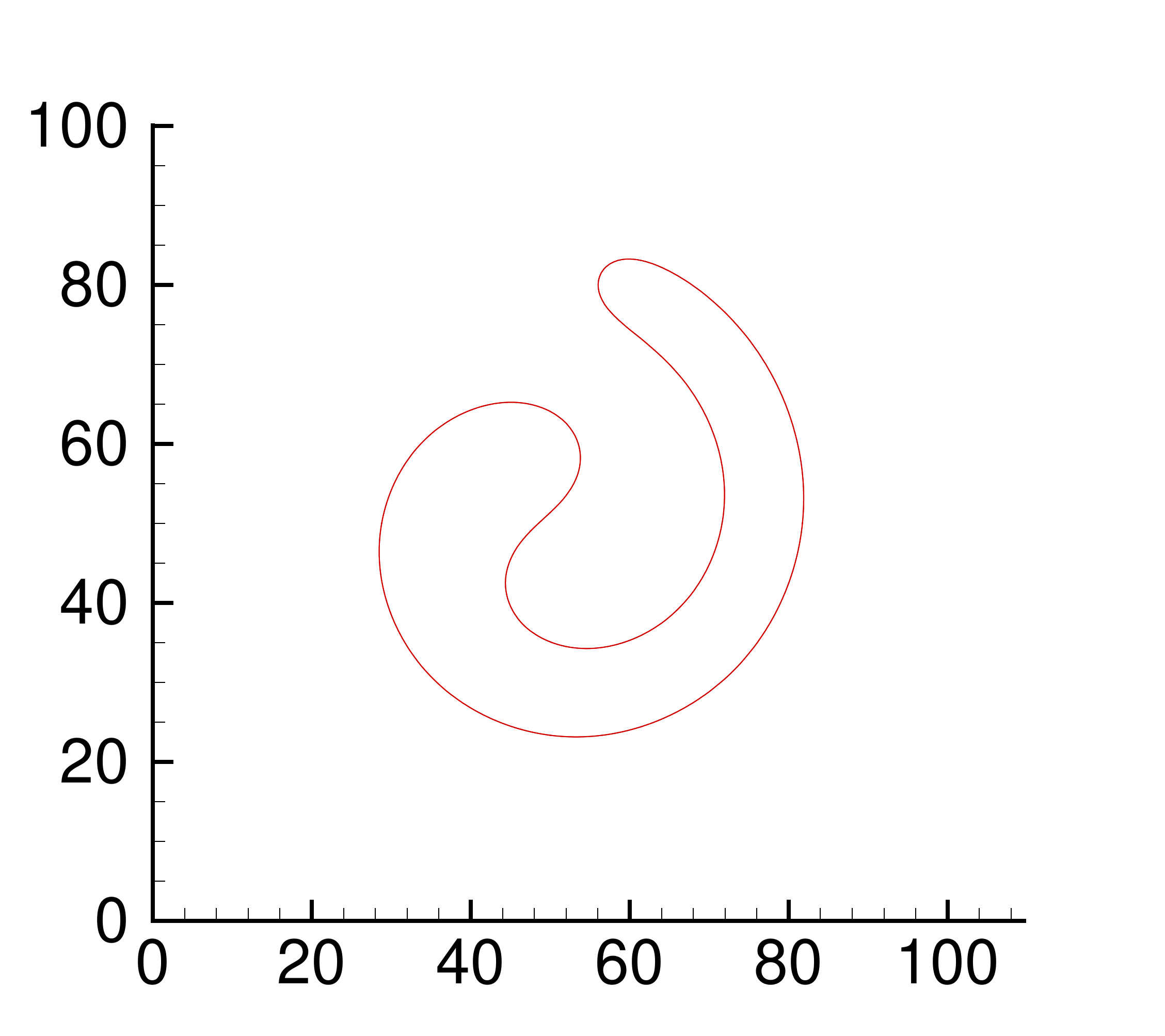}
\end{figure}

\begin{figure}[H]
\centering
\begin{subfigure}[c]{0.45\textwidth}
\includegraphics[width=\textwidth]{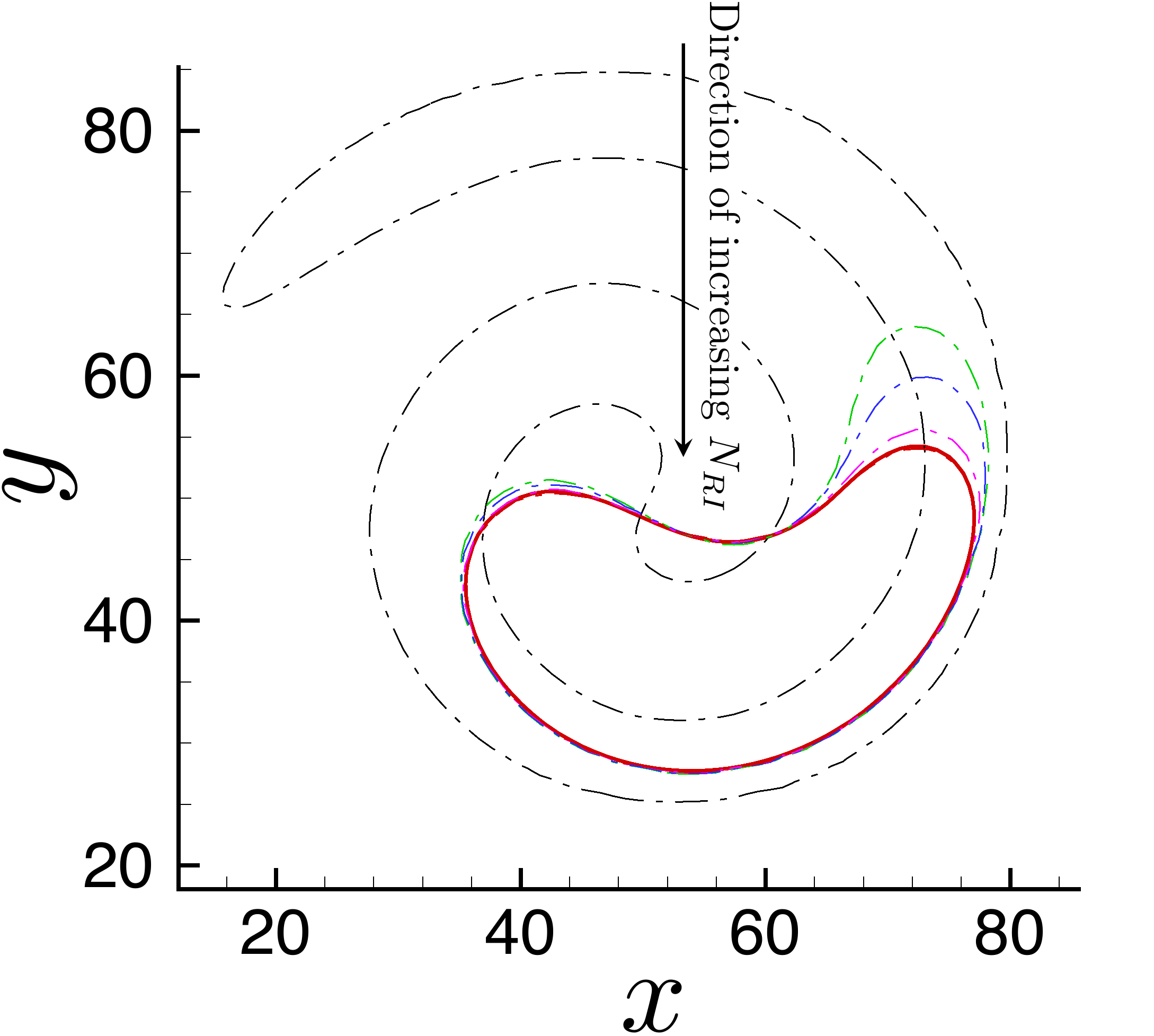}
\end{subfigure}
\qquad
\begin{subfigure}[c]{0.4\textwidth}
\includegraphics[width=\textwidth]{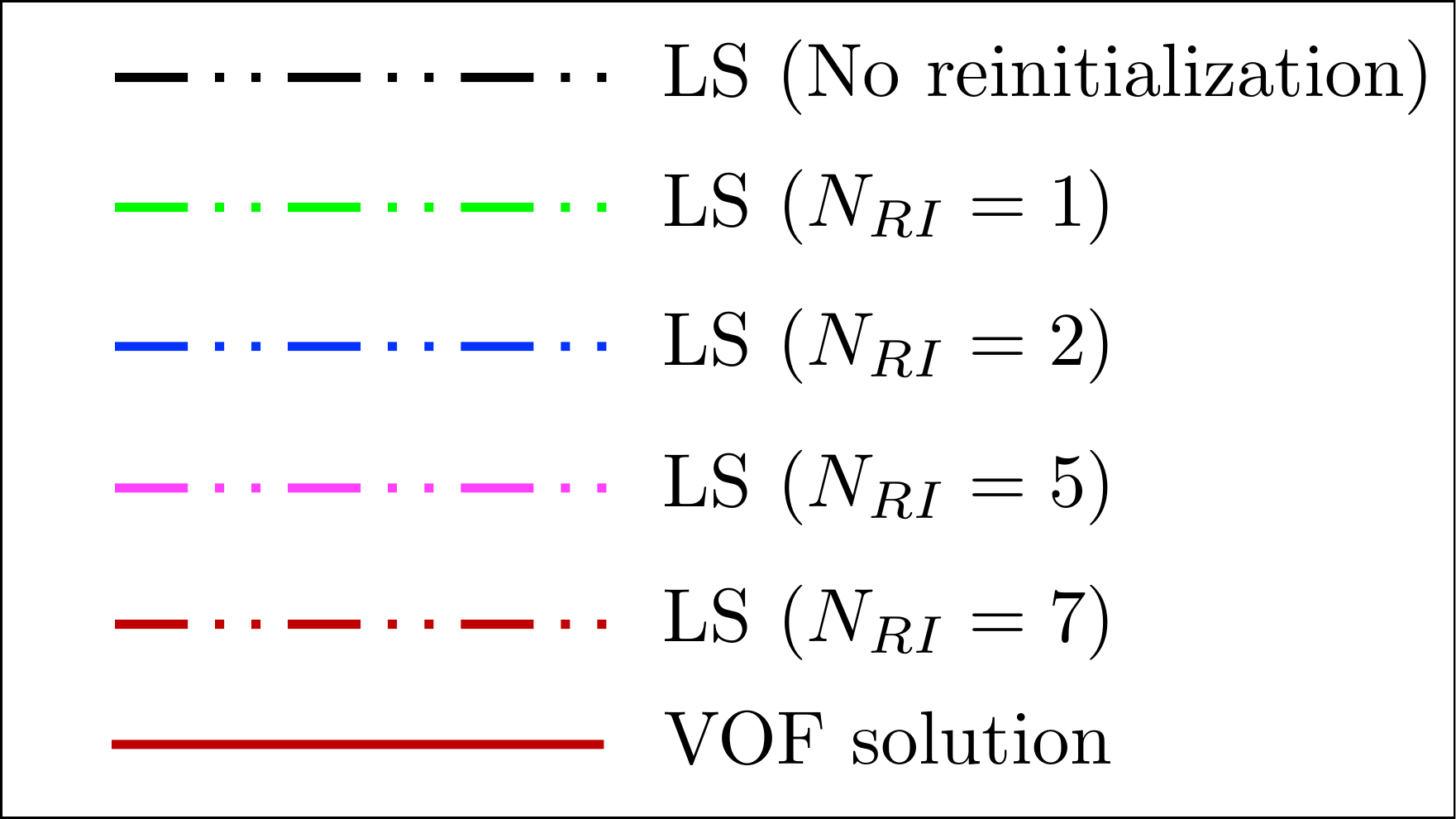}
\end{subfigure}
\caption{Comparison of final contour shape between VVOF and LS with and without reinitialization for the case of the wound spiral.}
\label{spiral_comp}
\end{figure}

Similar to the pointed star case, the shape of the final contour of the wound spiral is compared between VVOF and LS as shown in Fig. \ref{spiral_comp}. The effect of no reinitialization is immediately observable for the LS method where the contour deviates completely from the VVOF solution. As the number of reinitializations increases, the contour obtained from LS converges to the VVOF solution at $N_{RI}=7$. What is worth noting is that the result obtained with VVOF closely agrees with that obtained with DRLSE in \cite{Alame2020} for the same problem setup. Moreover, there is no clear criterion that determines the number of reinitialization steps needed to converge to the expected solution, this as a result makes VVOF advantageous and more robust.

\subsubsection{Dumbbell}\label{dumbbell_sec}
\begin{figure}[H]
\centering
\includegraphics[width=.33\textwidth]{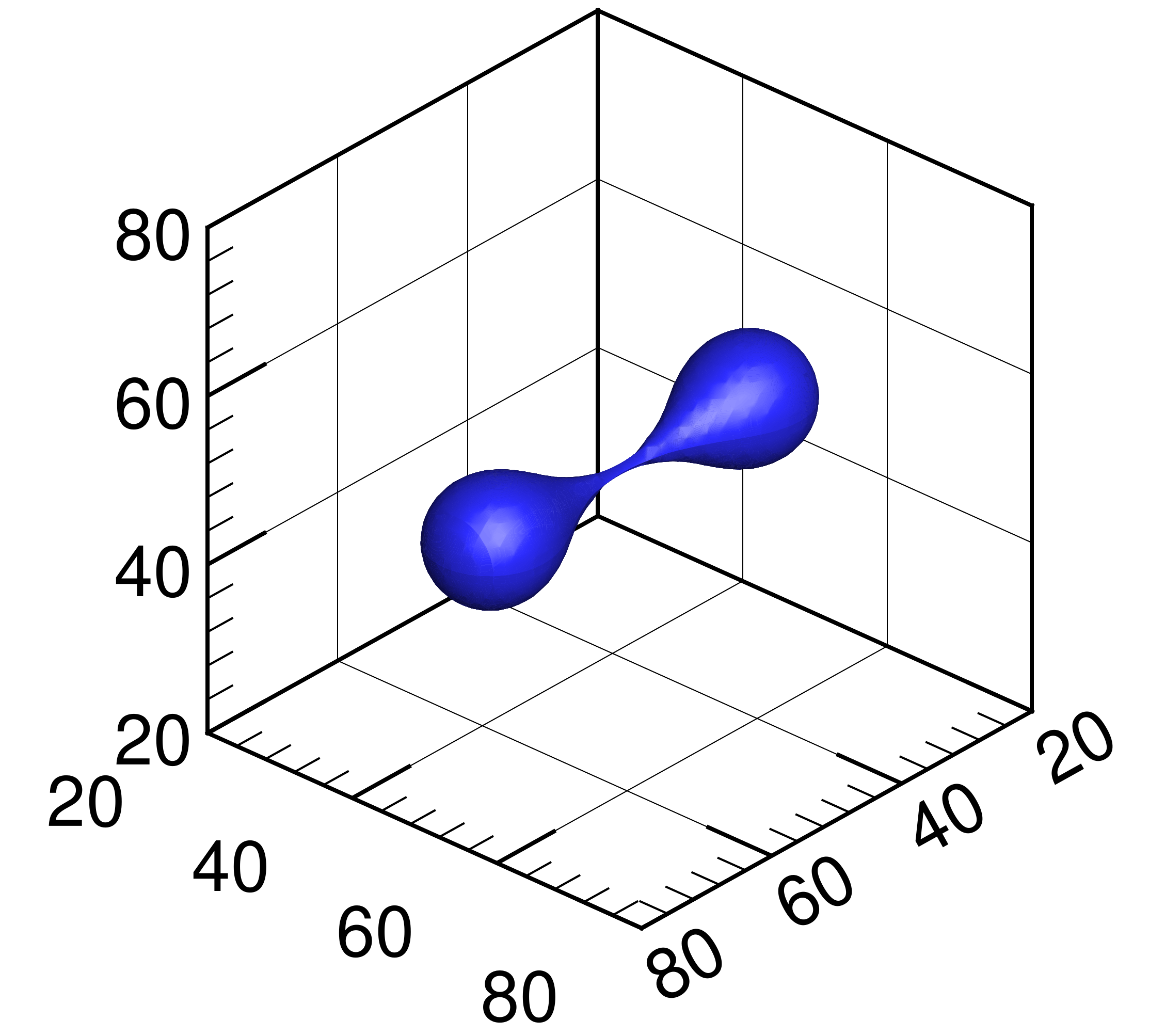}
\includegraphics[width=.33\textwidth]{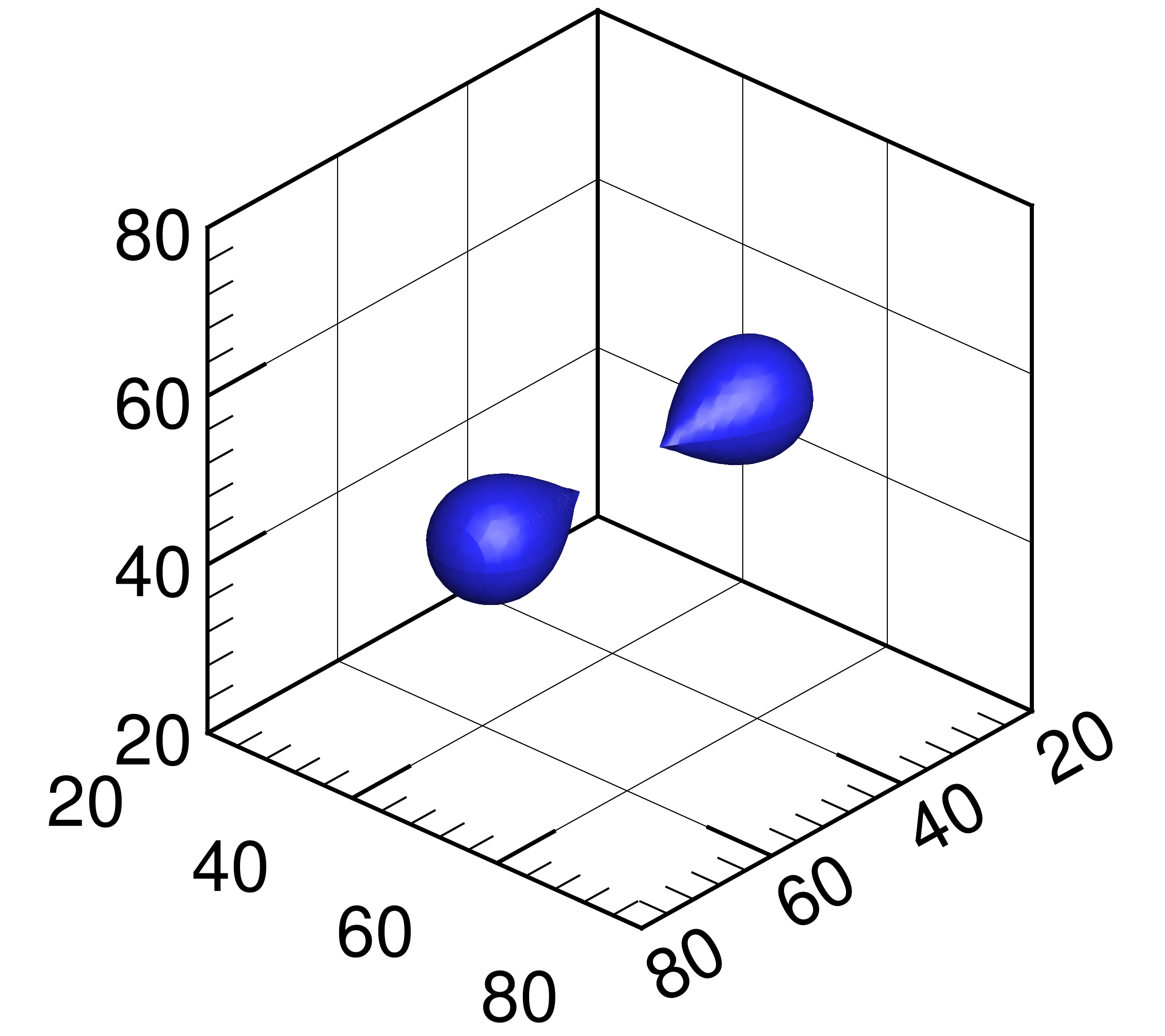} 
\includegraphics[width=.33\textwidth]{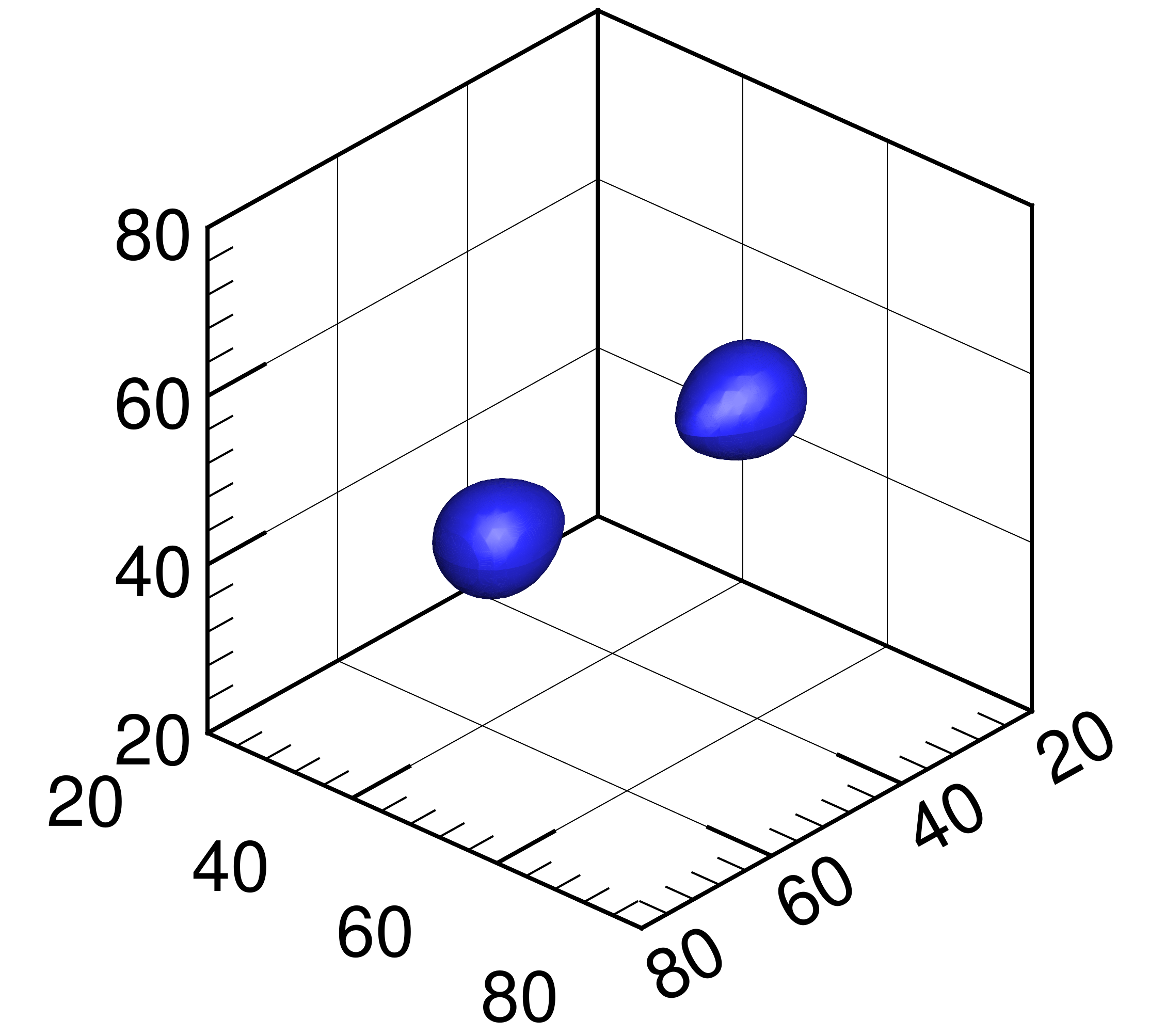}
\caption{Evolution of a dumbbell-shaped surface under curvature-driven flow. The images are snapshots of the interface in time from left to right, top to bottom where the handle of the dumbbell shrinks faster than the spherical shells due to its high curvature leading to a pinch-off for t = 0.0, 2.0, 4.0, 6.0, 8.0, 10.0, 12.0, 14.0, 16.0.}
\label{dumbbell}
\end{figure}

\noindent In this example, we demonstrate the method's capability to handle three-dimensional geometries. A dumbbell shape is initialized via the union of the level curves of two spheres and a cylinder such that 
\begin{equation}
\begin{cases}
    \phi_{sphere,R} = \sqrt{(x-x_c+o)^2+(y-y_c)^2+(z-z_c)^2}-r \\
    \phi_{sphere,L} = \sqrt{(x-x_c-o)^2+(y-y_c)^2+(z-z_c)^2}-r \\
    \phi_{cylinder} = max\bigg{[}(|x-x_c|-o),\bigg{(}\sqrt{(y-y_c)^2+(z-z_c)^2}-w\bigg{)}\bigg{]}
\end{cases}
\end{equation}
and the zeroth level set is given by
\begin{equation}
    \phi(\mathbf{x},0) = min[\phi_{sphere,R},\phi_{sphere,L},\phi_{cylinder}]
\end{equation}
The coordinates of the domain center are $\mathbf{x_c}$, $o$ is the center-to-center distance between the cylinder and each of the two spheres, $r$ is the radius of the spherical shells, and $w$ is the radius of the cylinder. The computational domain is a cube of side length equal to 100 units and grid size of $200 \times 200 \times 200$. The time-step $\Delta t$ is 0.01, the center of the dumbbell level curve is $\mathbf{x_c}=(50,50,50)$, $r=10$, $w=5$, and $o=20$. The total number of iterations is 1600 ($t_{final}=16.0$) and snapshots of interface evolution are shown in Fig \ref{dumbbell} for $t=0.0, 2.0, 4.0, 6.0, 8.0, 10.0, 12.0, 14.0,$ and $ 16.0$. The problem parameters describing dumbbell geometry are again chosen to match the reported values in \cite{Alame2020} for better comparison. Fig \ref{dumbbell} shows that the VOF methodology presented is capable of capturing all three unique features of the problem (pinching, merging, separation) without any special algorithmic treatment of the interface. The dumbbell handle initially collapses faster than the spheres on both ends due to higher curvature. The higher shrinkage rate at the handle then causes pinch-off whilst the spherical shells are shrinking simultaneously at a slower rate. The end result is two teardrop surfaces that will continue to collapse due to their curvature (see Fig. \ref{dumbbell}).

\par Upon running the problem with a coarser grid size of $100\times 100\times 100$, parasitic ``wisps" appear in the pinch-off region. Arrufat et al. \cite{Arrufat2021} define these wisps as cells with tiny values of $1-C_{i,j,k}$ in fluid 1 or $C_{i,j,k}$ in fluid 2. Hence, their appearance for the dumbbell problem is not surprising. \cite{Arrufat2021} suggest using a higher clipping tolerance to bypass that problem, however, there is still no clear criterion for choosing that tolerance as they had reported. The wisps did not decrease upon varying the clipping tolerance, nonetheless, they disappeared upon doubling the grid size. Note that our clipping tolerance ranged between $\epsilon=10^{-13}$ and $\epsilon=10^{-3}$ for the trial runs.
\par
\begin{figure}[H]
\centering
\begin{subfigure}[c]{0.45\textwidth}
\includegraphics[width=\textwidth]{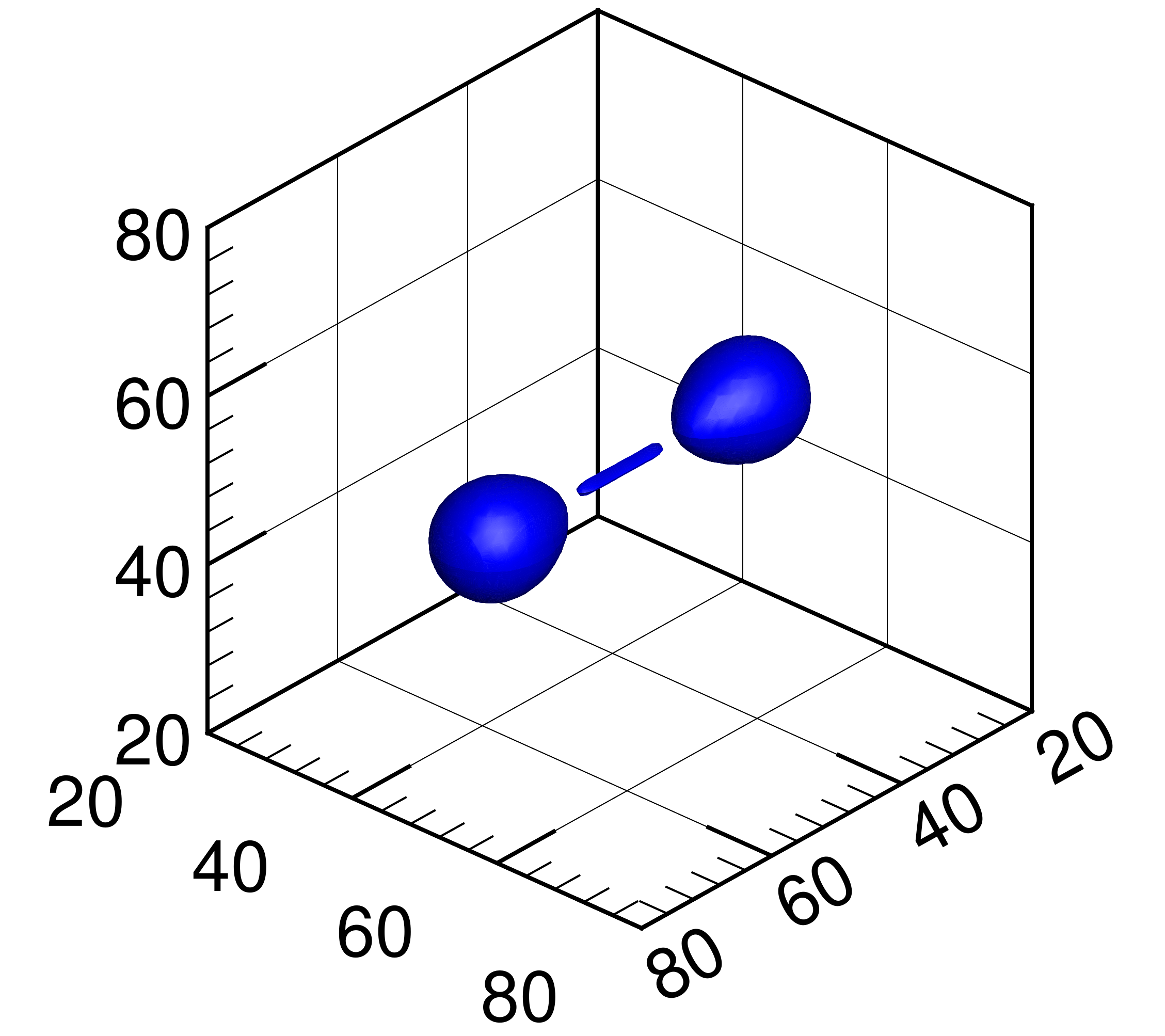}
\end{subfigure}
\qquad
\begin{subfigure}[c]{0.45\textwidth}
\includegraphics[width=\textwidth]{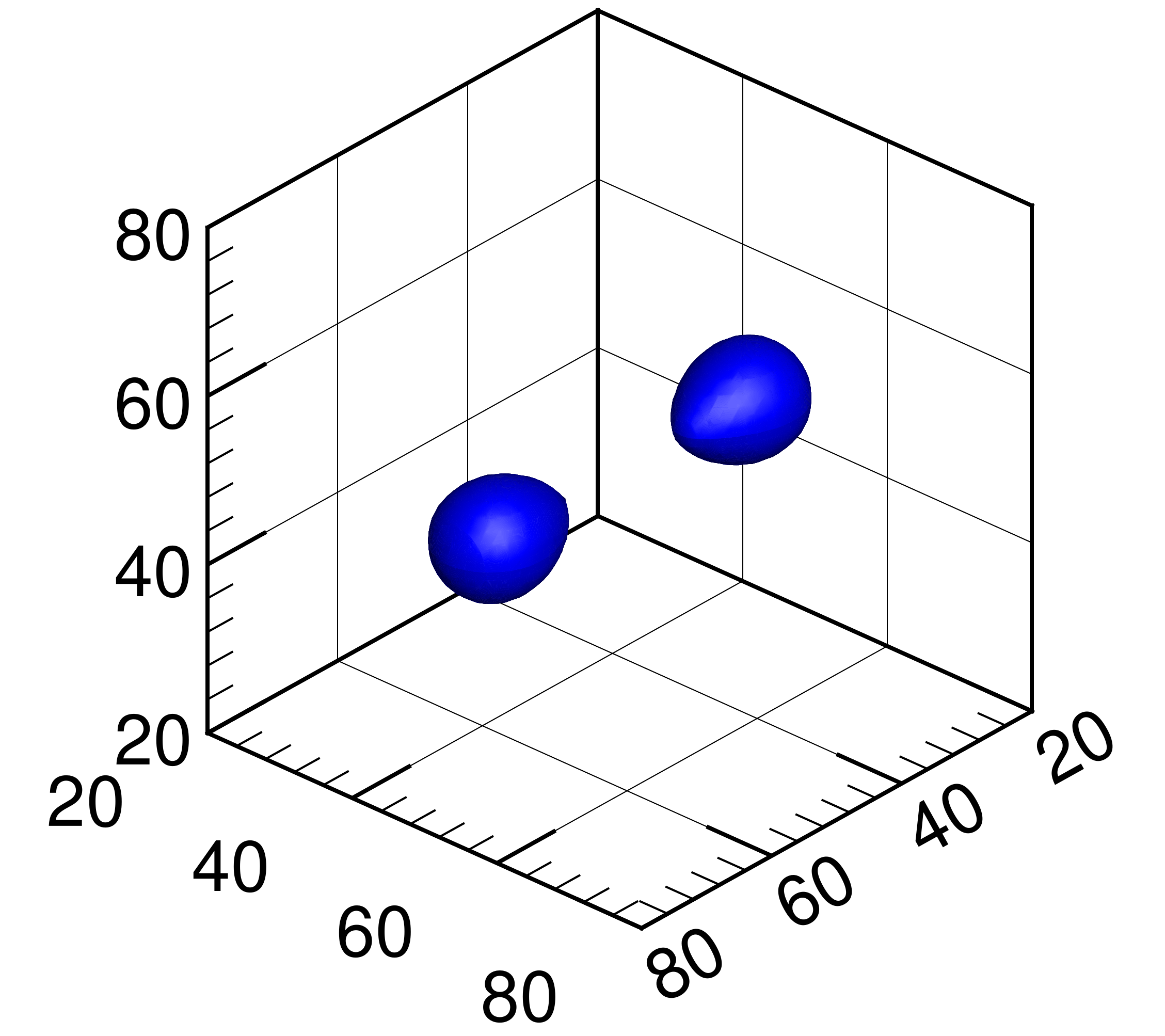}
\end{subfigure}
\caption{Dumbbell pinch-off into two tear drop shapes using WY method (left) and CIAM method (right).}
\label{dumb_comp}
\end{figure}

 Before doubling the grid size and varying the clipping tolerance, the advection of the volume fraction was also performed with CIAM (Calcul d'Interface Affine par Morceaux) \cite{Li1995,Arrufat2021} alongside the WY method to compare the pinch-off region. As shown in Fig. \ref{dumb_comp}, the pinch-off wisp does not appear for the $100\times 100\times 100$ grid using CIAM. The WY algorithm is still the advection method that we employ for all the presented work due to better conservation properties.

\vspace*{0.5cm}
\begin{figure}[H]
\centering
\includegraphics[width=.33\textwidth]{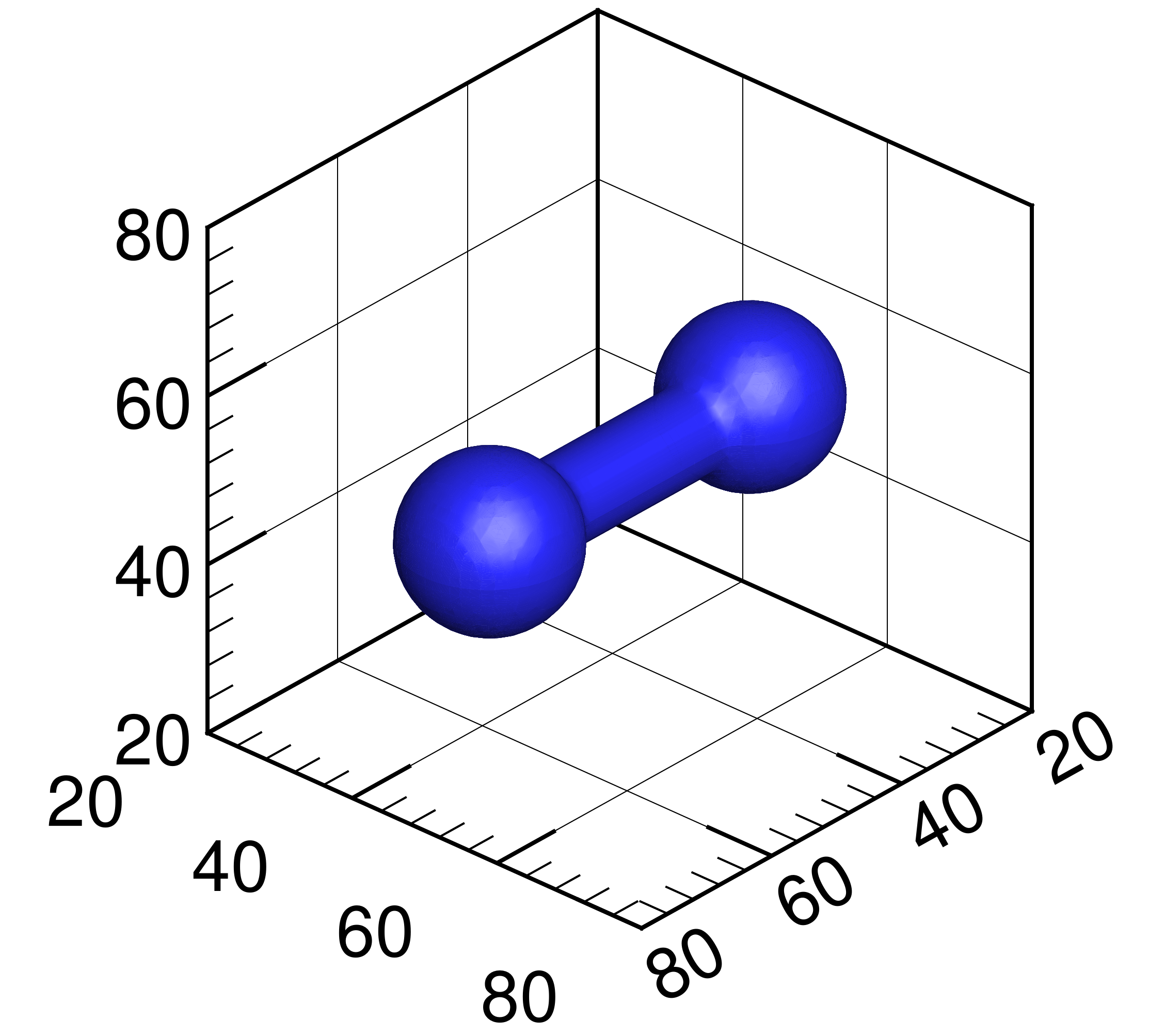}
\includegraphics[width=.33\textwidth]{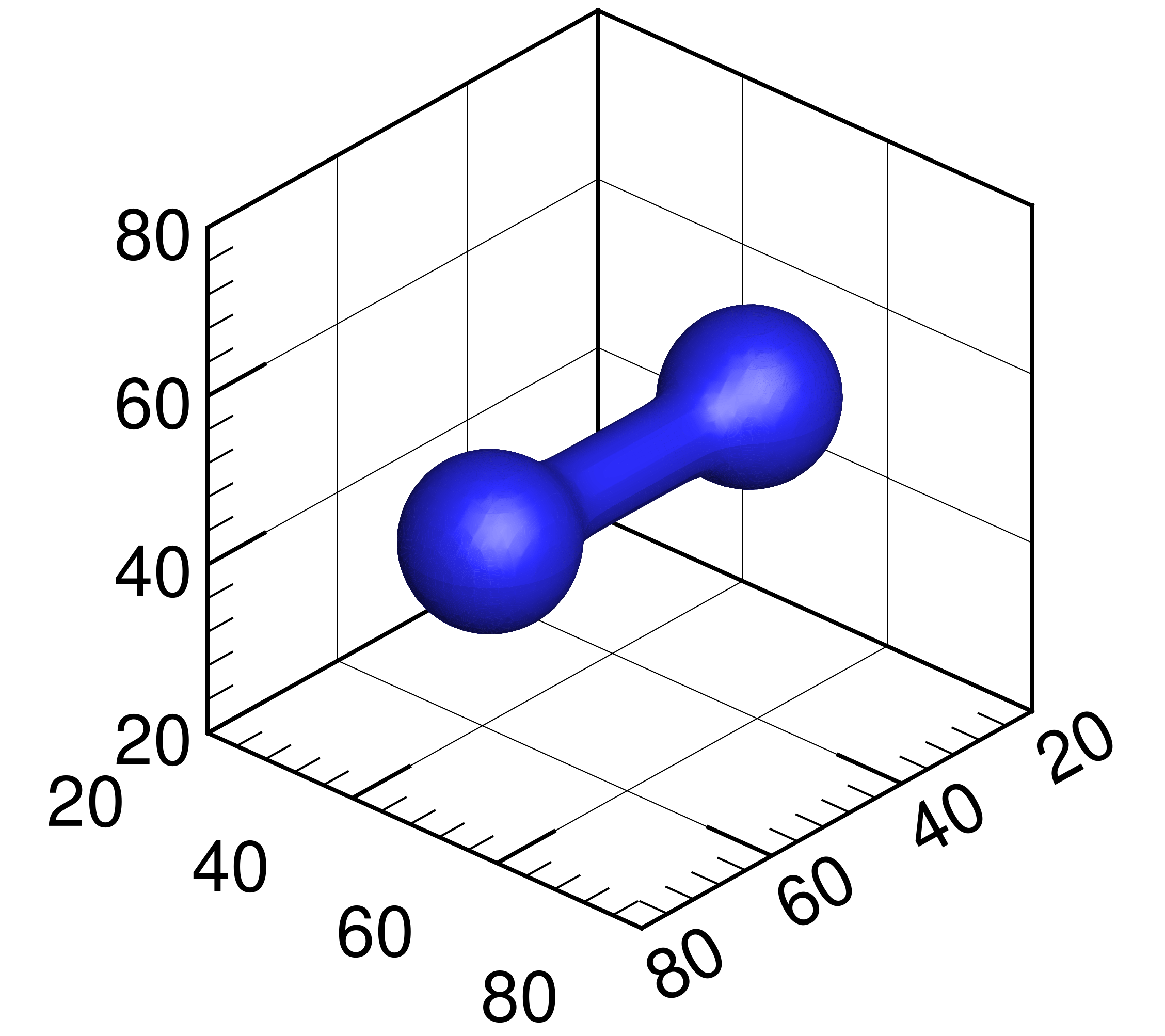} 
\includegraphics[width=.33\textwidth]{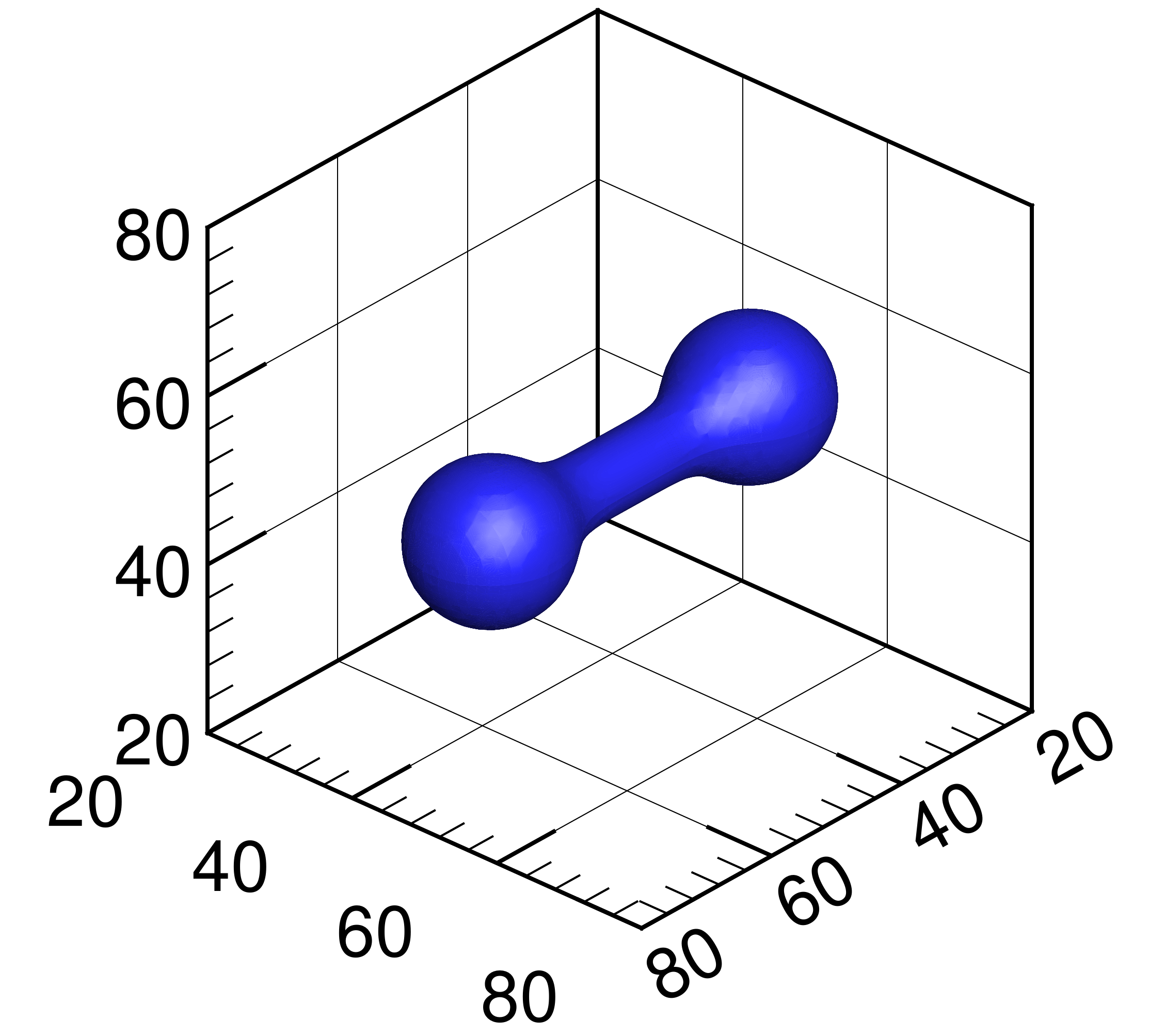}
\end{figure}
\begin{figure}[H]
\centering
\includegraphics[width=.33\textwidth]{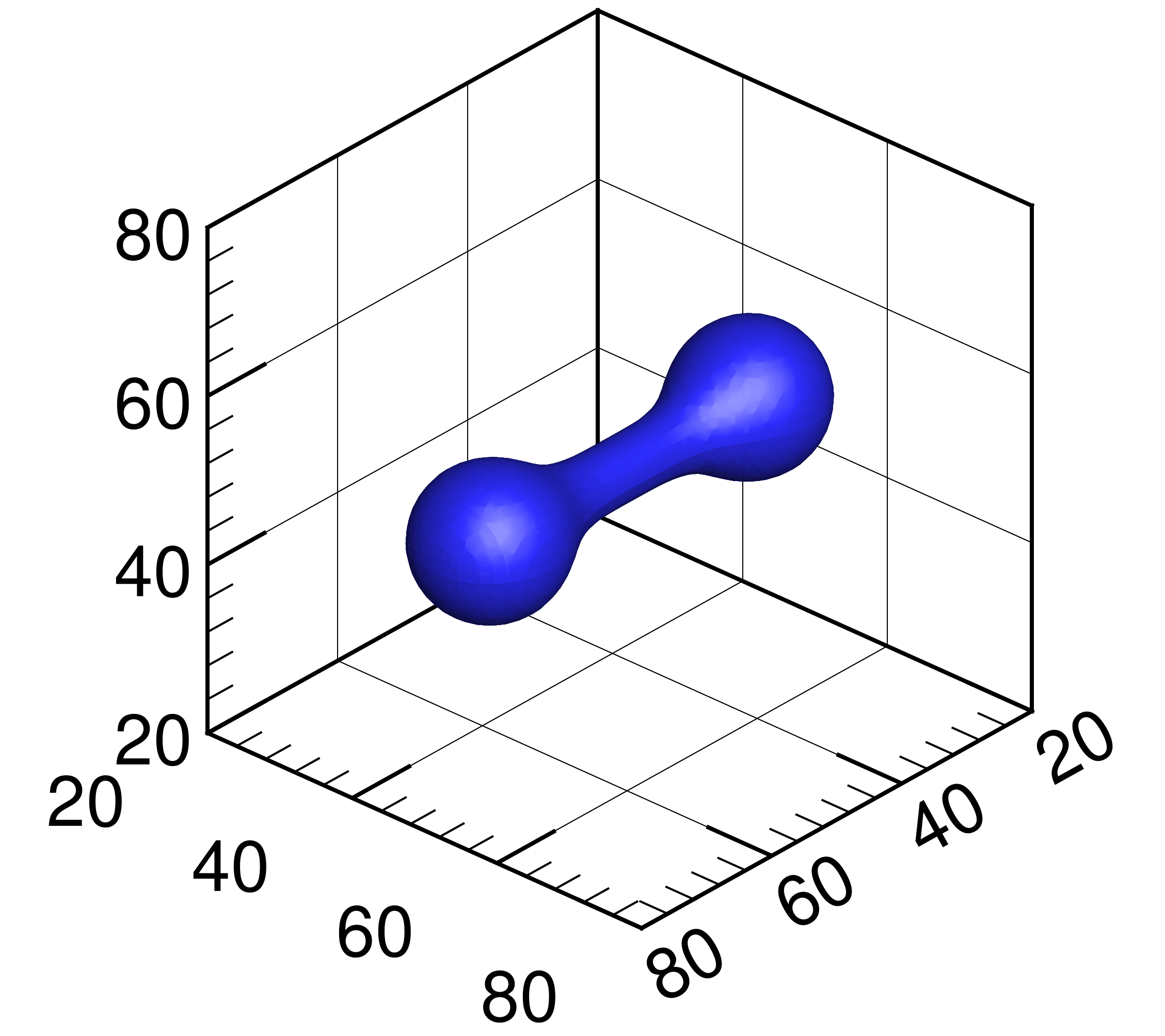}
\includegraphics[width=.33\textwidth]{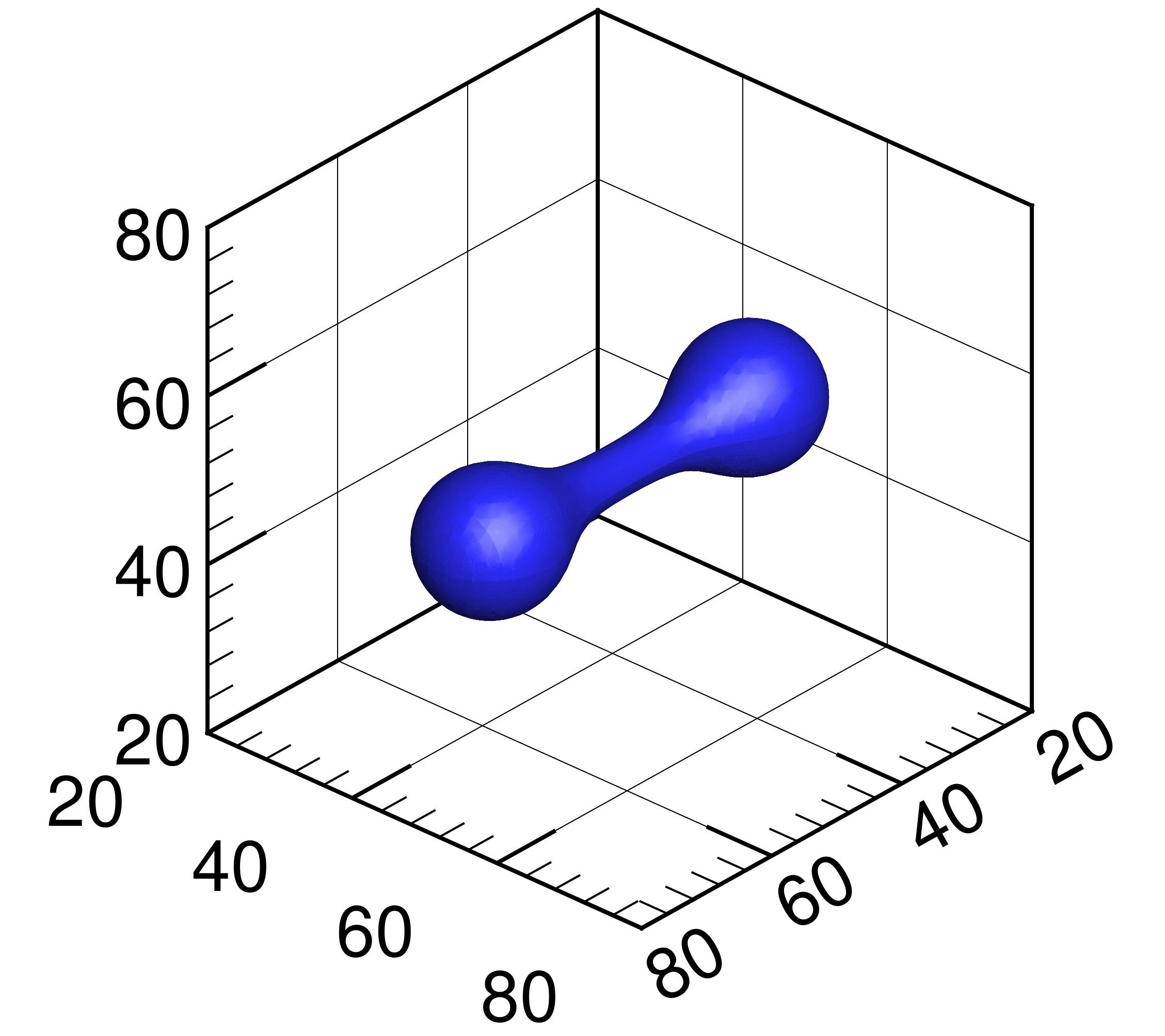} 
\includegraphics[width=.33\textwidth]{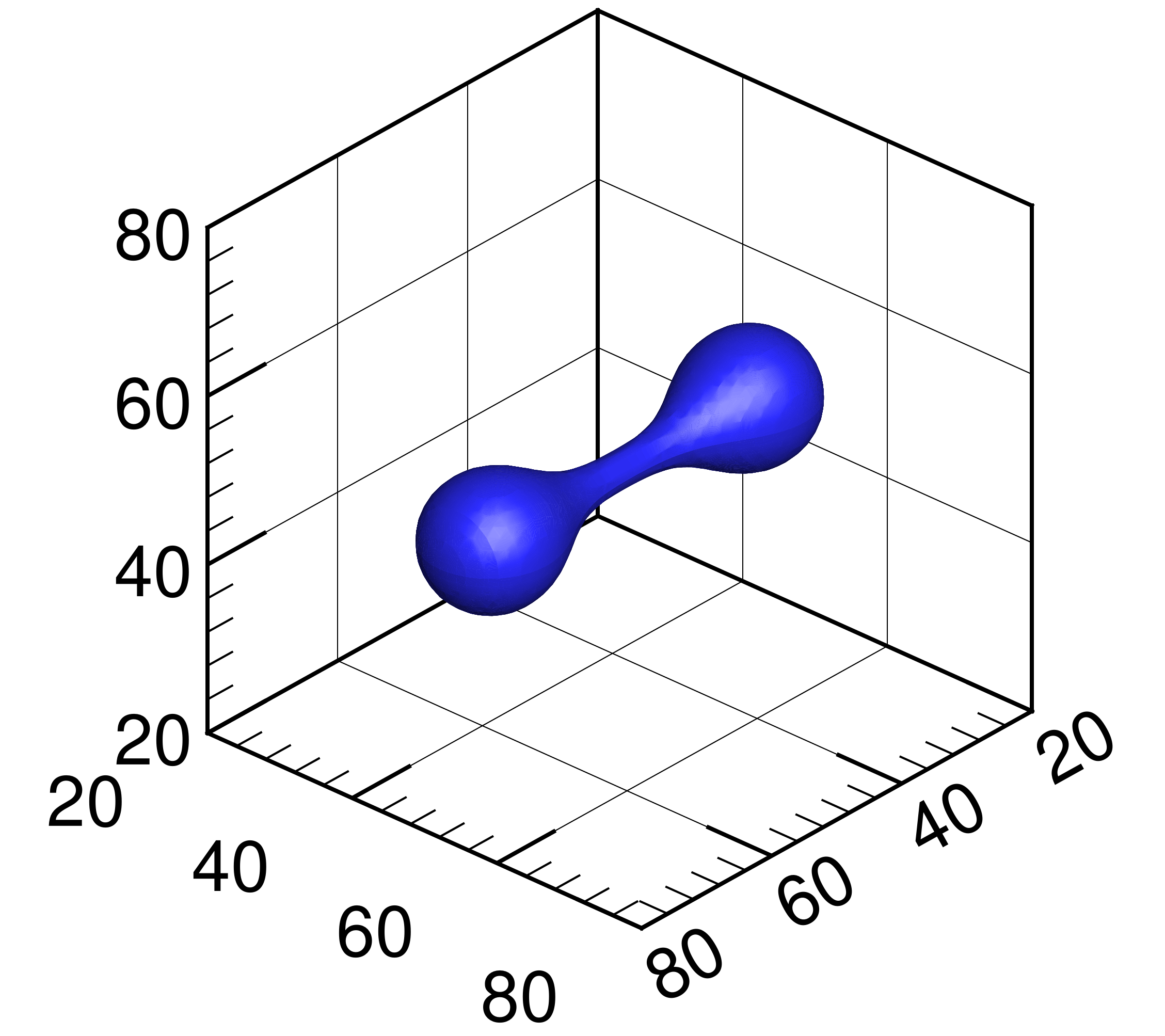}
\end{figure}

\begin{figure}[H]
\centering
\includegraphics[width=0.45\textwidth]{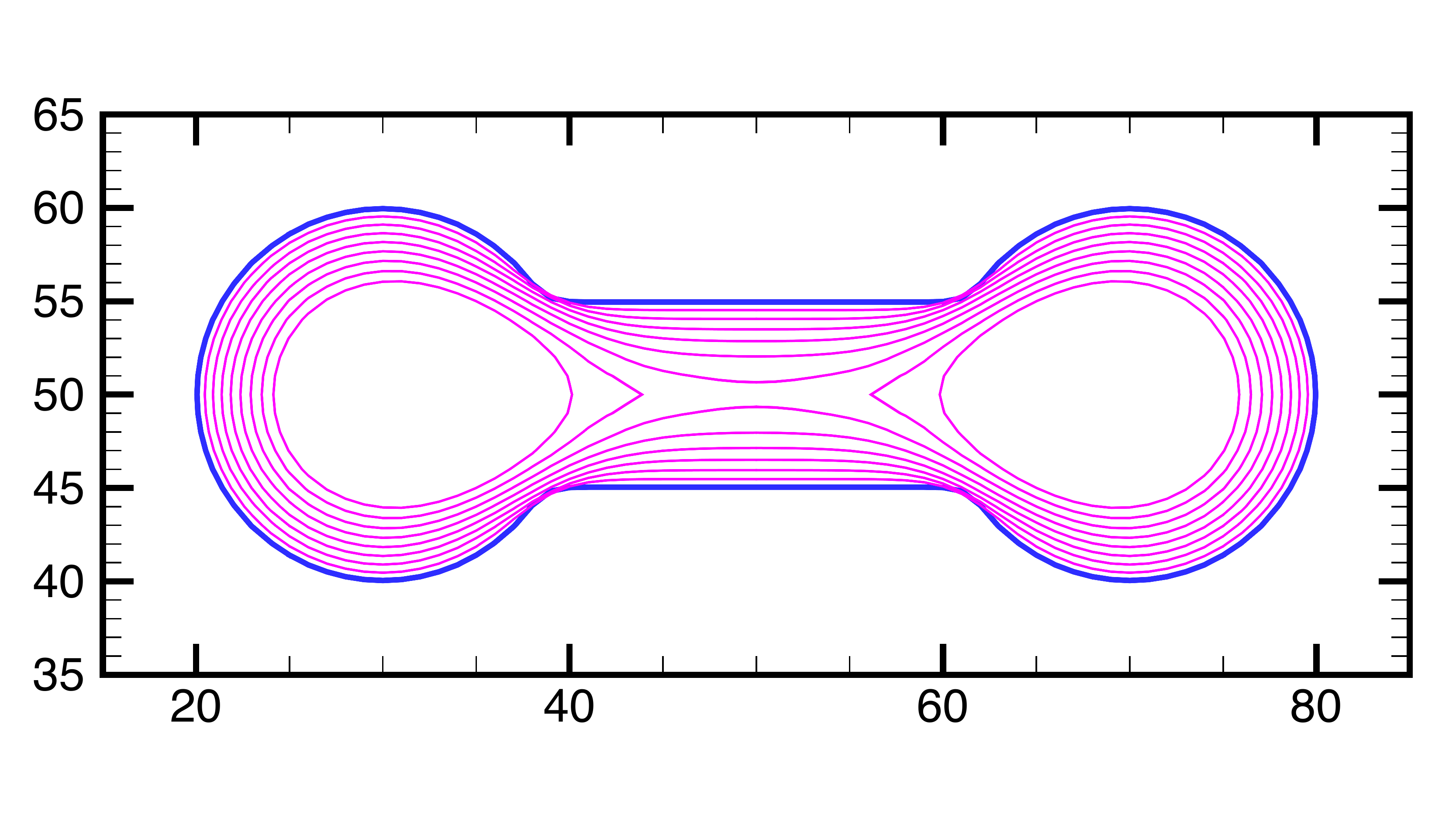}
\includegraphics[width=0.45\textwidth]{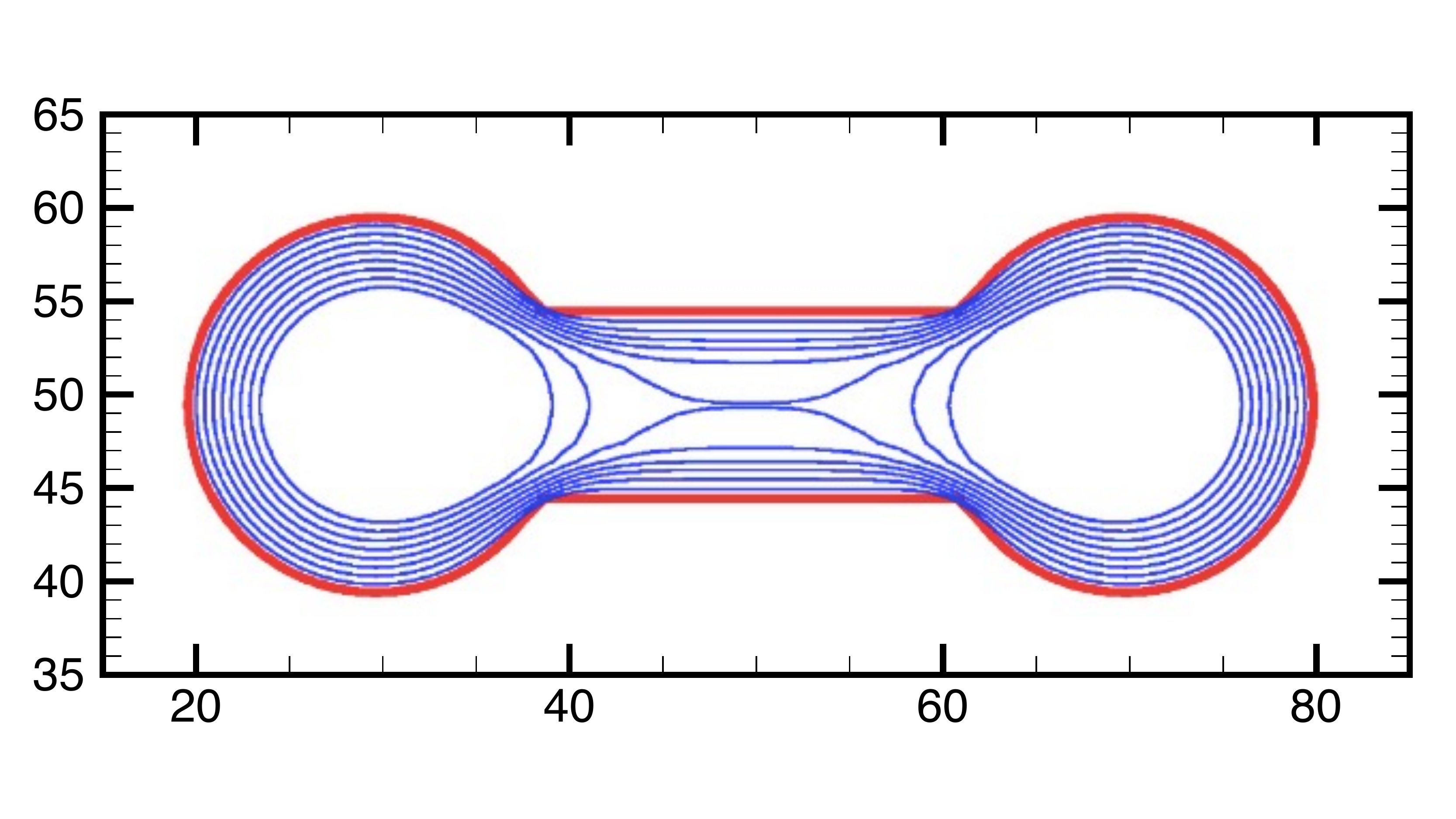}
\caption{Z-plane slice showing the evolution of dumbbell interface in time with VVOF (left) and DRLSE \cite{Alame2020} (right).}
\label{isocontours_alame}
\end{figure}

\begin{figure}[H]
\centering
\includegraphics[width=.32\textwidth]{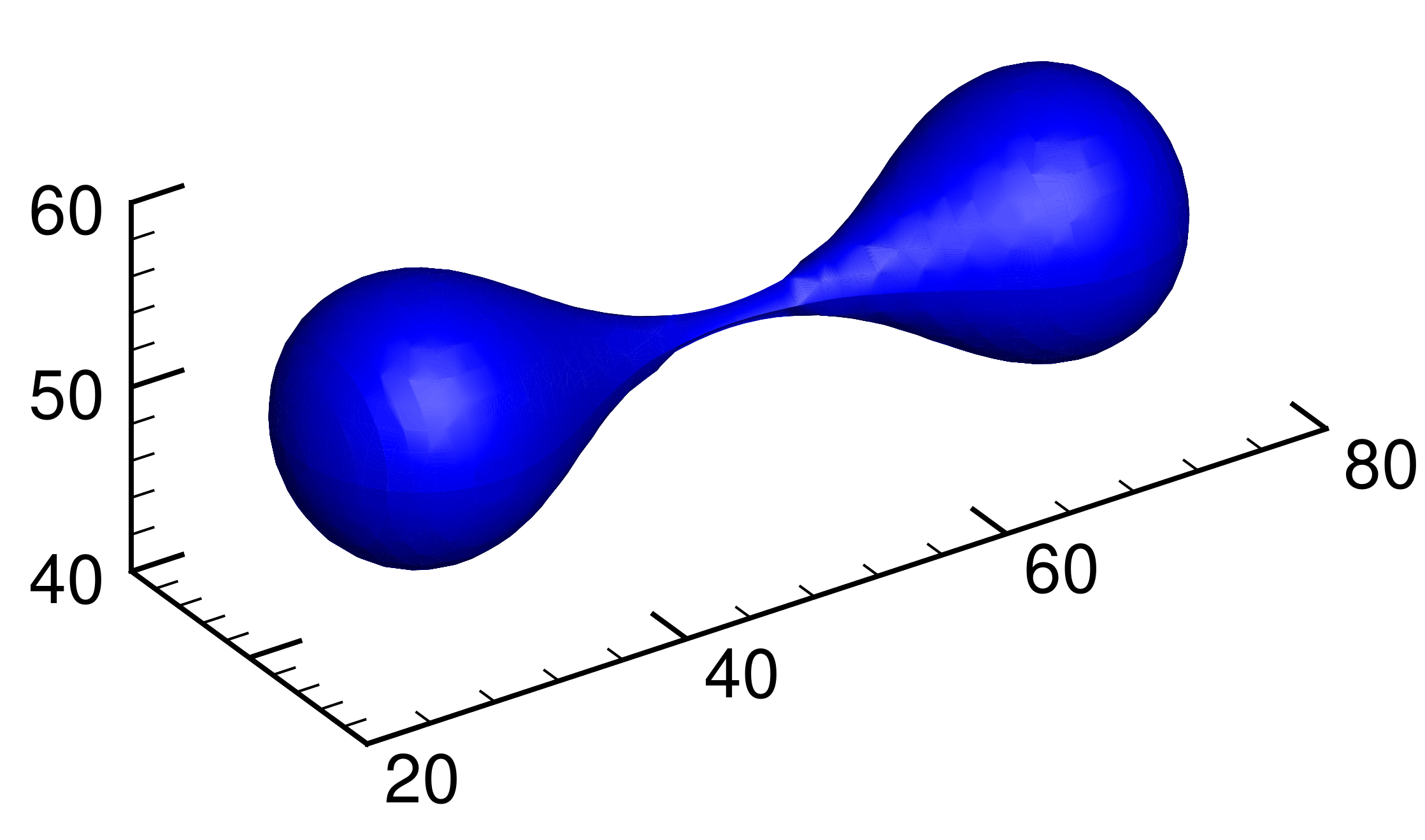}
\includegraphics[width=.32\textwidth]{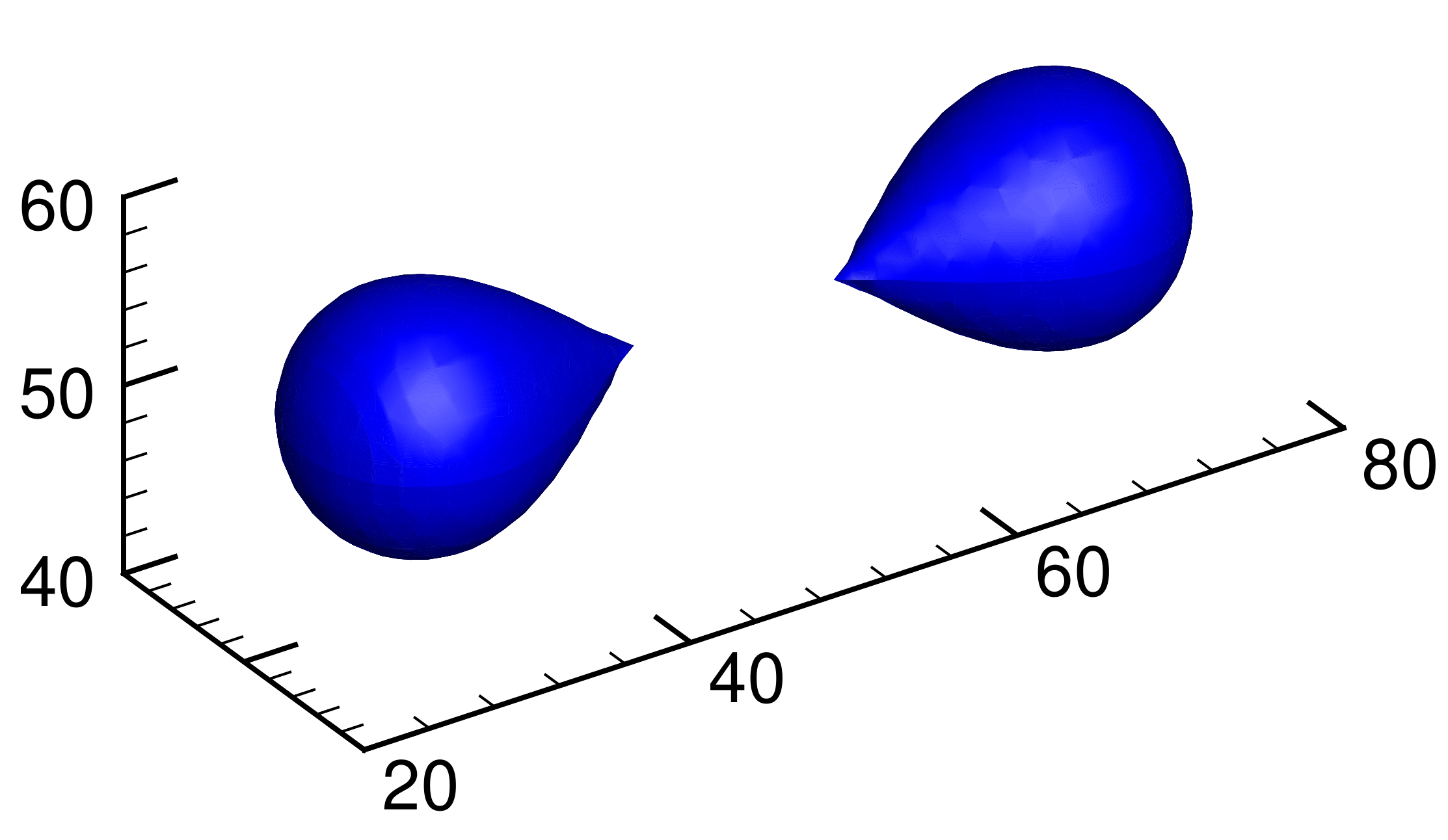} 
\includegraphics[width=.32\textwidth]{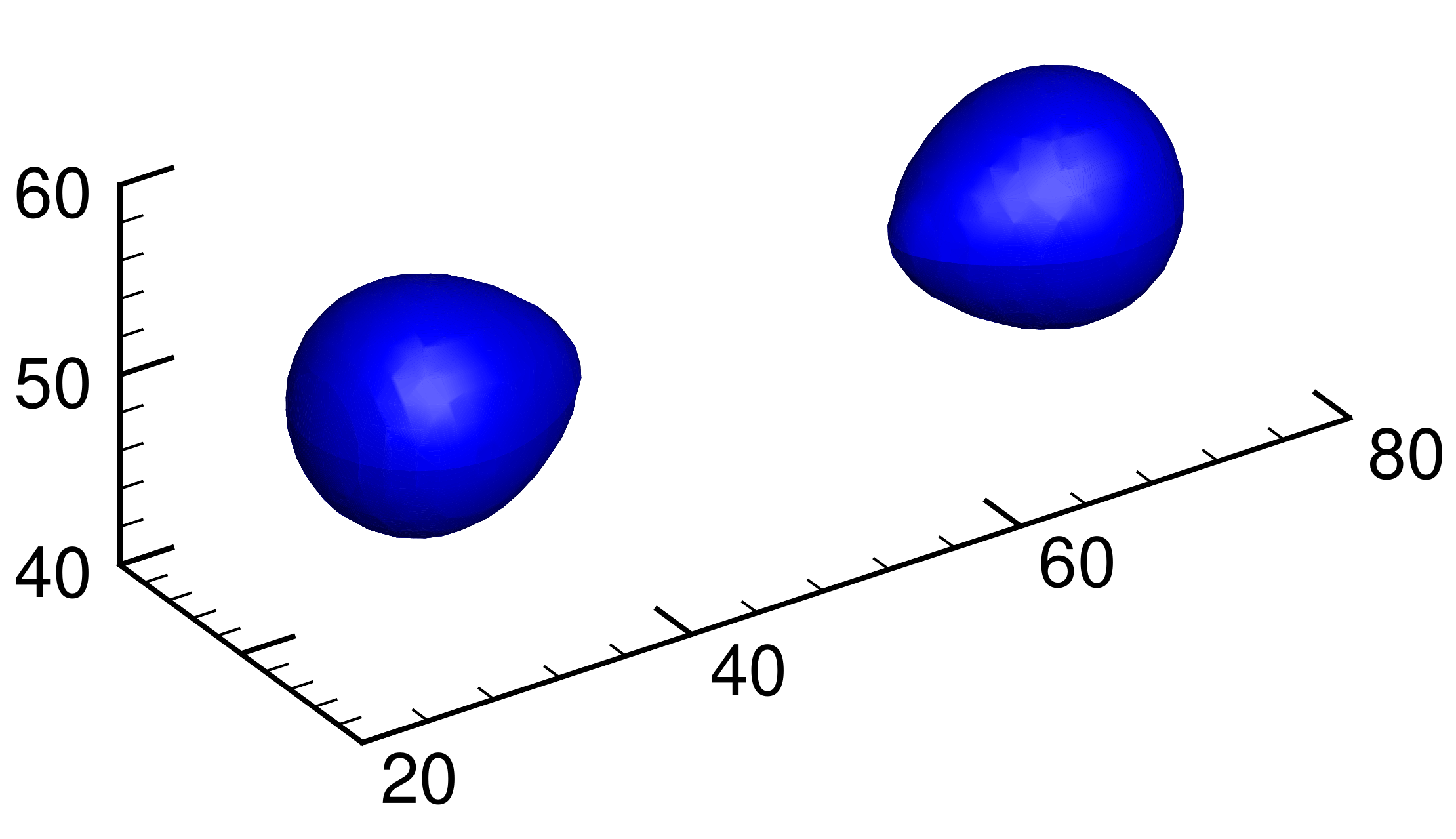} 
\caption{Snapshots of final stages of dumbbell pinch-off using VVOF.}
\label{isosurfaces_comp}
\end{figure}

The results for the dumbbell shape were compared with those obtained with the DRLSE method of \cite{Alame2020}. Fig. \ref{isocontours_alame} shows isocontours of the interface as the dumbbell evolves in time using VVOF (left) and the variational LS method (right). VVOF shows good agreement with variational LS with a minor difference in the shape of the pinch-off region. Isosurfaces of both methods at the last stages of pinch-off were also compared as shown in Fig. \ref{isosurfaces_comp}. The three main topological phenomena that characterize the dumbbell problem are pinching, merging, and separation, as mentioned earlier. VVOF is shown to capture all three features with good agreement compared to the variational LS method. 

\begin{figure}[H]
\centering
\includegraphics[width=.32\textwidth]{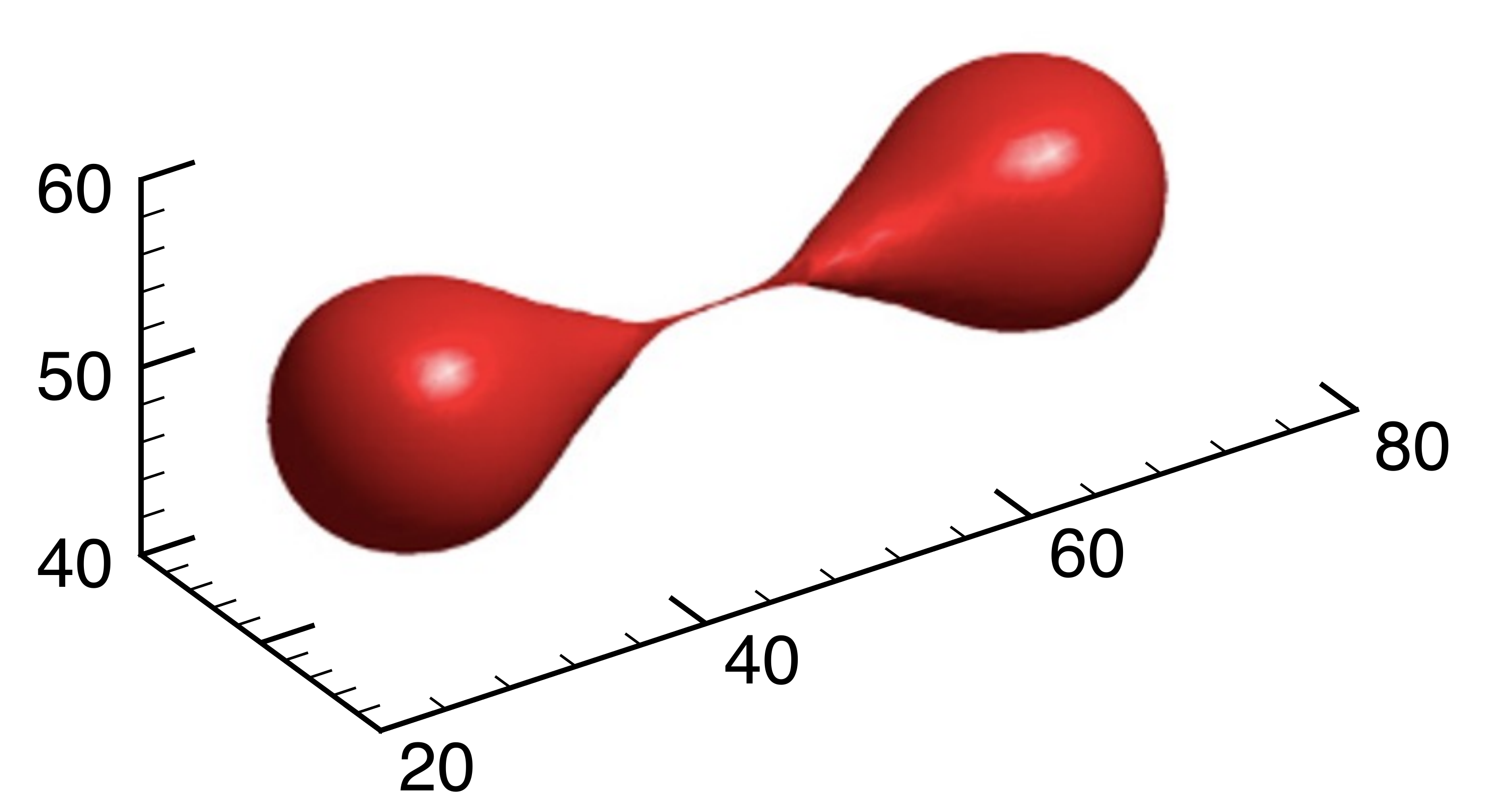}
\includegraphics[width=.32\textwidth]{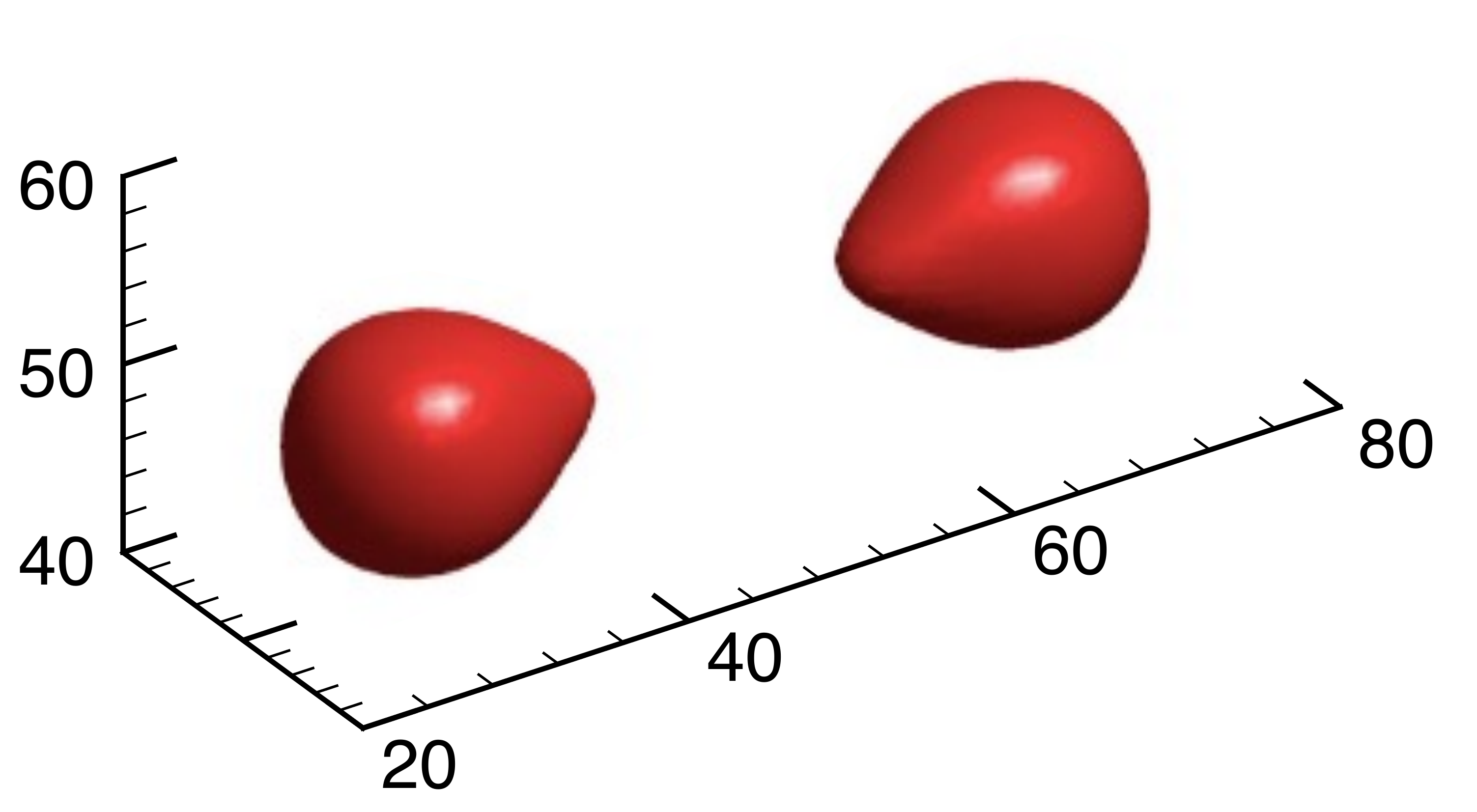}
\includegraphics[width=.32\textwidth]{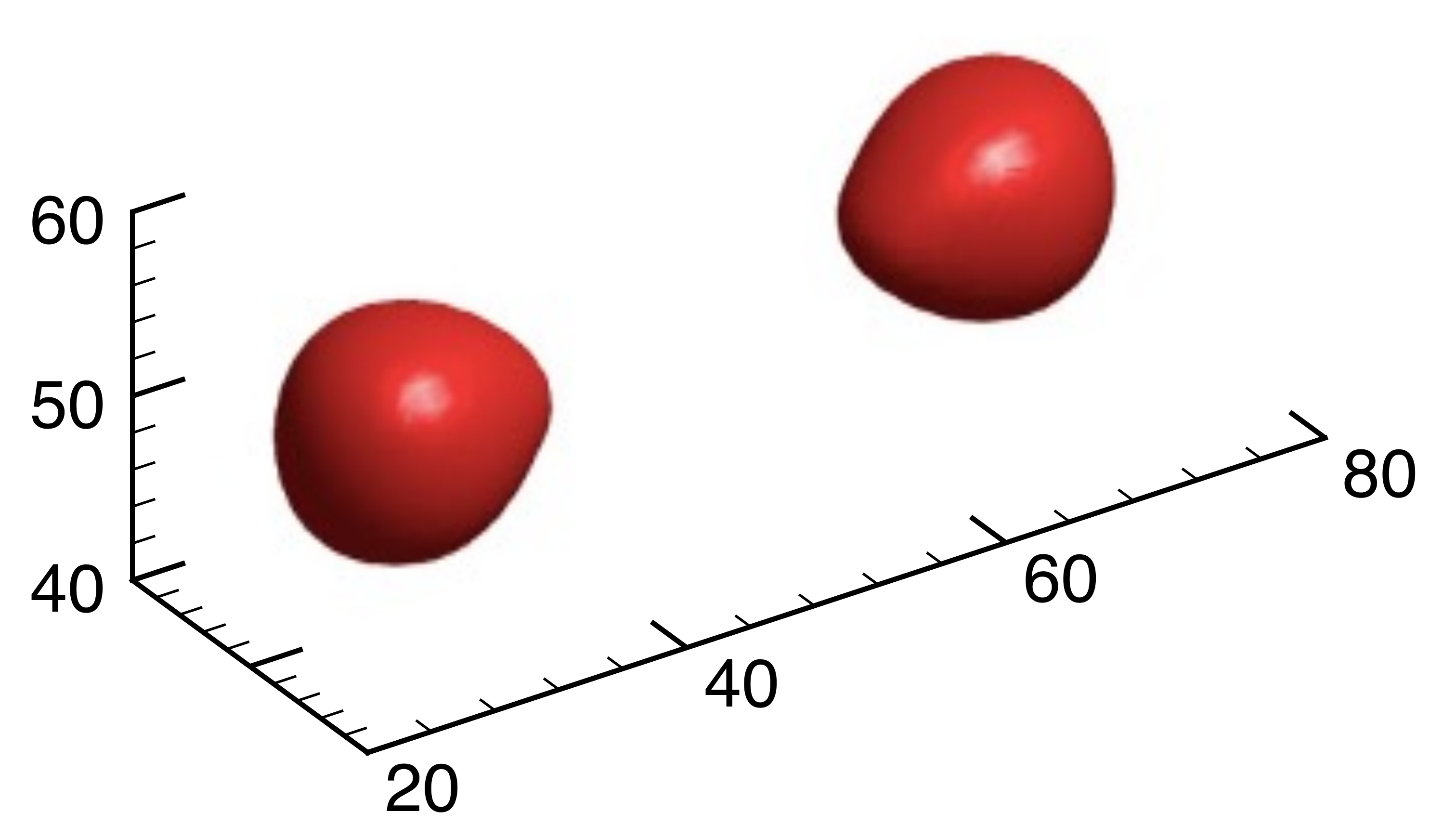}
\caption{Snapshots of final stages of dumbbell pinch-off using DRLSE \cite{Alame2020}.}
\end{figure}

\subsubsection{Irregular geometry}
\begin{figure}[H]
\centering
\includegraphics[width=.33\textwidth]{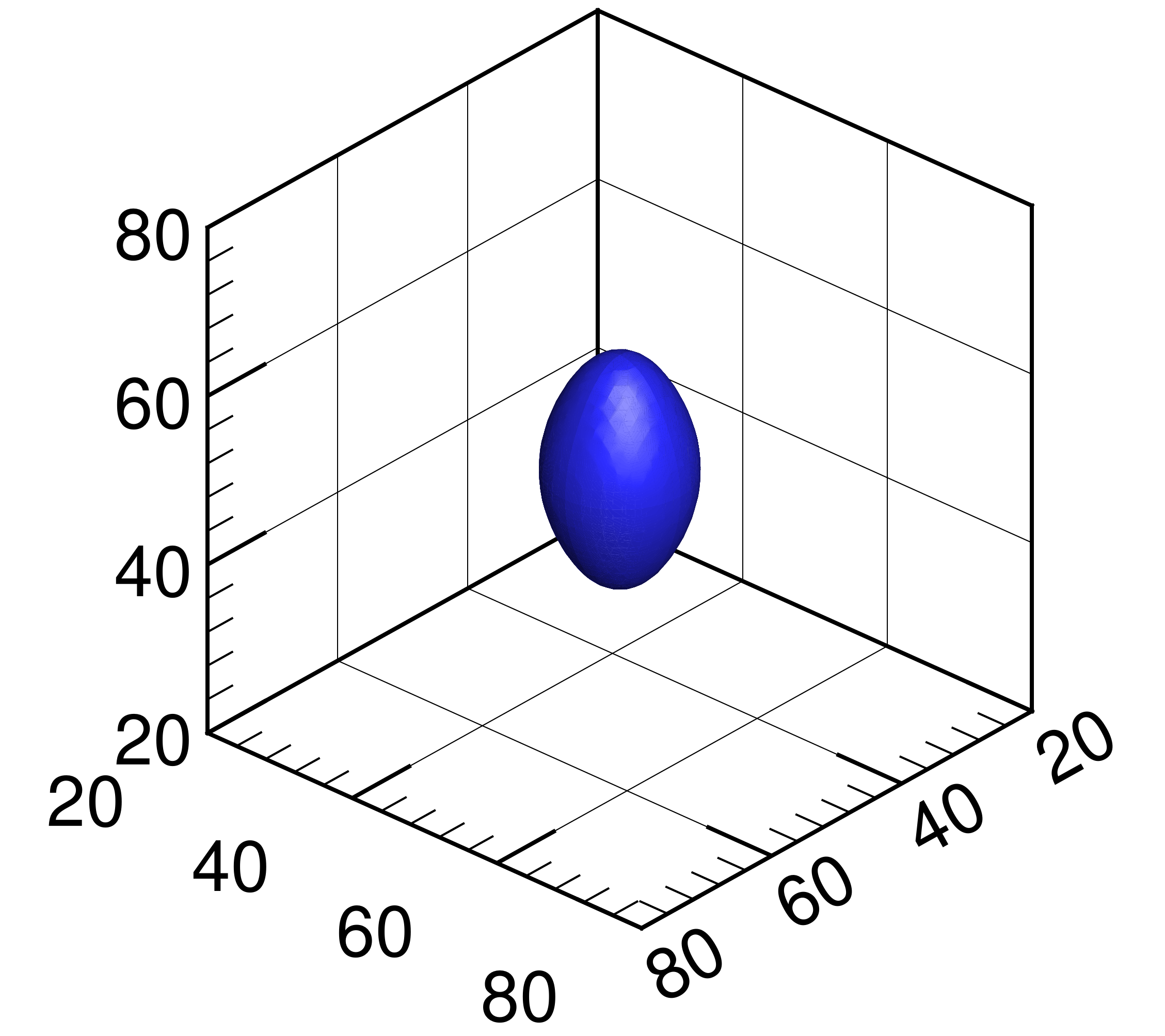}
\includegraphics[width=.33\textwidth]{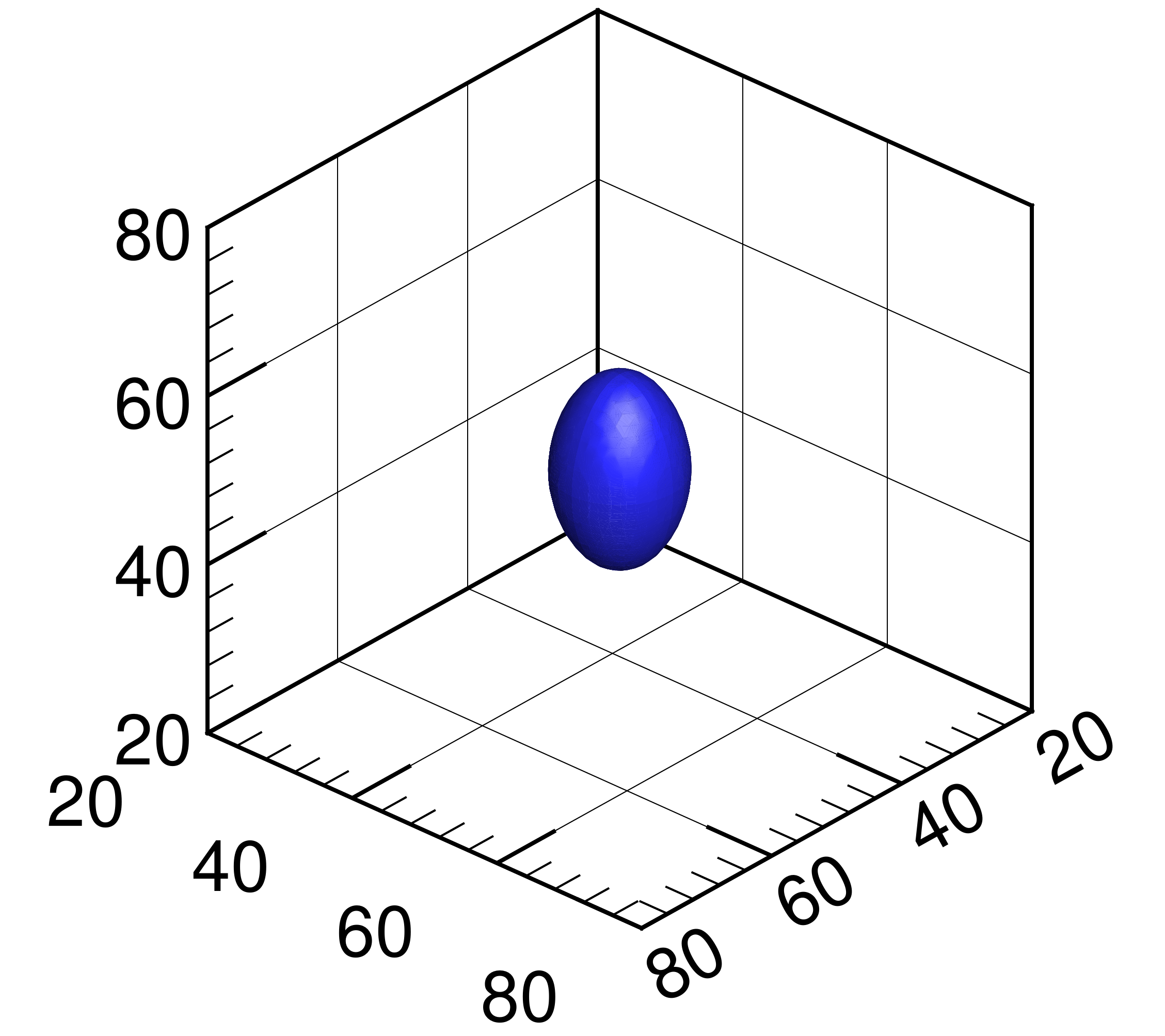} 
\includegraphics[width=.33\textwidth]{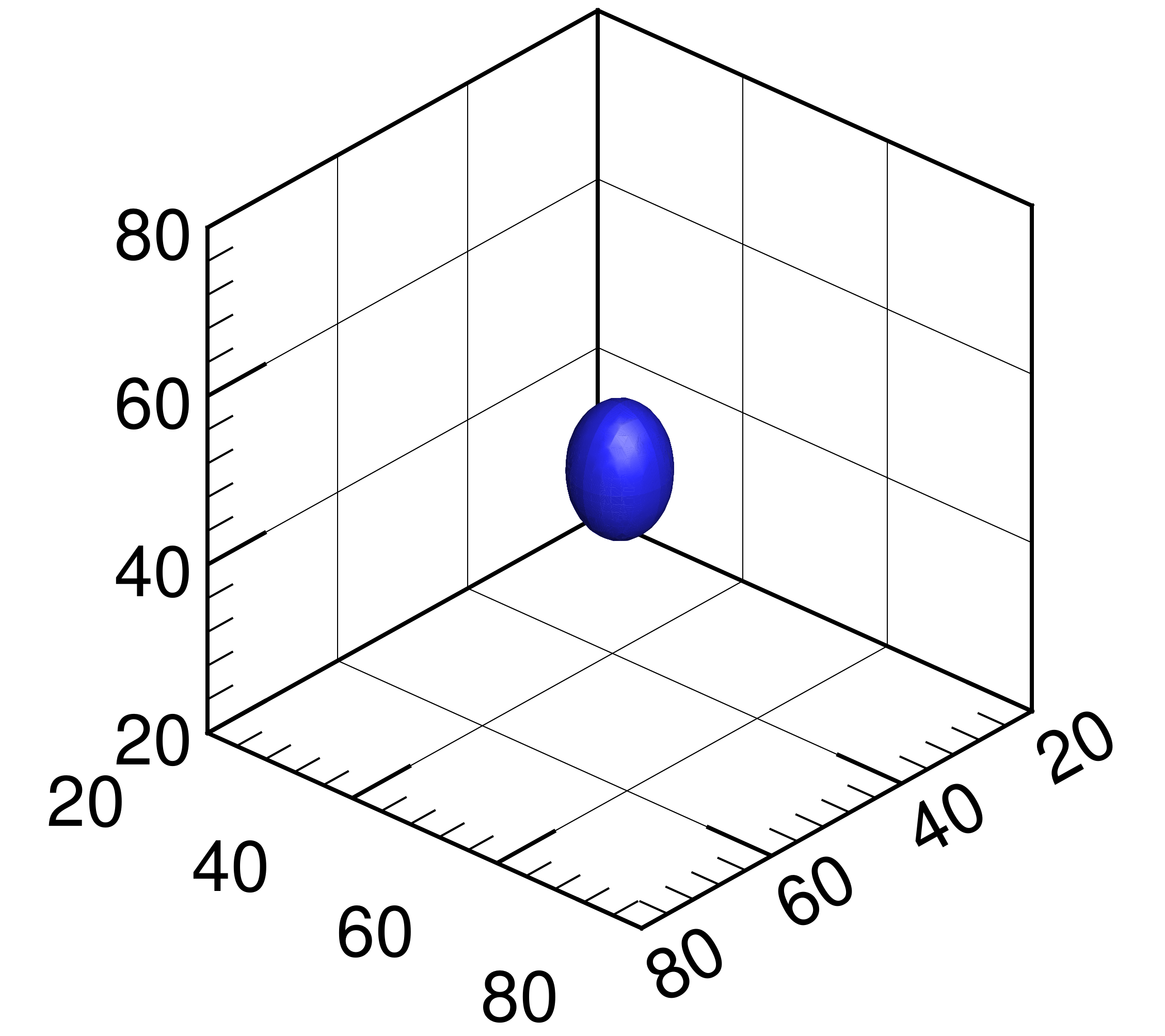} 
\caption{Evolution of an irregular shape under curvature-driven flow. The images are snapshots of the interface in time from left to right, top to bottom where the horizontal heads shrink faster than the vertical spheroid due to higher curvature for t = 0.0, 6.25, 12.5, 18.75, 25.0, 31.25, 37.5, 43.75, 50.0.}
\label{spheroid}
\end{figure}

Similar to Sec \ref{dumbbell_sec}, we demonstrate the method's capability to handle three-dimensional geometries. The shape under consideration is more irregular and naturally more challenging from the numerical point of view. An irregular shape is initialized via the union of the level curves of three prolate spheroids; a main stem spheroid oriented in the $z$-direction and two secondary spheroids oriented in the $x$ and $y$-directions, respectively. In general, the level curve of a tri-axial ellipsoid centered at point $\mathcal{C}$ with semi-axes $a$, $b$, and $c$ aligned along the coordinate axes is given by
\begin{equation}\label{ellipsoid}
\phi_{ellipsoid} = \frac{(x-x_c)^2}{a^2} + \frac{(y-y_c)^2}{b^2} + \frac{(z-z_c)^2}{c^2} -1 \\
\end{equation}
Using the definition in equation \ref{ellipsoid}, the level curves defining the geometry of this problem can be expressed as
\begin{equation}
\begin{cases}
    \displaystyle{\phi_{spheroid,x} = \frac{(x-x_c)^2}{a_1^2} + \frac{(y-y_c)^2+(z-z_c)^2}{c_1^2} - 1}\\
    \displaystyle{\phi_{spheroid,y} = \frac{(x-x_c)^2+(z-z_c)^2}{a_2^2} + \frac{(y-y_c)^2}{c_2^2} - 1} \\
    \displaystyle{\phi_{spheroid,z} = \frac{(x-x_c)^2+(y-y_c)^2}{a_3^2} + \frac{(z-z_c)^2}{c_3^2} - 1}
    \end{cases}
\end{equation}
and the zeroth level set is given by
\begin{equation}
    \phi(\mathbf{x},0) = min[\phi_{spheroid,x},\phi_{spheroid,y},\phi_{spheroid,z}]
\end{equation}
where the subscripts $x$, $y$, and $z$ refer to the direction of the symmetry axis of the spheroid, $a_i$ is the equatorial radius of the spheroid, and $c_i$ is the distance from center to pole along the symmetry axis. For $c_i>a_i$ the resulting geometry is a prolate spheroid. \\
The computational domain is a cube of side length equal to 100 units and grid size of $100\times 100\times 100$. The time-step $\Delta t=0.05$, the center of the geometry is at $\mathbf{x_c}=(50,50,50)$, $a_1=25$, $a_2=7.5$, $a_3=10$, $c_1=7.5$, $c_2=25$, and $c_3=35$. The total number of iterations is 1200 (t_{final} = 60.0) and snapshots of interface evolution are shown in Fig. \ref{spheroid} for  $t = 0.0, 6.25, 12.5, 18.75, 25.0, 31.25, 37.5, 43.75$, and $50.0$. The four arms parallel to the bottom plane propagate faster towards the stem of the geometry due to higher curvature as the junctions at the stem evolve outwards simultaneously. The difference in front propagation rate between the symmetry axes eventually leads to an egg-like shape that continues to shrink until it disappears at $t=60.0$ (not shown in Fig. \ref{spheroid}).
\begin{figure}[H]
\centering
\includegraphics[width=.52\textwidth]{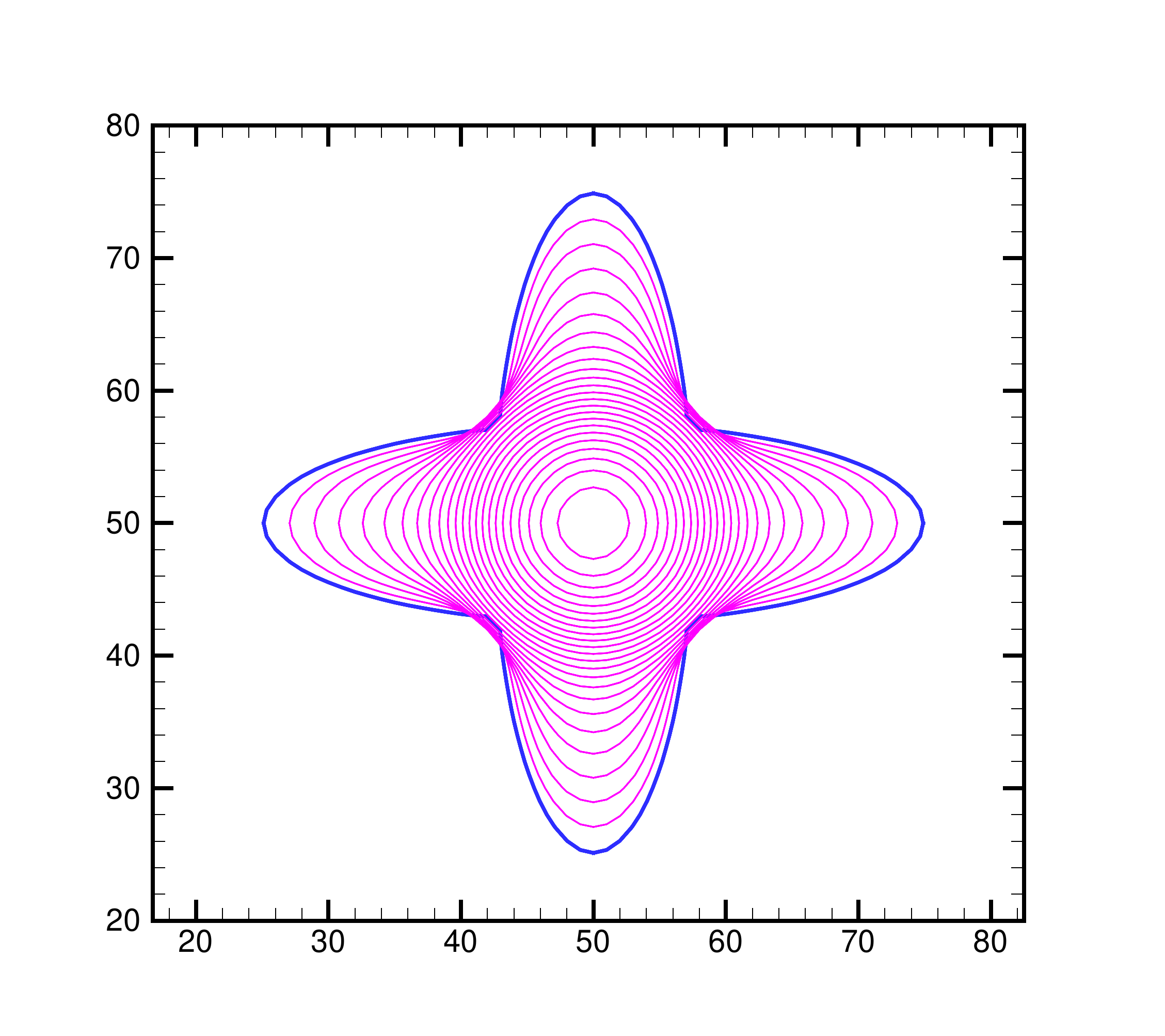}
\includegraphics[width=.45\textwidth]{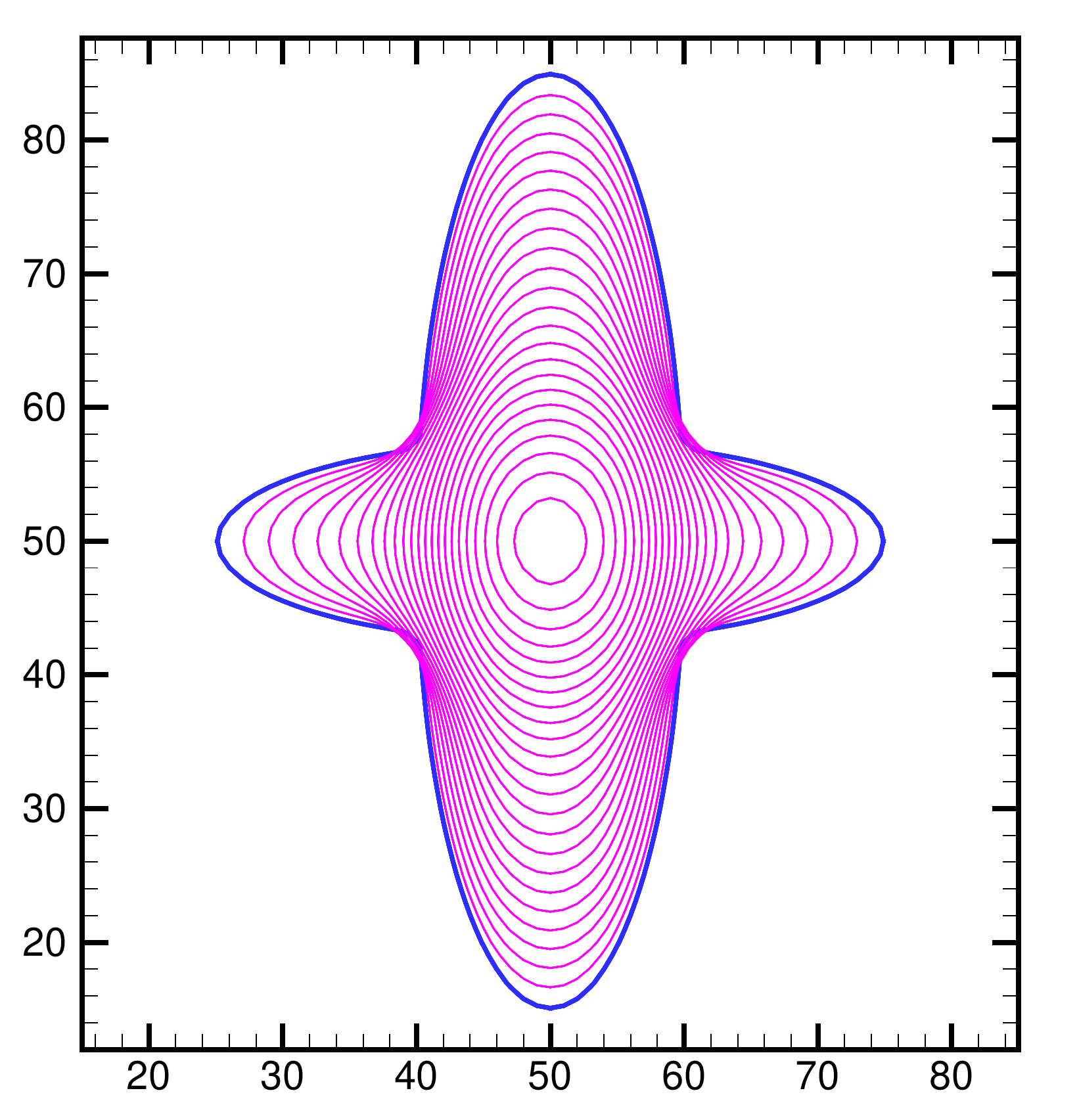} 
\caption{Y-plane slice (left) and Z-plane slice (right) showing the evolution of the interface of the irregular geometry in time.}
\end{figure}

\begin{figure}[H]
\centering
\includegraphics[width=.33\textwidth]{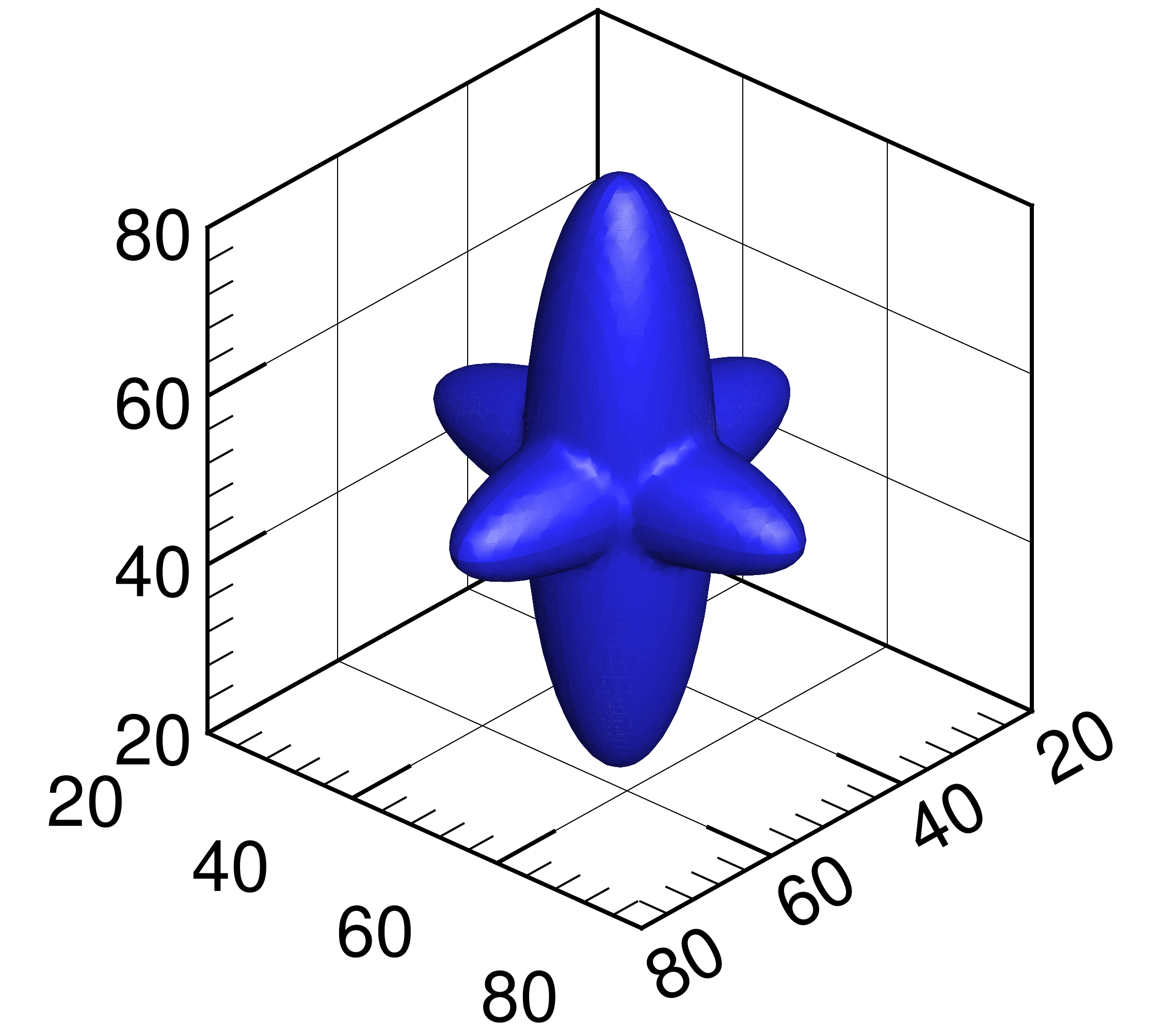}
\includegraphics[width=.33\textwidth]{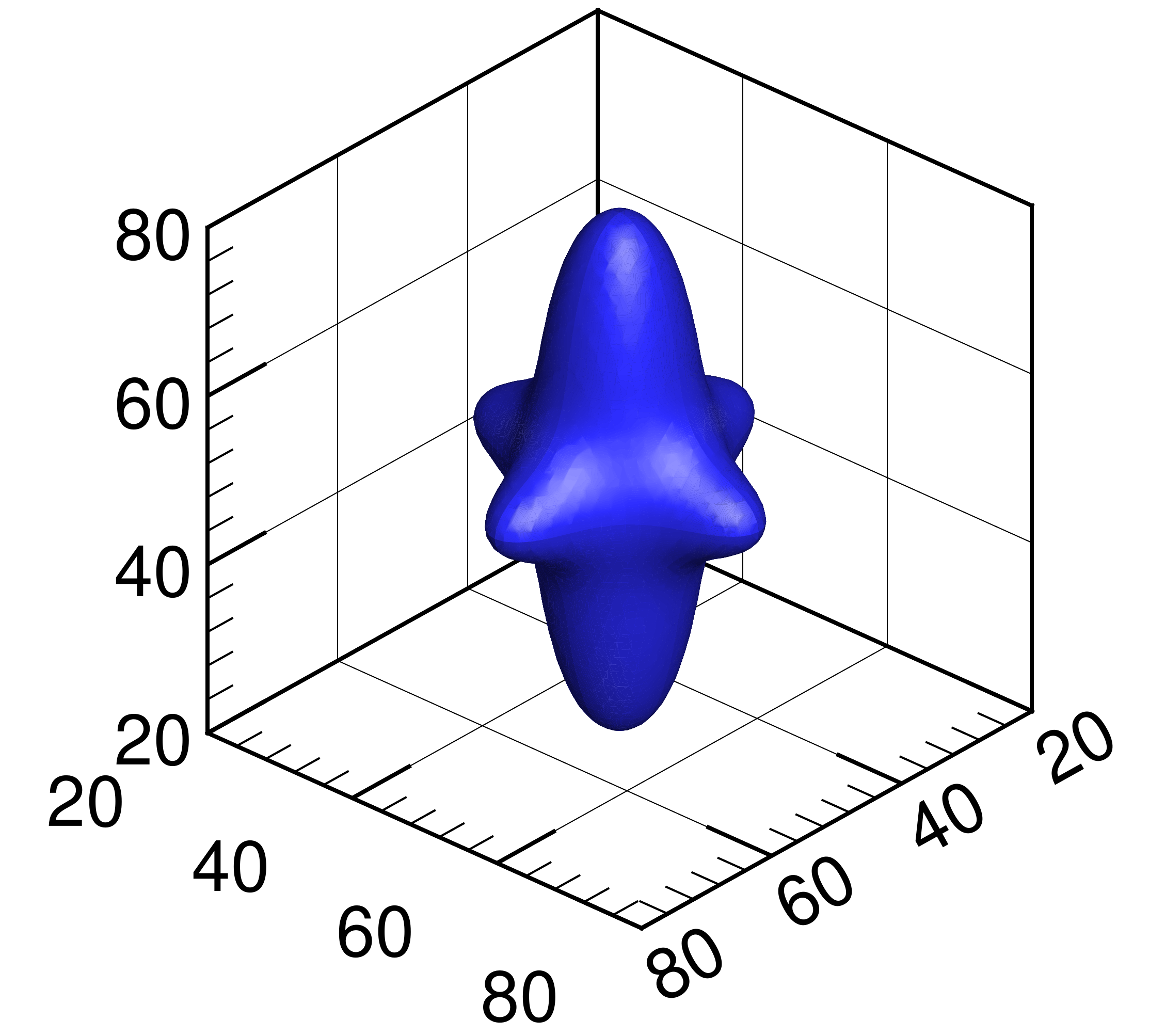} 
\includegraphics[width=.33\textwidth]{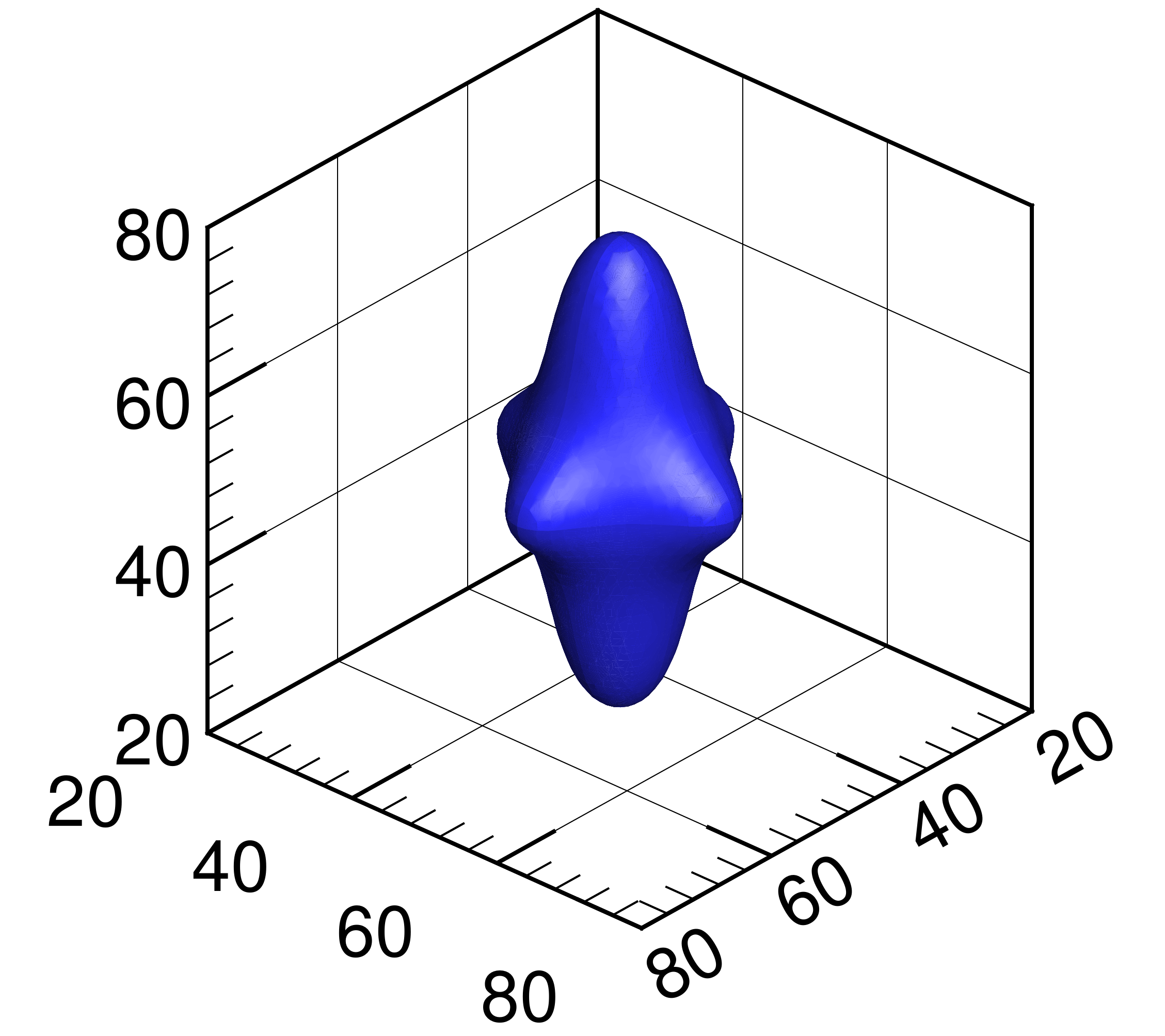}
\end{figure}
\begin{figure}[H]
\centering
\includegraphics[width=.33\textwidth]{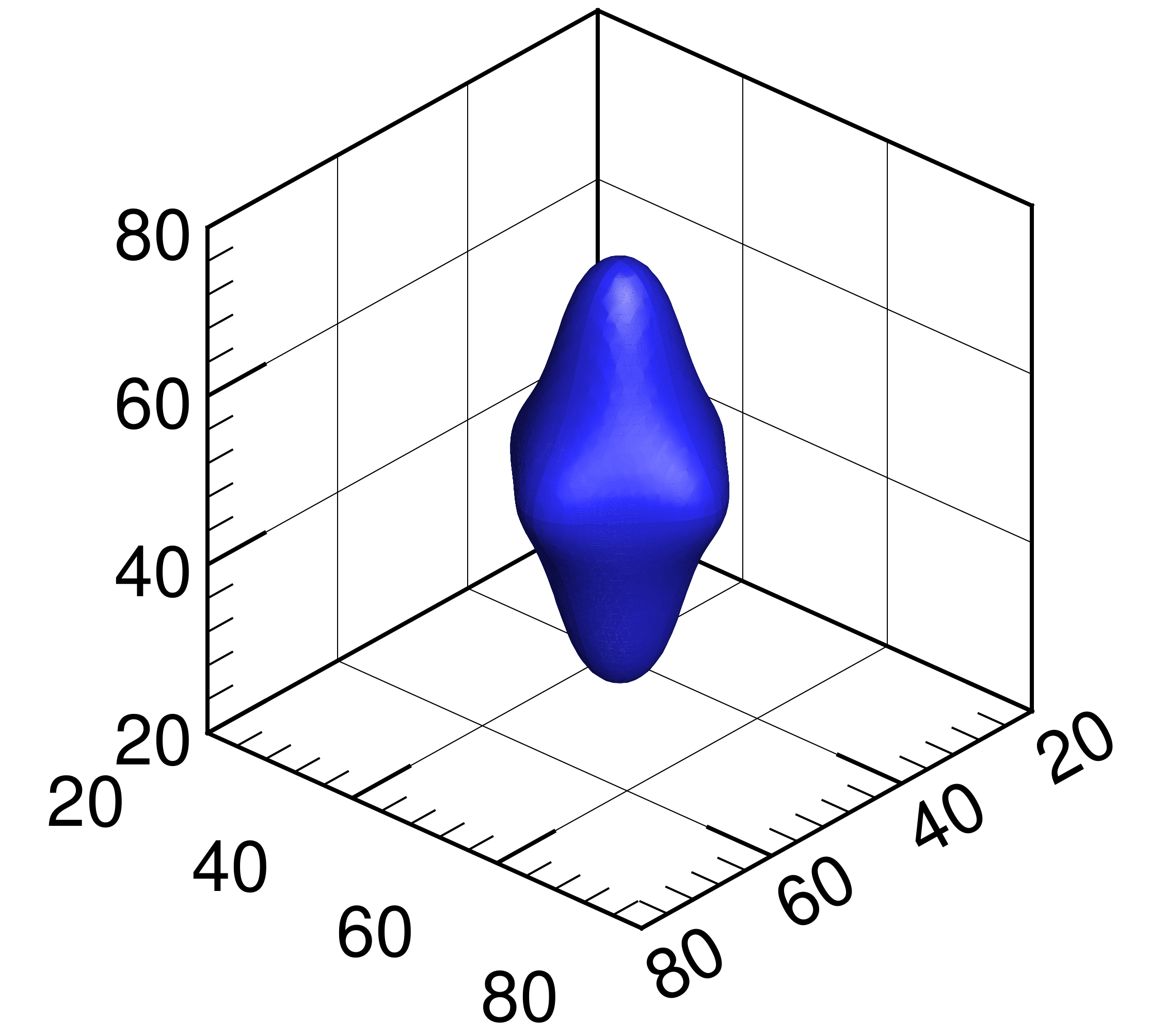}
\includegraphics[width=.33\textwidth]{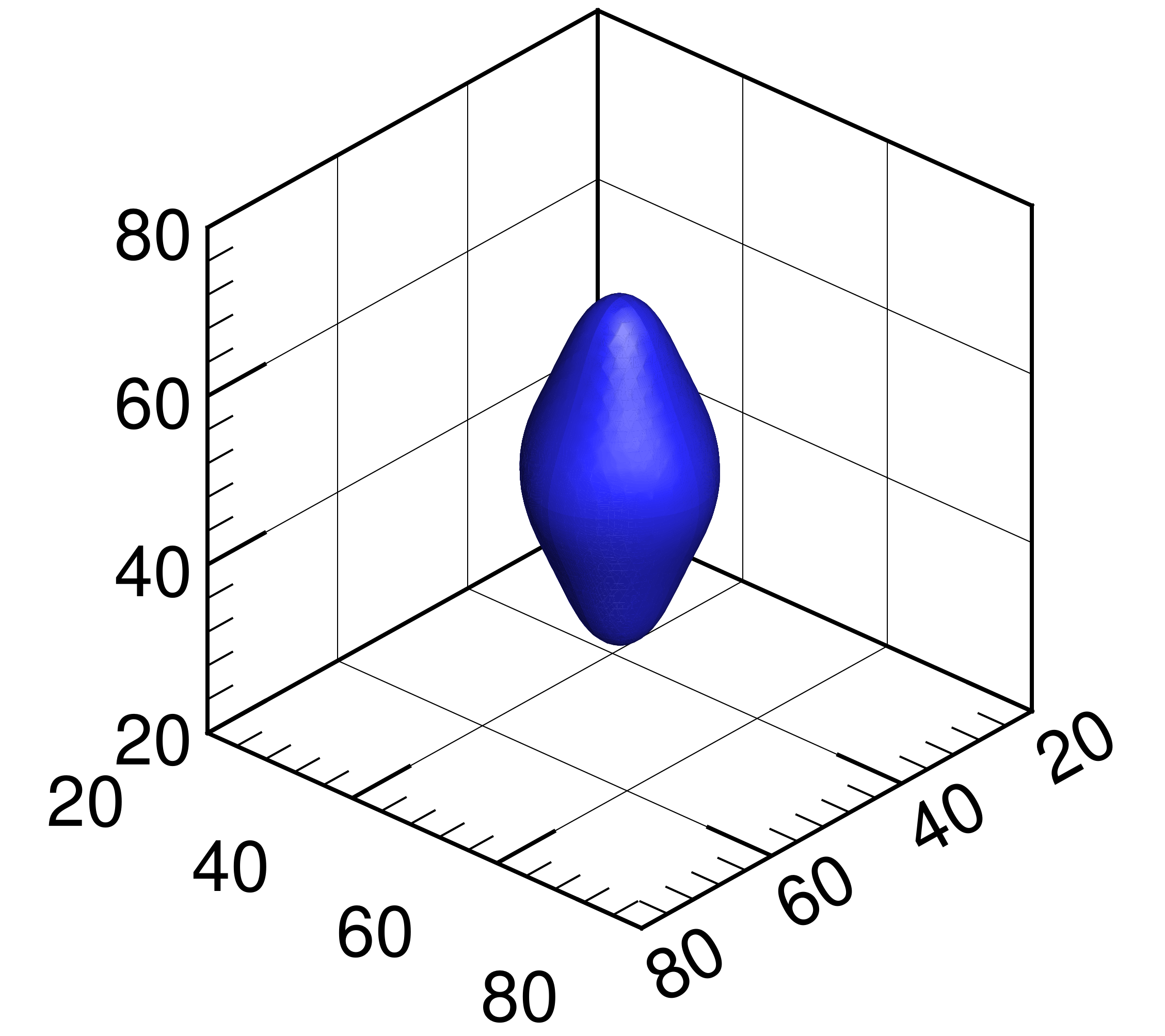} 
\includegraphics[width=.33\textwidth]{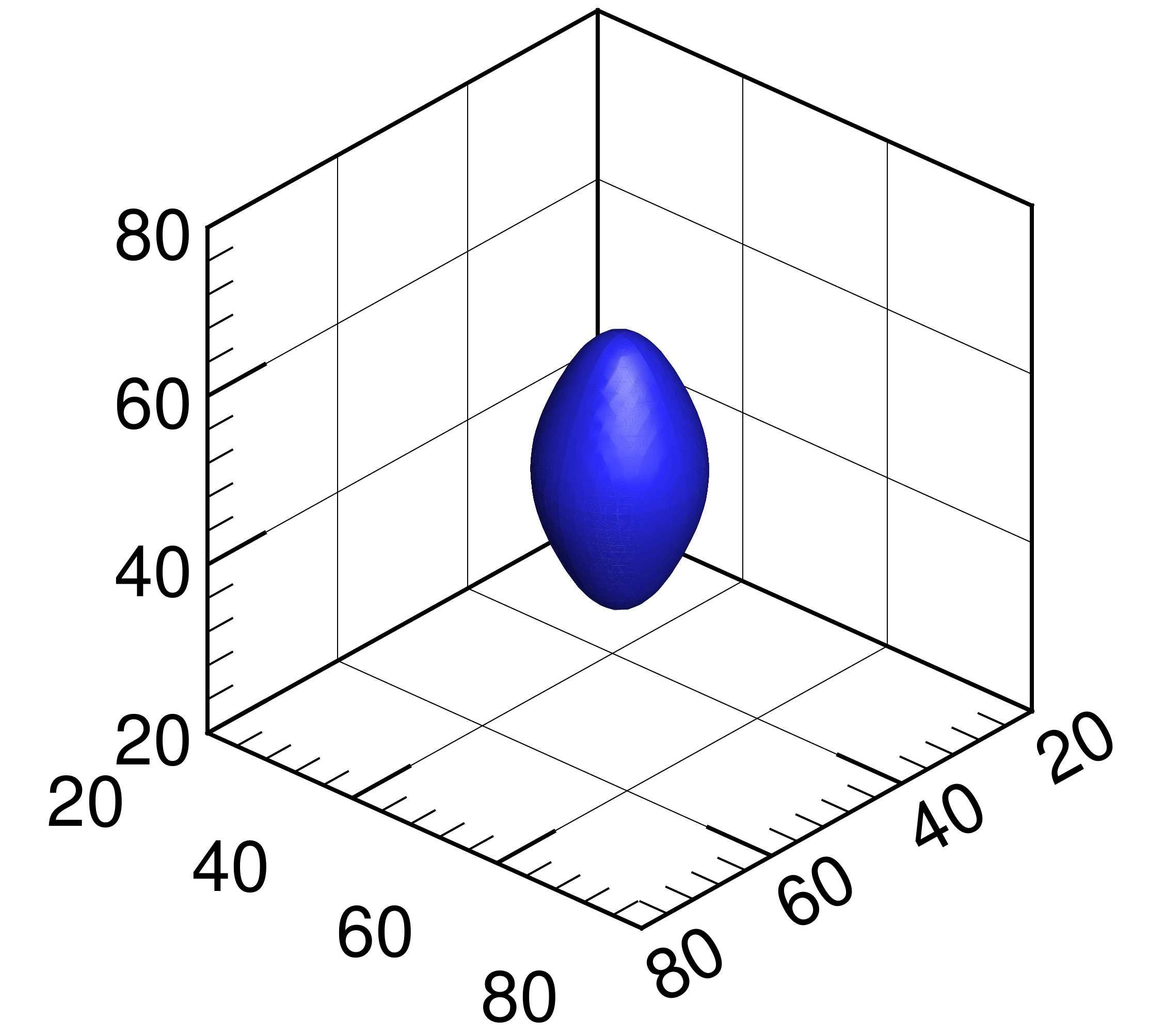}
\end{figure}

\subsection{Constrained curvature-driven motion}
In the results that follow, the constrained version of the velocity field is utilized. We investigate the interface evolution of different geometries under constrained curvature-driven motion whereby volume conservation is enforced. Each problem builds on the complexity of the problem that precedes it by adding an additional challenge via geometric features. \\
Not that volume conservation in the context of the proposed VOF approach should not be confused with volume conservation in VOF methods in general. The nature of the imposed conservation is purely geometric and not related to any alternation at the level of the flux computation algorithm. In other words, the reported errors in volume conservation quantify deviation from an imposed penalty on the volume of the entire initialized geometry.

\subsubsection{Ellipsoid}\label{sec:ell}
\noindent Consider the constrained curvature-driven motion of an ellipsoid whose zeroth level curve is defined as 
\begin{equation}
\phi(\mathbf{x},0) = \frac{(x-x_c)^2}{a^2} + \frac{(y-y_c)^2}{b^2} + \frac{(z-z_c)^2}{c^2} - 1\quad ,
\end{equation}
where $a$, $b$, and $c$ are the lengths of the semi-axes whose values are taken to be $35.0$, $15.625$, and $15.625$, respectively. This geometry represents a straightforward initial test of the methodology due to the smoothness of the interface and the absence of junctions. The computational domain is a unit cube and two resolutions were chosen for comparison, $50\times 50\times 50$ and $100\times 100\times 100$. The time-step is $\Delta t=10^{-5}$ and the total number of iterations for both grids is 10000 $(t_{final}=0.1)$. Fig. \ref{ellipsoid} shows snapshots of interface evolution for $t=0.0, 0.01, 0.02, 0.03, 0.04$, and $0.05$.

\begin{figure}[H]
\centering
\includegraphics[width=.32\textwidth]{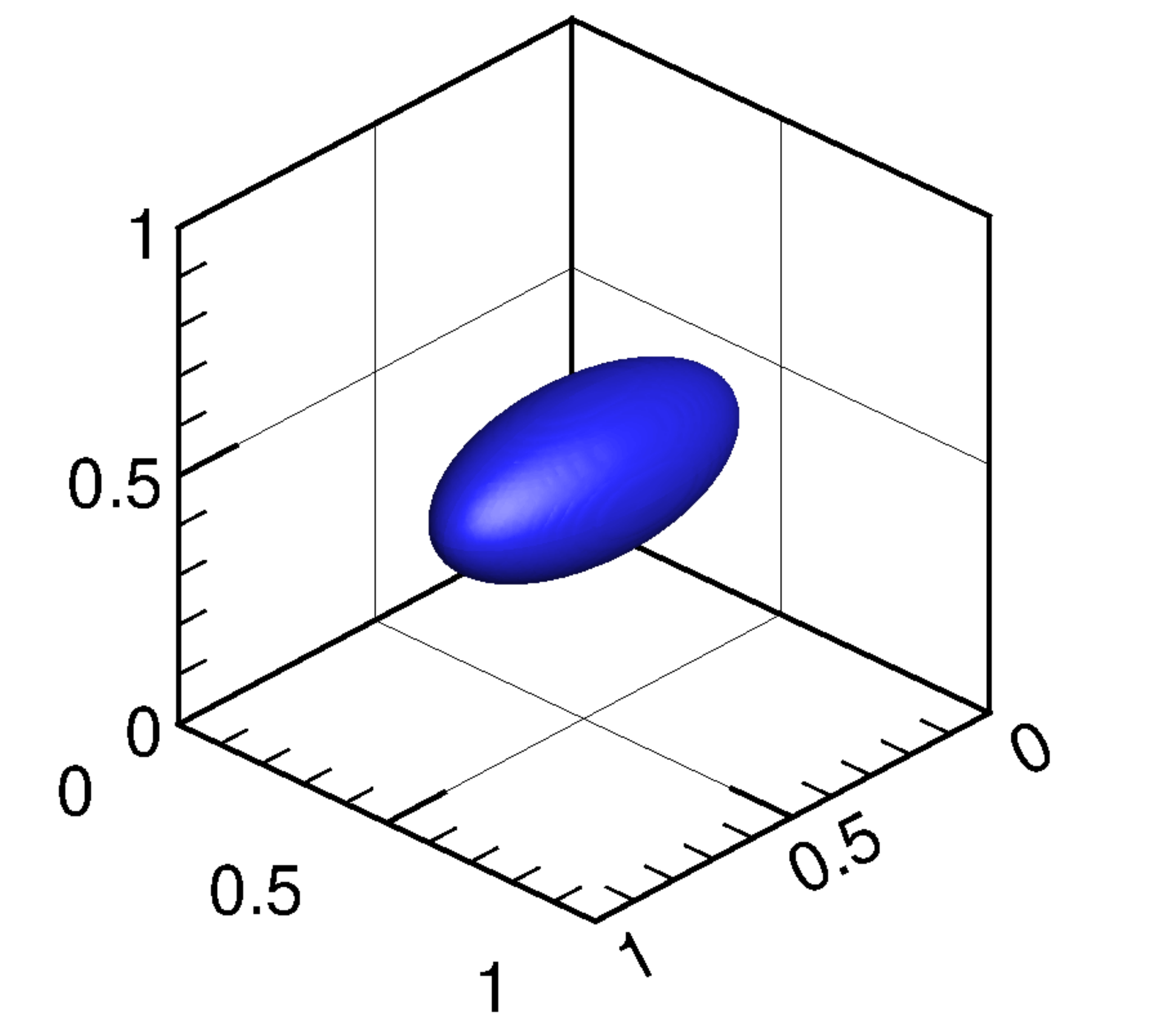}
\includegraphics[width=.32\textwidth]{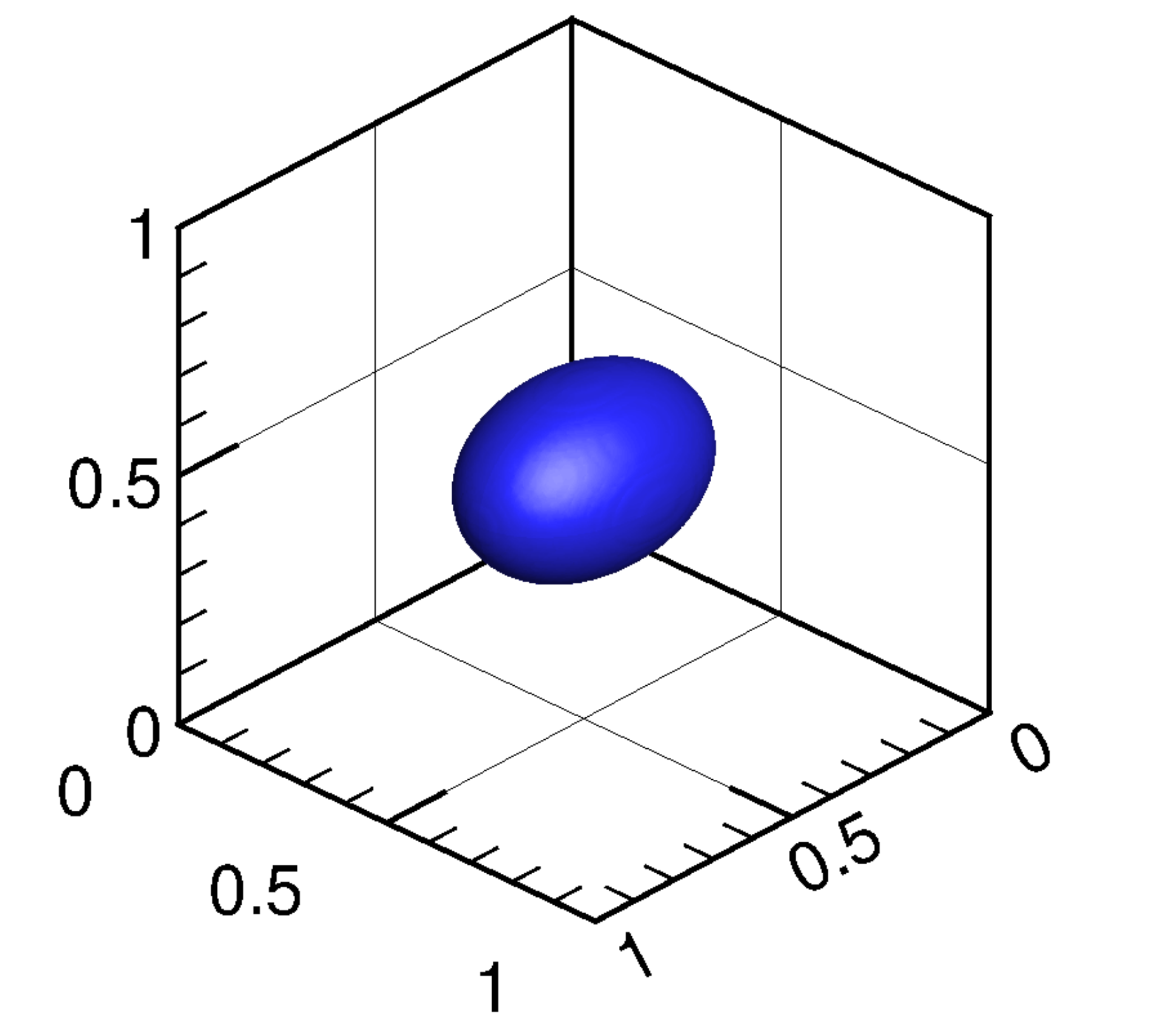} 
\includegraphics[width=.32\textwidth]{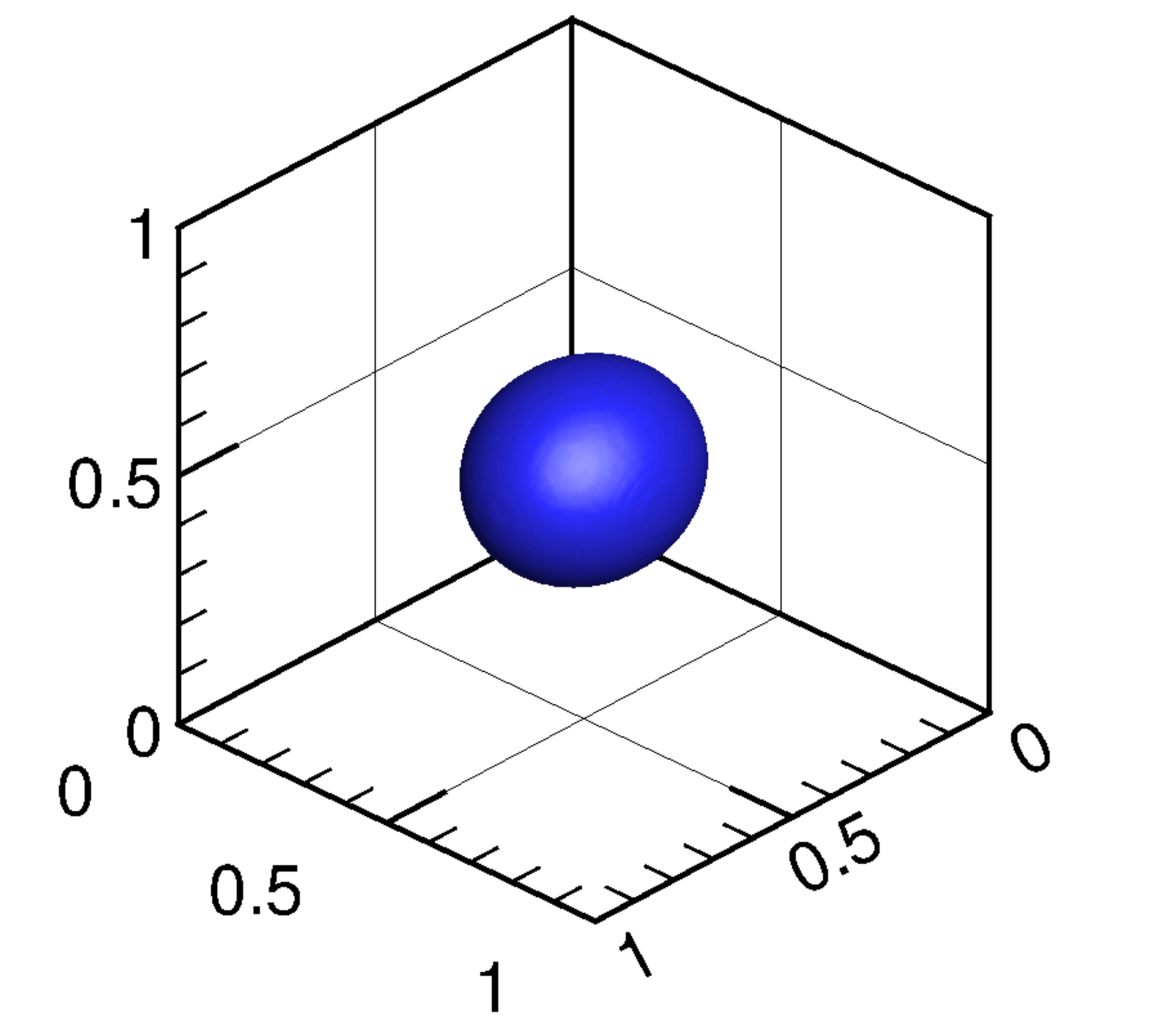} \\
\includegraphics[width=.32\textwidth]{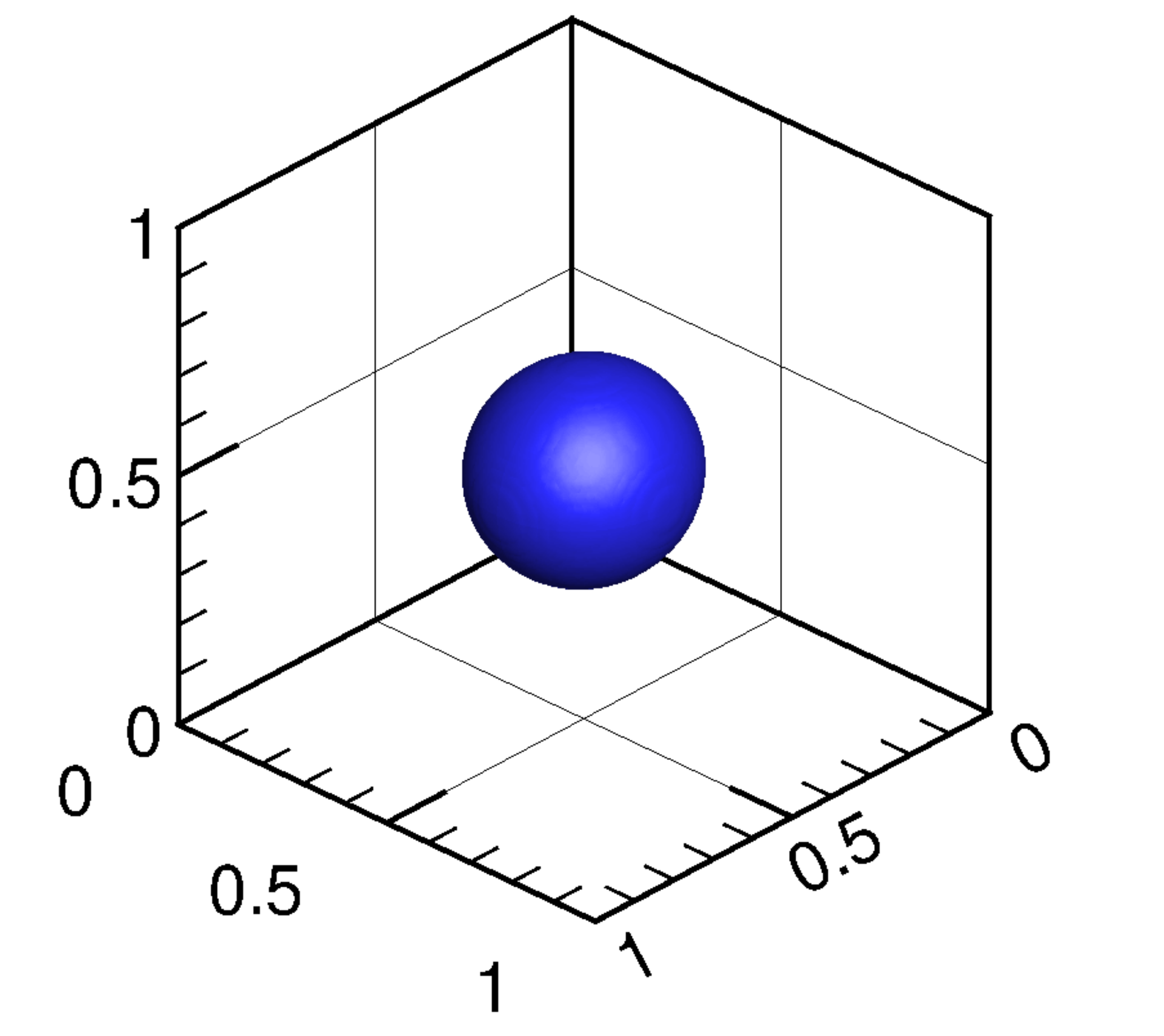}
\includegraphics[width=.32\textwidth]{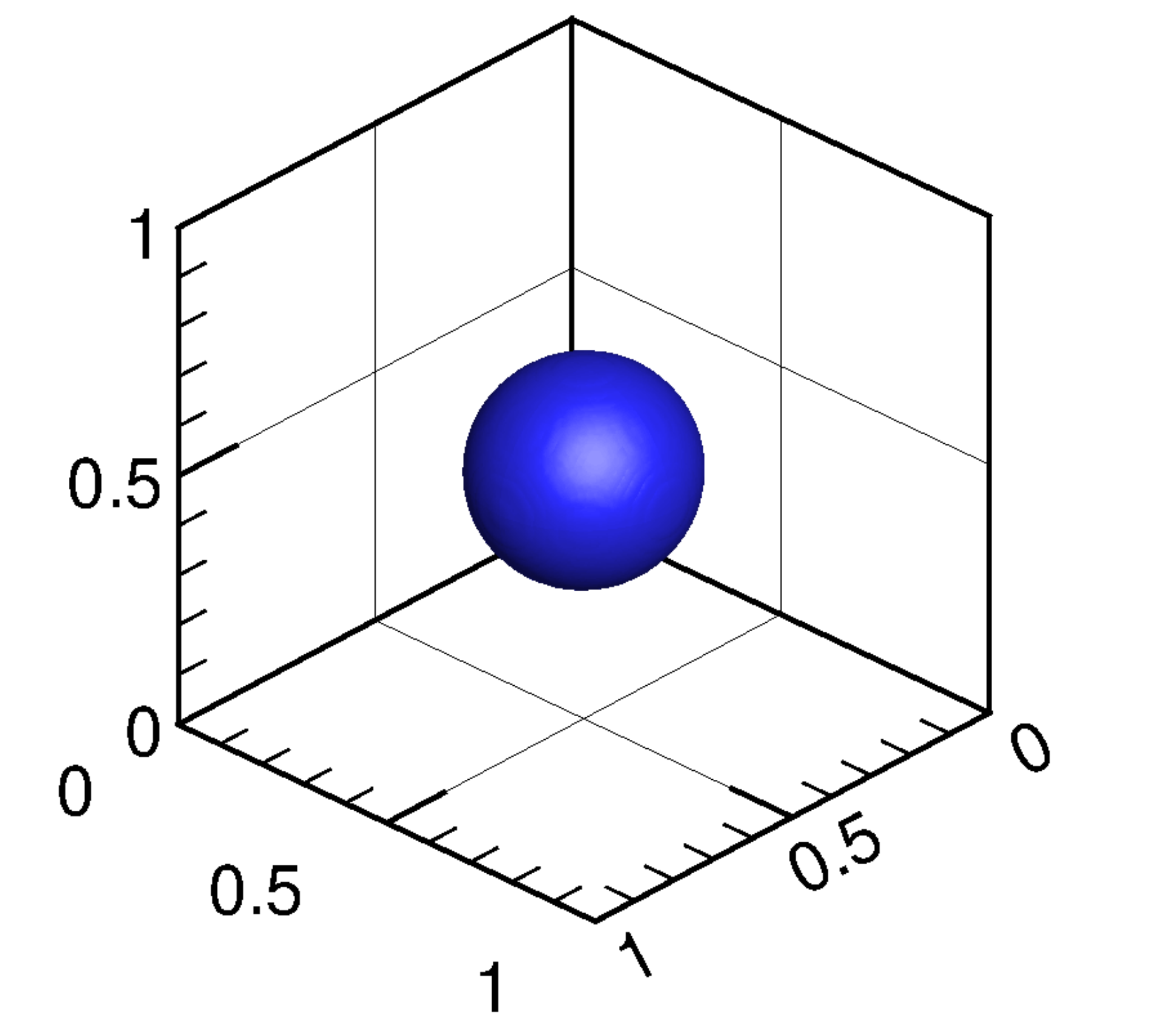} 
\includegraphics[width=.32\textwidth]{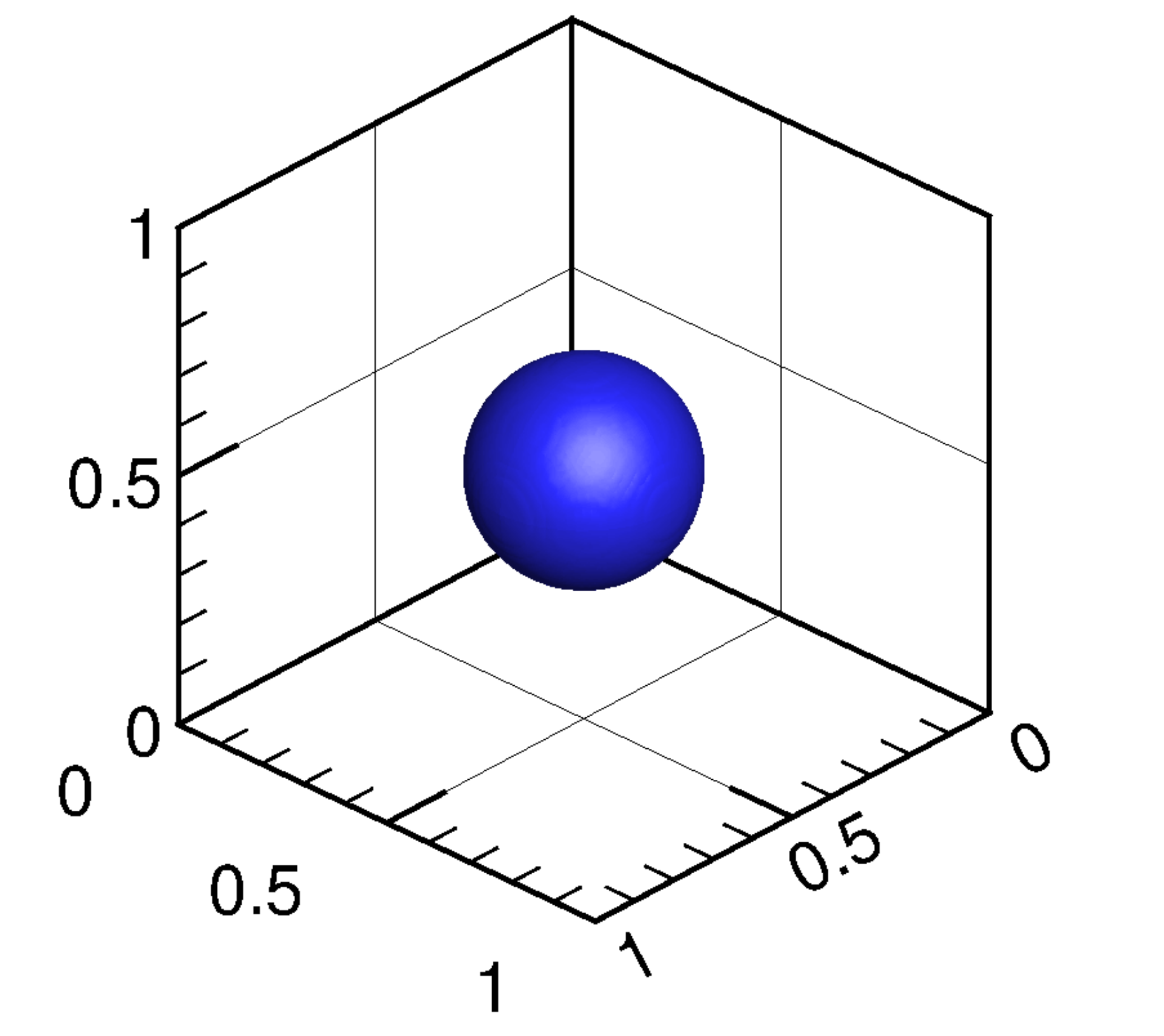} \\
\caption{Evolution of an ellipsoid under constrained curvature-driven flow. The images are snapshots of the interface in time from left to right, top to bottom $t = 0.0,0.01,0.02,0.03,0.04$ and $0.05$.}
\label{ellipsoid}
\end{figure}
As expected, the ellipsoid evolves into a sphere with time. When a spherical shape is reached, the volume no longer varies due to the balance between local curvature and mean curvature of the surface. The simulations were run 5000 iterations beyond the equilibrium point to ensure that built-up round-off errors will not cause the solution to drastically deviate from equilibrium or cause the rupture of the interface at any point in time. 
\begin{figure}[H]
    \centering
    \includegraphics[width=.45\textwidth]{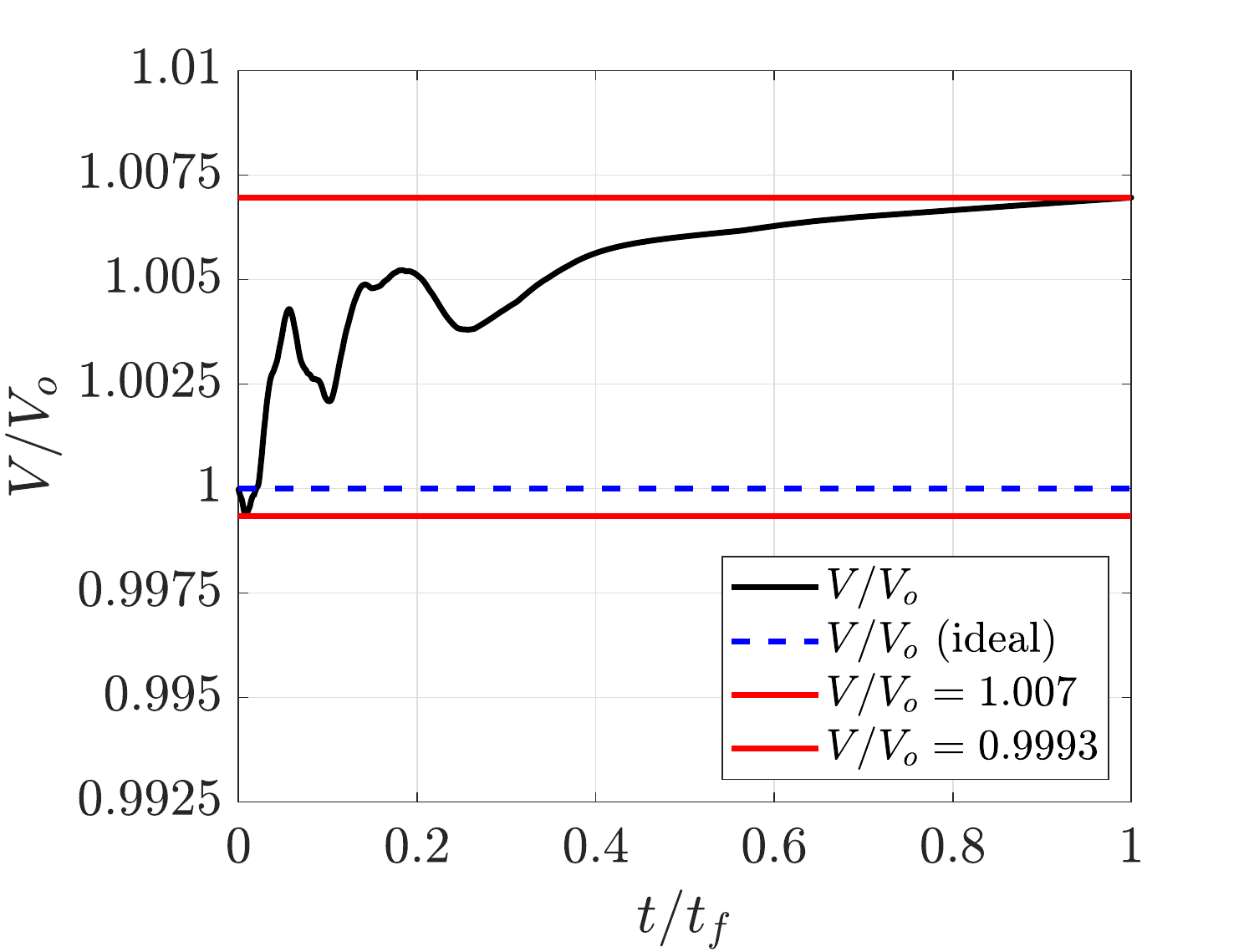}
    \includegraphics[width=.45\textwidth]{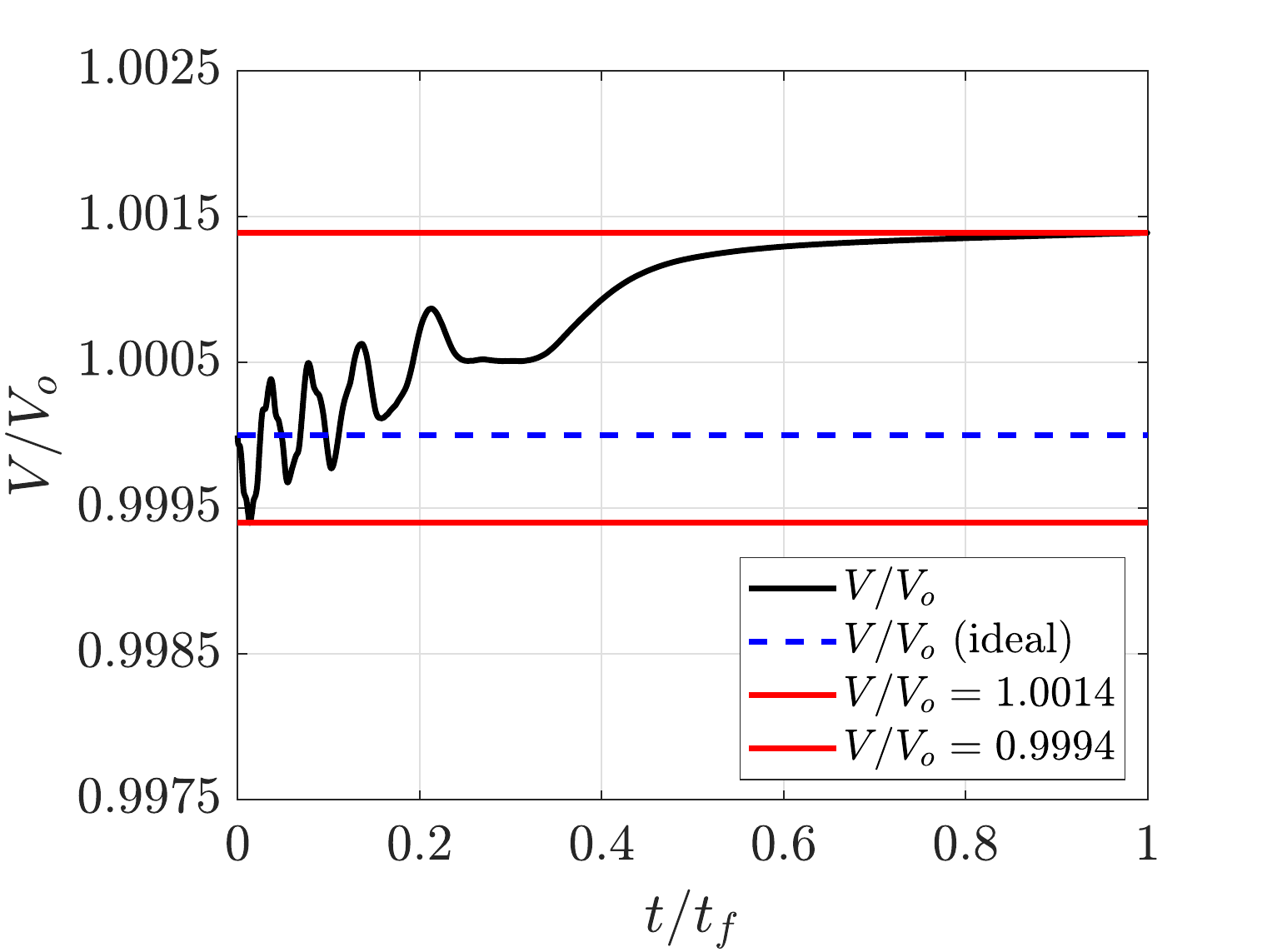}
    \caption{Variation of normalized volume of an ellipsoid with time under constrained curvature-driven motion.}
    \label{ellipsoid_conservation}
\end{figure}
Fig. \ref{ellipsoid_conservation} shows the variation of non-dimensionalized volume with respect to non-dimensional time for both the coarse grid (left) and finer grid (right). Error in volume conservation is approximately $0.7\%$ for the $50\times 50\times 50$ grid and $0.14\%$ for the $100\times 100\times 100$ indicating that higher grid resolution reduces errors in constraint violation.

\subsubsection{Superellipsoid}\label{sec:squircle}
\noindent In a Cartesian coordinate system, the superellipsoid is defined by
\begin{equation}
\displaystyle{\bigg|\frac{x-x_c}{a}\bigg|^n + \bigg|\frac{y-y_c}{b}\bigg|^n + \bigg|\frac{z-z_c}{c}\bigg|^n =1} 
\end{equation}
where $a$, $b$, and $c$ are the lengths of the semi-axes, $x_c$, $y_c$ and $z_c$ are the coordinates of the center of the ellipsoid, and $n$ is a positive integer. The choice of $n=4$ and equal semi-axes leads to a Squircle -- a circle with rounded edges. As $n$ increases, the resulting geometry attains sharper edges, this in turn helps in approximating a cubic shape without the singularities at the corners ($\kappa $ is not defined at a point). The idea behind simulating this geometry is to investigate the limitations of curvature-driven motion with VOF where $\kappa$ is not properly defined. \par
Consider the constrained curvature-driven motion of a round-edge cube whose zeroth level curve is represented by
\begin{equation}
\phi(\mathbf{x},0) = \bigg(\frac{x-x_c}{r}\bigg)^n + \bigg(\frac{y-y_c}{r}\bigg)^n + \bigg(\frac{z-z_c}{r}\bigg)^n - 1\quad ,
\end{equation}
where $n=12$, $r=0.25$, and $\mathbf{x_c}$ is $(0.5,0.5,0.5)$. The computational domain is a unit cube and the grid resolutions used are $50\times 50\times 50$ and $100\times 100\times 100$, respectively. The time-step chosen for both grids is $\Delta t=10^{-5}$ and the total number of iterations is 10000 ($t_{final}=0.1$). Fig. \ref{squircle} shows snapshots of interface evolution for $t=0.0, 0.002, 0.004, 0.006, 0.012$, and $0.014$.

\begin{figure}[H]
\centering
\includegraphics[width=.32\textwidth]{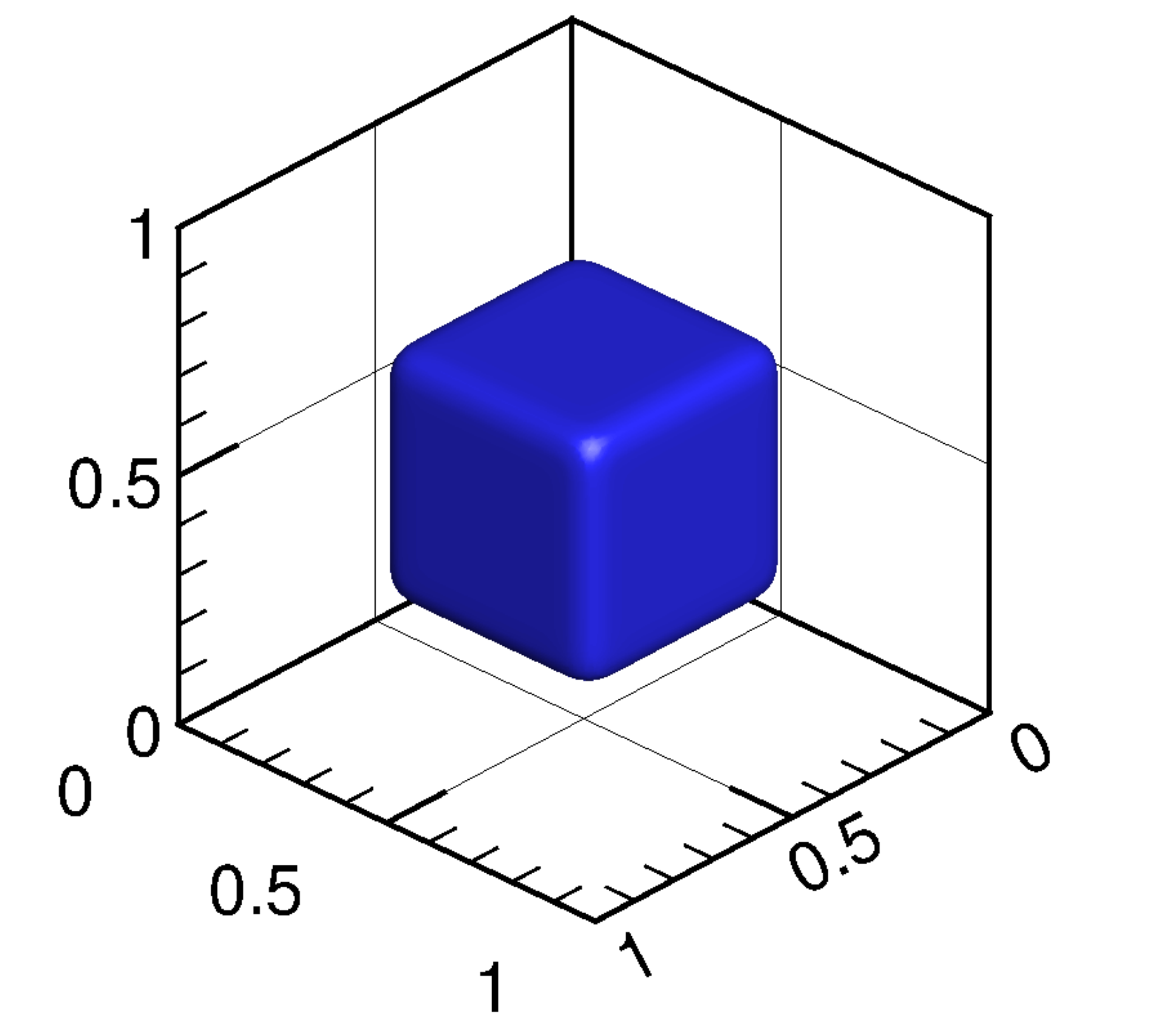}
\includegraphics[width=.32\textwidth]{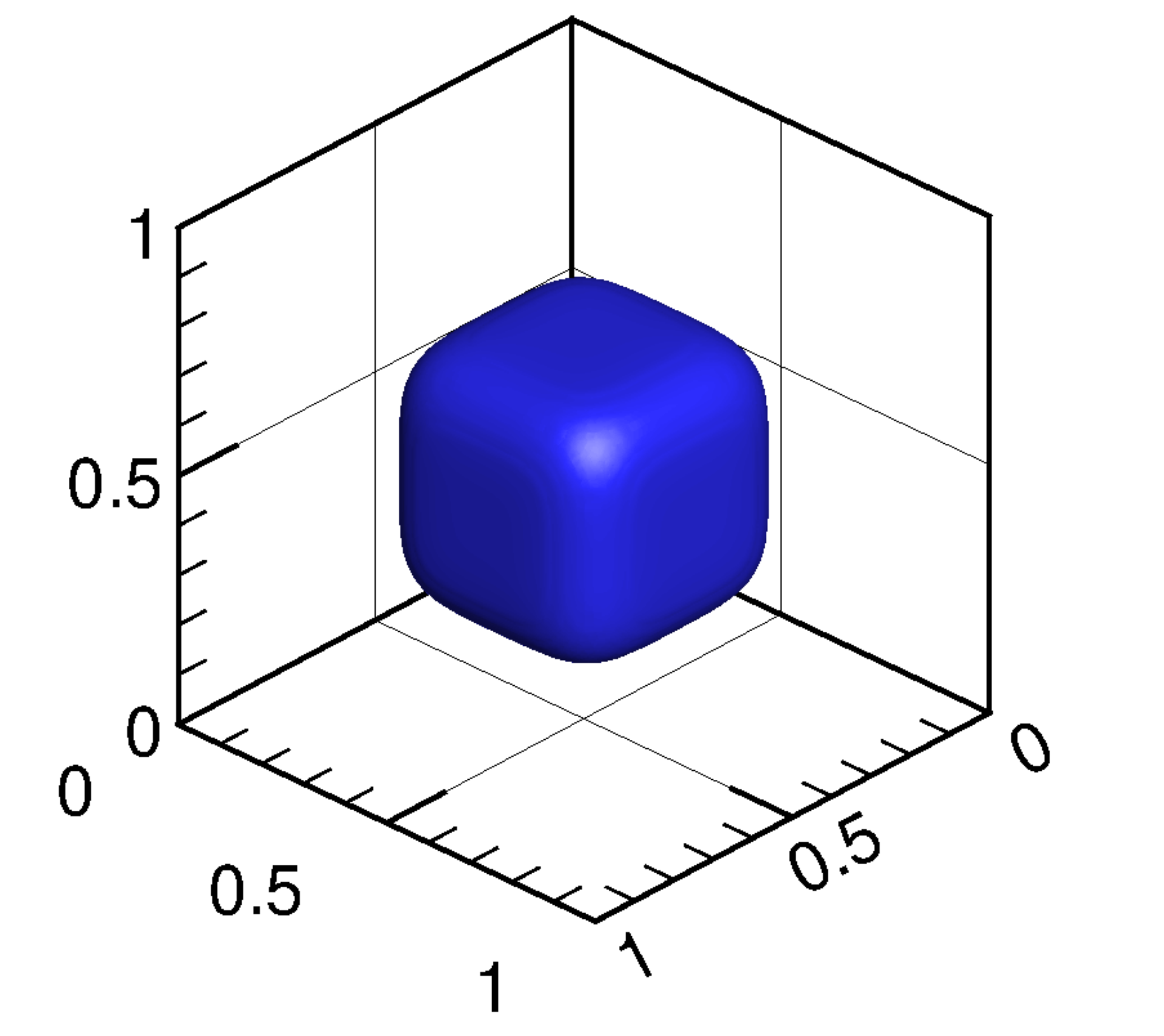} 
\includegraphics[width=.32\textwidth]{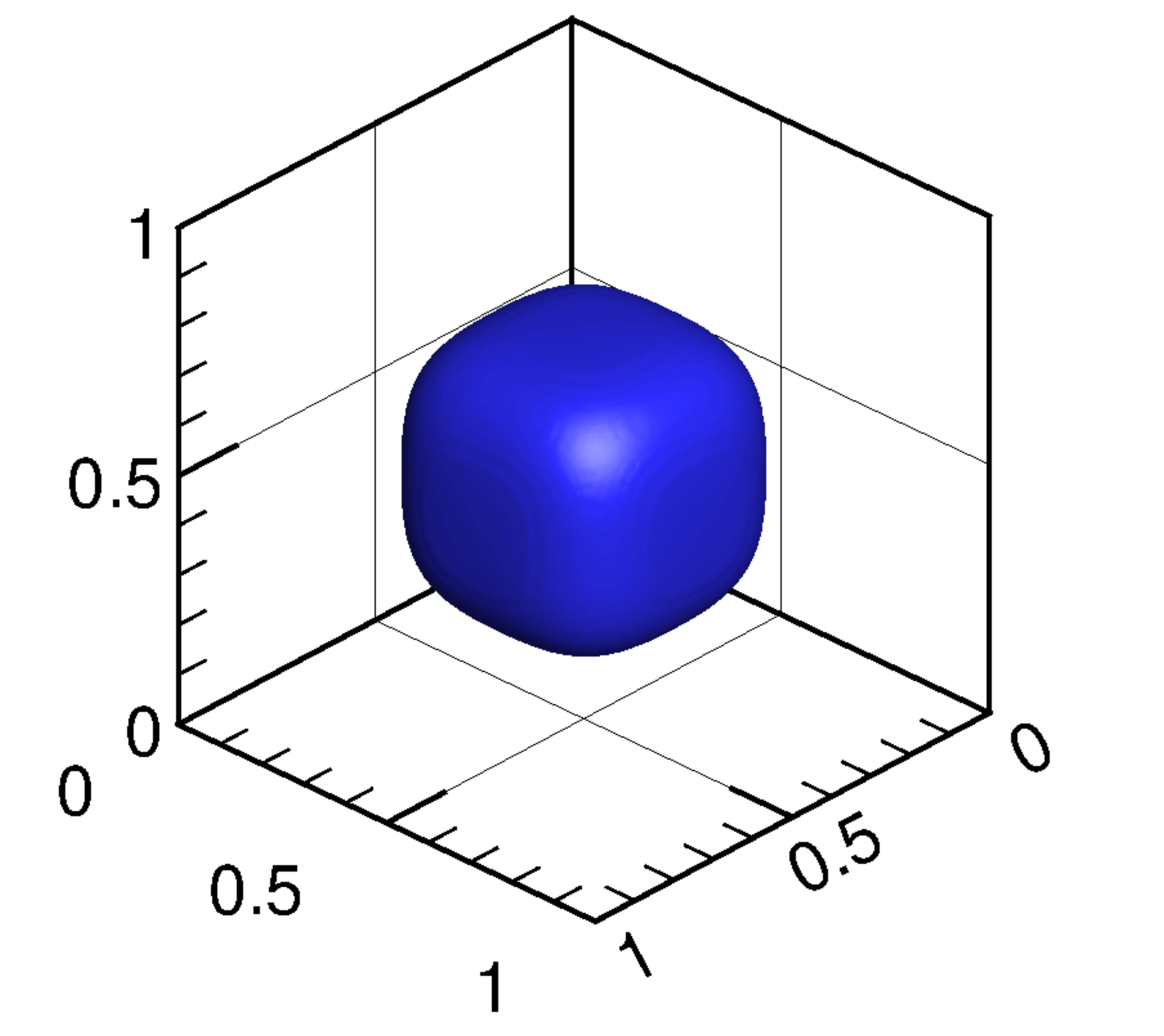} \\
\includegraphics[width=.32\textwidth]{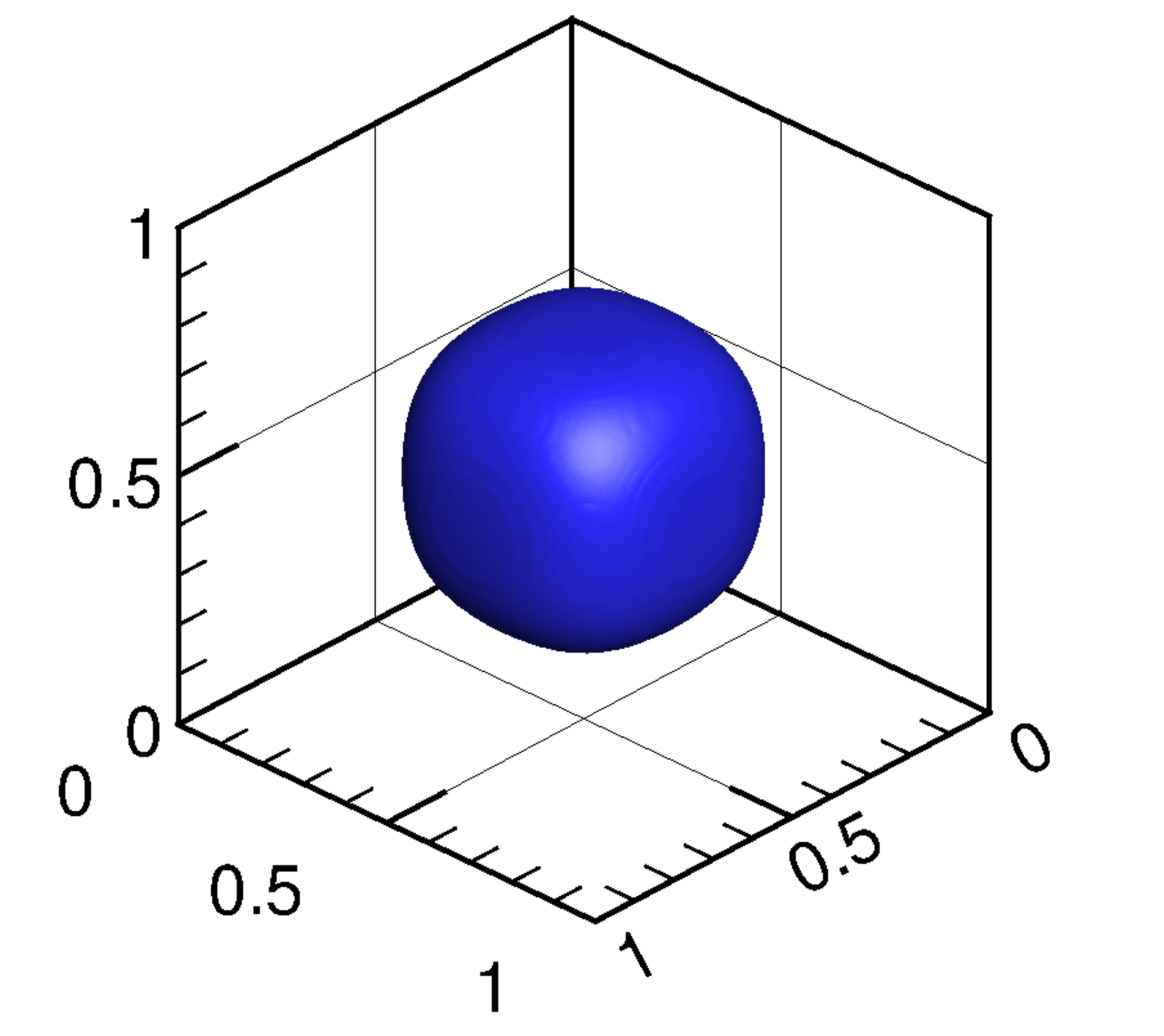}
\includegraphics[width=.32\textwidth]{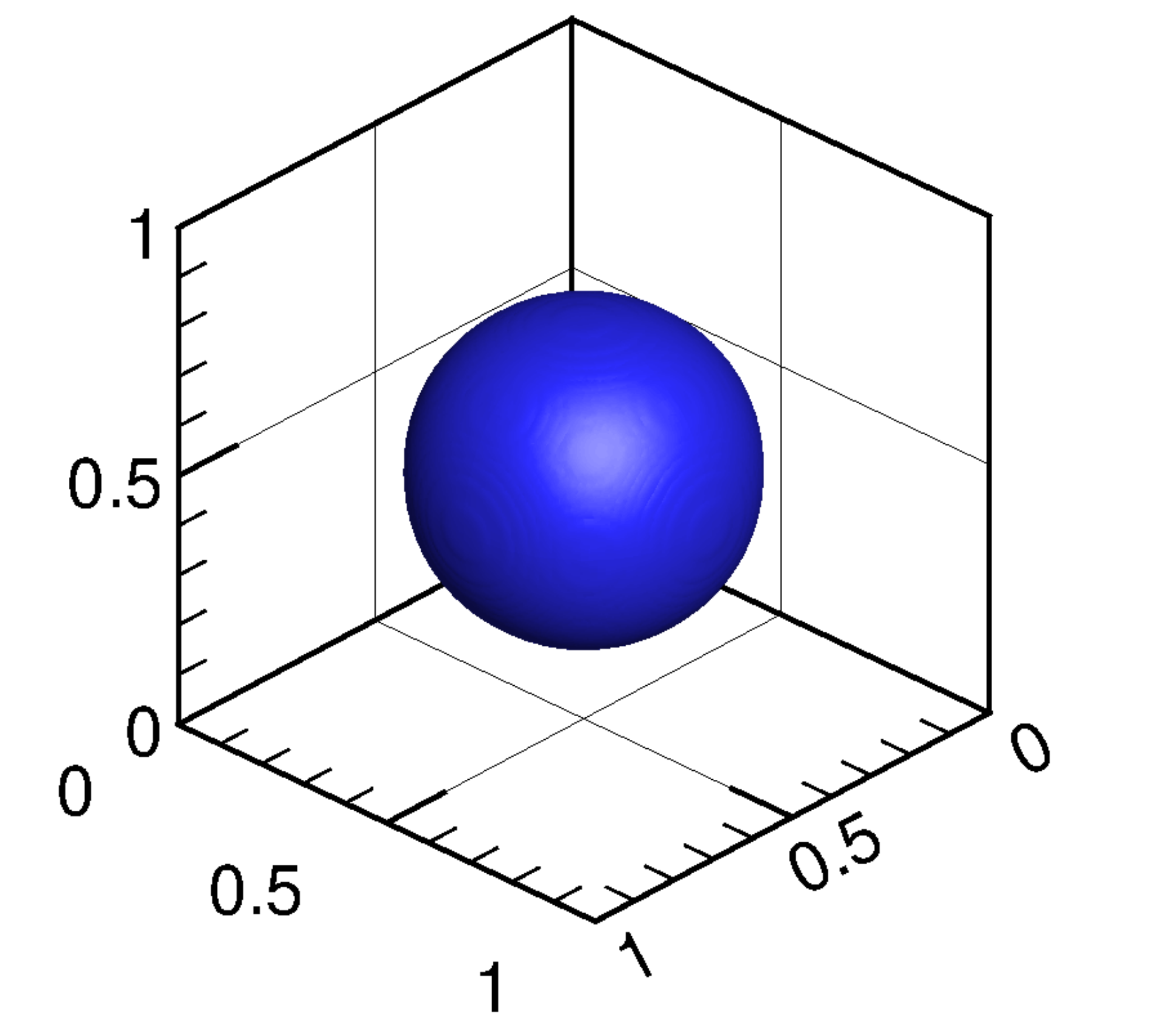} 
\includegraphics[width=.32\textwidth]{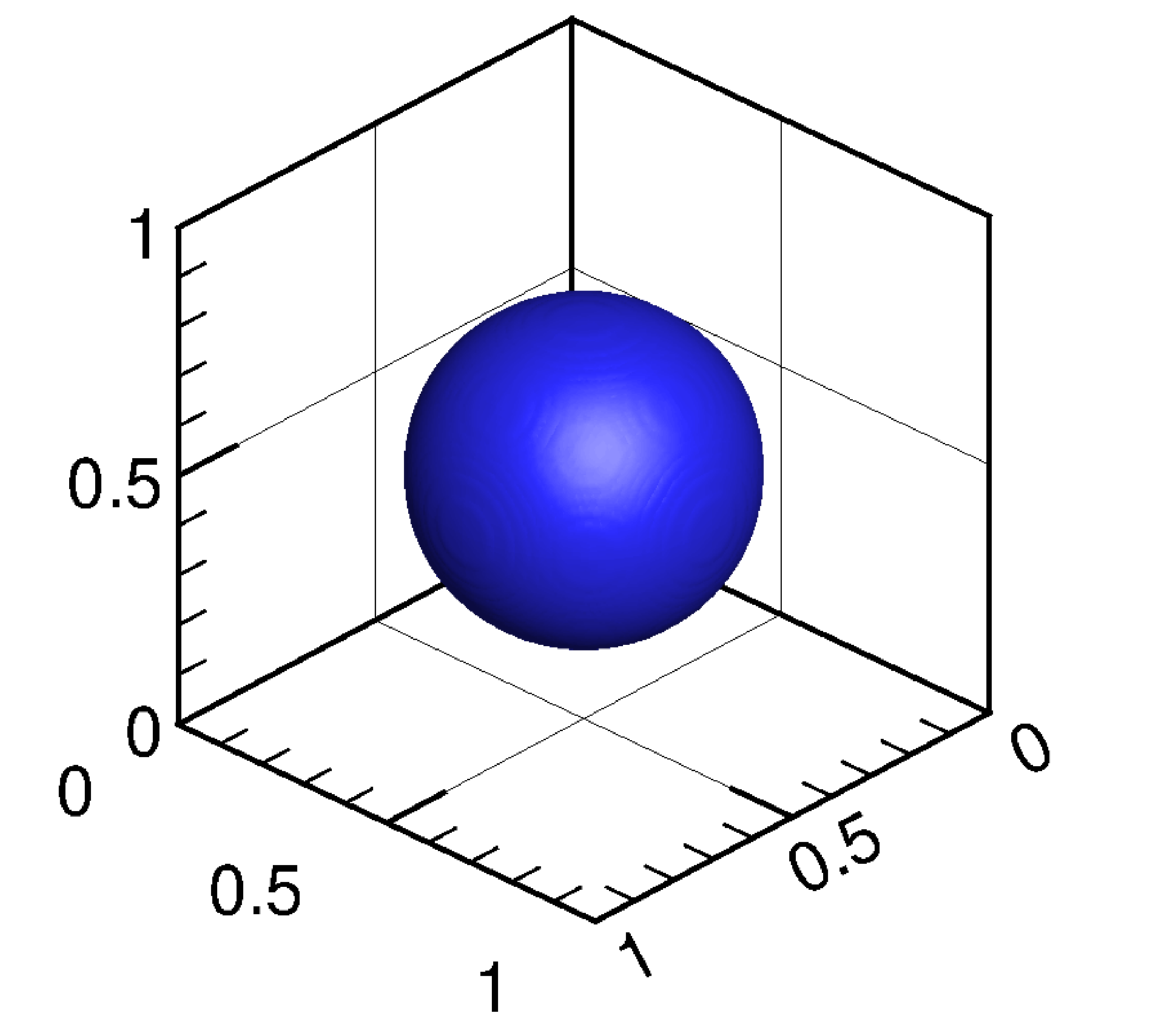} \\
\caption{Evolution of a round-edge cube under constrained curvature-driven flow. The images are snapshots of the interface in time from left to right, top to bottom where the cubic geometry evolves into a sphere for $t=0.0,0.002,0.004,0.006,0.012$ and $0.014$.}
\label{squircle}
\end{figure}

What is of particular importance in this problem, is the calculation of $\bar{\kappa}$. Through the numerical experiments, we have found that the definition of $\delta(C)$ directly affects the numerical stability of the method. Upon calculating $\delta(C)$ as $|\nabla C|$, the value of $\bar{\kappa}$ would fluctuate between consecutive time steps which eventually lead to parasitic velocities that ruptured the interface. This necessitated an alternative definition for the Dirac delta function that is bounded and smooth. A definition that satisfies the preceding conditions is $\displaystyle{\delta (C) = 4C(1-C)}$.

\begin{figure}[H]
\centering
\begin{subfigure}{0.47\textwidth}
    \includegraphics[width=\textwidth]{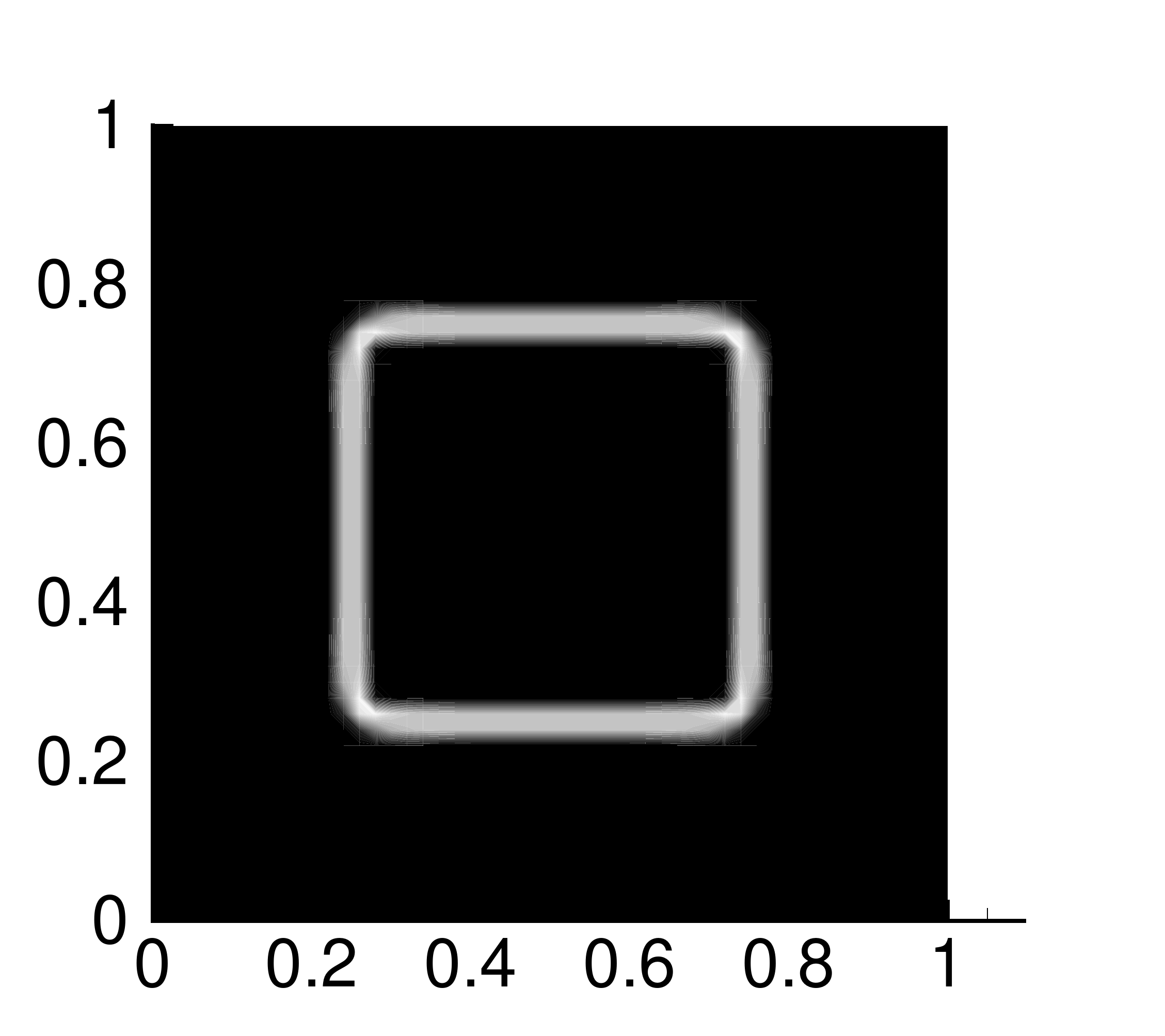}
    \caption{$\displaystyle{\delta_1(C)=4C(1-C)}$}
    \label{fig:first}
\end{subfigure}
\quad
\begin{subfigure}{0.47\textwidth}
    \includegraphics[width=\textwidth]{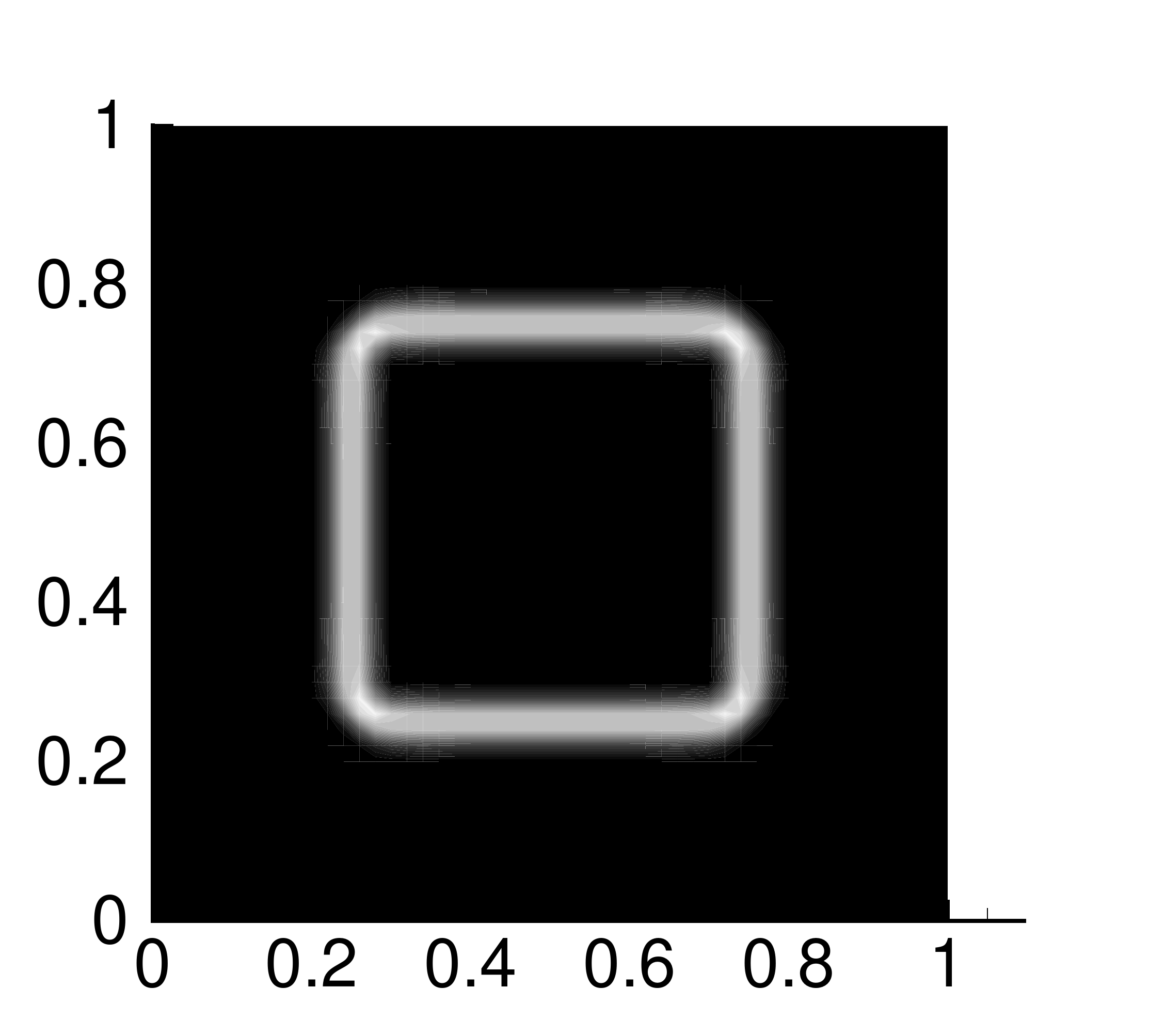}
    \caption{$\delta_2(C)=|\nabla C|$}
    \label{fig:second}
\end{subfigure}
        
\caption{Comparison of different definitions of the Dirac delta function.}
\label{dirac_comp}
\end{figure}

Fig. \ref{dirac_comp} presents a comparison between the two definitions of $\delta(C)$. It is evident that using the definition in Fig. \ref{fig:second} is not ideal, fluctuations are visible at the corners and the band defining the interface is diffuse. In contrast, the interface region in Fig. \ref{fig:first} is sharper and spurious features at the edges are less pronounced. Hence, $\delta_1(C)$ will be used throughout the rest of the paper.

\begin{figure}[H]
    \centering
    \includegraphics[width=.45\textwidth]{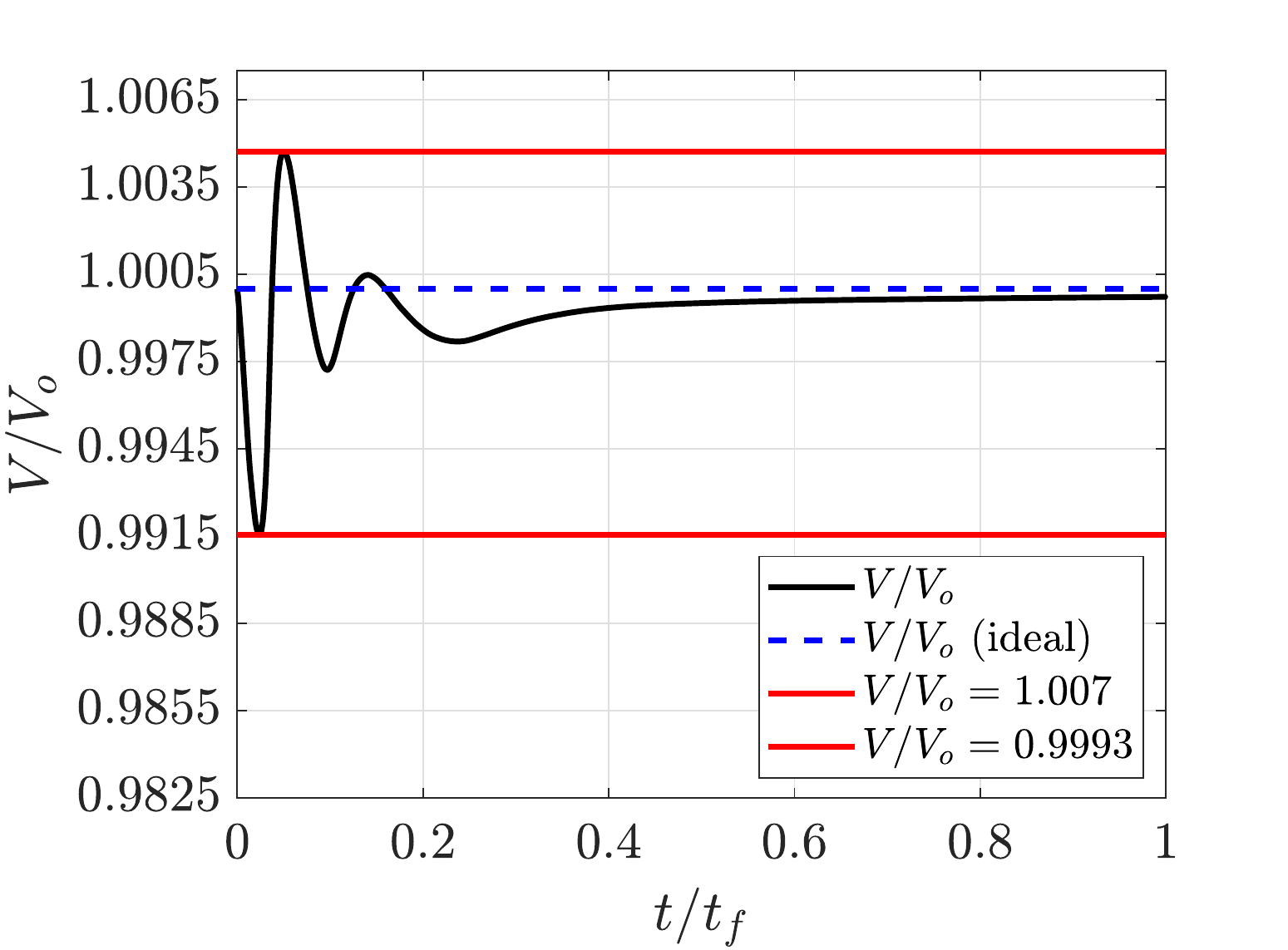}
    \includegraphics[width=.45\textwidth]{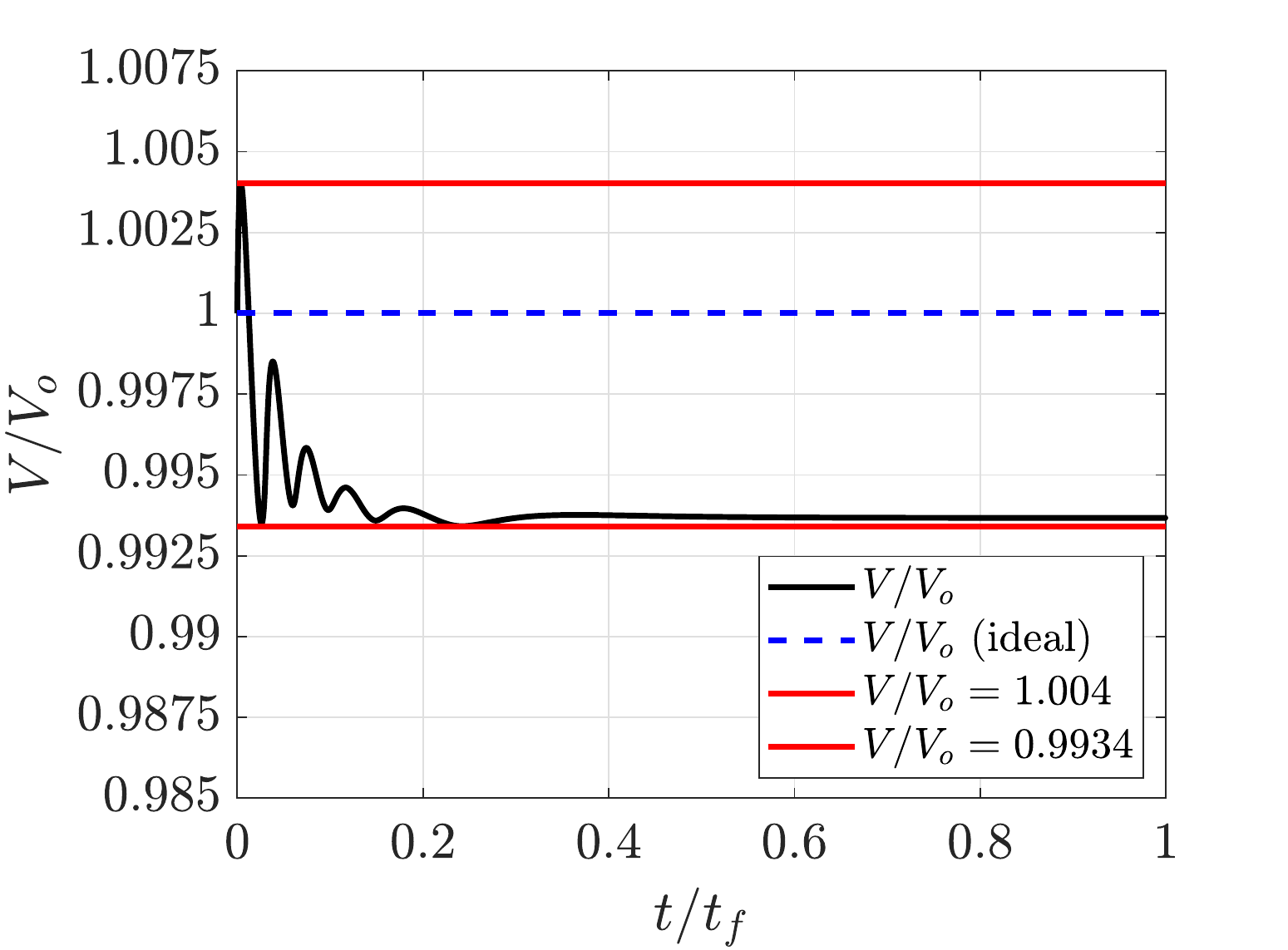}
    \caption{Variation of the normalized volume of a squircle with time under constrained curvature-driven motion.}
    \label{squircle_conservation}
\end{figure}
Fig. \ref{squircle_conservation} shows the variation of non-dimensionalized volume with respect to non-dimensional time for both the coarse grid (left) and finer grid (right). Error in volume conservation is approximately $0.7\%$ for the $50\times 50\times 50$ grid and $0.4\%$ for the $100\times 100\times 100$ indicating that higher grid resolution reduces errors in constraint violation which is consistent with the results in Sec. \ref{sec:ell}

\subsubsection{Octahedron}
Consider an Octahedron -- a polyhedron with eight faces. This geometry involves symmetry, sharp edges, and faces oriented at different angles. There were no attempts to smooth the Octahedron edges since using $\delta_1(C)$ proved to be sufficient to guarantee accuracy and stability. The zeroth level curve of the Octahedron is defined as
\begin{equation}
\phi(\mathbf{x},0)=|x-x_c|+|y-y_c|+|z-z_c|-r
\end{equation}
where $r$ is the radius. The computational domain is a unit cube, $r=0.3$, $\mathbf{x_c}$ is $(0.5,0.5,0.5)$ and the two grid resolutions used are the same as the preceding sections. The time-step is $\Delta t=0.0001$ and $\Delta t=0.000025$ for grid sizes $50\times 50\times 50$ and $100\times 100\times 100$, respectively, and $t_{final}=0.03$.
\begin{figure}[H]
\centering
\includegraphics[width=.32\textwidth]{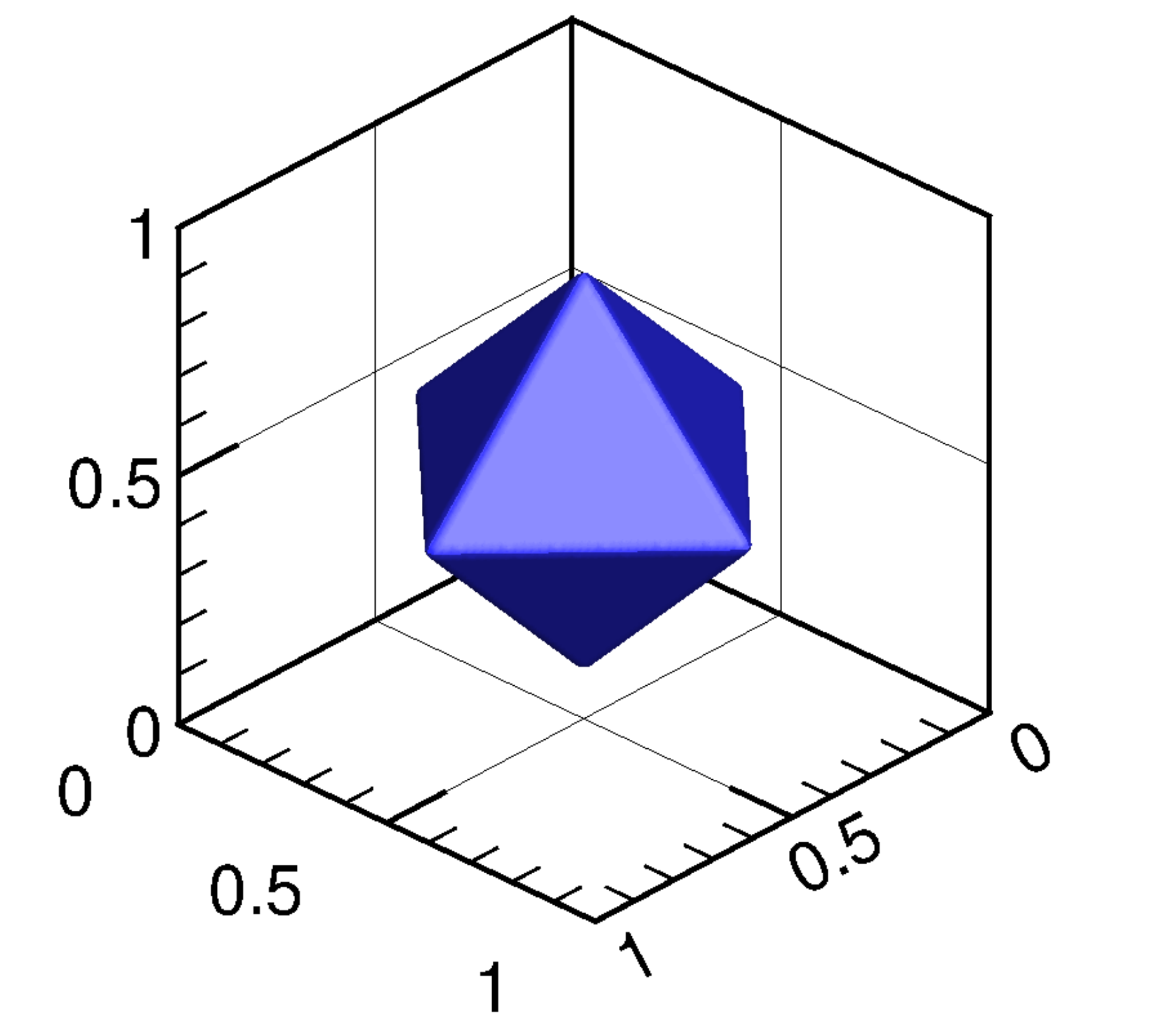}
\includegraphics[width=.32\textwidth]{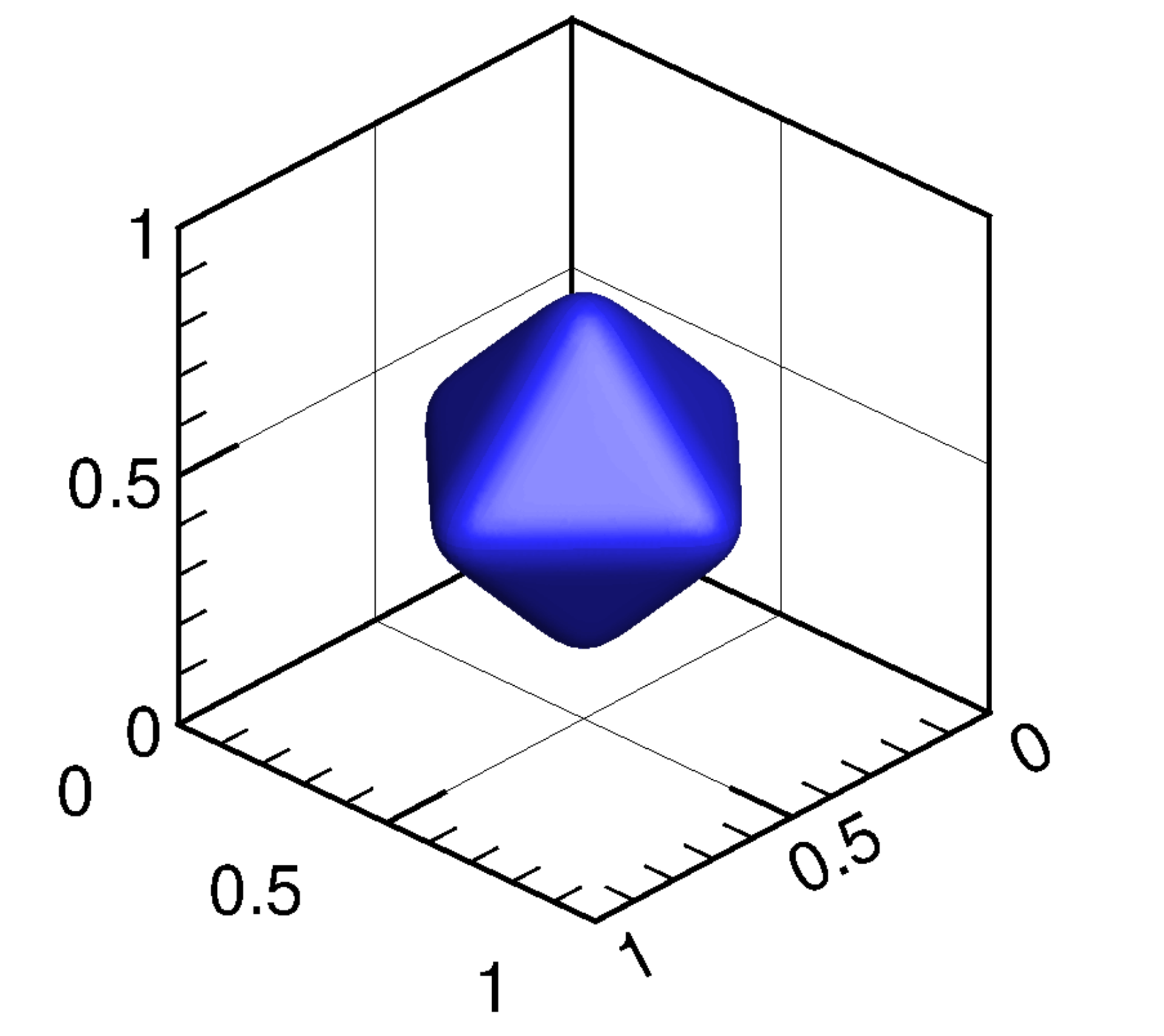} 
\includegraphics[width=.32\textwidth]{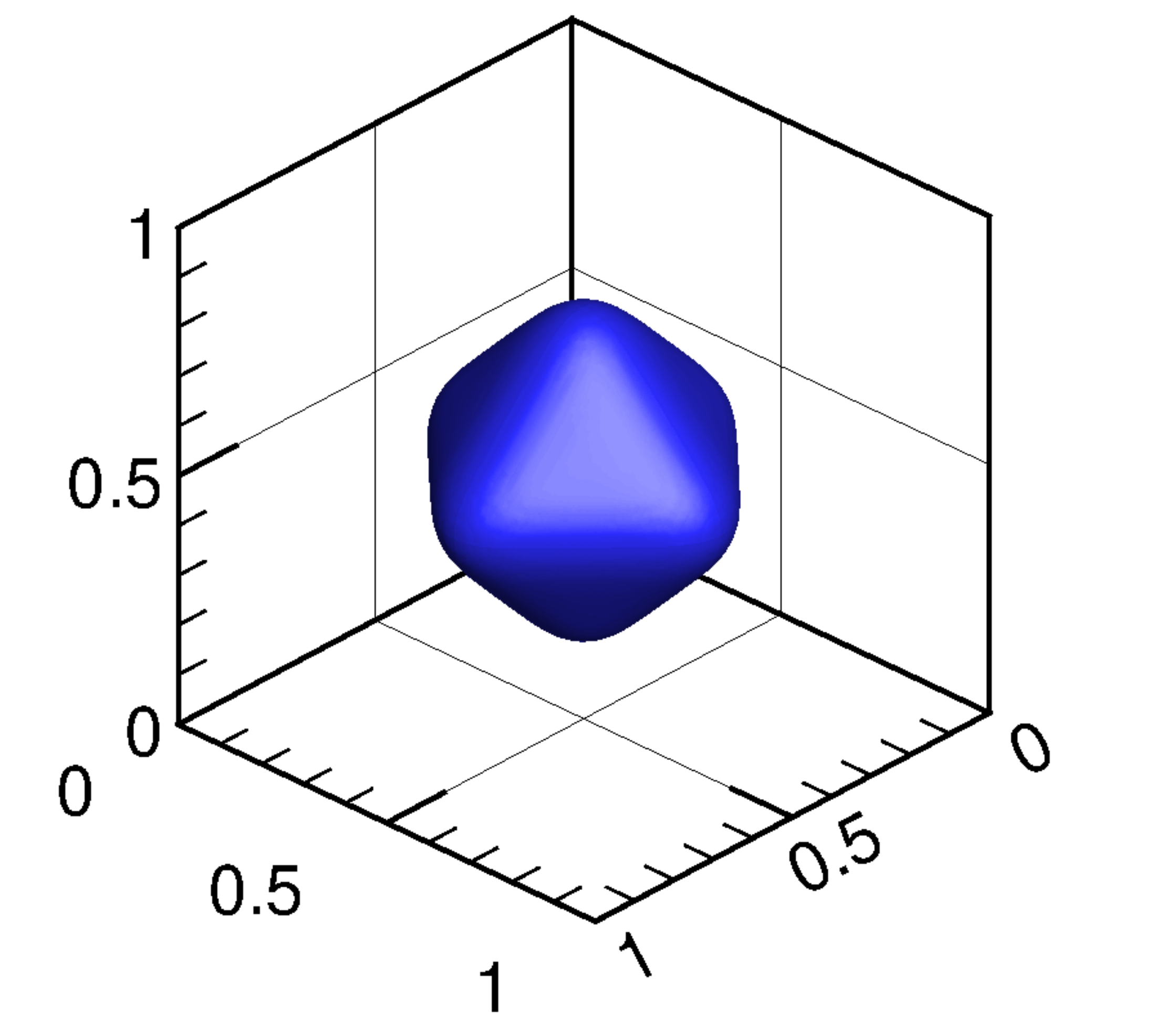} \\
\includegraphics[width=.32\textwidth]{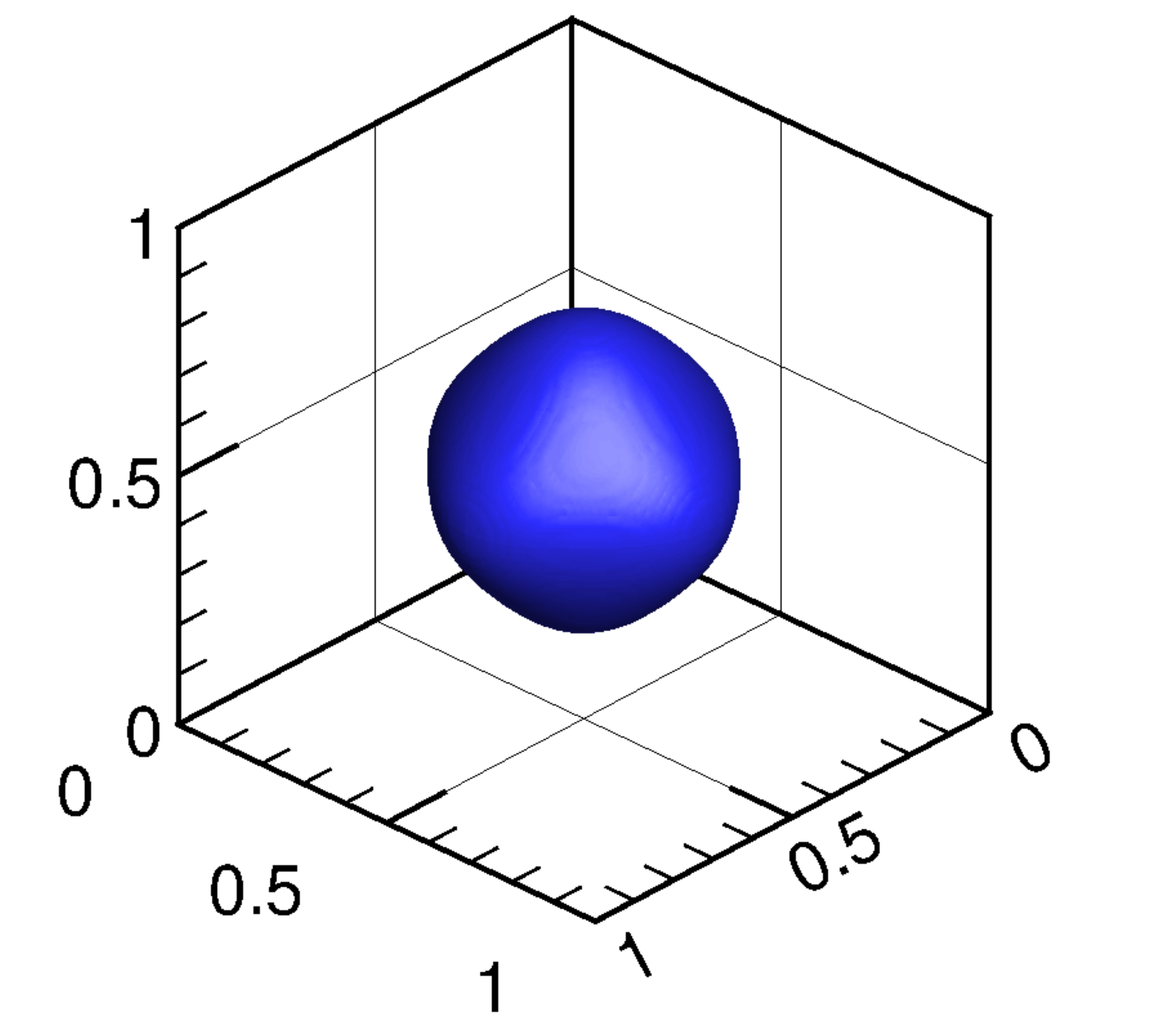}
\includegraphics[width=.32\textwidth]{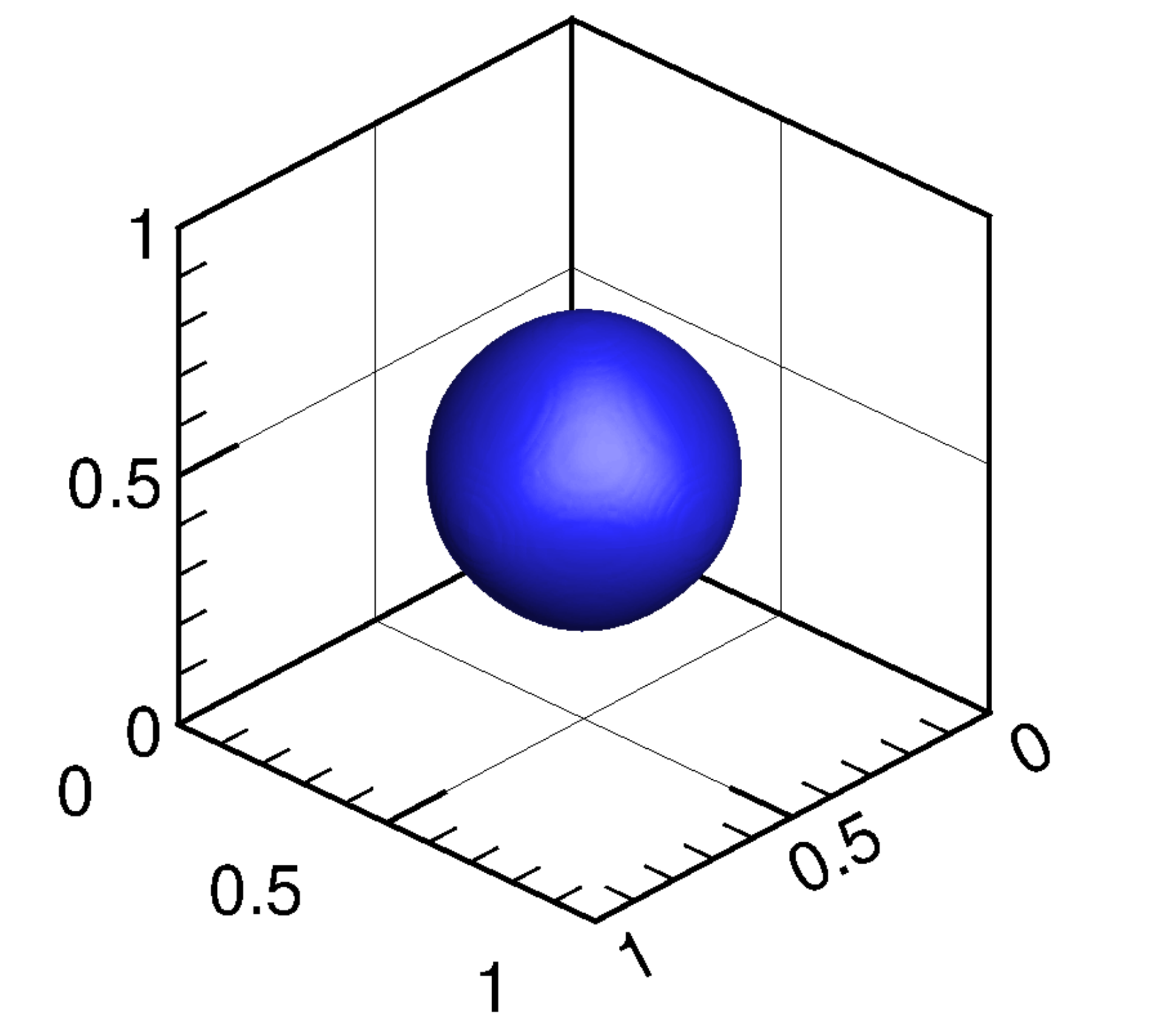} 
\includegraphics[width=.32\textwidth]{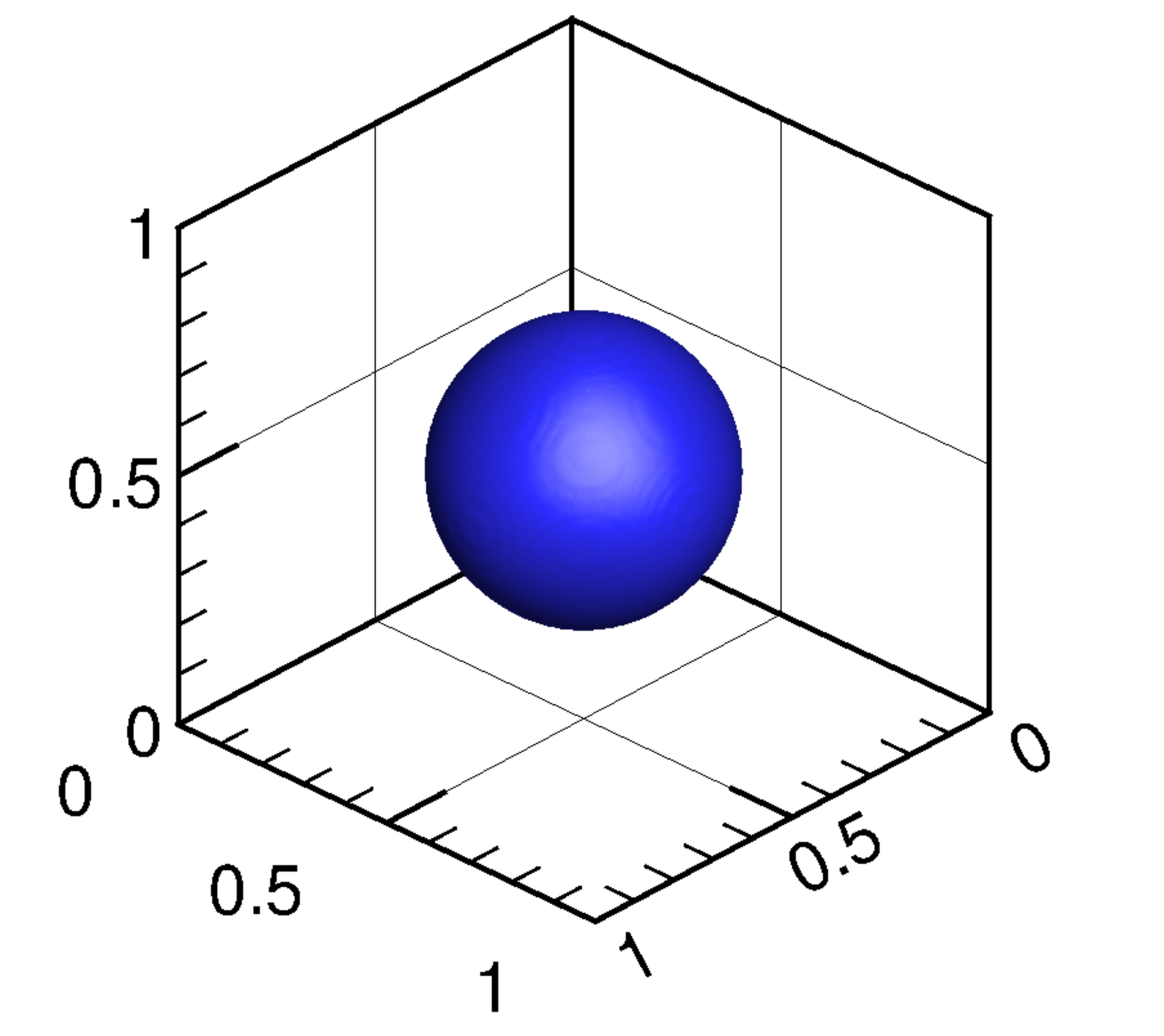} \\
\caption{Evolution of an Octahedron under constrained curvature-driven flow. The images are snapshots of the interface in time from left to right, top to bottom for  $t = 0.0, 0.001,0.002,0.005,0.007$ and $0.01$}
\label{squircle}
\end{figure}

\begin{figure}[H]
    \centering
    \includegraphics[width=.45\textwidth]{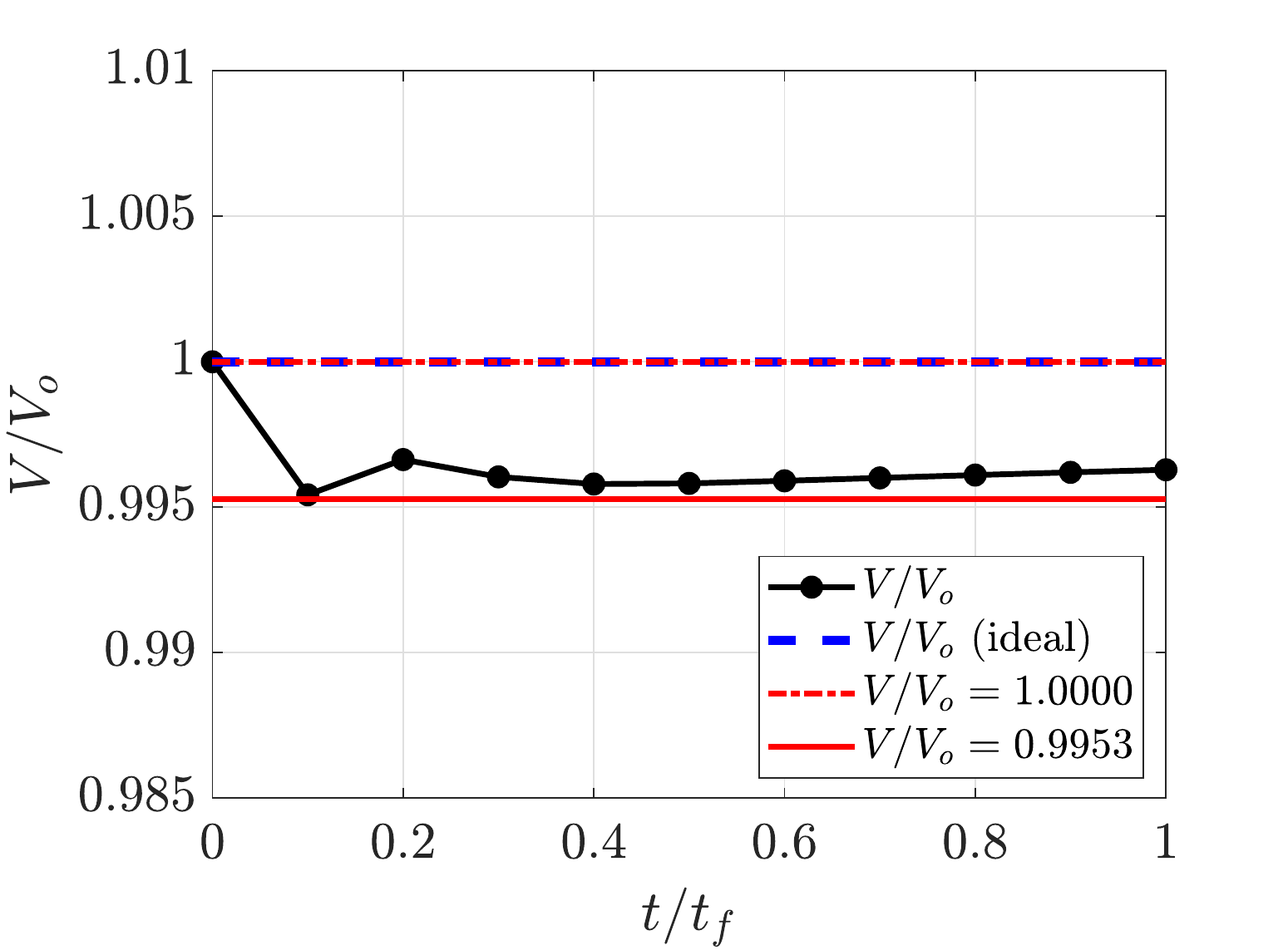}
    \includegraphics[width=.45\textwidth]{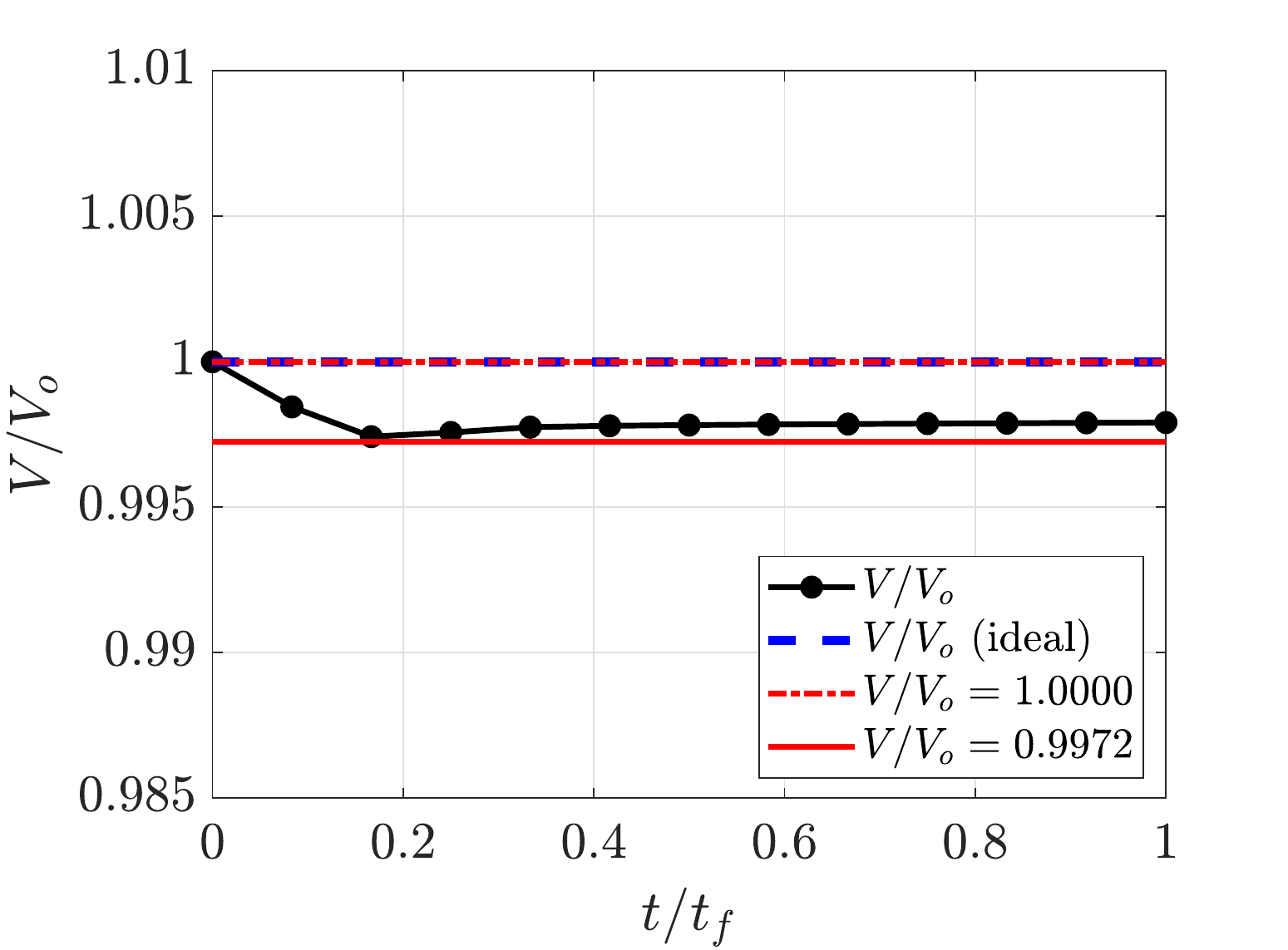}
    \caption{Variation of the normalized volume of an ellipsoid with time under constrained curvature-driven motion.}
    \label{octahedron_conservation}
\end{figure}
Fig. \ref{octahedron_conservation} shows the variation of non-dimensionalized volume with respect to non-dimensional time for both the coarse grid (left) and finer grid (right). Error in volume conservation is approximately $0.37\%$ for $50\times 50\times 50$ and $0.2\%$ for $100\times 100\times 100$. As the grid is refined the singularities at the face junctions of the octahedron become more pronounced, however, the method remains stable for the finer grid with a reduced mass error. Further refinement could possibly lead to ill-defined curvature at the corners; this can be remedied via smoothing as needed. 

\subsubsection{Sphere in a helical velocity field}
\begin{figure}[H]
\centering
\includegraphics[width=.33\textwidth]{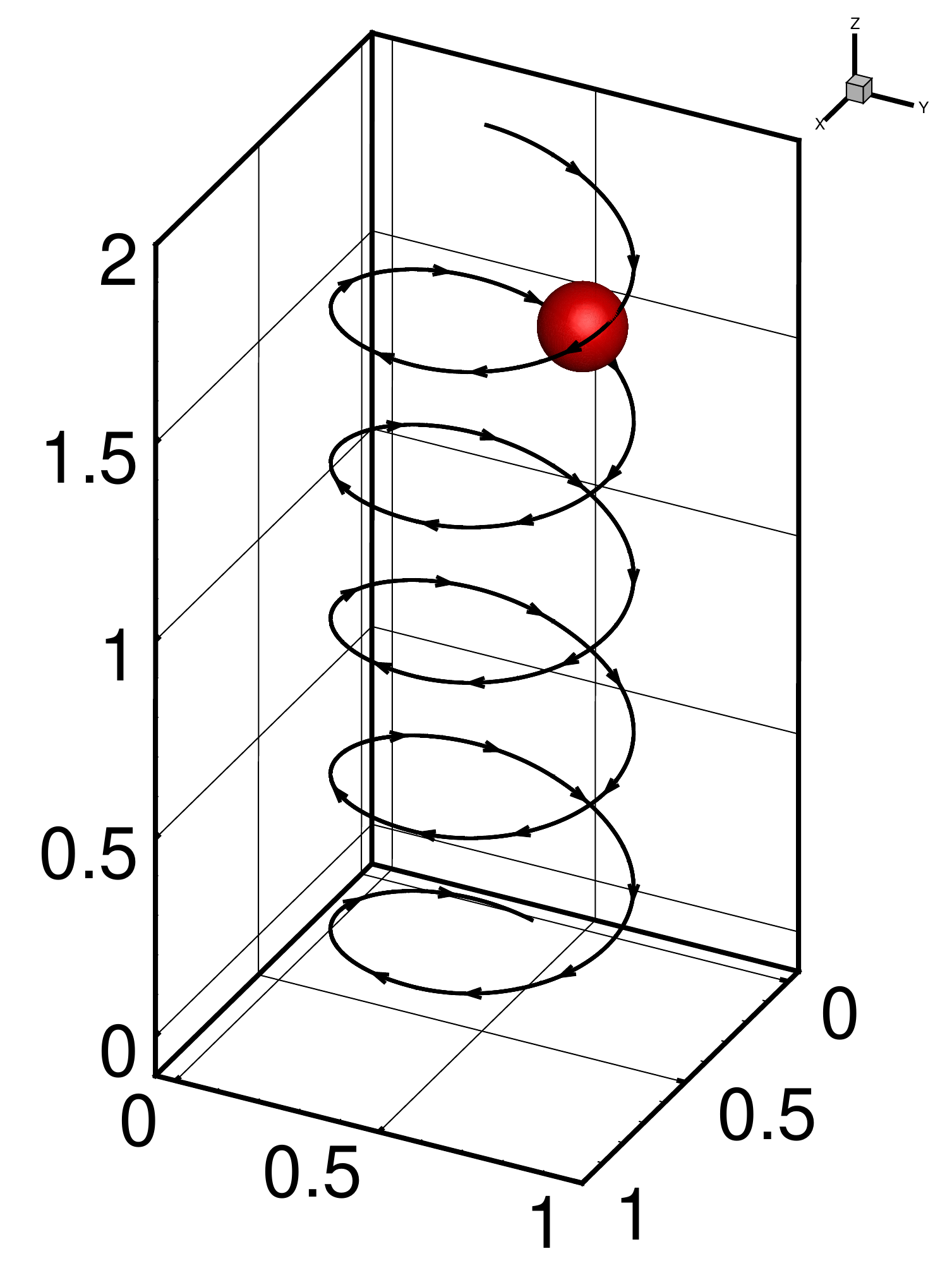}
\includegraphics[width=.33\textwidth]{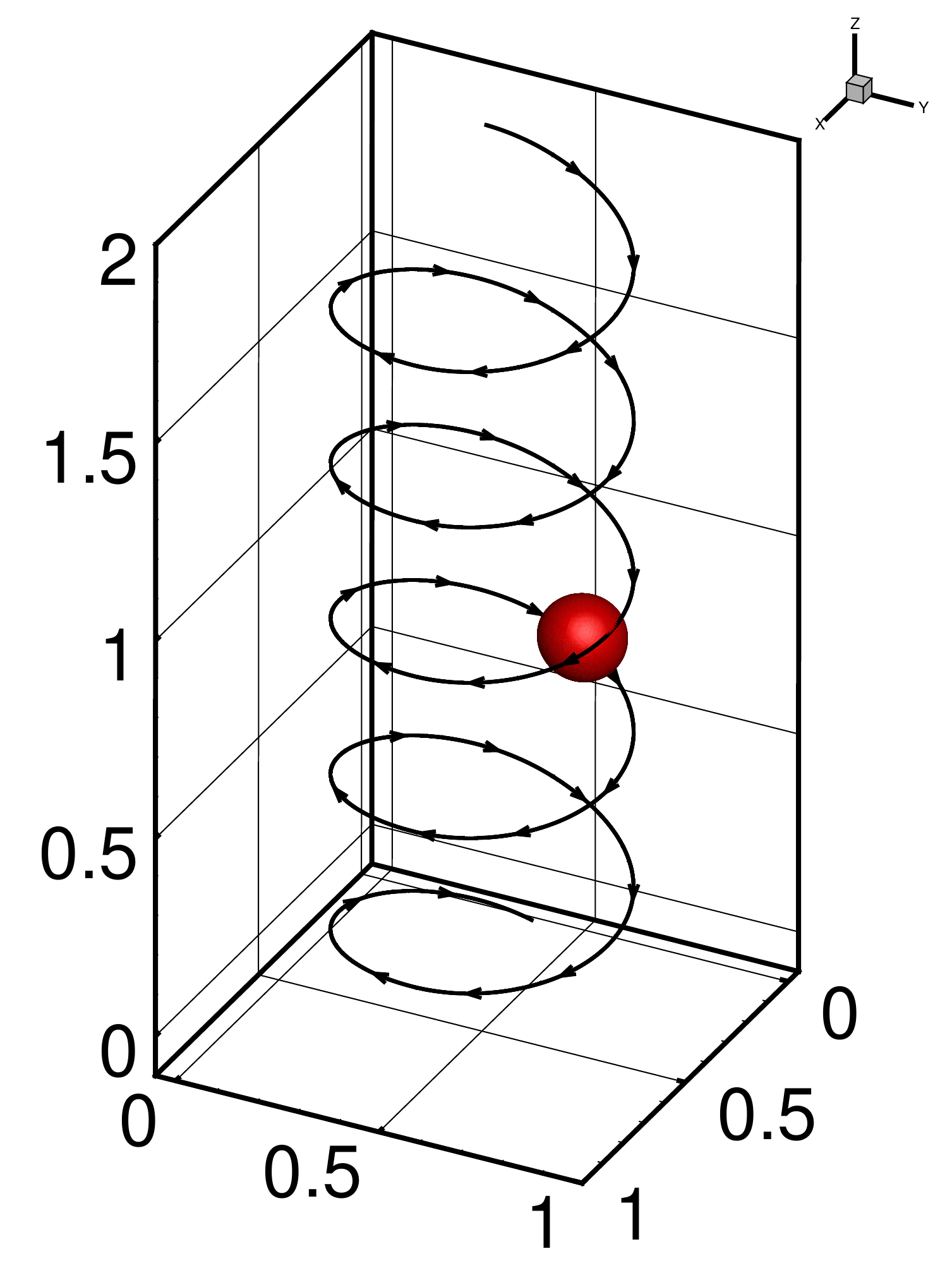}
\includegraphics[width=.33\textwidth]{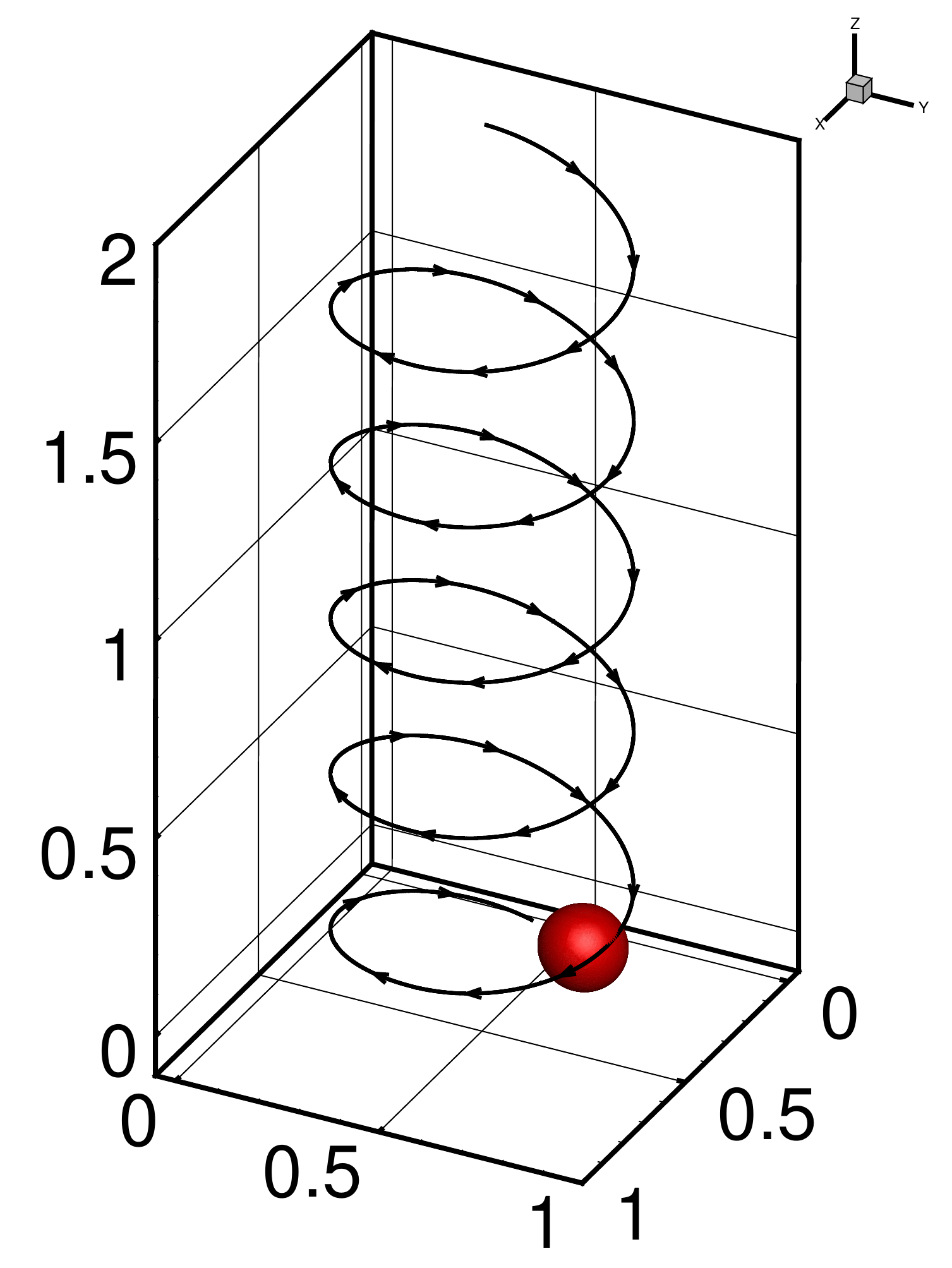}
\caption{Evolution of a sphere under constrained curvature-driven motion and Eulerian advection. The images are snapshots of the interface in time from left to right, top to bottom for $t=0.0,0.0125,$ and $0.025$.}
\label{sphere_adv}
\end{figure}

All the problems simulated in previous sections investigated pure interfacial motion in the normal direction. However, VOF is Eulerian in nature and a comprehensive test should reflect that property. Therefore, a sphere of $r=0.1$ was placed in a helical velocity field and allowed to evolve under constrained curvature-driven motion whilst being advected by the background velocity. The resulting advection equation is
\begin{equation}
\frac{\partial C}{\partial t} + \mathbf{u_{tot}}\cdot\nabla C=0 \quad \text{such that} \quad \mathbf{u_{tot}}=\mathbf{u}+\mathbf{u_{\Gamma}}
\end{equation}
where the helical velocity field $\mathbf{u}$ is defined as 
\begin{align}
u &= 2\pi \ U_{max} \ (y-y_c) \\
v &= 2\pi \ V_{max} \ (x_c-x) \\
w &= W_{max} \ \text{cos}^{-1}\phi
\end{align}
with $\phi=(x-x_c)/\sqrt{(x-x_c)^2+(y-y_c)^2}$. The computational domain is 2 units long in the $z$-direction and 1 unit long in the other two directions. The time-step is $\Delta t=0.000025$ for $75\times 75\times 150$ grid points and $0.0000125$ for $150\times 150\times 300$, $U_{max}=V_{max}=160$, $W_{max}=-40$, $\mathbf{x_c}=(0.5,0.75,0.25)$ and the total number of iterations is 1000 for the coarse grid and 2000 for the finer one. The boundary conditions are periodic in all three directions and the sphere is transported by the helical velocity field from its initial position to the bottom of the computational domain. The images in Fig. \ref{sphere_adv} are snapshots of the sphere as it advects to the bottom of the domain for $t=0.0,0.0125,$ and $0.025$.
\vspace*{0.5cm}

\begin{figure}[H]
    \centering
    \includegraphics[width=.45\textwidth]{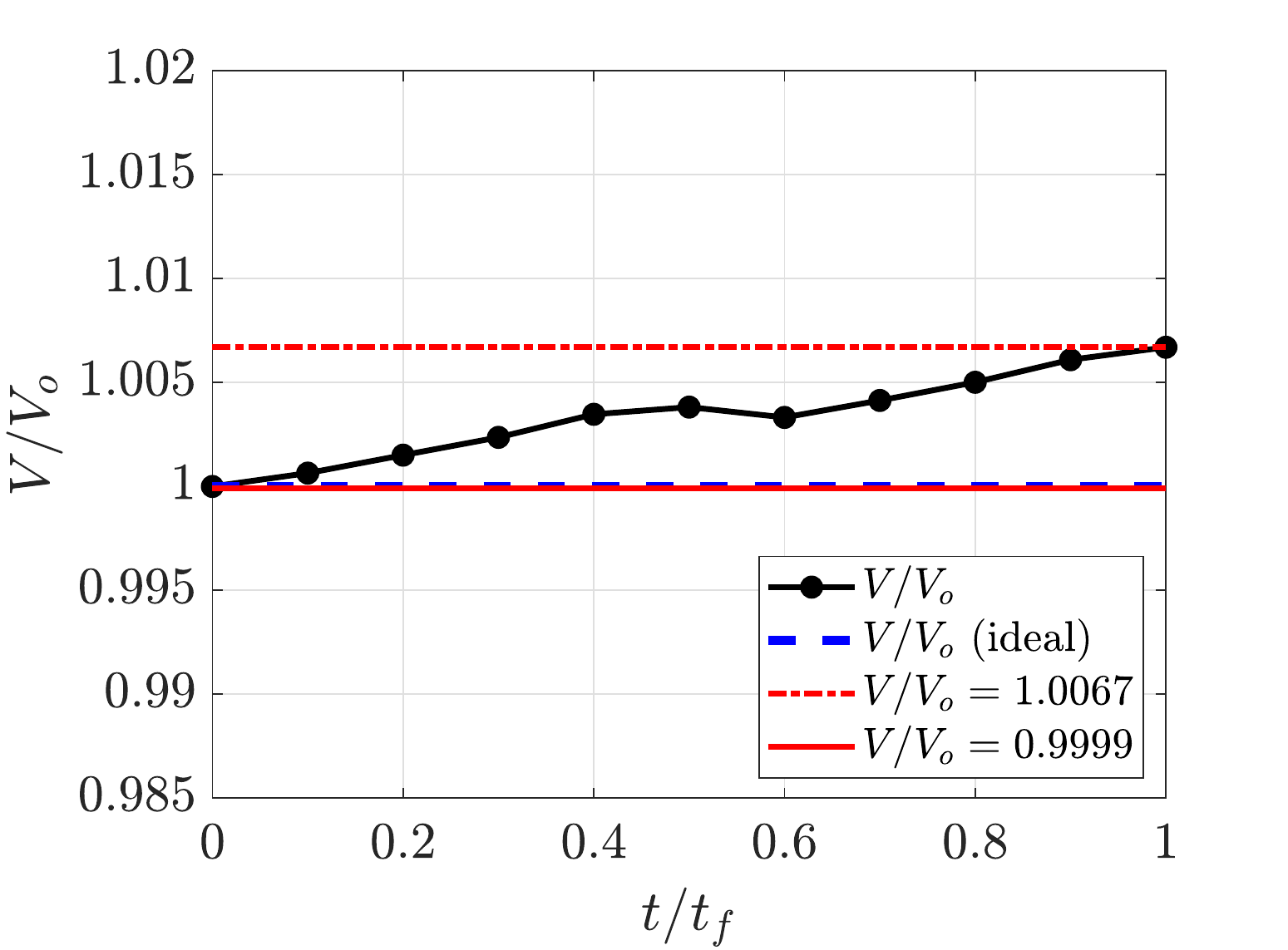}
    \includegraphics[width=.45\textwidth]{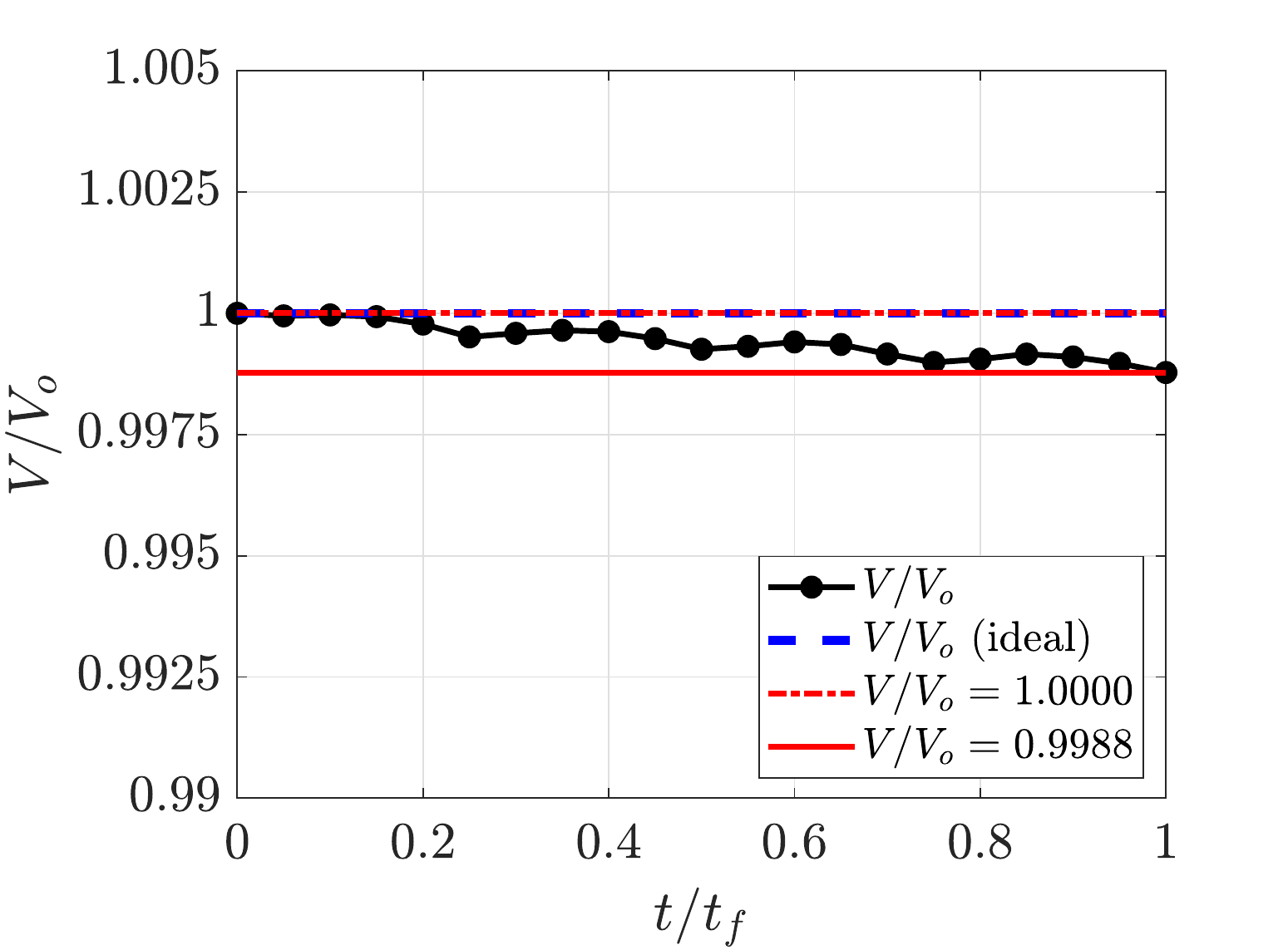}
    \caption{Variation of normalized volume of sphere with time under constrained curvature-driven motion with background velocity.}
    \label{sphere_conservation_advection}
\end{figure}

The sphere is shown to preserve its shape throughout the advection as shown in Fig. \ref{sphere_adv}.  The time scale of local interface equilibrium is relatively small compared to the time scale of the overall advection, this is particularly advantageous if the variational approach is applied to bubble-laden flows where the bubbles are used as tracers without the need for a surface tension model to maintain sphericity. Fig \ref{sphere_conservation_advection} shows the evolution of normalized volume versus non-dimensional time where the highest reported mass error is $0.67\%$ for the coarse grid and $0.12\%$ for the finer one.
\section{Summary and Concluding Remarks}
A variational volume-of-fluid approach whose kernel is the geometric VOF method is developed. The method is able to reproduce canonical interface evolution cases such as Rayleigh-Plesset (RP) bubble collapse as well traditional level set evolution cases in the context of curvature-driven flow. Key findings in this work can be summarized as follows:
\begin{enumerate}
\item Although algebraic advection algorithms promote flux conservation, they are not favorable at low grid resolution for problems involving pinch-off under curvature-driven motion due to the possible appearance of ``wisps".
\item The definition of the Dirac delta is of paramount importance when calculating $\bar{\kappa}$ for constrained curvature-driven motion. Using $\delta(C)=|\nabla C|$ for geometries with possible singularities in $\kappa$ promotes ``bad" velocities that end up rupturing the interface.
\item The degree to which enforced constraints are respected is an interplay between the complexity of interface topology and grid resolution.
\end{enumerate}
An intuitive extension of the current work is in the context of the Gibbs free energy framework whereby Gibbs free energy of a system is minimized. This moves the method from a numerical exercise to a physical model that can be used to solve real equilibrium problems. In that case, the energy functional $\mathscr{L}(C)$ includes additional terms such that
\begin{equation}
\mathscr{L}(C) = \mathcal{E}(C) + \mu G_1(C) + \lambda G_2(C)
\end{equation}
where $\mathcal{E}(C)$ contains the bulk energy due to the pressure difference at the solid-liquid interface and the interfacial energy, $G_1(C)$ represents the no penetration constraint between the interfaces, and $G_2(C)$ is the volume conservation constraint that we have discussed in this paper. Details of the extension of the current approach are left for future investigation.

\section*{Declaration of competing interest}
The authors declare that they have no known competing financial interests or personal relationships that could have appeared to influence the work reported in this paper.

\section*{Acknowledgments}
This work is supported by the United States Office of Naval Research (ONR) under ONR Grant N00014-17-1-2676 with Dr. Ki-Han Kim and Dr. Julie Young as the technical monitors. The computations were made possible through the Minnesota Supercomputing Institute (MSI) at the University of Minnesota.
\appendix
\section{Derivations of the energy functionals}
\subsection{Derivative of interfacial surface energy functional}
\noindent The interfacial surface energy functional is given by
\begin{equation}
	\mathcal{E}(C) = \int_{\Omega}\delta(C) \ |\nabla C| \ d\mathbf{x}
\end{equation}
The Fr\'echet derivative is defined as
\begin{equation}
	\Bigg{<}\frac{\partial \mathcal{E}}{\partial C},\chi\Bigg{>} = \lim_{\epsilon \to 0} \frac{\mathcal{E}(C+\epsilon\chi)-\mathcal{E}(C)}{\epsilon}
\end{equation}
where \\ \\
$\displaystyle{\lim_{\epsilon \to 0}\frac{\mathcal{E}(C+\epsilon\chi)-\mathcal{E}(C)}{\epsilon}} = \lim_{\epsilon \to 0}\frac{1}{\epsilon}\bigg{[}\int_{\Omega}\delta(C+\epsilon\chi) \ |\nabla(C+\epsilon\chi)| \ d\mathbf{x}-\int_{\Omega}\delta(C) \ |\nabla C| \ d\mathbf{x}\bigg{]}$ \\ \\
then add and subtract $\displaystyle{\int_{\Omega}\delta(C+\epsilon\chi) \ |\nabla C| \ d\mathbf{x}}$ from the right-hand side (RHS) such that \\
\begin{align*}
\lim_{\epsilon \to 0}\frac{\mathcal{E}(C+\epsilon\chi)-\mathcal{E}(C)}{\epsilon} &= \lim_{\epsilon \to 0}\frac{1}{\epsilon}\bigg{[}\int_{\Omega}\delta(C+\epsilon\chi) \ |\nabla(C+\epsilon\chi)| \ d\mathbf{x}-\int_{\Omega}\delta(C) \ |\nabla C| \ d\mathbf{x}+{\int_{\Omega}\delta(C+\epsilon\chi) \ |\nabla C| \ d\mathbf{x}} \ ... \\ & \mspace{50mu} - {\int_{\Omega}\delta(C+\epsilon\chi) \ |\nabla C| \ d\mathbf{x}}\bigg{]} \\ \\ &= \lim_{\epsilon \to 0}\frac{1}{\epsilon}\bigg{\{}\int_{\Omega}\delta(C+\epsilon\chi)\bigg{[}|\nabla(C+\epsilon\chi)|-|\nabla C|\bigg{]} \ d\mathbf{x}+\int_{\Omega}|\nabla C|\bigg{[}\delta(C+\epsilon\chi)-\delta(C)\bigg{]} \ d\mathbf{x} \bigg{\}} \\ \\
&= \int_{\Omega}\lim_{\epsilon \to 0}\frac{1}{\epsilon}\bigg{[}|\nabla(C+\epsilon\chi)|-|\nabla C|\bigg{]}\delta(C+\epsilon\chi) \ d\mathbf{x}+\int_{\Omega}\lim_{\epsilon \to 0}\frac{1}{\epsilon}\bigg{[}\delta(C+\epsilon\chi)-\delta(C)\bigg{]}|\nabla C| \ d\mathbf{x} \\ \\ &= 
\underbrace{\int_{\Omega}\lim_{\epsilon \to 0}\frac{1}{\epsilon}\bigg{[}[\nabla(C+\epsilon\chi)\cdot\nabla(C+\epsilon\chi)]^{1/2}-[\nabla C\cdot \nabla C]^{1/2}\bigg{]}\delta(C+\epsilon\chi) \ d\mathbf{x}}_{\text{Term I.}} + \underbrace{\int_{\Omega}\delta^{'}(C) \ |\nabla C| \ \chi \ d\mathbf{x}}_{\text{Term II.}}
\end{align*} \\
multiply and divide Term I. by $[\nabla(C+\epsilon\chi)\cdot\nabla(C+\epsilon\chi)]^{1/2} + [\nabla C\cdot \nabla C]^{1/2}$ such that \\
\begin{align*}
\lim_{\epsilon \to 0}\frac{\mathcal{E}(C+\epsilon\chi)-\mathcal{E}(C)}{\epsilon} &= \int_{\Omega}\lim_{\epsilon \to 0}\frac{1}{\epsilon} \bigg{[}[\nabla(C+\epsilon\chi)\cdot\nabla(C+\epsilon\chi)]^{1/2}-[\nabla C\cdot \nabla C]^{1/2}\bigg{]}  \ \times \\ & \mspace{50mu} \frac{\bigg{[}[\nabla(C+\epsilon\chi)\cdot\nabla(C+\epsilon\chi)]^{1/2}+[\nabla C\cdot \nabla C]^{1/2}\bigg{]}}{\bigg{[}[\nabla(C+\epsilon\chi)\cdot\nabla(C+\epsilon\chi)]^{1/2}+[\nabla C\cdot \nabla C]^{1/2}\bigg{]}} \ \delta(C+\epsilon\chi) \ d\mathbf{x} + \int_{\Omega}\delta^{'}(C) \ |\nabla C| \ \chi \ d\mathbf{x} \\ \\ & = \int_{\Omega}\lim_{\epsilon \to 0}\frac{1}{\epsilon} \frac{[\nabla(C+\epsilon\chi)\cdot\nabla(C+\epsilon\chi)]-|\nabla C|^2}{{[\nabla(C+\epsilon\chi)\cdot\nabla(C+\epsilon\chi)]^{1/2}+[\nabla C\cdot \nabla C]^{1/2}} } \ \delta(C+\epsilon\chi) \ d\mathbf{x} +  \int_{\Omega}\delta^{'}(C) \ |\nabla C| \ \chi \ d\mathbf{x} \\ \\ &= \int_{\Omega}\lim_{\epsilon \to 0}\frac{1}{\epsilon} \frac{(\nabla C \cdot \nabla C + 2 \ \epsilon \nabla \phi \cdot \nabla \chi+\epsilon^2\nabla\chi\cdot\nabla\chi)- |\nabla C|} {(\nabla C \cdot \nabla C + 2 \ \epsilon \nabla \phi \cdot \nabla \chi+\epsilon^2\nabla\chi\cdot\nabla\chi)^{1/2}+|\nabla C|} \ \delta(C+\epsilon\chi) \ d\mathbf{x} +  \int_{\Omega}\delta^{'}(C) \ |\nabla C| \ \chi \ d\mathbf{x} \\ \\ &= \int_{\Omega}\lim_{\epsilon \to 0}\frac{(\nabla C \cdot \nabla C + 2 \ \nabla \phi \cdot \nabla \chi+\epsilon\nabla\chi\cdot\nabla\chi)- |\nabla C|} {(\nabla C \cdot \nabla C + 2 \ \nabla \phi \cdot \nabla \chi+\epsilon\nabla\chi\cdot\nabla\chi)^{1/2}+|\nabla C|} \ \delta(C+\epsilon\chi) \ d\mathbf{x} +  \int_{\Omega}\delta^{'}(C) \ |\nabla C| \ \chi \ d\mathbf{x} \\ \\ &= \int_{\Omega}\delta(C) \ \frac{\nabla C}{|\nabla C|}\cdot\chi \ d\mathbf{x} +  \int_{\Omega}\delta^{'}(C) \ |\nabla C| \ \chi \ d\mathbf{x}
 \end{align*} \\
 Therefore, 
 \begin{align*}
 \lim_{\epsilon \to 0}\frac{\mathcal{E}(C+\epsilon\chi)-\mathcal{E}(C)}{\epsilon} = \int_{\Omega}\delta(C) \ \frac{\nabla C}{|\nabla C|}\cdot\chi \ d\mathbf{x} +  \int_{\Omega}\delta^{'}(C) \ |\nabla c| \ \chi \ d\mathbf{x}
 \end{align*}
 Using the following relation,
 \begin{align*}
 \int_{\Omega}\nabla\cdot\Bigg{(}\delta(C) \ \frac{\nabla C}{|\nabla C|}\ \chi\Bigg{)} \ d\mathbf{x} = \int_{\Omega}\delta(C) \ \frac{\nabla C}{|\nabla C|}\cdot\nabla \chi \ d\mathbf{x} + \int_{\Omega}\nabla\cdot\Bigg{(}\delta(C)\ \frac{\nabla C}{|\nabla C|}\Bigg{)} \ \chi \ d\mathbf{x} 
 \end{align*}
we can rewrite 
\begin{align*}
\int_{\Omega}\delta(C) \ \frac{\nabla C}{|\nabla C|}\cdot\nabla\chi \ d\mathbf{x} &= \int_{\Omega}\nabla\cdot\Bigg{(}\delta(C) \ \frac{\nabla C}{|\nabla C|}\ \chi\Bigg{)} \ d\mathbf{x} - \int_{\Omega}\nabla\cdot\Bigg{(}\delta(C)\ \frac{\nabla C}{|\nabla C|}\Bigg{)} \ \chi \ d\mathbf{x} \\ \\
&= \int_{\partial\Omega}\delta(C) \ \frac{\nabla C}{|\nabla C|}\cdot\mathbf{n} \ \chi \ d\mathbf{s} - \int_{\Omega}\nabla\delta(C)\cdot\Bigg{(}\frac{\nabla C}{|\nabla C|}\Bigg{)} \ \chi \ d\mathbf{x} - \int_{\Omega}\delta(C) \ \nabla\cdot\Bigg{(}\frac{\nabla C}{|\nabla C|}\Bigg{)} \ \chi \ d\mathbf{x} \\ \\ 
&=  \int_{\partial\Omega}\delta(C) \ \frac{1}{|\nabla C|}\frac{\partial C}{\partial n} \ \chi \ d\mathbf{s} - \int_{\Omega}\delta^{'}(C)\nabla C\cdot\Bigg{(}\frac{\nabla C}{|\nabla C|}\Bigg{)} \ \chi \ d\mathbf{x} - \int_{\Omega}\delta(C) \ \nabla\cdot\Bigg{(}\frac{\nabla C}{|\nabla C|}\Bigg{)} \ \chi \ d\mathbf{x} \\ \\ &=  \int_{\partial\Omega}\delta(C) \ \frac{1}{|\nabla C|}\frac{\partial C}{\partial n} \ \chi \ d\mathbf{s} - \int_{\Omega}\delta^{'}(C)\ |\nabla C| \ \chi \ d\mathbf{x} - \int_{\Omega}\delta(C) \ \nabla\cdot\Bigg{(}\frac{\nabla C}{|\nabla C|}\Bigg{)} \ \chi \ d\mathbf{x}
\end{align*} 
Henceforth,
\begin{align*}
\Bigg{<}\frac{\partial \mathcal{E}}{\partial C},\chi\Bigg{>} &= \int_{\partial\Omega}\delta(C) \ \frac{1}{|\nabla C|}\frac{\partial C}{\partial n} \ \chi \ d\mathbf{s} - \int_{\Omega}\delta(C) \ \nabla\cdot\Bigg{(}\frac{\nabla C}{|\nabla C|}\Bigg{)} \ \chi \ d\mathbf{x}
\end{align*}
Assuming zero-Neumann conditions on $\Omega$ implies $\displaystyle{\frac{\partial C}{\partial n}=0}$, then
\begin{equation}
\boxed{\Bigg{<}\frac{\partial \mathcal{E}}{\partial C},\chi\Bigg{>} = - \int_{\Omega}\delta(C) \ \nabla\cdot\Bigg{(}\frac{\nabla C}{|\nabla C|}\Bigg{)} \ \chi \ d\mathbf{x}}
\end{equation}

\subsection{Derivative of volume conservation constraint}
\noindent The volume conservation constraint is given by
\begin{equation}
G(C) = \int_{\Omega} \ C \ d\mathbf{x} = V_o
\end{equation}
where $V_o$ is the initial volume of the tracked phase. The Fr\'echet derivative is defined as
\begin{equation}
	\Bigg{<}\frac{\partial G}{\partial C},\chi\Bigg{>} = \lim_{\epsilon \to 0} \frac{G(C+\epsilon\chi)-G(C)}{\epsilon}
\end{equation}
where
\begin{align*}
\lim_{\epsilon \to 0} \frac{G(C+\epsilon\chi)-G(C)}{\epsilon} &= \lim_{\epsilon \to 0}\frac{1}{\epsilon}\Bigg{[}\int_{\Omega} \ (C+\epsilon\chi) \ d\mathbf{x}-\int_{\Omega} \ C \ d\mathbf{x}\Bigg{]} \\ \\ &= \int_{\Omega}\lim_{\epsilon \to 0}\frac{(C+\epsilon\chi)-C}{\epsilon} \ d\mathbf{x} \\ \\ &= \int_{\Omega}\delta(C) \ \chi \ d\mathbf{x}
\end{align*}
Therefore,
\begin{equation}
	\boxed{\Bigg{<}\frac{\partial G}{\partial C},\chi\Bigg{>} = \int_{\Omega}\delta(C) \ \chi \ d\mathbf{x}}
\end{equation}

\bibliographystyle{model1-num-names}
\bibliography{jcp_vvof}

\end{document}